\documentclass{cimart}
\usepackage{longtable}
 
\newtheorem*{problem1}{Problem 1}
\newtheorem*{problem2}{Problem 2}

\title[The geometric classification of non-associative algebras]{%
  The geometric classification of non-associative algebras: \\ a survey
    }

\author{%
    Ivan Kaygorodov, Mykola Khrypchenko and Pilar P\'{a}ez-Guill\'{a}n
    }

\authorinfo[%
    Ivan Kaygorodov]{
    CMA-UBI, University of Beira  Interior,  Portugal}{%
    kaygorodov.ivan@gmail.com
    }

\authorinfo[%
    Mykola Khrypchenko]{
      Federal University of Santa Catarina,   Brazil}{%
    nskhripchenko@gmail.com
    }

\authorinfo[%
    Pilar P\'{a}ez-Guill\'{a}n]{
      University of Santiago de Compostela,  Spain}{%
    pilar.paez@usc.es
    }

\abstract{%
This is a survey on the geometric classification of different varieties of algebras  
(nilpotent, 
nil-,
associative, 
commutative associative, 
cyclic  associative, 
Jordan, 
Kokoris, 
standard,
noncommutative Jordan, 
commutative power-associative 
weakly associative, 
terminal, 
Lie, 
Malcev, 
binary Lie, 
Tortkara, 
dual mock Lie, 
$\mathfrak{CD}$-,
commutative $\mathfrak{CD}$-,
anticommutative $\mathfrak{CD}$-,
symmetric Leibniz, 
Leibniz, 
Zinbiel, 
Novikov, 
bicommutative,  
assosymmetric, 
antiassociative, 
left-symmetric,
right alternative, 
and right commutative), 
$n$-ary algebras (Fillipov ($n$-Lie),  
Lie triple systems and anticommutative ternary),
superalgebras 
(Lie and Jordan),
and Poisson-type algebras
(Poisson, 
transposed Poisson, 
Leibniz-Poisson,
generic Poisson,
generic Poisson-Jordan,
transposed Leibniz-Poisson,
Novikov-Poisson, 
pre-Lie Poisson,
commutative pre-Lie, 
anti-pre-Lie Poisson,
pre-Poisson, 
compatible commutative associative,
compatible  associative,
compatible Novikov,
compatible pre-Lie).
We also discuss the degeneration level classification.
    }

\keywords{%
geometric classification,  
non-associative algebra,  superalgebra, $n$-ary algebra.
    }

\msc{%
17A30 (primary); 14L30,  17A40, 17A42, 17A70 (secondary).
    }

\VOLUME{32}
\YEAR{2024}
\ISSUE{2}
\firstpage{185}
\DOI{https://doi.org/10.46298/cm.14458}

\begin{document}

\tableofcontents

\section*{Introduction}

Given a complex $n$-dimensional vector space $\mathbb V$, the set ${\rm Hom}(\mathbb V \otimes \mathbb V,\mathbb V) \cong \mathbb V^* \otimes \mathbb V^* \otimes \mathbb V$
is a complex vector space of dimension $n^3$. This space has the structure of the affine variety $\mathbb{C}^{n^3}$. Indeed, if we fix a basis $\{e_1,\dots,e_n\}$ of $\mathbb V$, then any $\mu\in {\rm Hom}(\mathbb V \otimes \mathbb V,\mathbb V)$ is determined by $n^3$ structure constants $c_{ij}^k\in\mathbb{C}$ such that
$\mu(e_i\otimes e_j)=\sum_{k=1}^nc_{ij}^ke_k$. A subset of ${\rm Hom}(\mathbb V \otimes \mathbb V,\mathbb V)$ is {\it Zariski-closed} if it is the set of solutions of a system of polynomial equations in the variables $c_{ij}^k$ ($1\le i,j,k\le n$).

Let $T$ be a set of polynomial identities.
Every algebra structure on $\mathbb V$ satisfying polynomial identities from $T$ forms a Zariski-closed subset of the variety ${\rm Hom}(\mathbb V \otimes \mathbb V,\mathbb V)$. We denote this subset by $\mathbb{L}(T)$.
The general linear group ${\rm GL}(\mathbb V)$ acts on $\mathbb{L}(T)$ by conjugation:
$$ (g * \mu )(x\otimes y) = g\mu(g^{-1}x\otimes g^{-1}y)$$
for $x,y\in \mathbb V$, $\mu\in \mathbb{L}(T)\subset {\rm Hom}(\mathbb V \otimes\mathbb V, \mathbb V)$ and $g\in {\rm GL}(\mathbb V)$.
Thus, $\mathbb{L}(T)$ decomposes into ${\rm GL}(\mathbb V)$-orbits that correspond to the isomorphism classes of the algebras.
We shall denote by ${\mathcal O}(\mu)$ the orbit of $\mu\in\mathbb{L}(T)$ under the action of ${\rm GL}(\mathbb V)$ and by $\overline{O(\mu)}$ its Zariski closure.

Let $\mathcal A$ and $\mathcal B$ be two $n$-dimensional algebras satisfying the identities from $T$, and let $\mu,\lambda \in \mathbb{L}(T)$ represent $\mathcal A$ and $\mathcal B$, respectively.
We say that $\mathcal A$ {\it degenerates to} $\mathcal B$,
and write $\mathcal A\to \mathcal B$,
if $\lambda\in\overline{{\mathcal O}(\mu)}$.
Note that this implies $\overline{{\mathcal O}(\lambda)}\subset\overline{{\mathcal O}(\mu)}$. Hence, the definition of a degeneration does not depend on the choice of $\mu$ and $\lambda$. If $\mathcal A\not\cong \mathcal B$, then the assertion $\mathcal A\to \mathcal B$ is called a {\it proper degeneration}.
Following Gorbatsevich~\cite{gorb91}, we say that $\mathcal A$ has {\it level} $m$ if there exists a chain of proper degenerations of length $m$ starting in $\mathfrak{A}$ and there is no such chain of length $m+1$. Also, in~\cite{gorb93} it was introduced the notion of \emph{infinite level} of an algebra $\mathcal A$ as the limit of the usual levels of $\mathcal A \oplus \mathbb{C}^m$.

Let $\mathcal A$ be represented by $\mu\in\mathbb{L}(T)$. Then  $\mathcal A$ is  {\it rigid} in $\mathbb{L}(T)$ if ${\mathcal O}(\mu)$ is an open subset of $\mathbb{L}(T)$.
 Recall that a subset of a variety is called {\it irreducible} if it cannot be represented as a union of two non-trivial closed subsets.
 A maximal irreducible closed subset of a variety is called an {\it irreducible component}.
It is well known that any affine variety can be uniquely represented as a finite union of its irreducible components. Note that the algebra $\mathcal A$ is rigid in $\mathbb{L}(T)$ if and only if $\overline{{\mathcal O}(\mu)}$ is an irreducible component of $\mathbb{L}(T)$.

In this survey, we discuss the following problems:

\begin{problem1}[Geometric classification]
Let $\mathfrak{V}^n$ be a variety of $n$-dimen\-sio\-nal algebras defined by a family of identities $T$.
What are the irreducible components of $\mathfrak{V}^n$?
\end{problem1}

\begin{problem2}[Level classification]
Let $\mathfrak{V}^n$ be a variety of $n$-dimensional algebras defined by a family of identities $T$.
Which algebras from $\mathfrak{V}^n$ have level $m$?
\end{problem2}

Although we will not deal with the infinite level of algebras in this survey, some references that address this issue are~\cite{gorb93,gorb98,kv17,wolf1}.

We shall use a lot of papers on algebraic classification throughout the text without giving a precise reference. In such a situation we mention only the authors and the year of publication --- the precise references can be easily found in the articles on geometric classification, on MathSciNet, or in the references of \cite{ksurv}.

In all the multiplication tables,  the first column contains our notations for algebras, the second column is reserved for the original notation (from the original paper), and in the last column we have the multiplication law.  In particular, some authors prefer to denote the direct product of a finite number of algebras using the symbol $\times$, while some other algebraists use the direct sum operation $\oplus$. In such a case, in the second column, we follow the original authors' notation. 

\section{The geometric classification of algebras}
Throughout this section, we summarize the geometric classification of different varieties of (not necessarily associative) algebras over the field $\mathbb{C}$. In what follows, we will not refer to the base field anymore.

We will use the following notations:
\[\begin{array}{rcl}
[x,y]&:=& xy -yx \\
 x \circ  y&:=& xy+yx \\
(x,y,z)_{*} &:=& (x*y)*z-x*(y*z),\\
J(x,y,z)_{*} &:=& (x*y)*z+(y*z)*x+(z*x)*y.\\
\end{array}\]
We will define any $n$-dimensional algebra via its multiplication table in a fixed basis $\{e_1,\dots,e_n\}$, omitting the products that are zero. Moreover, in the commutative case we will write only the products $e_i e_j$ with $i\leq j$, and in the anticommutative case only the products $e_ie_j$  with $i<j$. On the left of the multiplication table we will write the name given to the algebra in the paper where the corresponding geometric classification was established.

Let $\mathfrak{V}$ be the class of algebras defined by a family of polynomial identities.
We denote:
\begin{itemize}
    \item the variety of $n$-dimensional $\mathfrak{V}$-algebras by  $\mathfrak{V}^n;$  
    \item the variety of $n$-dimensional nilpotent $\mathfrak{V}$-algebras by  $\mathfrak{NV}^n;$  
    \item the variety of $n$-dimensional commutative $\mathfrak{V}$-algebras by $\mathfrak{CV}^n;$
    \item the variety of $n$-dimensional anticommutative $\mathfrak{V}$-algebras by $\mathfrak{AV}^n.$
    
\end{itemize}

 \subsection{Non-associative algebras} 
The first case to be considered is the variety of all algebras. 
It is easy to prove that the variety of all $n$-dimensional algebras (as well as the varieties of $n$-dimensional commutative and anticommutative algebras) has only one irreducible component defined by an algebra or a family of algebras.
For example, in \cite{kv16} such a family of algebras 

{\small

\begin{longtable}{|l|l|llll|}
\caption*{ }\\

\hline \multicolumn{1}{|c|}{${\mathcal A}$} & \multicolumn{1}{c|}{ } & \multicolumn{4}{c|}{Multiplication table} \\ \hline 
\endfirsthead

 \multicolumn{6}{l}%
{{\bfseries  continued from previous page}} \\
\hline \multicolumn{1}{|c|}{${\mathcal A}$} & \multicolumn{1}{c|}{ } & \multicolumn{4}{c|}{Multiplication table} \\ \hline 
\endhead

\hline \multicolumn{6}{|r|}{{Continued on next page}} \\ \hline
\endfoot

\hline 
\endlastfoot

$\mathfrak{A}_{1}^2(\alpha,\beta,\gamma,\delta)$ & ${\bf E}_1(\alpha,\beta,\gamma,\delta)$ & $ e_1e_1=e_1$ & $e_1e_2=\alpha e_1+ \beta  e_2$ & $e_2e_1=\gamma  e_1 + \delta  e_2$ &  $e_2e_2=e_2$  \\ \hline

\end{longtable}
}\noindent was found in the variety of all $2$-dimensional algebras, whose algebraic classification has been obtained in different ways by Ananin and Mironov (2000), Petersson (2000), Goze and Remm (2011), and also in~\cite{kv16}.
In the case of $2$-dimensional commutative algebras (see \cite{kv16} for both an explicit list and the geometric classification), the irreducible component is defined by 

{\small

\begin{longtable}{|l|l|lll|}
\caption*{ }   \\

\hline \multicolumn{1}{|c|}{${\mathcal A}$} & \multicolumn{1}{c|}{ } & \multicolumn{3}{c|}{Multiplication table} \\ \hline 
\endfirsthead

 \multicolumn{5}{l}%
{{\bfseries  continued from previous page}} \\
\hline \multicolumn{1}{|c|}{${\mathcal A}$} & \multicolumn{1}{c|}{ } & \multicolumn{3}{c|}{Multiplication table} \\ \hline 
\endhead

\hline \multicolumn{5}{|r|}{{Continued on next page}} \\ \hline
\endfoot

\hline 
\endlastfoot

$\mathfrak{C A}_{1}^2(\alpha,\beta)$ & ${\bf E}_1(\alpha,\beta,\alpha,\beta)$  & $e_1e_1=e_1$ &  $e_1e_2=\alpha e_1+ \beta  e_2$ & $e_2e_2=e_2$ \\ \hline

\end{longtable}
}

The variety of $2$-dimensional anticommutative algebras has   one non-zero algebra   \cite{kv16}:
{\small

\begin{longtable}{|l|l|l|}
\caption*{ }   \\

\hline \multicolumn{1}{|c|}{${\mathcal A}$} & \multicolumn{1}{c|}{ } & \multicolumn{1}{c|}{Multiplication table} \\ \hline 
\endfirsthead

 \multicolumn{3}{l}%
{{\bfseries  continued from previous page}} \\
\hline \multicolumn{1}{|c|}{${\mathcal A}$} & \multicolumn{1}{c|}{ } & \multicolumn{1}{c|}{Multiplication table} \\ \hline 
\endhead

\hline \multicolumn{3}{|r|}{{Continued on next page}} \\ \hline
\endfoot

\hline 
\endlastfoot

$\mathfrak{AA}_{1}^2$ & ${\bf B}_3$ & $e_1e_2=e_2$  \\ \hline

\end{longtable}
}

\sloppy
The variety of $3$-dimensional anticommutative algebras has been classified algebraically  and geometrically in~\cite{ikv18}, and its irreducible component is defined by

{\small

\begin{longtable}{|l|l|lll|}
\caption*{ }  \\

\hline \multicolumn{1}{|c|}{${\mathcal A}$} & \multicolumn{1}{c|}{ } & \multicolumn{3}{c|}{Multiplication table} \\ \hline 
\endfirsthead

 \multicolumn{5}{l}%
{{\bfseries  continued from previous page}} \\
\hline \multicolumn{1}{|c|}{${\mathcal A}$} & \multicolumn{1}{c|}{ } & \multicolumn{3}{c|}{Multiplication table} \\ \hline 
\endhead

\hline \multicolumn{5}{|r|}{{Continued on next page}} \\ \hline
\endfoot

\hline 
\endlastfoot

$\mathfrak{AA}_{1}^3(\alpha)$ & $\mathcal{A}_1^{\alpha}$  & $e_1e_2 = e_3$ & $e_1e_3=e_1+e_3$ & $e_2e_3=\alpha e_2$ \\ \hline

\end{longtable}
}

 \subsubsection{Nilpotent algebras} 

First of all, let us recall the definition of nilpotent algebras. Given an arbitrary algebra $\mathfrak{N}$, we consider the series 
\[
\begin{array}{ccl} 
\mathfrak{N}^1=\mathfrak{N}, & \mathfrak{N}^{i+1}=\sum_{k=1}^{i} \mathfrak{N}^k \mathfrak{N}^{i+1-k}, & i\geq 1.
\end{array} \]
We say that $\mathfrak{N}$ is {\it nilpotent} if $\mathfrak{N}^i=0$ for some $i\geq 1$.

There is only one non-trivial $2$-dimensional nilpotent algebra.
The algebraic classification of all nilpotent algebras of dimension $3$ was given in a paper by Calder\'on Mart\'in,
Fern\'andez Ouaridi and Kaygorodov (2022). Using this result, in~\cite{fkkv} the authors constructed all the degenerations in the variety $\mathfrak{Nil}^3$ of nilpotent algebras of dimension $3$, showing that it has only one irreducible component defined by

{\small

\begin{longtable}{|l|l|lll|}
\caption*{ }  \\

\hline \multicolumn{1}{|c|}{${\mathcal A}$} & \multicolumn{1}{c|}{ } & \multicolumn{3}{c|}{Multiplication table} \\ \hline 
\endfirsthead

 \multicolumn{5}{l}%
{{\bfseries  continued from previous page}} \\
\hline \multicolumn{1}{|c|}{${\mathcal A}$} & \multicolumn{1}{c|}{ } & \multicolumn{3}{c|}{Multiplication table} \\ \hline 
\endhead

\hline \multicolumn{5}{|r|}{{Continued on next page}} \\ \hline
\endfoot

\hline 
\endlastfoot

$\mathfrak{N}_{1}^3$ & $\mathcal{N}_2$ & $e_1 e_1 = e_2$ & $e_2 e_1= e_3$ & $e_2 e_2=e_3$ \\ \hline

\end{longtable}
}

{ 
The present result was generalized in \cite{kkl22}.

\begin{theorem}[Theorem A, \cite{kkl22}]
For any $n\ge 2$, the variety of all $n$-dimensional nilpotent algebras is irreducible and has dimension $\frac{n(n-1)(n+1)}{3}$.
\end{theorem}

Let $n\geq 3$. Denote by $\mathfrak{R}_n$ the family of nilpotent algebras with basis $(e_i)_{i=1}^{n}$, whose structure constants $(c_{ij}^{k})_{i,j,k=1}^n,$ satisfy 
$c_{ij}^{k}=0,\ \forall k \leq \max\{i,j\},$
and  
$e_i^2=e_{i+1}$, for all $1\le i\le n-1$, $c_{21}^3=1$, $c_{1i}^{i+1}=0$, for all $2\le i\le n-1$, and with the remaining structure constants $c_{ij}^{k}$ being arbitrary independent complex parameters, for all $k>\max\{i,j\}$ and $1\le i\ne j\le n$.
It was shown that the family $\mathfrak{R}_n$   is generic in the variety of $n$-dimensional nilpotent algebras and inductively gives an algorithmic procedure to obtain any $n$-dimensional nilpotent algebra as a degeneration from $\mathfrak{R}_n$ \cite{kkl22}.

}

\subsubsection{Nilpotent commutative algebras} 
Let $\mathfrak{NCom}$ be the variety of nilpotent commutative algebras. 
Thanks to \cite{fkkv} we have the description of the geometry of the varieties $\mathfrak{NCom}^3$  and $\mathfrak{NCom}^4$. While the complete list of $3$-dimensional nilpotent commutative algebras can be extracted from Calder\'on Mart\'in, Fern\'andez Ouaridi and Kaygorodov (2018), in dimension $4$ the algebraic classification was made in~\cite{fkkv}.

The irreducible component of the variety $\mathfrak{NCom}^3$ is defined by 

{\small

\begin{longtable}{|l|l|ll|}
\caption*{ }   \\

\hline \multicolumn{1}{|c|}{${\mathcal A}$} & \multicolumn{1}{c|}{ } & \multicolumn{2}{c|}{Multiplication table} \\ \hline 
\endfirsthead

 \multicolumn{4}{l}%
{{\bfseries  continued from previous page}} \\
\hline \multicolumn{1}{|c|}{${\mathcal A}$} & \multicolumn{1}{c|}{ } & \multicolumn{2}{c|}{Multiplication table} \\ \hline 
\endhead

\hline \multicolumn{4}{|r|}{{Continued on next page}} \\ \hline
\endfoot

\hline 
\endlastfoot

$\mathfrak{NC}_{1}^3$ & $\mathcal{C}_{02}$ & $e_1e_1 = e_2$ & $e_2e_2=e_3$\\ \hline

\end{longtable}
}

The irreducible component of the variety $\mathfrak{NCom}^4$ is defined by 

{\small

\begin{longtable}{|l|l|lllll|}
\caption*{ }   \\

\hline \multicolumn{1}{|c|}{${\mathcal A}$} & \multicolumn{1}{c|}{ } & \multicolumn{5}{c|}{Multiplication table} \\ \hline 
\endfirsthead

 \multicolumn{7}{l}%
{{\bfseries  continued from previous page}} \\
\hline \multicolumn{1}{|c|}{${\mathcal A}$} & \multicolumn{1}{c|}{ } & \multicolumn{5}{c|}{Multiplication table} \\ \hline 
\endhead

\hline \multicolumn{7}{|r|}{{Continued on next page}} \\ \hline
\endfoot

\hline 
\endlastfoot

$\mathfrak{NC}_{1}^4(\alpha)$ & $\mathcal{C}_{19}(\alpha)$ &  $e_1 e_1 = e_2$ & $e_1e_3=\alpha e_4$ & $e_2 e_2=e_3$ & $e_2e_3= e_4$ & $e_3e_3=e_4$  \\ \hline

\end{longtable}
}

The complete graph of degenerations can be found in~\cite{fkkv}.{  
The present result was generalized in \cite{kkl22}.

\begin{theorem}[Theorem B, \cite{kkl22}]
For any $n\ge 2$, the variety of all $n$-dimensional commutative nilpotent algebras is irreducible and has dimension $\frac{n(n-1)(n+4)}{6}$.
\end{theorem}

Let $n\ge 4$. Denote by $\mathfrak{S}_n$ the family of commutative nilpotent  algebras with basis $(e_i)^n_{i=1},$
whose structure constants $(c_{ij}^{k})_{i,j,k=1}^n,$ satisfy  
$c_{ij}^{k}=0,\ \forall k \leq \max\{i,j\},$
and  
$e_i^2=e_{i+1}$ for all $1\le i\le n-1$, $c_{23}^4=1$, $c_{12}^4\neq 0$ and $c_{1i}^{i+1}=0$ for all $2\le i\le n-1$. The remaining structure constants $c_{ij}^{k}$ are arbitrary.
As above, it was shown that the family $\mathfrak{S}_n$  is generic in the variety of $n$-dimensional commutative nilpotent algebras and inductively gives an algorithmic procedure to obtain any $n$-dimensional nilpotent commutative algebra as a degeneration from $\mathfrak{S}_n$ \cite{kkl22}.

}

\subsubsection{Nilpotent anticommutative algebras} 
Let $\mathfrak{NACom}$ be the variety of nilpotent anticommutative algebras.  

The irreducible component of the variety $\mathfrak{NACom}^3$ is defined by the unique nilpotent anticommutative algebra of dimension $3$: 

{\small

\begin{longtable}{|l|l|l|}
\caption*{ }   \\

\hline \multicolumn{1}{|c|}{${\mathcal A}$} & \multicolumn{1}{c|}{ } & \multicolumn{1}{c|}{Multiplication table} \\ \hline 
\endfirsthead

 \multicolumn{3}{l}%
{{\bfseries  continued from previous page}} \\
\hline \multicolumn{1}{|c|}{${\mathcal A}$} & \multicolumn{1}{c|}{ } & \multicolumn{1}{c|}{Multiplication table} \\ \hline 
\endhead

\hline \multicolumn{3}{|r|}{{Continued on next page}} \\ \hline
\endfoot

\hline 
\endlastfoot

$\mathfrak{NAC}_{1}^3$ & $\mathfrak{A}_{01}$  & $e_1e_2=e_3$\\ \hline

\end{longtable}
}

The classifications, up to isomorphism, of all $4$- and $5$-dimensional nilpotent anticommutative algebras can be found in Calder\'on Mart\'in, Fern\'andez Ouaridi and Kaygorodov (2019) and in~\cite{fkkv}, respectively; their geometric description was studied in~\cite{fkkv, kkl19}.

The irreducible component of the variety $\mathfrak{NACom}^4$ is defined by 

{\small

\begin{longtable}{|l|l|ll|}
\caption*{ }\\

\hline \multicolumn{1}{|c|}{${\mathcal A}$} & \multicolumn{1}{c|}{ } & \multicolumn{2}{c|}{Multiplication table} \\ \hline 
\endfirsthead

 \multicolumn{4}{l}%
{{\bfseries  continued from previous page}} \\
\hline \multicolumn{1}{|c|}{${\mathcal A}$} & \multicolumn{1}{c|}{ } & \multicolumn{2}{c|}{Multiplication table} \\ \hline 
\endhead

\hline \multicolumn{4}{|r|}{{Continued on next page}} \\ \hline
\endfoot

\hline 
\endlastfoot

$\mathfrak{NAC}_{1}^4$ & $\mathfrak{A}_{02}$ & $e_1e_2=e_3$ & $e_1e_3=e_4$  \\ \hline

\end{longtable}
}

The irreducible component of the variety $\mathfrak{NACom}^5$ is defined by 

{\small

\begin{longtable}{|l|l|lll|}
\caption*{ }   \\

\hline \multicolumn{1}{|c|}{${\mathcal A}$} & \multicolumn{1}{c|}{ } & \multicolumn{3}{c|}{Multiplication table} \\ \hline 
\endfirsthead

 \multicolumn{5}{l}%
{{\bfseries  continued from previous page}} \\
\hline \multicolumn{1}{|c|}{${\mathcal A}$} & \multicolumn{1}{c|}{ } & \multicolumn{3}{c|}{Multiplication table} \\ \hline 
\endhead

\hline \multicolumn{5}{|r|}{{Continued on next page}} \\ \hline
\endfoot

\hline 
\endlastfoot

$\mathfrak{NAC}_{1}^5$ & $\mathfrak{A}_{11}$ & $e_1e_2=e_3$ & $e_1e_3=e_4$ & $e_3e_4=e_5$ \\ \hline

\end{longtable}
}

The complete graph of degenerations can be found in~\cite{fkkv}.
Dimension $6$ was studied in~\cite{kkl19} both algebraically and geometrically.
This paper yields that the irreducible component of $\mathfrak{NACom}^6$ is defined by 

{\small

\begin{longtable}{|l|l|llllll|}
\caption*{ }   \\

\hline \multicolumn{1}{|c|}{${\mathcal A}$} & \multicolumn{1}{c|}{ } & \multicolumn{6}{c|}{Multiplication table} \\ \hline 
\endfirsthead

 \multicolumn{8}{l}%
{{\bfseries  continued from previous page}} \\
\hline \multicolumn{1}{|c|}{${\mathcal A}$} & \multicolumn{1}{c|}{ } & \multicolumn{6}{c|}{Multiplication table} \\ \hline 
\endhead

\hline \multicolumn{8}{|r|}{{Continued on next page}} \\ \hline
\endfoot

\hline 
\endlastfoot

$\mathfrak{NAC}_{1}^6(\alpha)$ & $\mathfrak{A}_{82}(\alpha)$  & $e_1e_2=e_3$ &  $e_1e_3=e_4$ &  $e_2e_5=\alpha e_6$ & $e_3e_4=e_5$ & $e_3e_5=e_6$ &  $e_4e_5=e_6$ \\ \hline

\end{longtable}
}

{ 
The present result was generalized in \cite{kkl22}.

\begin{theorem}[Theorem C, \cite{kkl22}]
For any $n\ge 2$, the variety of all $n$-dimensional anticommutative nilpotent algebras is irreducible and has dimension $\frac{(n-2)(n^2+2n+3)}{6}$.
\end{theorem}

Let $n\ge 6$ (in case $n=6$, the condition $c_{46}^7=1$ is to be ignored). Denote by   $\mathfrak{T}_n$ is the family of
$n$-dimensional complex anticommutative algebras whose structure constants $(c_{ij}^k)^n_{i,j,k}$ relative to the basis $(e_i)_{i=1}^{n}$ satisfy 
$c_{i,j}^{k}=0,\ \forall k \leq \max\{i,j\}$ 
and such that:
\begin{itemize}
\item $e_i e_{i+1}=e_{i+2}$ for all $1\le i\le n-2$;
\item $c_{1i}^{i+2}=0=c_{2i}^{i+2}$, for all $4\leq i\leq n-2$;
\item $c_{13}^{4}=c_{14}^{5}=c_{24}^5=c_{15}^{6}=c_{25}^6=c_{13}^6=0$;
\item $c_{13}^5 \neq 0$;
\item $c_{35}^6=c_{46}^7=1$.
\end{itemize}
The remaining structure constants $c_{ij}^{k}$ are arbitrary, subject only to the anticommutativity constraint.
As above, it was shown that the family $\mathfrak{T}_n$  is generic in the variety of $n$-dimensional anticommutative nilpotent algebras and inductively gives an algorithmic procedure to obtain any $n$-dimensional nilpotent anticommutative algebra as a degeneration from $\mathfrak{T}_n$ \cite{kkl22}.

}

\subsubsection{2-step nilpotent  algebras}
Among the nilpotent algebras, those satisfying $\mathfrak{N}^3=0$ have been studied more in detail. They are called $2$-step nilpotent algebras and will be denoted by $\mathfrak{2N}$.

Selecting the $2$-step nilpotent algebras from the $3$-dimensional nilpotent algebras listed in \cite{fkkv}, we have the following list:

{\small
\vspace{-5mm}

\begin{longtable}{|l|l|lll|}
\caption*{ }   \\

\hline \multicolumn{1}{|c|}{${\mathcal A}$} & \multicolumn{1}{c|}{ } & \multicolumn{3}{c|}{Multiplication table} \\ \hline 
\endfirsthead

 \multicolumn{5}{l}%
{{\bfseries  continued from previous page}} \\
\hline \multicolumn{1}{|c|}{${\mathcal A}$} & \multicolumn{1}{c|}{ } & \multicolumn{3}{c|}{Multiplication table} \\ \hline 
\endhead

\hline \multicolumn{5}{|r|}{{Continued on next page}} \\ \hline
\endfoot

\hline 
\endlastfoot

$\mathfrak{2N}_{1}^3$ & $\mathcal{N}_{6}$  & $e_1e_1=e_3$ &  $e_2e_2=e_3$ & \\ \hline
$\mathfrak{2N}_{2}^3$ & $\mathcal{N}_{5}$  & $e_1e_1=e_2$ &  &\\ \hline
$\mathfrak{2N}_{3}^3$ & $\mathcal{N}_{7}$  & $e_1e_2=e_3$ &  $e_2e_1=-e_3$ & \\ \hline
$\mathfrak{2N}_{4}^3(\alpha)$ & $\mathcal{N}_{8}(\alpha)$  & $e_1e_1=\alpha e_3$ & $e_2e_1=e_3$ & $e_2e_2=e_3$  \\ \hline
\end{longtable}
}

The geometric classification of the variety $\mathfrak{2N}^3$ follows from the graph of degenerations of $3$-dimensional nilpotent algebras of \cite{fkkv}. 

\[{\rm Irr}(\mathfrak{2N}^3)=
\left\{\overline{{\mathcal O}\left(\mathfrak{2N}_{1}^{3}\right)}\right\}\cup\left\{
\overline{\bigcup {\mathcal O}\left(\mathfrak{2N}_{4}^3(\alpha)\right)}\right\}.\]

The list of all $4$-dimensional $2$-step nilpotent algebras altogether appeared in~\cite{kppv}. 
In the same paper, it was proved that the variety $\mathfrak{2N}^4$ has two irreducible components:
\[{\rm Irr}(\mathfrak{2N}^4)=
\left\{\overline{\bigcup {\mathcal O}\big(\mathfrak{2N}_{i}^4(\alpha)\big)}\right\}_{i=1}^2,\]
where

{\small

\begin{longtable}{|l|l|lllll|}
\caption*{ }   \\

\hline \multicolumn{1}{|c|}{${\mathcal A}$} & \multicolumn{1}{c|}{ } & \multicolumn{5}{c|}{Multiplication table} \\ \hline 
\endfirsthead

 \multicolumn{7}{l}%
{{\bfseries  continued from previous page}} \\
\hline \multicolumn{1}{|c|}{${\mathcal A}$} & \multicolumn{1}{c|}{ } & \multicolumn{5}{c|}{Multiplication table} \\ \hline 
\endhead

\hline \multicolumn{7}{|r|}{{Continued on next page}} \\ \hline
\endfoot

\hline 
\endlastfoot

$\mathfrak{2N}_{1}^4(\alpha)$ & $\mathfrak{N}_2(\alpha)$  & $e_1e_1 = e_3$ & $e_1e_2 = e_4$ &  $e_2e_1 = -\alpha e_3$ & $e_2e_2 = -e_4$ &  \\ \hline
$\mathfrak{2N}_{2}^4(\alpha)$ & $\mathfrak{N}_3(\alpha)$  & $e_1e_1 = e_4$ & $e_1e_2 = \alpha e_4$ & $e_2e_1 = -\alpha e_4$ & $e_2e_2 = e_4$ & $e_3e_3 = e_4$  \\ \hline

\end{longtable}
}

{ 
The present result was generalized in \cite{ikp22}.
For $k \leq n$ consider the (algebraic) subset $\mathfrak{2N}_{n,k}$ of the variety $\mathfrak{2N}^n$ of $2$-step nilpotent $n$-dimensional algebras defined by 
\[
\mathfrak{2N}_{n,k} = \{A \in \mathfrak{2N}^n : {\rm dim} \ A^2 \ \leq k, \ {\rm dim} \ {\rm Ann} \ A \geq k \}.
\]
It is easy to see that $\mathfrak{2N}^n = \cup_{k=1}^n\mathfrak{2N}_{n,k}.$

\begin{theorem}[Theorem A, \cite{ikp22}]
The sets $\mathfrak{2N}_{n,k}$ are irreducible and 
\begin{equation*}
\mathfrak{2N}^n = \bigcup_{k}\mathfrak{2N}_{n,k}, \quad \text{for} \quad 1 \le k \le \left\lfloor n + \frac{1 - \sqrt{4n+1}}{2}\right\rfloor
\end{equation*}
is the decomposition of $\mathfrak{2N}^n$ into irreducible components.
Moreover,
\begin{longtable}{lcl}
$\dim \mathfrak{2N}_{n,k}$ &$=$ &$(n-k)^2k + (n-k)k.$\\
\end{longtable}
\end{theorem}

}
\subsubsection{2-step nilpotent commutative algebras} 
Let $\mathfrak{2NC}$ denote the variety of $2$-step nilpotent commutative algebras. Checking the list of \cite{fkkv}, we obtain that there are only two $3$-dimensional $2$-step nilpotent algebras:

{\small
\vspace{-4mm}

\begin{longtable}{|l|l|ll|}
\caption*{ }   \\

\hline \multicolumn{1}{|c|}{${\mathcal A}$} & \multicolumn{1}{c|}{ } & \multicolumn{2}{c|}{Multiplication table} \\ \hline 
\endfirsthead

 \multicolumn{4}{l}%
{{\bfseries  continued from previous page}} \\
\hline \multicolumn{1}{|c|}{${\mathcal A}$} & \multicolumn{1}{c|}{ } & \multicolumn{2}{c|}{Multiplication table} \\ \hline 
\endhead

\hline \multicolumn{4}{|r|}{{Continued on next page}} \\ \hline
\endfoot

\hline 
\endlastfoot

$\mathfrak{2NC}_{1}^3$ & $\mathcal{N}_{6}$  & $e_1e_1=e_3$ &  $e_2e_2=e_3$  \\ \hline
$\mathfrak{2NC}_{2}^3$ & $\mathcal{N}_{5}$  & $e_1e_1=e_2$ &   \\ \hline
\end{longtable}
}

Again from \cite{fkkv}, it follows that the variety $\mathfrak{2NC}^3$ is irreducible and it is defined by the rigid algebra $\mathfrak{2NC}_{1}^3$.

The following list of $2$-step nilpotent commutative algebras of dimension $4$ is based on the classification of $4$-dimensional nilpotent commutative algebras from \cite{fkkv}.

{\small
\begin{longtable}{|l|l|ll|}
\caption*{ }   \\

\hline \multicolumn{1}{|c|}{${\mathcal A}$} & \multicolumn{1}{c|}{ } & \multicolumn{2}{c|}{Multiplication table} \\ \hline 
\endfirsthead

 \multicolumn{4}{l}%
{{\bfseries  continued from previous page}} \\
\hline \multicolumn{1}{|c|}{${\mathcal A}$} & \multicolumn{1}{c|}{ } & \multicolumn{2}{c|}{Multiplication table} \\ \hline 
\endhead

\hline \multicolumn{4}{|r|}{{Continued on next page}} \\ \hline
\endfoot

\hline 
\endlastfoot

$\mathfrak{2NC}_{1}^4$ & $\mathcal{C}_{06}$  & $e_1e_1=e_3$ &  $e_2e_2=e_4$  \\ \hline
$\mathfrak{2NC}_{2}^4$ & $\mathcal{C}_{08}$  & $e_1e_1=e_4$ & $e_2e_3=e_4$  \\ \hline
$\mathfrak{2NC}_{3}^4$ & $\mathcal{C}_{01}$  & $e_1e_1=e_2$ &   \\ \hline
$\mathfrak{2NC}_{4}^4$ & $\mathcal{C}_{04}$  & $e_1e_2= e_3$ & \\ \hline
$\mathfrak{2NC}_{5}^4$ & $\mathcal{C}_{07}$  & $e_1e_1= e_3$ & $e_1e_2=e_4$  \\ \hline
\end{longtable}
}

Analyzing the graph of degenerations of $4$-dimensional nilpotent algebras from \cite{fkkv}, we obtain that $\mathfrak{2NC}^4$ has two irreducible components:

\[{\rm Irr}(\mathfrak{2NC}^4)=
\left\{\overline{{\mathcal O}\left(\mathfrak{2NC}_{i}^{4}\right)}\right\}_{i=1}^2.\]

Regarding dimension $5$, the variety $\mathfrak{2NC}^5$ was algebraically and geometrically classified in~\cite{klp} based on the classifications of $5$-dimensional nilpotent associative commutative algebras.

In particular, the variety  $\mathfrak{2NC}^5$ has three irreducible components:
\[{\rm Irr}(\mathfrak{2NC}^5)= \left\{  \overline{{\mathcal O}\left(\mathfrak{2NC}_{i}^5\right)} \right\}_{i=1}^{3},\]
where

{\small

\begin{longtable}{|l|l|lll|}
\caption*{ }   \\

\hline \multicolumn{1}{|c|}{${\mathcal A}$} & \multicolumn{1}{c|}{ } & \multicolumn{3}{c|}{Multiplication table} \\ \hline 
\endfirsthead

 \multicolumn{5}{l}%
{{\bfseries  continued from previous page}} \\
\hline \multicolumn{1}{|c|}{${\mathcal A}$} & \multicolumn{1}{c|}{ } & \multicolumn{3}{c|}{Multiplication table} \\ \hline 
\endhead

\hline \multicolumn{5}{|r|}{{Continued on next page}} \\ \hline
\endfoot

\hline 
\endlastfoot

$\mathfrak{2NC}_{1}^5$ & $\mathbf{A}_{07}$ & $e_{1}e_{2}=e_{4}+e_{5}$ & $e_{1}e_{3}=e_{4}$ & $e_{2}e_{3}=e_{5}$  \\ \hline
$\mathfrak{2NC}_{2}^5$ & $\mathbf{A}_{16}$ & $e_{1}e_1=e_{3}$ & $e_{1}e_{2}=e_{4}$ & $e_{2}e_2=e_{5}$ \\ \hline
$\mathfrak{2NC}_{3}^5$ & $\mathbf{A}_{21}$  & $e_1e_4 = e_5$ & $e_2e_3 =e_5$ & \\ \hline
\end{longtable}
}

{ 
The present result was generalized in \cite{ikp22}.
For $k \leq n$ consider the (algebraic) subset $\mathfrak{2NC}_{n,k}$ of the variety $\mathfrak{2NC}^n$ of $2$-step  nilpotent commutative $n$-dimensional algebras defined by 
\[
\mathfrak{2NC}_{n,k} = \{A \in \mathfrak{2NC}^n : {\rm dim} \ A^2 \ \leq k, \ {\rm dim} \ {\rm Ann} \ A \geq k \}.
\]
It is easy to see that $\mathfrak{2NC}^n = \cup_{k=1}^n\mathfrak{2NC}_{n,k}.$

\begin{theorem}[Theorem A, \cite{ikp22}]
The sets $\mathfrak{2NC}_{n,k}$ are irreducible and 
\begin{equation*}
\mathfrak{2NC}^n \ = \ \bigcup_{k}\mathfrak{2NC}_{n,k}, \mbox{  for } 1 \le k \le \left\lfloor n + \frac{3 - \sqrt{8n+9}}{2}\right\rfloor,\end{equation*}
is the decomposition of $\mathfrak{2NC}^n$ into irreducible components.
Moreover,
\begin{longtable}{lcl}
$\dim \mathfrak{2NC}_{n,k}$ &$=$ & $\frac{(n-k)(n-k+1)}{2}k + (n-k)k.$\\
\end{longtable}
\end{theorem}

}

\subsubsection{2-step nilpotent anticommutative algebras} 
The variety of $2$-step nilpotent anticommutative algebras will be denoted by $\mathfrak{2NA}$. The unique $2$-step nilpotent anticommutative algebra of dimension $3$ is
\vspace{-5mm}
{\small

\begin{longtable}{|l|l|l|}
\caption*{ }   \\

\hline \multicolumn{1}{|c|}{${\mathcal A}$} & \multicolumn{1}{c|}{ } & \multicolumn{1}{c|}{Multiplication table} \\ \hline 
\endfirsthead

 \multicolumn{3}{l}%
{{\bfseries  continued from previous page}} \\
\hline \multicolumn{1}{|c|}{${\mathcal A}$} & \multicolumn{1}{c|}{ } & \multicolumn{1}{c|}{Multiplication table} \\ \hline 
\endhead

\hline \multicolumn{3}{|r|}{{Continued on next page}} \\ \hline
\endfoot

\hline
\endlastfoot

$\mathfrak{2NA}_{1}^3$ & $\mathcal{N}_{7}$  & $e_1e_2=e_3$  \\ \hline
\end{longtable}
}\noindent and in dimension $4$ there is also only one: 
\vspace{-5mm}
{\small

\begin{longtable}{|l|l|l|}
\caption*{ }   \\

\hline \multicolumn{1}{|c|}{${\mathcal A}$} & \multicolumn{1}{c|}{ } & \multicolumn{1}{c|}{Multiplication table} \\ \hline 
\endfirsthead

 \multicolumn{3}{l}%
{{\bfseries  continued from previous page}} \\
\hline \multicolumn{1}{|c|}{${\mathcal A}$} & \multicolumn{1}{c|}{ } & \multicolumn{1}{c|}{Multiplication table} \\ \hline 
\endhead

\hline \multicolumn{3}{|r|}{{Continued on next page}} \\ \hline
\endfoot

\hline
\endlastfoot

$\mathfrak{2NA}_{1}^4$ & $\mathcal{N}_{7}$  & $e_1e_2=e_3$  \\ \hline
\end{longtable}
}\noindent The notation is taken from \cite{fkkv}.

As for dimension $5$, the list can be extracted from the general list of the $5$-dimensional nilpotent anticommutative algebras of \cite{fkkv}. 
{\small
\vspace{-4mm}
\begin{longtable}{|l|l|ll|}
\caption*{ }   \\

\hline \multicolumn{1}{|c|}{${\mathcal A}$} & \multicolumn{1}{c|}{ } & \multicolumn{2}{c|}{Multiplication table} \\ \hline 
\endfirsthead

 \multicolumn{4}{l}%
{{\bfseries  continued from previous page}} \\
\hline \multicolumn{1}{|c|}{${\mathcal A}$} & \multicolumn{1}{c|}{ } & \multicolumn{2}{c|}{Multiplication table} \\ \hline 
\endhead

\hline \multicolumn{4}{|r|}{{Continued on next page}} \\ \hline
\endfoot

\hline 
\endlastfoot

$\mathfrak{2NA}_{1}^5$ & $\mathfrak{A}_{03}$  & $e_1e_2=e_4$ &  $e_1e_3=e_5$ \\ \hline
$\mathfrak{2NA}_{2}^5$ & $\mathfrak{A}_{05}$  & $e_1e_2=e_5$ & $e_3e_4=e_5$  \\ \hline
$\mathfrak{2NA}_{3}^5$ & $\mathfrak{A}_{01}$  & $e_1e_2=e_3$ & \\ \hline
\end{longtable}
}

The irreducible components of $\mathfrak{2NA}^5$ are deduced, again, from \cite{fkkv}.

\[{\rm Irr}(\mathfrak{2NA}^5)=
\left\{\overline{{\mathcal O}\left(\mathfrak{2NA}_{i}^{5}\right)}\right\}_{i=1}^2.\]

To determine all the $6$-dimensional $2$-step nilpotent anticommutative algebras, we select them from the list of \cite{S90} of $6$-dimensional nilpotent Lie algebras. Note that every $2$-step nilpotent anticommutative algebra verifies Jacobi identity and is therefore a Lie algebra.

{\small
\vspace{-6mm}
\begin{longtable}{|l|l|lll|}
\caption*{ }   \\

\hline \multicolumn{1}{|c|}{${\mathcal A}$} & \multicolumn{1}{c|}{ } & \multicolumn{3}{c|}{Multiplication table} \\ \hline 
\endfirsthead

 \multicolumn{5}{l}%
{{\bfseries  continued from previous page}} \\
\hline \multicolumn{1}{|c|}{${\mathcal A}$} & \multicolumn{1}{c|}{ } & \multicolumn{3}{c|}{Multiplication table} \\ \hline 
\endhead

\hline \multicolumn{5}{|r|}{{Continued on next page}} \\ \hline
\endfoot

\hline 
\endlastfoot

$\mathfrak{2NA}_{1}^6$ & $g_3\times g_3$  & $e_1e_3=e_5$ &  $e_2e_4=e_6$ & \\ \hline
$\mathfrak{2NA}_{2}^6$ & $g_{6,24}$  & $e_1e_2=e_4$ & $e_1e_3=e_5$ & $e_2e_3=e_6$  \\ \hline
$\mathfrak{2NA}_{3}^6$ & $g_{6,21}$  & $e_1e_2=e_5$ & $e_1e_3=e_6$  & $e_3e_4=e_5$  \\ \hline
$\mathfrak{2NA}_{4}^6$ & $g_{5,2}\times \mathbb{C}$  & $e_1e_2=e_5$   & $e_3e_4=e_5$ & \\ \hline
$\mathfrak{2NA}_{5}^6$ & $g_{5,5}\times \mathbb{C}$  & $e_1e_2=e_4$ & $e_1e_3=e_5$  &  \\ \hline
$\mathfrak{2NA}_{6}^6$ & $g_{3}\times \mathbb{C}^3$  & $e_1e_2=e_3$  & &  \\ \hline
\end{longtable}
}

It follows from the general graph of degenerations of \cite{S90} that $\mathfrak{2NA}^6$ has two irreducible components:

\[{\rm Irr}(\mathfrak{2NA}^6)=
\left\{\overline{{\mathcal O}\left(\mathfrak{2NA}_{i}^{6}\right)}\right\}_{i=1}^2.\]

Dimensions $7$ and $8$ have been studied with the aim of contributing to the knowledge of the varieties of $7$- and $8$-dimensional nilpotent Lie algebras. The irreducible components of $\mathfrak{2NA}^7$ were determined in \cite{ale2}, employing the algebraic classification of all $7$-dimensional nilpotent Lie algebras by Gong (1998). Note that the rigid algebras had already been identified in~\cite{ale}.

The variety $\mathfrak{2NA}^7$  has three irreducible components:
\[{\rm Irr}\left(\mathfrak{2NA}^7\right)=\left\{\overline{{\mathcal O}\big(\mathfrak{2NA}_{i}^7\big)}\right\}_{i=1}^{3},\]
where

{\small

\begin{longtable}{|l|l|llll|}
\caption*{ }   \\

\hline \multicolumn{1}{|c|}{${\mathcal A}$} & \multicolumn{1}{c|}{ } & \multicolumn{4}{c|}{Multiplication table} \\ \hline 
\endfirsthead

 \multicolumn{6}{l}%
{{\bfseries  continued from previous page}} \\
\hline \multicolumn{1}{|c|}{${\mathcal A}$} & \multicolumn{1}{c|}{ } & \multicolumn{4}{c|}{Multiplication table} \\ \hline 
\endhead

\hline \multicolumn{6}{|r|}{{Continued on next page}} \\ \hline
\endfoot

\hline 
\endlastfoot

$\mathfrak{2NA}_{1}^7$ & $(17)$ & $e_1e_2=e_7$ & $e_3e_4=e_7$ & $e_5e_6=e_7$ & \\ \hline
$\mathfrak{2NA}_{2}^7$ & $(27B)$ & $e_1e_2=e_6$ & $e_1e_5=e_7$ & $e_2e_3= e_7$ & $e_3e_4=e_6$  \\ \hline
$\mathfrak{2NA}_{3}^7$ & $(37D)$ & $e_1e_2=e_5$ & $e_1e_3=e_6$ & $e_2e_4= e_7$ & $e_3e_4=e_5$  \\ \hline
\end{longtable}
}

The complete graph of degenerations can be found in~\cite{ale2}.

The algebraic classification of the $2$-step nilpotent Lie algebras of dimension $8$ was made by Yan and Deng in 2013. Their graph of degenerations was constructed in~\cite{ale3}, although~\cite{ale} had already shown that the variety $\mathfrak{2NA}^8$ has three irreducible components:
\[{\rm Irr}\left(\mathfrak{2NA}^8\right)=
\left\{\overline{{\mathcal O}\big(\mathfrak{2NA}_{i}^8\big)}\right\}_{i=1}^{3},\]
where

{\small

\begin{longtable}{|l|l|llllll|}
\caption*{ }   \\

\hline \multicolumn{1}{|c|}{${\mathcal A}$} & \multicolumn{1}{c|}{ } & \multicolumn{6}{c|}{Multiplication table} \\ \hline 
\endfirsthead

 \multicolumn{8}{l}%
{{\bfseries  continued from previous page}} \\
\hline \multicolumn{1}{|c|}{${\mathcal A}$} & \multicolumn{1}{c|}{ } & \multicolumn{6}{c|}{Multiplication table} \\ \hline 
\endhead

\hline \multicolumn{8}{|r|}{{Continued on next page}} \\ \hline
\endfoot

\hline 
\endlastfoot

$\mathfrak{2NA}_{1}^8$ & $N_1^{8,2}$ & $e_1e_2=e_7$ &  $e_3e_4=e_8$ & $e_5e_6=e_7+e_8$ &&& \\ \hline
$\mathfrak{2NA}_{2}^8$ & $N_1^{8,4}$ & $e_1e_2=e_5$ & $e_2e_3=e_6$ & $e_3e_4=e_7$ & $e_4e_1=e_8$ && \\ \hline
$\mathfrak{2NA}_{3}^8$ & $N_9^{8,3}$ & $e_1e_2=e_6$ & $e_1e_3=e_7$ & $e_1e_4=e_8$ & $e_2e_3=e_8$  & $e_2e_5=e_7$ & $e_4e_5=e_6$ \\ \hline
\end{longtable}
}

{ 
The present result was generalized in \cite{ikp22}.
For $k \leq n$ consider the (algebraic) subset $\mathfrak{2NA}_{n,k}$ of the variety $\mathfrak{2NA}^n$ of $2$-step  nilpotent commutative $n$-dimensional algebras defined by 
\[
\mathfrak{2NA}_{n,k} = \{A \in \mathfrak{2NA}^n : {\rm dim} \ A^2 \ \leq k, \ {\rm dim} \ {\rm Ann} \ A \geq k \}.
\]
It is easy to see that $\mathfrak{2NA}^n = \cup_{k=1}^n\mathfrak{2NA}_{n,k}.$

\begin{theorem}[Theorem A, \cite{ikp22}]
The sets $\mathfrak{2NA}_{n,k}$ are irreducible and 
\begin{equation*}
\mathfrak{2NA}^n \ = \ \bigcup_{k}\mathfrak{2NA}_{n,k}, \mbox{  for } 1 + (n + 1) \operatorname{mod} 2 \le k \le \left\lfloor n + \frac{1 - \sqrt{8n+1}}{2}\right\rfloor \text{ for } n \ge 3,\end{equation*}
is the decomposition of $\mathfrak{2NA}^n$ into irreducible components.
Moreover,
\begin{longtable}{lcl}
$\dim \mathfrak{2NA}_{n,k}$ &$=$ & $\frac{(n-k)(n-k-1)}{2}k + (n-k)k.$\\
\end{longtable}
\end{theorem}

}

\subsubsection{Noncommutative Heisenberg algebras} 
An algebra $\mathfrak{A}$ is a noncommutative Heisenberg algebra if ${\rm dim} \ \mathfrak{A}^2 \leq  1$ and $\mathfrak{A}^2\mathfrak{A}+\mathfrak{A}\mathfrak{A}^2= 0.$ Let $\mathfrak{NH}$ denote the variety of noncommutative Heisenberg algebras. 
$\mathfrak{NH}$ is a special subvariety in  the variety of $2$-step nilpotent algebras. 
It is easy to see that $\mathfrak{NH}^n$  is irreducible.
The full graph of degenerations of the variety $\mathfrak{NH}^5$  was obtained in \cite{kv23}.
The variety $\mathfrak{NH}^5$  is determined by following family of algebras

{\small

\begin{longtable}{|l|l|llll|}
\caption*{ }   \\

\hline \multicolumn{1}{|c|}{${\mathcal A}$} & \multicolumn{1}{c|}{ } & \multicolumn{4}{c|}{Multiplication table} \\ \hline 
\endfirsthead

 \multicolumn{6}{l}%
{{\bfseries  continued from previous page}} \\
\hline \multicolumn{1}{|c|}{${\mathcal A}$} & \multicolumn{1}{c|}{ } & \multicolumn{4}{c|}{Multiplication table} \\ \hline 
\endhead

\hline \multicolumn{6}{|r|}{{Continued on next page}} \\ \hline
\endfoot

\hline 
\endlastfoot

$\mathfrak{NH}_{1}^5(\alpha,\beta)$ & 
$\mathfrak{H}_{13}^{\alpha,\beta}$ & 
$e_1e_2=e_5$ & $e_2e_1=\beta e_5$ & $e_3e_4=e_5$ & $e_4e_3=\alpha e_5$\\  \hline

\end{longtable}
}

 \subsection{Nilalgebras} 
An element $x$ is nil with nilindex $n$, if for each $k\geq n$ we have $x^k=0$\footnote{By $x^k$ we mean all possible arrangements of non-associative products.}.
An algebra is called a nilalgebra with nilindex $n$ if each element is nil
and $n$ is the maximal nilindex of elements from the algebra. 
The variety of nilalgebras with nilindex $n$ 
will be denoted by $\mathfrak{Nil(n)}$.

\subsubsection{3-dimensional nilalgebras with nilindex 3} 
The algebraic and geometric classification of $3$-dimensional nilalgebras with nilindex $3$ can be found in \cite{ks25}. 
In particular, it is proven that the variety  $\mathfrak{Nil(3)}^3$ has two irreducible components:
\[{\rm Irr}(\mathfrak{Nil(3)}^3)=\left\{
\overline{{\mathcal O}\big(\mathfrak{Nil(3)}_{1}^3\big)}\right\} \cup 
\left\{
\overline{\bigcup {\mathcal O}\big(\mathfrak{Nil(3)}_{2}^3(\alpha)\big)}\right\} ,\]
where algebras are defined as follows:

{\small
\vspace{-4mm}
\begin{longtable}{|l|l|llll|}
\caption*{ }   \\

\hline \multicolumn{1}{|c|}{${\mathcal A}$} & \multicolumn{1}{c|}{ } & \multicolumn{4}{c|}{Multiplication table} \\ \hline 
\endfirsthead

 \multicolumn{6}{l}%
{{\bfseries  continued from previous page}} \\
\hline \multicolumn{1}{|c|}{${\mathcal A}$} & \multicolumn{1}{c|}{ } & \multicolumn{4}{c|}{Multiplication table} \\ \hline 
\endhead

\hline \multicolumn{6}{|r|}{{Continued on next page}} \\ \hline
\endfoot

\hline 
\endlastfoot

$\mathfrak{Nil(3)}_{1}^3$ & ${\mathcal N}_{5}$ &$ e_1 e_1 = e_2$  & $e_1 e_3=e_3$ & $e_3 e_1=-e_3$ & $e_3 e_3=e_2$ \\ \hline
$\mathfrak{Nil(3)}_{2}^3(\alpha)$ & $\mathcal{A}_1^{\alpha}$ &  
$e_1e_2=e_3$ & $e_1e_3 =e_1+e_3$ & $e_2e_3=\alpha e_2$ &\\&&
$e_2e_1=-e_3$ & $e_3e_1 =-e_1-e_3$ & $e_3e_2=-\alpha e_2$ & \\ \hline 
\end{longtable}
}

\subsubsection{3-dimensional nilalgebras with nilindex 4} 
The algebraic and geometric classification of $3$-dimensional nilalgebras with nilindex $4$ can be found in \cite{ks25}. 
In particular, it is proven that the variety  $\mathfrak{Nil(4)}^3$ has three irreducible components:
\[{\rm Irr}(\mathfrak{Nil(4)}^3)=\left\{
\overline{{\mathcal O}\big(\mathfrak{Nil(4)}_{1}^3\big)}\right\} \cup 
\left\{
\overline{\bigcup {\mathcal O}\big(\mathfrak{Nil(4)}_{i}^3(\alpha)\big)}\right\}_{i=2}^3 ,\]
where algebras are defined as follows:

\vspace{-5mm}

{\small

\begin{longtable}{|l|l|llll|}
\caption*{ }   \\

\hline \multicolumn{1}{|c|}{${\mathcal A}$} & \multicolumn{1}{c|}{ } & \multicolumn{4}{c|}{Multiplication table} \\ \hline 
\endfirsthead

 \multicolumn{6}{l}%
{{\bfseries  continued from previous page}} \\
\hline \multicolumn{1}{|c|}{${\mathcal A}$} & \multicolumn{1}{c|}{ } & \multicolumn{4}{c|}{Multiplication table} \\ \hline 
\endhead

\hline \multicolumn{6}{|r|}{{Continued on next page}} \\ \hline
\endfoot

\hline 
\endlastfoot
$\mathfrak{Nil(4)}_{1}^3$ & ${\mathcal N}_{5}$ &$ e_1 e_1 = e_2$  & $e_1 e_3=e_3$ & $e_3 e_1=-e_3$ & $e_3 e_3=e_2$ \\ \hline
$\mathfrak{Nil(4)}_{2}^3(\alpha)$ & $\mathcal{A}_1^{\alpha}$ &  
$e_1e_2=e_3$ & $e_1e_3 =e_1+e_3$ & $e_2e_3=\alpha e_2$ &\\&&
$e_2e_1=-e_3$ & $e_3e_1 =-e_1-e_3$ & $e_3e_2=-\alpha e_2$ & \\ \hline 

$\mathfrak{Nil(4)}_{2}^3(\alpha)$ & $\rm{N}_2^\alpha$&  $e_1 e_1 = e_2$ & $e_1 e_2=e_3$ & $e_2 e_1=\alpha e_3$  & \\

\end{longtable}
}

\subsubsection{3-dimensional nilalgebras with nilindex 5} 
The algebraic and geometric classification of $3$-dimensional nilalgebras with nilindex $5$ can be found in \cite{ks25}. 
In particular, it is proven that the variety  $\mathfrak{Nil(5)}^3$ has three irreducible components:
\[{\rm Irr}(\mathfrak{Nil(5)}^3)=\left\{
\overline{{\mathcal O}\big(\mathfrak{Nil(5)}_{i}^3\big)}\right\}_{i=1}^2 \cup 
\left\{
\overline{\bigcup {\mathcal O}\big(\mathfrak{Nil(5)}_{3}^3(\alpha)\big)}\right\},\]
where algebras are defined as follows:

\vspace{-5mm}
{\small

\begin{longtable}{|l|l|llll|}
\caption*{ }   \\

\hline \multicolumn{1}{|c|}{${\mathcal A}$} & \multicolumn{1}{c|}{ } & \multicolumn{4}{c|}{Multiplication table} \\ \hline 
\endfirsthead

 \multicolumn{6}{l}%
{{\bfseries  continued from previous page}} \\
\hline \multicolumn{1}{|c|}{${\mathcal A}$} & \multicolumn{1}{c|}{ } & \multicolumn{4}{c|}{Multiplication table} \\ \hline 
\endhead

\hline \multicolumn{6}{|r|}{{Continued on next page}} \\ \hline
\endfoot

\hline 
\endlastfoot
$\mathfrak{Nil(5)}_{1}^3$ & ${\mathcal N}_{5}$ &$ e_1 e_1 = e_2$  & $e_1 e_3=e_3$ & $e_3 e_1=-e_3$ & $e_3 e_3=e_2$ \\ \hline

$\mathfrak{Nil(5)}_{1}^3$ & $\bf{N}_2$ & $e_1 e_1 = e_2$ & $e_2 e_1= e_3$ & $e_2 e_2=e_3$ & \\ \hline
$\mathfrak{Nil(5)}_{2}^3(\alpha)$ & $\mathcal{A}_1^{\alpha}$ &  
$e_1e_2=e_3$ & $e_1e_3 =e_1+e_3$ & $e_2e_3=\alpha e_2$ &\\&&
$e_2e_1=-e_3$ & $e_3e_1 =-e_1-e_3$ & $e_3e_2=-\alpha e_2$ & \\ \hline

\end{longtable}
}

 \subsection{Associative algebras}

An  algebra $\mathfrak{A}$ is called  {\it associative} if it satisfies the identity 
\[(xy)z=x(yz).\] The variety of associative algebras will be denoted by $\mathfrak{Ass}$.

\subsubsection{2-dimensional associative algebras}
The variety of $2$-dimensional associative  algebras has three irreducible components:

\[{\rm Irr}(\mathfrak{Ass}^2)=\left\{\overline{{\mathcal O}\left(\mathfrak{A}_{i}^2\right)} \right\}_{i=1}^{3},\]
where  

{\small

\begin{longtable}{|l|l|ll|}
\caption*{ }   \\

\hline \multicolumn{1}{|c|}{${\mathcal A}$} & \multicolumn{1}{c|}{ } & \multicolumn{2}{c|}{Multiplication table} \\ \hline 
\endfirsthead

 \multicolumn{4}{l}%
{{\bfseries  continued from previous page}} \\
\hline \multicolumn{1}{|c|}{${\mathcal A}$} & \multicolumn{1}{c|}{ } & \multicolumn{2}{c|}{Multiplication table} \\ \hline 
\endhead

\hline \multicolumn{4}{|r|}{{Continued on next page}} \\ \hline
\endfoot

\hline 
\endlastfoot

$\mathfrak{A}_{1}^2$ & & $e_1e_1 = e_1$ & $e_2e_2 = e_2$  \\ \hline
$\mathfrak{A}_{2}^2$ & & $e_1e_1 = e_1$ & $e_1e_2 =e_2$ \\ \hline
$\mathfrak{A}_{3}^2$ & & $e_1e_1 = e_1$ & $e_2e_1 = e_2$ \\ \hline

\end{longtable}
}

\subsubsection{3-dimensional nilpotent associative algebras}
The algebraic and geometric classification of $3$-dimensional nilpotent associative algebras can be obtained from the classification and description of degenerations of $3$-dimensional nilpotent algebras given in \cite{fkkv}.
Hence, we have  that the variety   $\mathfrak{NAss}^3$ has two irreducible components:
\[{\rm Irr}(\mathfrak{NAss}^3)=
\left\{\overline{{\mathcal O}\left(\mathfrak{NA}_{1}^{3}\right)}\right\}\cup\left\{
\overline{\bigcup {\mathcal O}\left(\mathfrak{NA}_{2}^3(\alpha)\right)}\right\},\]
where

{\small

\begin{longtable}{|l|l|llllll|}
\caption*{ }   \\

\hline \multicolumn{1}{|c|}{${\mathcal A}$} & \multicolumn{1}{c|}{ } & \multicolumn{6}{c|}{Multiplication table} \\ \hline 
\endfirsthead

 \multicolumn{8}{l}%
{{\bfseries  continued from previous page}} \\
\hline \multicolumn{1}{|c|}{${\mathcal A}$} & \multicolumn{1}{c|}{ } & \multicolumn{6}{c|}{Multiplication table} \\ \hline 
\endhead

\hline \multicolumn{8}{|r|}{{Continued on next page}} \\ \hline
\endfoot

\hline 
\endlastfoot
$\mathfrak{NA}_{1}^3$ & $\mathcal{N}_{4}(1)$  & $e_1e_1=e_3$ &  $e_1e_2=e_3$ & $e_2e_1=e_3$ &&&\\ \hline

$\mathfrak{NA}_{2}^3(\alpha)$ &  $\mathcal{N}_{8}(\alpha)$  & $e_1e_1=\alpha e_3$ & $e_2e_1=e_3$ & $e_2e_2=e_3$ && &\\ \hline
\end{longtable}
}

\subsubsection{4-dimensional nilpotent associative algebras}
  The algebraic classification of $4$-dimensional nilpotent associative algebras can be found in a paper by Karimjanov (2021) and the geometric classification was given in \cite{ikp22}. 
 
In particular, it is proven that the variety  $\mathfrak{NAss}^4$ has four irreducible components:
 
\[{\rm Irr}(\mathfrak{NAss}^4)=
\left\{\overline{{\mathcal O}\big(\mathfrak{NA}_{i}^4\big)}\right\}_{i=1}^2 
\cup
\left\{\overline{\bigcup {\mathcal O}\big(\mathfrak{NA}_{i}^4(\alpha)\big)}\right\}_{i=3}^4,\]
where

{\small

\begin{longtable}{|l|l|llllll|}
\caption*{ }   \\

\hline \multicolumn{1}{|c|}{${\mathcal A}$} & \multicolumn{1}{c|}{ } & \multicolumn{6}{c|}{Multiplication table} \\ \hline 
\endfirsthead

 \multicolumn{8}{l}%
{{\bfseries  continued from previous page}} \\
\hline \multicolumn{1}{|c|}{${\mathcal A}$} & \multicolumn{1}{c|}{ } & \multicolumn{6}{c|}{Multiplication table} \\ \hline 
\endhead

\hline \multicolumn{8}{|r|}{{Continued on next page}} \\ \hline
\endfoot

\hline 
\endlastfoot
$\mathfrak{NA}_{1}^4$ & $\mu_{0}^4$ & 
$ e_1e_1=e_2$ &$ e_1e_2=e_3$ &$ e_1e_3=e_4$ &  $ e_2e_1=e_3$ &$ e_2e_2=e_4$&$ e_3e_1=e_4$ \\ 
\hline
$\mathfrak{NA}_{2}^4$ & ${\mathcal A}^4_{05}$  &  
$e_1 e_1 = e_2$ &  $ e_1e_2=e_4$ & $e_1e_3= e_4$ & $e_2e_1= e_4$ & $e_3e_3=e_4$&\\ \hline
$\mathfrak{NA}_{3}^4(\alpha)$ & $\mathfrak{N}_2(\alpha)$ & $e_1e_1 = e_3$ & $e_1e_2 = e_4$ &  $e_2e_1 = -\alpha e_3$ & $e_2e_2 = -e_4$ & & \\ \hline
$\mathfrak{NA}_{4}^4(\alpha)$ & $\mathfrak{N}_3(\alpha)$ & $e_1e_1 = e_4$ & $e_1e_2 = \alpha e_4$ & $e_2e_1 = -\alpha e_4$ & $e_2e_2 = e_4$ & $e_3e_3 = e_4$ & \\ \hline
\end{longtable}
}

\subsubsection{5-dimensional nilpotent associative algebras }
 The algebraic classification of $5$-dimensional nilpotent associative algebras can be found in a paper by Karimjanov (2021) and the geometric classification was given in \cite{ikp22}. 
 
In particular, it is proven that the variety  $\mathfrak{NAss}^5$ has eleven irreducible components:

\begin{align*}
&{\rm Irr}(\mathfrak{NAss}^5) \\
&{\ }=
\left\{\overline{{\mathcal O}\big(\mathfrak{NA}_{i}^5\big)}\right\}_{i=1}^{7} 
\cup
\left\{\overline{\bigcup {\mathcal O}\big(\mathfrak{NA}_{i}^5(\alpha)\big)}\right\}_{i=8}^{9}
\cup
\left\{\overline{\bigcup {\mathcal O}\big(\mathfrak{NA}_{10}^5(\alpha,\beta)\big)}\right\} 
\cup
\left\{\overline{\bigcup {\mathcal O}\big(\mathfrak{NA}_{11}^5(\overline{\mu})\big)}\right\},
\end{align*}
where

{\small
\setlength{\tabcolsep}{5pt} 
\begin{longtable}{|l|l|lllll|}
\caption*{ }   \\

\hline \multicolumn{1}{|c|}{${\mathcal A}$} & \multicolumn{1}{c|}{ } & \multicolumn{5}{c|}{Multiplication table} \\ \hline 
\endfirsthead

 \multicolumn{7}{l}%
{{\bfseries  continued from previous page}} \\
\hline \multicolumn{1}{|c|}{${\mathcal A}$} & \multicolumn{1}{c|}{ } & \multicolumn{5}{c|}{Multiplication table} \\ \hline 
\endhead

\hline \multicolumn{7}{|r|}{{Continued on next page}} \\ \hline
\endfoot

\hline 
\endlastfoot
$\mathfrak{NA}_{1}^5$ &  
$\mu_{0}^5$ & 
$ e_1e_1=e_2$ &$ e_1e_2=e_3$ &$ e_1e_3=e_4$ & $e_1e_4=e_{5}$ &\\
&& $ e_2e_1=e_3$ &$ e_2e_2=e_4$& $ e_2e_3=e_5$&&\\
&&$ e_3e_1=e_4$  & $ e_3e_2=e_5$  & $e_4e_1=e_{5}$ && \\

 \hline

$\mathfrak{NA}_{2}^5$ & 
$\mu_{1,4}^5$ & 
$ e_1e_1=e_2$ &$ e_1e_2=e_3$ &$ e_1e_3=e_4$ & $e_1e_5=e_{4}$ &\\&    
&$ e_2e_1=e_3$ &$ e_2e_2=e_4$&$ e_3e_1=e_4$  & $e_5e_5=e_{4}$ & \\

 \hline

$\mathfrak{NA}_{3}^5$ &$\mu_ {11}$& $e_1e_1=e_2$&$e_1e_2=e_3$&$e_1e_4=e_5$&&\\ &&
 $e_2e_1=e_3$&$e_4e_4=e_3+e_5$ &&& \\
 \hline

$\mathfrak{NA}_{4}^5$ &$\mu_ {15}$& $e_1e_1=e_2$&$e_1e_2=e_3$&$e_1e_4=e_5$&$e_2e_1=e_3$&\\
&&$e_4e_1=e_2+e_5$&$e_4e_2=2e_3$&$e_4e_4=e_3+2e_5$&$e_5e_1=e_3$&\\ \hline

$\mathfrak{NA}_{5}^5$ & $\mu_ {17}$& $e_1e_1=e_2$&$e_1e_2=e_3$&$e_1e_4=e_5$&$e_2e_1=e_3$&\\
&&$e_4e_1=e_3+e_5$&$e_4e_4=e_2$&$e_4e_5=e_3$&$e_5e_4=e_3$&\\
\hline 

$\mathfrak{NA}_{6}^5$ &$\mu_ {18}$ &$e_1e_1=e_2$&$e_1e_2=e_3$&$ e_1e_4=e_5$&$e_2e_1=e_3$&\\
&&$e_4e_1=-e_5$&$e_4e_4=e_2$&$e_4e_5=-e_3$&$e_5e_4=e_3$ & \\

 \hline
$\mathfrak{NA}_{7}^5$ &$\mu_ {20}$& $e_1e_1=e_2$&$e_1e_2=e_3$&$e_1e_4=e_5$&$e_2e_1=e_3$&\\
&&$ e_4e_1=e_3+e_5$&$e_4e_4=-e_2+2e_5$&$e_4e_5=e_3$&$e_5e_4=-e_3$&\\
 \hline

$\mathfrak{NA}_{8}^5(\alpha)$ &$\lambda_6^{\alpha}$  & 
$e_1e_1=e_2$& $e_1e_2=e_3$ & $e_2e_1=e_3$ &&\\
&&$e_4e_5=e_3$ & $e_5e_4=\alpha e_3$&&& \\ \hline

$\mathfrak{NA}_{9}^5(\alpha)$&$\mu_ {22}^\alpha$& 
$ e_1e_1=e_2$&$  e_1e_2=e_3$&$e_1e_4=e_5$&&\\&&
$e_2e_1=e_3$& \multicolumn{2}{l}{$e_4e_1=(1-\alpha)e_2+\alpha e_5$}&&\\
&&$e_4e_2=(1-\alpha^2)e_3$&
\multicolumn{2}{l}{$e_4e_4=-\alpha e_2+(1+\alpha)e_5$}&&\\&&
$e_4e_5=-\alpha^2e_3$&
$e_5e_1=(1-\alpha)e_3$&
$e_5e_4=-\alpha e_3$&&\\
 \hline

$\mathfrak{NA}_{10}^5(\alpha, \beta)$ &${\mathfrak V}_{4+1}$ &    
$e_1e_2=e_5$& $e_2e_1=\alpha e_5$ &$e_3e_4=e_5$&$e_4e_3=\beta e_5$&\\

 \hline

$\mathfrak{NA}_{11}^5(\overline{\mu})$ &${\mathfrak V}_{3+2}$ & 
$e_1e_1 =  e_4$& $e_1e_2 = \mu_1 e_5$ & $e_1e_3 =\mu_2 e_5$&&\\
&&
$e_2e_1 = \mu_3 e_5$  & $e_2e_2 = \mu_4 e_5$   & $e_2e_3 = \mu_5 e_5$  &&\\
&&$e_3e_1 = \mu_6 e_5$  &$e_3e_2 = \mu_0 e_4+ \mu_7 e_5$  & $e_3e_3 =  e_5$  &&\\
 \hline
\end{longtable}
}

\subsection{Commutative associative  algebras} 

We consider now the associative algebras $\mathfrak{CA}$ which are also commutative. This variety will be denoted by $\mathfrak{CAss}$.
To study the varieties $\mathfrak{CAss}^n$, $n=2,3,4$, we rely on results for Jordan algebras, selecting the associative ones among them. 
The varieties $\mathfrak{CAss}^n$, $n=2,3,4$, are irreducible.

\subsubsection{2-dimensional commutative associative algebras}\label{2dasc}
For the variety $\mathfrak{CAss}^2$, we rely on the classification from~\cite{gkp}. It is determined by the rigid algebra:
{\small
\begin{longtable}{|l|l|ll|}
\caption*{ }   \\

\hline \multicolumn{1}{|c|}{${\mathcal A}$} & \multicolumn{1}{c|}{ } & \multicolumn{2}{c|}{Multiplication table} \\ \hline 
\endfirsthead

 \multicolumn{4}{l}%
{{\bfseries  continued from previous page}} \\
\hline \multicolumn{1}{|c|}{${\mathcal A}$} & \multicolumn{1}{c|}{ } & \multicolumn{2}{c|}{Multiplication table} \\ \hline 
\endhead

\hline \multicolumn{4}{|r|}{{Continued on next page}} \\ \hline
\endfoot

\hline 
\endlastfoot

$\mathfrak{CA}_{1}^2$ & $\mathfrak{B}_4$ & $e_1e_1=e_1$ & $e_2e_2=e_2$  \\ \hline

\end{longtable}
}

\subsubsection{3-dimensional nilpotent commutative associative algebras}
The algebraic and geometric classification of $3$-dimensional nilpotent commutative associative algebras can be obtained from the classification and description of degenerations of $3$-dimensional nilpotent algebras given in \cite{fkkv}.
Hence, we have  that the variety   $\mathfrak{NCAss}^3$ is irreducible, defined by the rigid algebra: 

{\small

\begin{longtable}{|l|l|llllll|}
\caption*{ }   \\

\hline \multicolumn{1}{|c|}{${\mathcal A}$} & \multicolumn{1}{c|}{ } & \multicolumn{6}{c|}{Multiplication table} \\ \hline 
\endfirsthead

 \multicolumn{8}{l}%
{{\bfseries  continued from previous page}} \\
\hline \multicolumn{1}{|c|}{${\mathcal A}$} & \multicolumn{1}{c|}{ } & \multicolumn{6}{c|}{Multiplication table} \\ \hline 
\endhead

\hline \multicolumn{8}{|r|}{{Continued on next page}} \\ \hline
\endfoot

\hline 
\endlastfoot
$\mathfrak{NCA}_{1}^3$ & $\mathfrak{N}_{4}(1)$  & $e_1e_1=e_3$ &  $e_1e_2=e_3$ &   &&&\\ \hline
 
\end{longtable}
}

\subsubsection{3-dimensional commutative associative algebras}
Again, we employ the classification of~\cite{gkp} to see that the variety $\mathfrak{CAss}^3$ is determined by the rigid algebra:

{\small

\begin{longtable}{|l|l|lll|}
\caption*{ }   \\

\hline \multicolumn{1}{|c|}{${\mathcal A}$} & \multicolumn{1}{c|}{ } & \multicolumn{3}{c|}{Multiplication table} \\ \hline 
\endfirsthead

 \multicolumn{5}{l}%
{{\bfseries  continued from previous page}} \\
\hline \multicolumn{1}{|c|}{${\mathcal A}$} & \multicolumn{1}{c|}{ } & \multicolumn{3}{c|}{Multiplication table} \\ \hline 
\endhead

\hline \multicolumn{5}{|r|}{{Continued on next page}} \\ \hline
\endfoot

\hline 
\endlastfoot

$\mathfrak{CA}_{1}^3$ & $\mathbb{T}_{01}$ & $e_1e_1=e_1$ & $e_2e_2=e_2$ & $e_3e_3=e_3$  \\ \hline

\end{longtable}
}

\subsubsection{4-dimensional nilpotent commutative associative algebras}
Also, from~\cite{esp} and~\cite{fkkv} (where some of the results of~\cite{esp} were corrected) we can extract that the variety $\mathfrak{NCAss}^4$ is irreducible, defined by the rigid algebra: 


{\small

\begin{longtable}{|l|l|llll|}
\caption*{ }   \\

\hline \multicolumn{1}{|c|}{${\mathcal A}$} & \multicolumn{1}{c|}{ } & \multicolumn{4}{c|}{Multiplication table} \\ \hline 
\endfirsthead

 \multicolumn{6}{l}%
{{\bfseries  continued from previous page}} \\
\hline \multicolumn{1}{|c|}{${\mathcal A}$} & \multicolumn{1}{c|}{ } & \multicolumn{4}{c|}{Multiplication table} \\ \hline 
\endhead

\hline \multicolumn{6}{|r|}{{Continued on next page}} \\ \hline
\endfoot

\hline 
\endlastfoot

$\mathfrak{NCA}_{1}^4$ & $\phi_1$ & $e_1e_1  = e_2$ & $e_1 e_2 = e_3$ & $e_1 e_3 = e_4$ & $e_2 e_2 = e_4$ \\ \hline

\end{longtable}
}

Note that we are employing the notation of~\cite{esp}.

\subsubsection{4-dimensional commutative associative algebras}
The degenerations of Jordan algebras of dimension $4$ were studied in~\cite{KM}, but the authors did not present a complete graph of degenerations. However, they did prove that every associative Jordan algebra degenerates from $\mathfrak{CA}_{1}^4$. We deduce that 
the variety $\mathfrak{CAss}^4$ is determined by the rigid algebra:

{\small

\begin{longtable}{|l|l|llll|}
\caption*{ }   \\

\hline \multicolumn{1}{|c|}{${\mathcal A}$} & \multicolumn{1}{c|}{ } & \multicolumn{4}{c|}{Multiplication table} \\ \hline 
\endfirsthead

 \multicolumn{6}{l}%
{{\bfseries  continued from previous page}} \\
\hline \multicolumn{1}{|c|}{${\mathcal A}$} & \multicolumn{1}{c|}{ } & \multicolumn{4}{c|}{Multiplication table} \\ \hline 
\endhead

\hline \multicolumn{6}{|r|}{{Continued on next page}} \\ \hline
\endfoot

\hline 
\endlastfoot

$\mathfrak{CA}_{1}^4$ & $\mathfrak{J}_3$ & $e_1e_1=e_1$ & $e_2e_2=e_2$ & $e_3e_3=e_3$ & $e_4e_4=e_4$  \\ \hline

\end{longtable}
}

\subsubsection{5-dimensional nilpotent commutative associative algebras}

In dimension $5$, we will first focus on nilpotent algebras, whose algebraic classification was given in~\cite{maz80}, and whose degenerations were established in~\cite{klp}. The variety $\mathfrak{NCAss}^5$ is determined by the rigid algebra:

{\small
\begin{longtable}{|l|l|llllll|}
\caption*{ }   \\

\hline \multicolumn{1}{|c|}{${\mathcal A}$} & \multicolumn{1}{c|}{ } & \multicolumn{6}{c|}{Multiplication table} \\ \hline 
\endfirsthead

 \multicolumn{8}{l}%
{{\bfseries  continued from previous page}} \\
\hline \multicolumn{1}{|c|}{${\mathcal A}$} & \multicolumn{1}{c|}{ } & \multicolumn{6}{c|}{Multiplication table} \\ \hline 
\endhead

\hline \multicolumn{8}{|r|}{{Continued on next page}} \\ \hline
\endfoot

\hline 
\endlastfoot

$\mathfrak{NCA}_{1}^5$ & $\mathbf{A}_{1}$ & $e_{1}e_1=e_{2}$ &  $e_{1}e_{2}=e_{3}$ &  $e_{1}e_{3}=e_{4}$  & $e_{1}e_{4}=e_5$ & $e_{2}e_{2}=e_4$ &  $e_{2}e_{3}=e_{5}$ \\ \hline

\end{longtable}
}

The complete graph of degenerations can be found in~\cite{klp}.

\subsection{Cyclic associative  algebras} 

We now consider the associative algebras which also satisfy the cyclic identity
\begin{center}
$(xy)z=(yz)x.$
\end{center}
 This variety will be denoted by $\mathfrak{CyAss}$.

\subsubsection{2-dimensional cyclic associative algebras} 
Thanks to \cite{aks24}, each $2$-dimensional cyclic associative algebra is commutative associative.
The variety $\mathfrak{CyAss}^2$  is determined by the rigid algebra:
{\small
\begin{longtable}{|l|l|ll|}
\caption*{ }   \\

\hline \multicolumn{1}{|c|}{${\mathcal A}$} & \multicolumn{1}{c|}{ } & \multicolumn{2}{c|}{Multiplication table} \\ \hline 
\endfirsthead

 \multicolumn{4}{l}%
{{\bfseries  continued from previous page}} \\
\hline \multicolumn{1}{|c|}{${\mathcal A}$} & \multicolumn{1}{c|}{ } & \multicolumn{2}{c|}{Multiplication table} \\ \hline 
\endhead

\hline \multicolumn{4}{|r|}{{Continued on next page}} \\ \hline
\endfoot

\hline 
\endlastfoot

$\mathfrak{CyAss}_{1}^2$ & $\mathcal{A}_{01}$ & $e_1e_1=e_1$ & $e_2e_2=e_2$  \\ \hline

\end{longtable}
}

\subsubsection{3-dimensional nilpotent cyclic associative algebras}
The algebraic and geometric classification of $3$-dimensional nilpotent cyclic associative algebras can be found in \cite{aks24}. 
In particular, it is proven that the variety  $\mathfrak{NCyAss}^3$ has two irreducible components:
\[{\rm Irr}(\mathfrak{NCyAss}^3)=\left\{
\overline{{\mathcal O}\big(\mathfrak{NCyAss}_{1}^3\big)}\right\} \cup 
\left\{
\overline{\bigcup {\mathcal O}\big(\mathfrak{NCyAss}_{2}^3(\alpha)\big)}\right\} ,\]
where the algebras $\mathfrak{NCyAss}_{1}^3$ and $\mathfrak{NCyAss}_{2}^3(\alpha)$ are defined as follows:

{\small

\begin{longtable}{|l|l|llll|}
\caption*{ }   \\

\hline \multicolumn{1}{|c|}{${\mathcal A}$} & \multicolumn{1}{c|}{ } & \multicolumn{4}{c|}{Multiplication table} \\ \hline 
\endfirsthead

 \multicolumn{6}{l}%
{{\bfseries  continued from previous page}} \\
\hline \multicolumn{1}{|c|}{${\mathcal A}$} & \multicolumn{1}{c|}{ } & \multicolumn{4}{c|}{Multiplication table} \\ \hline 
\endhead

\hline \multicolumn{6}{|r|}{{Continued on next page}} \\ \hline
\endfoot

\hline 
\endlastfoot

$\mathfrak{NCyAss}_{1}^3$ & $\mathcal{A}_{06}$&$ e_{1}  e_{1}=e_{2}$ & $e_{1} e_{2}=e_{3}$ & $e_{2} e_{1}=e_{3}$ & \\ \hline
$\mathfrak{NCyAss}_{2}^3(\alpha)$ & $\mathfrak{a}_{02}^{\alpha } $&$ e_{1}  e_{1}=e_{3}$&$e_{1} e_{2}=e_{3}$&$e_{2}  e_{1}=-e_{3}$&$e_{2}  e_{2}=\alpha e_{3}$\\  \hline 
\end{longtable}
}

\subsubsection{3-dimensional cyclic associative algebras}
The algebraic and geometric classification of $3$-dimensional  cyclic associative algebras can be found in \cite{aks24}. 
In particular, it is proven that the variety  $\mathfrak{CyAss}^3$ has two irreducible components:
\[{\rm Irr}(\mathfrak{CyAss}^3)=\left\{
\overline{{\mathcal O}\big(\mathfrak{CyAss}_{1}^3\big)}\right\} \cup 
\left\{
\overline{\bigcup {\mathcal O}\big(\mathfrak{CyAss}_{2}^3(\alpha)\big)}\right\} ,\]
where the algebras $\mathfrak{CyAss}_{1}^3$ and $\mathfrak{CyAss}_{2}^3(\alpha)$ are defined as follows:

{\small
\begin{longtable}{|l|l|llll|}
\caption*{ }   \\

\hline \multicolumn{1}{|c|}{${\mathcal A}$} & \multicolumn{1}{c|}{ } & \multicolumn{4}{c|}{Multiplication table} \\ \hline 
\endfirsthead

 \multicolumn{6}{l}%
{{\bfseries  continued from previous page}} \\
\hline \multicolumn{1}{|c|}{${\mathcal A}$} & \multicolumn{1}{c|}{ } & \multicolumn{4}{c|}{Multiplication table} \\ \hline 
\endhead

\hline \multicolumn{6}{|r|}{{Continued on next page}} \\ \hline
\endfoot

\hline 
\endlastfoot

$\mathfrak{CyAss}_{1}^3$ & $\mathcal{A}_{08}$&$ e_{1}  e_{1}=e_{1}$&$e_{2}  e_{2}=e_{2}$&$e_{2}   e_{3}=e_{3}$ & \\ \hline
$\mathfrak{CyAss}_{2}^3(\alpha)$ & $\mathfrak{a}_{02}^{\alpha } $&$ e_{1}  e_{1}=e_{3}$&$e_{1} e_{2}=e_{3}$&$e_{2}  e_{1}=-e_{3}$&$e_{2}  e_{2}=\alpha e_{3}$\\  \hline 
\end{longtable}
}

\subsubsection{4-dimensional nilpotent cyclic associative algebras}
The algebraic and geometric classification of $4$-dimensional nilpotent cyclic associative algebras can be found in \cite{aks24}. 
In particular, it is proven that the variety  $\mathfrak{NCyAss}^4$ has four irreducible components:
\[{\rm Irr}(\mathfrak{NCyAss}^4)=\left\{
\overline{{\mathcal O}\big(\mathfrak{NCyAss}_{1}^4\big)}\right\} \cup 
\left\{
\overline{\bigcup {\mathcal O}\big(\mathfrak{NCyAss}_{i}^4(\alpha)\big)}\right\}_{i=2}^4 ,\]
where the algebras $\mathfrak{NCyAss}_{1}^4$ and $\mathfrak{NCyAss}_{i}^4(\alpha)$ are defined as follows:

{\small
\begin{longtable}{|l|l|llll|}
\caption*{ }   \\

\hline \multicolumn{1}{|c|}{${\mathcal A}$} & \multicolumn{1}{c|}{ } & \multicolumn{4}{c|}{Multiplication table} \\ \hline 
\endfirsthead

 \multicolumn{6}{l}%
{{\bfseries  continued from previous page}} \\
\hline \multicolumn{1}{|c|}{${\mathcal A}$} & \multicolumn{1}{c|}{ } & \multicolumn{4}{c|}{Multiplication table} \\ \hline 
\endhead

\hline \multicolumn{6}{|r|}{{Continued on next page}} \\ \hline
\endfoot

\hline 
\endlastfoot

$\mathfrak{NCyAss}_{1}^4$ & $\mathcal{A}_{16}$&
$e_{1}  e_1 =e_{2} $&$e_{1}  e_{2}=e_{3}$&  $e_{1}   e_{3}=e_{4}$& \\
&&$e_{2}  e_{1}=e_{3}$& $e_{2} e_{2}=e_{4}$&
 $e_{3}   e_{1}=e_{4}$&
\\ \hline
$\mathfrak{NCyAss}_{2}^4(\alpha)$ & $\mathfrak{a}_{02}^{\alpha } $&$ e_{1}  e_{1}=e_{3}$&$e_{1} e_{2}=e_{3}$&$e_{2}  e_{1}=-e_{3}$&$e_{2}  e_{2}=\alpha e_{3}$\\  \hline 

$\mathfrak{NCyAss}_{3}^4(\alpha)$ & $\mathfrak{a}_{10}^{\alpha }$& $e_{1}  e_{1}=e_{3}$&$e_{1}  e_{2}=e_{3}+e_{4}$&$e_{2}  e_{1}=e_{4}-e_{3}$&$e_{2}  e_{2}=\alpha e_{3}$\\ \hline 

$\mathfrak{NCyAss}_{4}^4(\alpha)$ & $\mathfrak{a}_{13}$&$ e_{1}  e_{1}=e_{4}$&$e_{1}  e_{2}=e_{3}$&$e_{2}  e_{1}=-e_{3}$&\\&& 
 $e_{2}  e_{2}=e_{3}$&$e_{1}  e_{4}=e_{3}$&$e_{4}  e_{1}=e_{3}$&\\

\end{longtable}}
\subsubsection{4-dimensional cyclic associative algebras}
The algebraic and geometric classification of $4$-dimensional   cyclic associative algebras can be found in \cite{aks24}. 
In particular, it is proven that the variety  $\mathfrak{CyAss}^4$ has three irreducible components:
\[{\rm Irr}(\mathfrak{CyAss}^4)=\left\{
\overline{{\mathcal O}\big(\mathfrak{CyAss}_{1}^4\big)}\right\} \cup 
\left\{
\overline{\bigcup {\mathcal O}\big(\mathfrak{CyAss}_{i}^4(\alpha)\big)}\right\}_{i=2}^3 ,\]
where the algebras $\mathfrak{CyAss}_{1}^4$ and $\mathfrak{CyAss}_{i}^4(\alpha)$ are defined as follows:

{\small

\begin{longtable}{|l|l|lllll|}
\caption*{ }   \\

\hline \multicolumn{1}{|c|}{${\mathcal A}$} & \multicolumn{1}{c|}{ } & \multicolumn{5}{c|}{Multiplication table} \\ \hline 
\endfirsthead

 \multicolumn{7}{l}%
{{\bfseries  continued from previous page}} \\
\hline \multicolumn{1}{|c|}{${\mathcal A}$} & \multicolumn{1}{c|}{ } & \multicolumn{5}{c|}{Multiplication table} \\ \hline 
\endhead

\hline \multicolumn{7}{|r|}{{Continued on next page}} \\ \hline
\endfoot

\hline 
\endlastfoot

$\mathfrak{NCyAss}_{1}^4$ &  $\mathcal{A}_{17}$& 
$e_{1} e_1=e_{1}$&$e_{2}  e_2=e_{2}$&$ e_{3}  e_3=e_{3}$&$ e_{4} 
 e_{4}=e_{4} $ & 
\\ \hline

$\mathfrak{NCyAss}_{2}^4(\alpha)$ & $\mathfrak{a}_{10}^{\alpha }$& $e_{1}  e_{1}=e_{3}$&$e_{1}  e_{2}=e_{3}+e_{4}$&$e_{2}  e_{1}=e_{4}-e_{3}$&$e_{2}  e_{2}=\alpha e_{3}$ &\\ \hline 

$\mathfrak{NCyAss}_{3}^4(\alpha)$ & $\mathfrak{a}_{12}^{\alpha }$& 
$e_{1}  e_{1}=e_{3}$&$e_{1}  e_{2}=e_{3} $&$ e_{2}  e_{1}=-e_{3}$&$e_{2}  e_{2}=\alpha e_{3}$&$e_{4}  e_{4}=e_{4}$\\

\end{longtable}
}

 \subsection{Jordan algebras} 
A commutative algebra $\mathfrak{J}$ is called a {\it Jordan} algebra if it satisfies the identity 
\[x^2(yx)=(x^2y)x.\] Let $\mathfrak{Jord}$ be the variety of Jordan algebras.

\subsubsection{2-dimensional Jordan algebras}
The algebraic classification of $2$-dimensional Jordan algebras was made between 1975 (a result by Gabriel, who described the associative ones) and 1989 (a result by Sherkulov for the non-associative ones).
The graph of degenerations can be deduced from~\cite{kv16} and is explicitly given in~\cite{gkp}. In particular, it is proven that the variety  $\mathfrak{Jord}^2$ has two irreducible components:
\[{\rm Irr}(\mathfrak{Jord}^2)=\left\{
\overline{{\mathcal O}\big(\mathfrak{J}_{i}^2\big)}\right\}_{i=1}^{2},\]
where the algebras $\mathfrak{J}_{1}^2$ and $\mathfrak{J}_{2}^2$ are defined as follows:

{\small
\vspace{-4mm}
\begin{longtable}{|l|l|ll|}
\caption*{ }   \\

\hline \multicolumn{1}{|c|}{${\mathcal A}$} & \multicolumn{1}{c|}{ } & \multicolumn{2}{c|}{Multiplication table} \\ \hline 
\endfirsthead

 \multicolumn{4}{l}%
{{\bfseries  continued from previous page}} \\
\hline \multicolumn{1}{|c|}{${\mathcal A}$} & \multicolumn{1}{c|}{ } & \multicolumn{2}{c|}{Multiplication table} \\ \hline 
\endhead

\hline \multicolumn{4}{|r|}{{Continued on next page}} \\ \hline
\endfoot

\hline 
\endlastfoot

$\mathfrak{J}_{1}^2$ & $\mathfrak{B}_2$ & $e_1e_1 = e_1$ & $e_1e_2 =\frac{1}{2}e_2$  \\ \hline
$\mathfrak{J}_{2}^2$ & $\mathfrak{B}_4$ & $e_1e_1 = e_1$ & $e_2e_2 = e_2$  \\ \hline
\end{longtable}
}

\subsubsection{3-dimensional nilpotent  Jordan algebras}
The nilpotent Jordan algebras of dimension $3$ were classified algebraically and geometrically in~\cite{esp}. The variety $\mathfrak{NJord}^3$ is irreducible, determined by the rigid algebra:

{\small
\vspace{-4mm}
\begin{longtable}{|l|l|ll|}
\caption*{ }   \\

\hline \multicolumn{1}{|c|}{${\mathcal A}$} & \multicolumn{1}{c|}{ } & \multicolumn{2}{c|}{Multiplication table} \\ \hline 
\endfirsthead

 \multicolumn{4}{l}%
{{\bfseries  continued from previous page}} \\
\hline \multicolumn{1}{|c|}{${\mathcal A}$} & \multicolumn{1}{c|}{ } & \multicolumn{2}{c|}{Multiplication table} \\ \hline 
\endhead

\hline \multicolumn{4}{|r|}{{Continued on next page}} \\ \hline
\endfoot

\hline 
\endlastfoot

$\mathfrak{NJ}_{1}^3$ & $\phi_1$ & $e_1e_1 = e_2$ &  $e_1e_2 = e_3$  \\ \hline

\end{longtable}
}

\subsubsection{3-dimensional Jordan algebras}
Also in~\cite{gkp}, the geometric classification in dimension $3$ (initiated in~\cite{KS07} together with the algebraic classification) was completed. We will refer to the notation of~\cite{gkp}. There exist five irreducible components in $\mathfrak{Jord}^3$:
\[{\rm Irr}(\mathfrak{Jord}^3)=\left\{
\overline{{\mathcal O}\big(\mathfrak{J}_{i}^3\big)}\right\}_{i=1}^{5},\]
where 

{\small

\begin{longtable}{|l|l|lllll|}
\caption*{ }   \\

\hline \multicolumn{1}{|c|}{${\mathcal A}$} & \multicolumn{1}{c|}{ } & \multicolumn{5}{c|}{Multiplication table} \\ \hline 
\endfirsthead

\multicolumn{7}{l}%
{{\bfseries  continued from previous page}} \\
\hline \multicolumn{1}{|c|}{${\mathcal A}$} & \multicolumn{1}{c|}{ } & \multicolumn{5}{c|}{Multiplication table} \\ \hline 
\endhead

\hline \multicolumn{7}{|r|}{{Continued on next page}} \\ \hline
\endfoot

\hline 
\endlastfoot

$\mathfrak{J}_{1}^3$ & $\mathbb{T}_{01}$ & $e_1e_1 = e_1$ & $e_2e_2 = e_2$ & $e_3e_3 = e_3$ & & \\ \hline
$\mathfrak{J}_{2}^3$ & $\mathbb{T}_{02}$ & $e_1e_1 = e_1$ & $e_1e_3 = \frac{1}{2}e_3$ & $e_2e_2 =e_2$ & $e_2e_3 = \frac{1}{2}e_3$ & $e_3e_3 = e_1 + e_2$  \\ \hline
$\mathfrak{J}_{3}^3$ & $\mathbb{T}_{05}$ & $e_1e_1 = e_1$ & $e_1e_3 = \frac{1}{2}e_3$ & $e_2e_2 = e_2$ & &  \\ \hline
$\mathfrak{J}_{4}^3$ & $\mathbb{T}_{10}$ & $e_1e_1 = e_1$ & $e_1e_2 = \frac{1}{2}e_2$ & $e_1e_3 = e_3$ & $e_2e_2 = e_3$ & \\ \hline
$\mathfrak{J}_{5}^3$ & $\mathbb{T}_{12}$ & $e_1e_1 = e_1$ & $e_1e_2 = \frac{1}{2}e_2$ & $e_1e_3 = \frac{1}{2}e_3$ && \\ \hline

\end{longtable}
}

The complete graph of degenerations can be found in~\cite{gkp}.

\subsubsection{4-dimensional nilpotent Jordan algebras}
The variety $\mathfrak{NJord}^4$  was studied in 
\cite{esp} both algebraically and geometrically. It has two irreducible components:
\[{\rm Irr}(\mathfrak{NJord}^4)=
\left\{\overline{{\mathcal O}\left(\mathfrak{NJ}_{i}^4\right)}\right\}_{i=1}^{2},\]
where

{\small

\begin{longtable}{|l|l|llll|}
\caption*{ }   \\

\hline \multicolumn{1}{|c|}{${\mathcal A}$} & \multicolumn{1}{c|}{ } & \multicolumn{4}{c|}{Multiplication table} \\ \hline 
\endfirsthead

 \multicolumn{6}{l}%
{{\bfseries  continued from previous page}} \\
\hline \multicolumn{1}{|c|}{${\mathcal A}$} & \multicolumn{1}{c|}{ } & \multicolumn{4}{c|}{Multiplication table} \\ \hline 
\endhead

\hline \multicolumn{6}{|r|}{{Continued on next page}} \\ \hline
\endfoot

\hline 
\endlastfoot

$\mathfrak{NJ}_{1}^4$ & $\phi_1$ & $e_1e_1  = e_2$ & $e_1 e_2 = e_3$ &  $e_1 e_3 = e_4$ & $e_2 e_2 = e_4$  \\ \hline 
$\mathfrak{NJ}_{2}^4$ & $\phi_2$ & $e_1 e_1 = e_3$ & $e_1 e_3 = e_4$ & $e_2 e_2 = e_3$ & \\ \hline

\end{longtable}
}

\subsubsection{4-dimensional Jordan algebras}
Regarding $\mathfrak{Jord}^4$, the algebraic classification was made by Martin in 2013, and the ten irreducible components were found later in~\cite{KM}:
\[{\rm Irr}(\mathfrak{Jord}^4)=\left\{\overline{{\mathcal O}\big(\mathfrak{J}_{i}^4\big)} \right\}_{i=1}^{10},\]
where

{\small

\begin{longtable}{|l|l|llllll|}
\caption*{ }   \\

\hline \multicolumn{1}{|c|}{${\mathcal A}$} & \multicolumn{1}{c|}{ } & \multicolumn{6}{c|}{Multiplication table} \\ \hline 
\endfirsthead

 \multicolumn{7}{l}%
{{\bfseries  continued from previous page}} \\
\hline \multicolumn{1}{|c|}{${\mathcal A}$} & \multicolumn{1}{c|}{ } & \multicolumn{6}{c|}{Multiplication table} \\ \hline 
\endhead

\hline \multicolumn{8}{|r|}{{Continued on next page}} \\ \hline
\endfoot

\hline 
\endlastfoot

$\mathfrak{J}_{1}^4$ & $\mathfrak{J}_{1}$ & $e_1e_1 = e_1$ & $e_1e_3 = \frac{1}{2}e_3$ & $e_2e_2 =e_2$& $e_2e_3 = \frac{1}{2}e_3$ & $e_3e_3 = e_1 + e_2$  & $e_4e_4=e_4$   \\ \hline
$\mathfrak{J}_{2}^4$ & $\mathfrak{J}_{2}$ & $e_1e_3 = \frac{1}{2}e_3$ & $e_1e_4 = \frac{1}{2}e_4$ & $e_2e_3 = \frac{1}{2}e_3$ & $e_2e_4 =\frac{1}{2}e_4$ & $e_3e_4 = \frac{1}{2}(e_1 + e_2)$  & \\ \hline
$\mathfrak{J}_{3}^4$ & $\mathfrak{J}_{3}$ & $e_1e_1=e_1$ & $e_2e_2=e_2$ & $e_3e_3=e_3$ & $e_4e_4=e_4$ && \\ \hline
$\mathfrak{J}_{4}^4$ & $\mathfrak{J}_{6}$ & $e_1e_1 = e_1$ & $e_1e_2 =\frac{1}{2}e_2$ & $e_3e_3=e_3$ & $e_4e_4=e_4$ & & \\ \hline
$\mathfrak{J}_{5}^4$ & $\mathfrak{J}_{12}$ & $e_1e_1 = e_1$ & $e_1e_2 = \frac{1}{2}e_2$ & $e_1e_3 = \frac{1}{2}e_3$ & $e_4e_4=e_4$ & &\\ \hline
$\mathfrak{J}_{6}^4$ & $\mathfrak{J}_{13}$ & $e_1e_1 = e_1$ & $e_1e_2 =\frac{1}{2}e_2$ & $e_3e_3=e_3$ & $e_3e_4 =\frac{1}{2}e_4$ && \\ \hline
$\mathfrak{J}_{7}^4$ & $\mathfrak{J}_{16}$ & $e_1e_3 = \frac{1}{2}e_3$ & $e_1e_4 = \frac{1}{2}e_4$ & $e_2e_3 = \frac{1}{2}e_3$ && &  \\ \hline
$\mathfrak{J}_{8}^4$ & $\mathfrak{J}_{24}$ & $e_1e_1 = e_1$ & $e_1e_2 = \frac{1}{2}e_2$ & $e_1e_3 = e_3$ & $e_2e_2 = e_3$ & $e_4e_4=e_4$ & \\ \hline
$\mathfrak{J}_{9}^4$ & $\mathfrak{J}_{33}$ & $e_1e_2 = \frac{1}{2}e_2$ & $e_1e_3 = \frac{1}{2}e_3$ & $e_1e_4 = \frac{1}{2}e_4$ && &\\ \hline
$\mathfrak{J}_{10}^4$ & $\mathfrak{J}_{59}$ & $e_1e_2 = e_2$ & $e_1e_3 = \frac{1}{2}e_3$ & $e_1e_4= \frac{1}{2}e_4$ & $e_3e_4=e_2$ & $e_4e_4=e_2$  &\\ \hline

\end{longtable}
}

\subsubsection{5-dimensional nilpotent Jordan algebras}
In dimension $5$ there is no complete algebraic classification of  Jordan algebras yet. However, nilpotent algebras were classified thanks to~\cite{maz80} and the work of Abdelwahab and Hegazi (2016). The geometric classification is given in~\cite{KM18}. The authors found that the variety $\mathfrak{NJord}^5$ has five irreducible components:
\[{\rm Irr}(\mathfrak{NJord}^5)=\left\{\overline{{\mathcal O}\left(\mathfrak{NJ}_{i}^5\right)}\right\}_{i=1}^{4}\cup\left\{\overline{\bigcup {\mathcal O}\left(\mathfrak{NJ}_{5}^5(\alpha,\beta)\right)}\right\},\]
where

{\small
\vspace{-5mm}
\begin{longtable}{|l|l|llllll|}
\caption*{ }   \\

\hline \multicolumn{1}{|c|}{${\mathcal A}$} & \multicolumn{1}{c|}{ } & \multicolumn{6}{c|}{Multiplication table} \\ \hline 
\endfirsthead

 \multicolumn{8}{l}%
{{\bfseries  continued from previous page}} \\
\hline \multicolumn{1}{|c|}{${\mathcal A}$} & \multicolumn{1}{c|}{ } & \multicolumn{6}{c|}{Multiplication table} \\ \hline 
\endhead

\hline \multicolumn{8}{|r|}{{Continued on next page}} \\ \hline
\endfoot

\hline 
\endlastfoot

$\mathfrak{NJ}_{1}^5$ & $\epsilon_1$ & $e_1e_1 = e_2$ & $e_1e_2=e_3$ & $e_1e_3=e_4$ & $e_1e_4=e_5$ & $e_2e_2 =e_4$ & $e_2e_3 = e_5$ \\ \hline
$\mathfrak{NJ}_{2}^5$ & $\mathfrak{J}_{21}$ & $e_1e_1 = e_5$ & $e_1e_2=e_4$ & $e_2e_2=e_5$ & $e_3e_3=e_4$ & $e_3e_4 = e_5$ & \\ \hline
$\mathfrak{NJ}_{3}^5$ & $\mathfrak{J}_{22}$ & $e_1e_1=e_2$ & $e_1e_2=e_4$ & $e_1e_4=e_5$ & $e_2e_2=e_2$ & $e_3e_3=e_4$ & \\ \hline
$\mathfrak{NJ}_{4}^5$ & $\mathfrak{J}_{40}$ & $e_1e_1 = e_5$ & $e_1e_2=e_3$ & $e_1e_3=e_4$ & $e_2e_2=e_4$ & $e_2e_3 = e_5$ & \\ \hline
$\mathfrak{NJ}_{5}^5(\alpha,\beta)$ & $\mathfrak{N}_{27}^{\#}$ & $e_1e_1 = e_3$ & $e_1e_3 = \alpha e_5$ & $e_1e_4 = e_5$& $e_2e_2=e_4$ & $e_2e_3=e_5$  & $e_2e_4=\beta e_5$ \\ \hline
\end{longtable}
}

\subsection{Kokoris algebras} 

An algebra $\mathfrak{K}$ is called a {\it Kokoris} algebra if it satisfies the identities 
\[(x , y, z)_{\circ}=0, \  (x,y,z)=-(z,y,x).\] 
Let $\mathfrak{K}$ be the variety of Kokoris algebras.

\subsubsection{2-dimensional Kokoris  algebras}
The algebraic and geometric classification of $2$-dimensional   Kokoris algebras can be found in \cite{aak24}.  In particular, it is proven that the variety  $\mathfrak{K}^2$ has two irreducible components:
\[{\rm Irr}(\mathfrak{K}^2)=\left\{
\overline{{\mathcal O}\big(\mathfrak{K}_{i}^2\big)}\right\}_{i=1}^2,\]
where  

{\small
\vspace{-5mm}
\begin{longtable}{|l|l|llll|}
\caption*{ }   \\

\hline \multicolumn{1}{|c|}{${\mathcal A}$} & \multicolumn{1}{c|}{ } & \multicolumn{4}{c|}{Multiplication table} \\ \hline 
\endfirsthead

 \multicolumn{6}{l}%
{{\bfseries  continued from previous page}} \\
\hline \multicolumn{1}{|c|}{${\mathcal A}$} & \multicolumn{1}{c|}{ } & \multicolumn{4}{c|}{Multiplication table} \\ \hline 
\endhead

\hline \multicolumn{6}{|r|}{{Continued on next page}} \\ \hline
\endfoot

\hline 
\endlastfoot
$\mathfrak{K}^2_{1}$ & ${\bf A}_{08}$ &    $e_1e_1=e_1$ & $e_2e_2=e_2$  &&\\ \hline

$\mathfrak{K}^2_{2} $ & ${\bf A}_{22}^{0}$ &      $ e_1e_2 =e_2$ &  $ e_2e_1 =-e_2$  && \\ \hline
\end{longtable}
}

\subsubsection{3-dimensional Kokoris  nilpotent algebras}
The algebraic and geometric classification of $3$-dimensional nilpotent Kokoris algebras can be obtained from the classification and description of degenerations of $3$-dimensional nilpotent algebras given in \cite{fkkv}.
Hence, we have  that the variety   $\mathfrak{NK}^3$ has two irreducible components:
\[{\rm Irr}(\mathfrak{NK}^3)=
\left\{\overline{{\mathcal O}\left(\mathfrak{NK}_{1}^{3}\right)}\right\}\cup\left\{
\overline{\bigcup {\mathcal O}\left(\mathfrak{NK}_{4}^3(\alpha)\right)}\right\},\]
where

{\small
\vspace{-5mm}
\begin{longtable}{|l|l|llllll|}
\caption*{ }   \\

\hline \multicolumn{1}{|c|}{${\mathcal A}$} & \multicolumn{1}{c|}{ } & \multicolumn{6}{c|}{Multiplication table} \\ \hline 
\endfirsthead

 \multicolumn{8}{l}%
{{\bfseries  continued from previous page}} \\
\hline \multicolumn{1}{|c|}{${\mathcal A}$} & \multicolumn{1}{c|}{ } & \multicolumn{6}{c|}{Multiplication table} \\ \hline 
\endhead

\hline \multicolumn{8}{|r|}{{Continued on next page}} \\ \hline
\endfoot

\hline 
\endlastfoot
$\mathfrak{NA}_{1}^3$ & $\mathcal{N}_{4}(1)$  & $e_1e_1=e_3$ &  $e_1e_2=e_3$ & $e_2e_1=e_3$ &&&\\ \hline

$\mathfrak{NA}_{2}^3(\alpha)$ &  $\mathcal{N}_{8}(\alpha)$  & $e_1e_1=\alpha e_3$ & $e_2e_1=e_3$ & $e_2e_2=e_3$ && &\\ \hline
\end{longtable}
}

\subsubsection{3-dimensional Kokoris    algebras}
The algebraic and geometric classification of $3$-dimensional   Kokoris algebras can be found in \cite{aak24}.  In particular, it is proven that the variety  $\mathfrak{K}^3$ has five irreducible components:
\[{\rm Irr}(\mathfrak{K}^3)=\left\{
\overline{{\mathcal O}\big(\mathfrak{K}_{i}^3\big)}\right\}_{i=1}^3 \cup \left\{
\overline{\bigcup {\mathcal O}\big(\mathfrak{K}_{i}^3(\alpha)\big)} \right\}_{i=4}^5 ,\]
where  

{\small
\vspace{-5mm}
\begin{longtable}{|l|l|llll|}
\caption*{ }   \\

\hline \multicolumn{1}{|c|}{${\mathcal A}$} & \multicolumn{1}{c|}{ } & \multicolumn{4}{c|}{Multiplication table} \\ \hline 
\endfirsthead

 \multicolumn{6}{l}%
{{\bfseries  continued from previous page}} \\
\hline \multicolumn{1}{|c|}{${\mathcal A}$} & \multicolumn{1}{c|}{ } & \multicolumn{4}{c|}{Multiplication table} \\ \hline 
\endhead

\hline \multicolumn{6}{|r|}{{Continued on next page}} \\ \hline
\endfoot

\hline 
\endlastfoot
$\mathfrak{K}_{1}^3$& ${\bf A}_{04} $ &   $e_{1} e_{1}=e_{1}$ & $ e_{2} e_{2}=e_{2}$ && \\ 
\hline

$\mathfrak{K}_{2}^3$&
${\bf A}_{29}$ &   $e_{1} e_{1}=e_{1}$ & $e_{1} e_{2}=e_{2}$ & $e_{2} e_{1}=e_{2}$ & $e_{1} e_{3}=e_{3}$ \\
&& $e_{3} e_{1}=e_{3}$ & $e_{2} e_{3}=e_{3}$ & $e_{3}\ e_{2}=-e_{3}$ & \\
\hline

$\mathfrak{K}_{3}^3$&${\bf A}_{30}$ &   
$e_{1} e_{1}=e_{1}$ & $e_{2} e_{3}=e_{3}$ & $e_{3} e_{2}=-e_{3}$ &\\
\hline

$\mathfrak{K}_{4}^3(\alpha)$&  ${\bf A}_{02}^{\alpha}$ &  $e_1 e_2=\left( 1+\alpha \right) e_3$ & $e_2 e_1=\left(1-\alpha \right) e_3$ && \\

\hline 
$\mathfrak{K}_{5}^3(\alpha)$&${\bf A}_{24}^{\alpha }$ &   
$  e_{1}e_{2}  =e_{3}$ & 
$e_{1}e_{3}  =e_{1}+e_{3}$ &
$e_{2}e_{1}   =-e_{3}$ &\\
&&$ e_{2}e_{3} =\alpha e_{2}$ &
$e_{3}e_{1}  =-e_{1}-e_{3}$
& $ e_{3}e_{2} =-\alpha e_{2}$&\\
\hline 

\end{longtable}
}
\subsubsection{4-dimensional     nilpotent Kokoris algebras}

The variety $\mathfrak{NK}^4$  was studied in 
\cite{aak24} both algebraically, and geometrically. It has five irreducible components:
\[{\rm Irr}(\mathfrak{NK}^4)=
\left\{\overline{{\mathcal O}\left(\mathfrak{NK}_{i}^4\right)}\right\}_{i=1}^{3}\cup\left\{\overline{\bigcup {\mathcal O}\left(\mathfrak{NK}_{i}^4(\alpha)\right)}\right\}_{i=4}^{5},\]
where

{\small
\vspace{-5mm}
\begin{longtable}{|l|l|lllll|}
\caption*{ }   \\

\hline \multicolumn{1}{|c|}{${\mathcal A}$} & \multicolumn{1}{c|}{ } & \multicolumn{5}{c|}{Multiplication table} \\ \hline 
\endfirsthead

 \multicolumn{7}{l}%
{{\bfseries  continued from previous page}} \\
\hline \multicolumn{1}{|c|}{${\mathcal A}$} & \multicolumn{1}{c|}{ } & \multicolumn{5}{c|}{Multiplication table} \\ \hline 
\endhead

\hline \multicolumn{7}{|r|}{{Continued on next page}} \\ \hline
\endfoot

\hline 
\endlastfoot

$\mathfrak{NK}_{1}^4$ &  ${\mathcal J}_{03}$ & $e_1 e_1 = e_2$ & $e_1 e_2=e_4$  & $e_2 e_1=e_4$ & $e_3 e_1=e_4$ & $e_3 e_3=e_4$\\\hline
$\mathfrak{NK}_{1}^4$ & $ {\mathcal J}_{17}$ & $e_1 e_1 = e_4$ & $e_1 e_2 = e_3$ & $e_2 e_1=-e_3$ && \\&& $e_1 e_3=e_4$ & $e_2 e_2=e_4$ & $e_3 e_1=-e_4$ &&\\ \hline
$\mathfrak{NK}_{1}^4$ & ${\mathcal J}_{18}$ & $e_1 e_1 = e_2$ & $e_1 e_2 = e_3$ & $e_1 e_3=e_4$ && \\&& $e_2 e_1=e_3$ & $e_2 e_2=e_4$  & $e_3 e_1=e_4$ &&\\ \hline
$\mathfrak{NK}_{4}^4(\alpha)$ & $\mathfrak{N}_2(\alpha)$ & $e_1e_1 = e_3$ & $e_1e_2 = e_4$  &$e_2e_1 = -\alpha e_3$ & $e_2e_2 = -e_4$ & \\ \hline
$\mathfrak{NK}_{5}^4(\alpha)$ & $\mathfrak{N}_3(\alpha)$ & $e_1e_1 = e_4$ & $e_1e_2 = \alpha e_4$  & $e_2e_1 = -\alpha e_4$ & $e_2e_2 = e_4$ & \\ \hline
\end{longtable}
}

\subsection{Standard  algebras}

An algebra $\mathfrak{S}$ is called a {\it standard} algebra if it satisfies the identities 
\[(x,y,z)+(z,x,y)=(x,z,y), \ \ (x,y,wz)+(w,y,xz)+(z,y,wx) =0.\] 
Let $\mathfrak{S}$ be the variety of standard algebras.

\subsubsection{2-dimensional standard  algebras}

The algebraic and geometric classification of $2$-dimensional standard algebras can be found in \cite{aak24}. 
It is proven that the variety  $\mathfrak{S}^2$ has four irreducible components:
\[{\rm Irr}(\mathfrak{S}^2)=\left\{
\overline{{\mathcal O}\big(\mathfrak{S}_{i}^2\big)}\right\}_{i=1}^{4},\]
where  

{\small
\vspace{-5mm}
\begin{longtable}{|l|l|lll|}
\caption*{ }   \\

\hline \multicolumn{1}{|c|}{${\mathcal A}$} & \multicolumn{1}{c|}{ } & \multicolumn{3}{c|}{Multiplication table} \\ \hline 
\endfirsthead

 \multicolumn{5}{l}%
{{\bfseries  continued from previous page}} \\
\hline \multicolumn{1}{|c|}{${\mathcal A}$} & \multicolumn{1}{c|}{ } & \multicolumn{3}{c|}{Multiplication table} \\ \hline 
\endhead

\hline \multicolumn{5}{|r|}{{Continued on next page}} \\ \hline
\endfoot

\hline 
\endlastfoot

$\mathfrak{S}_{1}^2$ & ${\bf A}_{08} $ & $e_1e_1 = e_1$ & $e_2e_2 = e_2$ &\\ \hline
$\mathfrak{S}_{2}^2$ & ${\bf A}_{18}^{0}$  & $e_1e_1 = e_1$ & $e_1e_2 =\frac{1}{2}e_2$ & $e_2e_1 =\frac{1}{2}e_2$  \\ \hline
$\mathfrak{S}_{3}^2$ &${\bf A}_{18}^{\frac{1}{2}}$ &   $e_{1} e_{1}=e_{1}$ & $e_{1}  e_{2}=e_{2}$& \\\hline

$\mathfrak{S}_{4}^2$ & ${\bf A}_{18}^{-\frac{1}{2}}$ &   $e_{1}  e_{1}=e_{1}$ & $e_{2}  e_{1}=e_{2}$&\\
\end{longtable}
}

\subsubsection{3-dimensional nilpotent standard  algebras}
Thanks to \cite{aak24},
the varieties of $3$-dimensional nilpotent standard and nilpotent noncommutative Jordan algebras coincide. 
The geometric classification of  $3$-dimensional nilpotent noncommutative Jordan algebras is given in \cite{jkk19}.   Hence,  $\mathfrak{NS}^3$ has two irreducible components:
\[{\rm Irr}(\mathfrak{NS}^3)=\left\{
\overline{{\mathcal O}\big(\mathfrak{NS}_{1}^3\big)}\right\} \cup \left\{
\overline{\bigcup {\mathcal O}\big(\mathfrak{NS}_{2}^3(\alpha)\big)}\right\} ,\]
where  

{\small
\vspace{-5mm}
\begin{longtable}{|l|l|lll|}
\caption*{ }   \\

\hline \multicolumn{1}{|c|}{${\mathcal A}$} & \multicolumn{1}{c|}{ } & \multicolumn{3}{c|}{Multiplication table} \\ \hline 
\endfirsthead

 \multicolumn{5}{l}%
{{\bfseries  continued from previous page}} \\
\hline \multicolumn{1}{|c|}{${\mathcal A}$} & \multicolumn{1}{c|}{ } & \multicolumn{3}{c|}{Multiplication table} \\ \hline 
\endhead

\hline \multicolumn{5}{|r|}{{Continued on next page}} \\ \hline
\endfoot

\hline 
\endlastfoot

$\mathfrak{NS}_{1}^3$ & ${\mathcal J}^3_{01}$ & $e_1 e_1 = e_2$ & $e_1 e_2=e_3$ & $e_2 e_1= e_3$ \\ \hline
$\mathfrak{NS}_{2}^3(\alpha)$ & ${\mathcal J}^{3*}_{04}(\alpha)$ & $e_1 e_1 = \alpha e_3$ & $e_2 e_1=e_3$ & $e_2 e_2=e_3$  \\ \hline 
\end{longtable}
}

\subsubsection{3-dimensional   standard  algebras}
The algebraic and geometric classification of $3$-dimensional   standard algebras can be found in \cite{aak24}.  In particular, it is proven that the variety  $\mathfrak{S}^3$ has fourteen irreducible components:
\[{\rm Irr}(\mathfrak{S}^3)=\left\{
\overline{{\mathcal O}\big(\mathfrak{S}_{i}^3\big)}\right\}_{i=1}^{13} \cup \left\{
\overline{\bigcup {\mathcal O}\big(\mathfrak{S}_{14}^3(\alpha)\big)} \right\},\]
where  
\vspace{-5mm}
{\small

\begin{longtable}{|l|l|llll|}
\caption*{ }   \\

\hline \multicolumn{1}{|c|}{${\mathcal A}$} & \multicolumn{1}{c|}{ } & \multicolumn{4}{c|}{Multiplication table} \\ \hline 
\endfirsthead

 \multicolumn{6}{l}%
{{\bfseries  continued from previous page}} \\
\hline \multicolumn{1}{|c|}{${\mathcal A}$} & \multicolumn{1}{c|}{ } & \multicolumn{4}{c|}{Multiplication table} \\ \hline 
\endhead

\hline \multicolumn{6}{|r|}{{Continued on next page}} \\ \hline
\endfoot

\hline 
\endlastfoot
$\mathfrak{S}_{1}^3$& ${\bf A}_{04} $ &   $e_{1} e_{1}=e_{1}$ & $ e_{2} e_{2}=e_{2}$ && \\ 
\hline

$\mathfrak{S}_{2}^3$& ${\bf A}_{12}$ &  
$e_{1} e_{1}=e_{1}$ & $e_{1} e_{3}= \frac{1}{2}e_{3}$ & $e_{2} e_{2}=e_{2}$ & 
  $e_{2} e_{3}=\frac{1}{2}e_{3}$ \\
&&
$e_{3} e_{1}= \frac{1}{2}e_{3}$ & $ e_{3} e_{3}=e_{1}+e_{2}$ & 
  $e_{3} e_{2}=\frac{1}{2}e_{3}$ &\\
\hline

$\mathfrak{S}_{3}^3$& ${\bf A}_{14}^{0,0}$ &  
$e_{1} e_{1}=e_{1}$ & $e_{1} e_{2}=\frac{1}{2}e_{2}$ & $e_{2} e_{1}=\frac{1}{2}e_{2}$  & \\
&& $e_{1} e_{3}=\frac{1}{2} e_{3}$ & $e_{3} e_{1}=\frac{1}{2}e_{3}$ && \\
\hline

$\mathfrak{S}_{4}^3$& ${\bf A}_{14}^{\frac{1}{2},\frac{1}{2}}$ & 
$e_{1} e_{1}=e_{1}$ & $e_{1} e_{2}=e_{2}$ &    $e_{1} e_{3}= e_{3}$ &\\
\hline

$\mathfrak{S}_{5}^3$& ${\bf A}_{14}^{-\frac{1}{2},-\frac{1}{2}}$   & 
$e_{1} e_{1}=e_{1}$   & $e_{2} e_{1}= e_{2}$  & $e_{3} e_{1}=e_{3}$ &\\
\hline

$\mathfrak{S}_{6}^3$&  ${\bf A}_{14}^{0,\frac{1}{2}}$ & 
$e_{1} e_{1}=e_{1}$ & $e_{1} e_{2}=\frac{1}{2}e_{2}$ & 
 $e_{2} e_{1}= \frac{1}{2} e_{2}$ & $e_{1} e_{3}= e_{3}$  \\
\hline

$\mathfrak{S}_{7}^3$&  ${\bf A}_{14}^{0,-\frac{1}{2}}$ & 
$e_{1} e_{1}=e_{1}$ & $e_{1} e_{2}=\frac{1}{2}e_{2}$ & $e_{2} e_{1}=\frac{1}{2}e_{2}$  
&  $e_{3} e_{1}= e_{3}$ \\
\hline

$\mathfrak{S}_{8}^3$&  ${\bf A}_{14}^{\frac{1}{2},-\frac{1}{2}}$ &   $e_{1} e_{1}=e_{1}$ & $e_{1} e_{2}= e_{2}$ & $e_{2} e_{1}=  e_{2}$  
&   $e_{3}  e_{1}= e_{3}$ \\
\hline

$\mathfrak{S}_{9}^3$& ${\bf A}_{16}$ &   
$e_{1} e_{1}=e_{1}$ & $e_{1} e_{2}= \frac{1}{2}e_{2}$ & $e_{1} e_{3}=e_{3}$   & \\&&$e_{2} e_{1}= \frac{1}{2}e_{2}$ & $e_{2}  e_{2}=e_{3}$ & $e_{3} e_{1}=e_{3}$&\\
\hline

$\mathfrak{S}_{10}^3$& 
${\bf A}_{17}^{\frac{1}{2}}$ &   $e_{1} e_{1}=e_{1}$ & $e_{2} e_{2}=e_{2}$ & $e_{1} e_{3}=e_{3}$ & $e_{3} e_{2}=e_{3}$\\
\hline

$\mathfrak{S}_{11}^3$& 
${\bf A}_{19}^{0}$ &   $e_{1} e_{1}=e_{1}$ & $e_{2} e_{2}=e_{2}$ & $e_{1} e_{3}=\frac{1}{2} e_{3}$ & $e_{3}  e_{1}=\frac{1}{2} e_{3}$   \\
\hline

$\mathfrak{S}_{12}^3$& 
${\bf A}_{19}^{\frac{1}{2}}$ &  $e_{1} e_{1}=e_{1}$ & $e_{1} e_{2}=e_{2}$ & $e_{3} e_{3}=e_{3}$ & \\
\hline
$\mathfrak{S}_{13}^3$& 
${\bf A}_{19}^{-\frac{1}{2}}$ &  $e_{1} e_{1}=e_{1}$ & $e_{2} e_{1}=e_{2}$ & $e_{3} e_{3}=e_{3}$& \\
\hline
 
$\mathfrak{S}_{14}^3(\alpha)$&${\bf A}_{02}^{\alpha }$ &   
$e_1  e_2=\left( 1+\alpha \right) e_3$ & $e_2 e_1=\left(1-\alpha \right) e_3$ &&\\
\hline 

\end{longtable}
}

\subsubsection{4-dimensional nilpotent standard   algebras}
Thanks to \cite{aak24},
the varieties of $4$-dimensional nilpotent standard and noncommutative Jordan algebras coincide. Hence,  $\mathfrak{NS}^4$ has five irreducible components:
\[{\rm Irr}(\mathfrak{NS}^4)=
\left\{\overline{{\mathcal O}\left(\mathfrak{NS}_{i}^4\right)}\right\}_{i=1}^{3}\cup\left\{\overline{\bigcup {\mathcal O}\left(\mathfrak{NS}_{i}^4(\alpha)\right)}\right\}_{i=4}^{5},\]
where

{\small
\vspace{-5mm}
\begin{longtable}{|l|l|lllll|}
\caption*{ }   \\

\hline \multicolumn{1}{|c|}{${\mathcal A}$} & \multicolumn{1}{c|}{ } & \multicolumn{5}{c|}{Multiplication table} \\ \hline 
\endfirsthead

 \multicolumn{7}{l}%
{{\bfseries  continued from previous page}} \\
\hline \multicolumn{1}{|c|}{${\mathcal A}$} & \multicolumn{1}{c|}{ } & \multicolumn{5}{c|}{Multiplication table} \\ \hline 
\endhead

\hline \multicolumn{7}{|r|}{{Continued on next page}} \\ \hline
\endfoot

\hline 
\endlastfoot

$\mathfrak{NS}_{1}^4$ & ${\mathcal J}^4_{07}$ & $e_1 e_2 = e_3$ & $e_1 e_3=e_4$  & $e_2 e_1=e_3+e_4$& & \\&& $e_2 e_3=e_4$ & $e_3 e_1=e_4$ & $e_3 e_2=e_4$ &&\\ \hline
$\mathfrak{NS}_{1}^4$ & $ {\mathcal J}^4_{17}$ & $e_1 e_1 = e_4$ & $e_1 e_2 = e_3$ & $e_2 e_1=-e_3$ && \\&& $e_1 e_3=e_4$ & $e_2 e_2=e_4$ & $e_3 e_1=-e_4$ &&\\ \hline
$\mathfrak{NS}_{1}^4$ & ${\mathcal J}^4_{18}$ & $e_1 e_1 = e_2$ & $e_1 e_2 = e_3$ & $e_1 e_3=e_4$ && \\&& $e_2 e_1=e_3$ & $e_2 e_2=e_4$  & $e_3 e_1=e_4$ &&\\ \hline
$\mathfrak{NS}_{4}^4(\alpha)$ & $\mathfrak{N}_2(\alpha)$ & $e_1e_1 = e_3$ & $e_1e_2 = e_4$  &$e_2e_1 = -\alpha e_3$ & $e_2e_2 = -e_4$ & \\ \hline
$\mathfrak{NS}_{5}^4(\alpha)$ & $\mathfrak{N}_3(\alpha)$ & $e_1e_1 = e_4$ & $e_1e_2 = \alpha e_4$  & $e_2e_1 = -\alpha e_4$ & $e_2e_2 = e_4$ & \\ \hline
\end{longtable}
}

 \subsection{Noncommutative Jordan algebras} 
 
An algebra $\mathfrak{J}$ is called a {\it noncommutative Jordan} algebra if it satisfies the identities 
\[(xy)x=x(yz), \  x^2(yx)=(x^2y)x.\] 
Let $\mathfrak{NCJord}$ be the variety of noncommutative  Jordan algebras.

\subsubsection{2-dimensional noncommutative Jordan algebras}
The algebraic and geometric classification of $2$-dimensional noncommutative Jordan algebras can be found in \cite{jkk19}.  In particular, it is proven that the variety  $\mathfrak{NCJord}^2$ has two irreducible components:
\[{\rm Irr}(\mathfrak{NCJord}^2)=\left\{
\overline{{\mathcal O}\big(\mathfrak{NCJ}_{1}^2\big)} \right\} \cup 
\left\{\overline{\bigcup {\mathcal O}\big(\mathfrak{NCJ}_{2}^2(\alpha)\big)}\right\} ,\]
where the algebras $\mathfrak{NCJ}_{1}^2$ and $\mathfrak{NCJ}_{2}^2(\alpha)$ are defined as follows:

{\small
\begin{longtable}{|l|l|lll|}
\caption*{ }   \\

\hline \multicolumn{1}{|c|}{${\mathcal A}$} & \multicolumn{1}{c|}{ } & \multicolumn{3}{c|}{Multiplication table} \\ \hline 
\endfirsthead

 \multicolumn{5}{l}%
{{\bfseries  continued from previous page}} \\
\hline \multicolumn{1}{|c|}{${\mathcal A}$} & \multicolumn{1}{c|}{ } & \multicolumn{3}{c|}{Multiplication table} \\ \hline 
\endhead

\hline \multicolumn{5}{|r|}{{Continued on next page}} \\ \hline
\endfoot

\hline 
\endlastfoot

$\mathfrak{NCJ}_{1}^2$ & ${\bf E}_{1}(0,0,0,0)$ & $e_1 e_1 = e_1$ & $e_2 e_2=e_2$ & \\ \hline
$\mathfrak{NCJ}_{2}^2(\alpha)$ & ${\bf E}_{5}(\alpha)$ & $e_1 e_1 = e_1$ & $e_1 e_2=(1-\alpha)e_1+\alpha e_2$ & \\&& $e_2 e_1=\alpha e_1+(1-\alpha) e_2$ &  $e_2 e_2 = e_2$ &   \\ \hline 
\end{longtable}
}

\subsubsection{3-dimensional nilpotent noncommutative Jordan algebras}
The algebraic and geometric classification of $3$-dimensional nilpotent noncommutative Jordan algebras can be found in \cite{jkk19}.  In particular, it is proven that the variety  $\mathfrak{NNCJord}^3$ has two irreducible components:
\[{\rm Irr}(\mathfrak{NNCJord}^3)=\left\{
\overline{{\mathcal O}\big(\mathfrak{NNCJ}_{1}^3\big)}\right\} \cup \left\{
\overline{\bigcup {\mathcal O}\big(\mathfrak{NNCJ}_{2}^3(\alpha)\big)}\right\} ,\]
where the algebras $\mathfrak{NNCJ}_{1}^3$ and $\mathfrak{NNCJ}_{2}^3(\alpha)$ are defined as follows:

{\small
\begin{longtable}{|l|l|lll|}
\caption*{ }   \\

\hline \multicolumn{1}{|c|}{${\mathcal A}$} & \multicolumn{1}{c|}{ } & \multicolumn{3}{c|}{Multiplication table} \\ \hline 
\endfirsthead

 \multicolumn{5}{l}%
{{\bfseries  continued from previous page}} \\
\hline \multicolumn{1}{|c|}{${\mathcal A}$} & \multicolumn{1}{c|}{ } & \multicolumn{3}{c|}{Multiplication table} \\ \hline 
\endhead

\hline \multicolumn{5}{|r|}{{Continued on next page}} \\ \hline
\endfoot

\hline 
\endlastfoot

$\mathfrak{NNCJ}_{1}^3$ & ${\mathcal J}^3_{01}$ & $e_1 e_1 = e_2$ & $e_1 e_2=e_3$ & $e_2 e_1= e_3$ \\ \hline
$\mathfrak{NNCJ}_{2}^3(\alpha)$ & ${\mathcal J}^{3*}_{04}(\alpha)$ & $e_1 e_1 = \alpha e_3$ & $e_2 e_1=e_3$ & $e_2 e_2=e_3$  \\ \hline 
\end{longtable}
}

\subsubsection{3-dimensional   noncommutative Jordan algebras}
The algebraic and geometric classification of $3$-dimensional   noncommutative Jordan algebras can be found in \cite{aak24}.  In particular, it is proven that the variety  $\mathfrak{NCJord}^3$ has eight irreducible components:
\[{\rm Irr}(\mathfrak{NCJord}^3)=\left\{
\overline{{\mathcal O}\big(\mathfrak{NCJ}_{i}^3\big)}\right\}_{i=1}^5 \cup \left\{
\overline{\bigcup {\mathcal O}\big(\mathfrak{NCJ}_{i}^3(\alpha)\big)} \right\}_{i=6}^8 ,\]
where  

{\small
\begin{longtable}{|l|l|llll|}
\caption*{ }   \\

\hline \multicolumn{1}{|c|}{${\mathcal A}$} & \multicolumn{1}{c|}{ } & \multicolumn{4}{c|}{Multiplication table} \\ \hline 
\endfirsthead

 \multicolumn{6}{l}%
{{\bfseries  continued from previous page}} \\
\hline \multicolumn{1}{|c|}{${\mathcal A}$} & \multicolumn{1}{c|}{ } & \multicolumn{4}{c|}{Multiplication table} \\ \hline 
\endhead

\hline \multicolumn{6}{|r|}{{Continued on next page}} \\ \hline
\endfoot

\hline 
\endlastfoot
$\mathfrak{NCJ}_{1}^3$& ${\bf A}_{04} $ &   $e_{1} e_{1}=e_{1}$ & $ e_{2} e_{2}=e_{2}$ && \\ 
\hline

$\mathfrak{NCJ}_{2}^3$& ${\bf A}_{12}$ &  
$e_{1} e_{1}=e_{1}$ & $e_{1} e_{3}= \frac{1}{2}e_{3}$ & $e_{2} e_{2}=e_{2}$ & 
  $e_{2} e_{3}=\frac{1}{2}e_{3}$ \\
&&
$e_{3} e_{1}= \frac{1}{2}e_{3}$ & $ e_{3} e_{3}=e_{1}+e_{2}$ & 
  $e_{3} e_{2}=\frac{1}{2}e_{3}$ &\\
\hline

$\mathfrak{NCJ}_{3}^3$& ${\bf A}_{16}$ &   
$e_{1} e_{1}=e_{1}$ & $e_{1} e_{2}= \frac{1}{2}e_{2}$ & $e_{1} e_{3}=e_{3}$   & \\&&$e_{2} e_{1}= \frac{1}{2}e_{2}$ & $e_{2}  e_{2}=e_{3}$ & $e_{3} e_{1}=e_{3}$&\\
\hline

$\mathfrak{NCJ}_{4}^3$&${\bf A}_{30}$ &   
$e_{1} e_{1}=e_{1}$ & $e_{2} e_{3}=e_{3}$ & $e_{3} e_{2}=-e_{3}$ &\\
\hline 

$\mathfrak{NCJ}_{5}^3$& ${\bf A}_{32}$ &  
$e_{1} e_{1}=e_{1}$ & $e_{1}  e_{2}=\frac{1}{2}e_{2}+e_{3}$ & $e_{1} e_{3}=\frac{1%
}{2}e_{3}$   & $e_{2} e_{1}=\frac{1}{2}e_{2}-e_{3}$    \\
&& $e_{2}  e_{3}=e_{2}$ &  $e_{3} e_{1}=\frac{1}{2}e_{3}$   & $e_{3} e_{2}=-e_{2}$& \\

\hline 

$\mathfrak{NCJ}_{6}^3(\alpha)$& ${\bf A}_{17}^{\alpha}$ &  
$e_{1} e_{1}=e_{1}$  & $e_{1} e_{3}=(\frac{1}{2}+\alpha )e_{3}$ & $e_{2} e_{2}=e_{2}$ & \\
&& $e_{2} e_{3}=(\frac{1}{2}-\alpha
)e_{3}$& $e_{3} e_{1}=(
\frac{1}{2}-\alpha )e_{3}$   &  $e_{3} e_{2}=(\frac{1}{2}+\alpha )e_{3}$ & \\
\hline 
$\mathfrak{NCJ}_{7}^3(\alpha)$&${\bf A}_{19}^{\alpha }$ &   
$e_{1} e_{1}=e_{1}$ & $e_{1} e_{3}=(\frac{1}{2}+\alpha )e_{3}$ & $e_{2} e_{2}=e_{2}$  & $e_{3} e_{1}=(\frac{1}{2}-\alpha )e_{3}$   \\
\hline 
$\mathfrak{NCJ}_{8}^3(\alpha)$&${\bf A}_{24}^{\alpha }$ &   
$  e_{1}e_{2}  =e_{3}$ & 
$e_{1}e_{3}  =e_{1}+e_{3}$ &
$e_{2}e_{1}   =-e_{3}$ &\\
&&$ e_{2}e_{3} =\alpha e_{2}$ &
$e_{3}e_{1}  =-e_{1}-e_{3}$
& $ e_{3}e_{2} =-\alpha e_{2}$&\\
\hline 

\end{longtable}
}

\subsubsection{4-dimensional nilpotent noncommutative Jordan algebras}

The variety $\mathfrak{NNCJord}^4$  was studied in 
\cite{jkk19} both algebraically, and geometrically. It has five irreducible components:
\[{\rm Irr}(\mathfrak{NNCJord}^4)=
\left\{\overline{{\mathcal O}\left(\mathfrak{NNCJ}_{i}^4\right)}\right\}_{i=1}^{3}\cup\left\{\overline{\bigcup {\mathcal O}\left(\mathfrak{NNCJ}_{i}^4(\alpha)\right)}\right\}_{i=4}^{5},\]
where

{\small
\vspace{-5mm}
\begin{longtable}{|l|l|lllll|}
\caption*{ }   \\

\hline \multicolumn{1}{|c|}{${\mathcal A}$} & \multicolumn{1}{c|}{ } & \multicolumn{5}{c|}{Multiplication table} \\ \hline 
\endfirsthead

 \multicolumn{7}{l}%
{{\bfseries  continued from previous page}} \\
\hline \multicolumn{1}{|c|}{${\mathcal A}$} & \multicolumn{1}{c|}{ } & \multicolumn{5}{c|}{Multiplication table} \\ \hline 
\endhead

\hline \multicolumn{7}{|r|}{{Continued on next page}} \\ \hline
\endfoot

\hline 
\endlastfoot

$\mathfrak{NNCJ}_{1}^4$ & ${\mathcal J}^4_{07}$ & $e_1 e_2 = e_3$ & $e_1 e_3=e_4$  & $e_2 e_1=e_3+e_4$& & \\&& $e_2 e_3=e_4$ & $e_3 e_1=e_4$ & $e_3 e_2=e_4$ &&\\ \hline
$\mathfrak{NNCJ}_{2}^4$ & $ {\mathcal J}^4_{17}$ & $e_1 e_1 = e_4$ & $e_1 e_2 = e_3$ & $e_2 e_1=-e_3$ && \\&& $e_1 e_3=e_4$ & $e_2 e_2=e_4$ & $e_3 e_1=-e_4$ &&\\ \hline
$\mathfrak{NNCJ}_{3}^4$ & ${\mathcal J}^4_{18}$ & $e_1 e_1 = e_2$ & $e_1 e_2 = e_3$ & $e_1 e_3=e_4$ && \\&& $e_2 e_1=e_3$ & $e_2 e_2=e_4$  & $e_3 e_1=e_4$ &&\\ \hline
$\mathfrak{NNCJ}_{4}^4(\alpha)$ & $\mathfrak{N}_2(\alpha)$ & $e_1e_1 = e_3$ & $e_1e_2 = e_4$  &$e_2e_1 = -\alpha e_3$ & $e_2e_2 = -e_4$ & \\ \hline
$\mathfrak{NNCJ}_{5}^4(\alpha)$ & $\mathfrak{N}_3(\alpha)$ & $e_1e_1 = e_4$ & $e_1e_2 = \alpha e_4$  & $e_2e_1 = -\alpha e_4$ & $e_2e_2 = e_4$ & \\ \hline
\end{longtable}
}

 {

\subsection{Commutative power-associative algebras}
 
A commutative algebra $\mathfrak{CPA}$ is called  
{\it commutative power-associative } if it satisfies the identity 
\[ x^2x^2=(x^2x)x. \]
 We will denote the variety of commutative power-associative algebras by $\mathfrak{CPA}$.
 
\subsubsection{2-dimensional commutative power-asso\-cia\-tive algebras}
The variety of Jordan algebras is a proper subvariety of the variety of commutative power-asso\-cia\-tive algebras.
In dimension $2$, these two varieties coincide \cite{rpe20}. 
The graph of degenerations can be deduced from~\cite{kv16} and is explicitly given in~\cite{gkp}. In particular, it is proven that the variety  $\mathfrak{CPA}^2$ has two irreducible components:

{\small
\vspace{-5mm}
\begin{longtable}{|l|l|ll|}
\caption*{ }   \\

\hline \multicolumn{1}{|c|}{${\mathcal A}$} & \multicolumn{1}{c|}{ } & \multicolumn{2}{c|}{Multiplication table} \\ \hline 
\endfirsthead

 \multicolumn{4}{l}%
{{\bfseries  continued from previous page}} \\
\hline \multicolumn{1}{|c|}{${\mathcal A}$} & \multicolumn{1}{c|}{ } & \multicolumn{2}{c|}{Multiplication table} \\ \hline 
\endhead

\hline \multicolumn{4}{|r|}{{Continued on next page}} \\ \hline
\endfoot

\hline 
\endlastfoot

$\mathfrak{CPA}_{1}^2$ & $\mathfrak{B}_2$ & $e_1e_1 = e_1$ & $e_1e_2 =\frac{1}{2}e_2$  \\ \hline
$\mathfrak{CPA}_{2}^2$ & $\mathfrak{B}_4$ & $e_1e_1 = e_1$ & $e_2e_2 = e_2$  \\ \hline
\end{longtable}
}

\subsubsection{3-dimensional nilpotent  commutative power-associative algebras}
The nilpotent commutative power-associative algebras of dimension $3$  
are coincides with the nilpotent Jordan algebras \cite{rpe20},  
that were classified algebraically and geometrically in~\cite{esp}. 
The variety $\mathfrak{NCPA}^3$ is irreducible, determined by the rigid algebra

{\small

\begin{longtable}{|l|l|ll|}
\caption*{ }   \\

\hline \multicolumn{1}{|c|}{${\mathcal A}$} & \multicolumn{1}{c|}{ } & \multicolumn{2}{c|}{Multiplication table} \\ \hline 
\endfirsthead

 \multicolumn{4}{l}%
{{\bfseries  continued from previous page}} \\
\hline \multicolumn{1}{|c|}{${\mathcal A}$} & \multicolumn{1}{c|}{ } & \multicolumn{2}{c|}{Multiplication table} \\ \hline 
\endhead

\hline \multicolumn{4}{|r|}{{Continued on next page}} \\ \hline
\endfoot

\hline 
\endlastfoot

$\mathfrak{NCPA}_{1}^3$ & $\phi_1$ & $e_1e_1 = e_2$ &  $e_1e_2 = e_3$  \\ \hline

\end{longtable}
}

\subsubsection{3-dimensional commutative power-associative algebras}
The variety of Jordan algebras is a proper subvariety of the variety of commutative power-asso\-cia\-tive algebras.
In dimension $3$, these two varieties coincide \cite{rpe20}. 
Hence, in~\cite{gkp}, the geometric classification in dimension $3$ (initiated in~\cite{KS07} together with the algebraic classification) was completed. We will refer to the notation of~\cite{gkp}. There exist five irreducible components in $\mathfrak{CPA}^3$:
\[{\rm Irr}(\mathfrak{CPA}^3)=\left\{
\overline{{\mathcal O}\big(\mathfrak{CPA}_{i}^3\big)}\right\}_{i=1}^{5},\]
where 

{\small

\begin{longtable}{|l|l|lllll|}
\caption*{ }   \\

\hline \multicolumn{1}{|c|}{${\mathcal A}$} & \multicolumn{1}{c|}{ } & \multicolumn{5}{c|}{Multiplication table} \\ \hline 
\endfirsthead

\multicolumn{7}{l}%
{{\bfseries  continued from previous page}} \\
\hline \multicolumn{1}{|c|}{${\mathcal A}$} & \multicolumn{1}{c|}{ } & \multicolumn{5}{c|}{Multiplication table} \\ \hline 
\endhead

\hline \multicolumn{7}{|r|}{{Continued on next page}} \\ \hline
\endfoot

\hline 
\endlastfoot

$\mathfrak{CPA}_{1}^3$ & $\mathbb{T}_{01}$ & $e_1e_1 = e_1$ & $e_2e_2 = e_2$ & $e_3e_3 = e_3$ & & \\ \hline
$\mathfrak{CPA}_{2}^3$ & $\mathbb{T}_{02}$ & $e_1e_1 = e_1$ & $e_1e_3 = \frac{1}{2}e_3$ & $e_2e_2 =e_2$ & $e_2e_3 = \frac{1}{2}e_3$ & $e_3e_3 = e_1 + e_2$  \\ \hline
$\mathfrak{CPA}_{3}^3$ & $\mathbb{T}_{05}$ & $e_1e_1 = e_1$ & $e_1e_3 = \frac{1}{2}e_3$ & $e_2e_2 = e_2$ & &  \\ \hline
$\mathfrak{CPA}_{4}^3$ & $\mathbb{T}_{10}$ & $e_1e_1 = e_1$ & $e_1e_2 = \frac{1}{2}e_2$ & $e_1e_3 = e_3$ & $e_2e_2 = e_3$ & \\ \hline
$\mathfrak{CPA}_{5}^3$ & $\mathbb{T}_{12}$ & $e_1e_1 = e_1$ & $e_1e_2 = \frac{1}{2}e_2$ & $e_1e_3 = \frac{1}{2}e_3$ && \\ \hline

\end{longtable}
}

The complete graph of degenerations can be found in~\cite{gkp}.

\subsubsection{4-dimensional nilpotent commutative power-associative algebras}
The variety of nilpotent Jordan algebras is a proper subvariety of the variety of nilpotent commutative power-associative algebras.
In dimension $4$, these two varieties coincide \cite{rpe20}. 
The variety $\mathfrak{NCPA}^4$  was studied in 
\cite{esp} both algebraically and geometrically. It has two irreducible components:
\[{\rm Irr}(\mathfrak{NCPA}^4)=
\left\{\overline{{\mathcal O}\left(\mathfrak{NCPA}_{i}^4\right)}\right\}_{i=1}^{2},\]
where

{\small

\begin{longtable}{|l|l|llll|}
\caption*{ }   \\

\hline \multicolumn{1}{|c|}{${\mathcal A}$} & \multicolumn{1}{c|}{ } & \multicolumn{4}{c|}{Multiplication table} \\ \hline 
\endfirsthead

 \multicolumn{6}{l}%
{{\bfseries  continued from previous page}} \\
\hline \multicolumn{1}{|c|}{${\mathcal A}$} & \multicolumn{1}{c|}{ } & \multicolumn{4}{c|}{Multiplication table} \\ \hline 
\endhead

\hline \multicolumn{6}{|r|}{{Continued on next page}} \\ \hline
\endfoot

\hline 
\endlastfoot

$\mathfrak{NCPA}_{1}^4$ & $\phi_1$ & $e_1e_1  = e_2$ & $e_1 e_2 = e_3$ &  $e_1 e_3 = e_4$ & $e_2 e_2 = e_4$  \\ \hline 
$\mathfrak{NCPA}_{2}^4$ & $\phi_2$ & $e_1 e_1 = e_3$ & $e_1 e_3 = e_4$ & $e_2 e_2 = e_3$ & \\ \hline

\end{longtable}
}

\subsubsection{4-dimensional commutative power-associative algebras}
Regarding $\mathfrak{CPA}^4$, the algebraic classification was made in \cite{rpe20}, 
and the twelve irreducible components were found  in~\cite{rpe20}:
\[{\rm Irr}(\mathfrak{CPA}^4)=\left\{\overline{{\mathcal O}\big(\mathfrak{CPA}_{i}^4\big)} \right\}_{i=1}^{12},\]
where

{\small

\begin{longtable}{|l|l|lllll|}
\caption*{ }   \\

\hline \multicolumn{1}{|c|}{${\mathcal A}$} & \multicolumn{1}{c|}{ } & \multicolumn{5}{c|}{Multiplication table} \\ \hline 
\endfirsthead

 \multicolumn{7}{l}%
{{\bfseries  continued from previous page}} \\
\hline \multicolumn{1}{|c|}{${\mathcal A}$} & \multicolumn{1}{c|}{ } & \multicolumn{5}{c|}{Multiplication table} \\ \hline 
\endhead

\hline \multicolumn{7}{|r|}{{Continued on next page}} \\ \hline
\endfoot

\hline 
\endlastfoot

$\mathfrak{CPA}_{1}^4$ & $\mathfrak{J}_{1}$ & $e_1e_1 = e_1$ & $e_1e_3 = \frac{1}{2}e_3$ & $e_2e_2 =e_2$& &
\\&&$e_2e_3 = \frac{1}{2}e_3$ & $e_3e_3 = e_1 + e_2$  & $e_4e_4=e_4$ &&  \\ \hline
$\mathfrak{CPA}_{2}^4$ & $\mathfrak{J}_{2}$ & $e_1e_3 = \frac{1}{2}e_3$ & $e_1e_4 = \frac{1}{2}e_4$ & $e_2e_3 = \frac{1}{2}e_3$ & $e_2e_4 =\frac{1}{2}e_4$ & $e_3e_4 = \frac{1}{2}(e_1 + e_2)$   \\ \hline
$\mathfrak{CPA}_{3}^4$ & $\mathfrak{J}_{3}$ & $e_1e_1=e_1$ & $e_2e_2=e_2$ & $e_3e_3=e_3$ & $e_4e_4=e_4$ & \\ \hline
$\mathfrak{CPA}_{4}^4$ & $\mathfrak{J}_{6}$ & $e_1e_1 = e_1$ & $e_1e_2 =\frac{1}{2}e_2$ & $e_3e_3=e_3$ & $e_4e_4=e_4$ &  \\ \hline
$\mathfrak{CPA}_{5}^4$ & $\mathfrak{J}_{12}$ & $e_1e_1 = e_1$ & $e_1e_2 = \frac{1}{2}e_2$ & $e_1e_3 = \frac{1}{2}e_3$ & $e_4e_4=e_4$ & \\ \hline
$\mathfrak{CPA}_{6}^4$ & $\mathfrak{J}_{13}$ & $e_1e_1 = e_1$ & $e_1e_2 =\frac{1}{2}e_2$ & $e_3e_3=e_3$ & $e_3e_4 =\frac{1}{2}e_4$ & \\ \hline
$\mathfrak{CPA}_{7}^4$ & $\mathfrak{J}_{16}$ & $e_1e_3 = \frac{1}{2}e_3$ & $e_1e_4 = \frac{1}{2}e_4$ & $e_2e_3 = \frac{1}{2}e_3$ &&   \\ \hline
$\mathfrak{CPA}_{8}^4$ & $\mathfrak{J}_{24}$ & $e_1e_1 = e_1$ & $e_1e_2 = \frac{1}{2}e_2$ & $e_1e_3 = e_3$ & $e_2e_2 = e_3$ & $e_4e_4=e_4$  \\ \hline
$\mathfrak{CPA}_{9}^4$ & $\mathfrak{J}_{33}$ & $e_1e_2 = \frac{1}{2}e_2$ & $e_1e_3 = \frac{1}{2}e_3$ & $e_1e_4 = \frac{1}{2}e_4$ && \\ \hline
$\mathfrak{CPA}_{10}^4$ & $\mathfrak{J}_{59}$ & $e_1e_2 = e_2$ & $e_1e_3 = \frac{1}{2}e_3$ & $e_1e_4= \frac{1}{2}e_4$ & $e_3e_4=e_2$ & $e_4e_4=e_2$  \\ \hline

$\mathfrak{CPA}_{11}^4$ & $A_{7}$ & $e_1e_1=e_1$ & $e_1e_2=\frac{1}{2}e_2$ &  $e_1e_3=e_3+e_4$ & 
   $e_2e_2=e_3$ & $e_2e_3=e_4$ \\ \hline

$\mathfrak{CPA}_{12}^4$ & $A_{8}$ &
$e_1e_1=e_1$ &  $e_1e_2=\frac{1}{2}e_2$ & $e_1e_3=e_3$& &\\
&&$e_1e_4=e_4$ & 
  $e_2e_2=e_3$ &  $e_2e_3=e_4$  &&\\
  \hline

\end{longtable}
}

 \subsection{Weakly associative algebras} 

An algebra $\mathfrak{WA}$ is called  {\it weakly associative} if it satisfies the identity 
\[ (xy)z-x(yz)+(yz)x-y(zx) = (yx)z- y(xz). \]
 We will denote the variety of weakly associative algebras by $\mathfrak{WA}$.
 
\subsubsection{2-dimensional weakly associative algebras}
The variety of weakly associative algebras is a proper subvariety of the variety of flexible algebras defined by the following identity $(xy)x=x(yx).$
In dimension $2$, these two varieties coincide. Therefore, the algebraic and geometric classification of $2$-dimensional weakly associative algebras follows from \cite[Section $7.1$]{kv16}.
Hence,  $\mathfrak{WA}^2$ has two irreducible components, namely:

\[{\rm Irr}(\mathfrak{WA}^2)=
\left\{ \overline{{\mathcal O}\big(\mathfrak{WA}_{1}^2(\alpha)\big)} \right\}
\cup
\left\{\overline{\bigcup {\mathcal O}\big(\mathfrak{WA}_{2}^2(\alpha,\beta)\big)}\right\},\]
where

{\small

\begin{longtable}{|l|l|llll|}
\caption*{ }   \\

\hline \multicolumn{1}{|c|}{${\mathcal A}$} & \multicolumn{1}{c|}{ } & \multicolumn{4}{c|}{Multiplication table} \\ \hline 
\endfirsthead

 \multicolumn{6}{l}%
{{\bfseries  continued from previous page}} \\
\hline \multicolumn{1}{|c|}{${\mathcal A}$} & \multicolumn{1}{c|}{ } & \multicolumn{4}{c|}{Multiplication table} \\ \hline 
\endhead

\hline \multicolumn{6}{|r|}{{Continued on next page}} \\ \hline
\endfoot

\hline 
\endlastfoot

$\mathfrak{WA}_{1}^2(\alpha)$ & ${\bf E}_5(\alpha) $& $e_1e_1=e_1$ & $ e_1e_2=(1-\alpha) e_1+ \alpha  e_2$ &&\\&& $e_2e_1=\alpha  e_1 + (1-\alpha)e_2$ & $ e_2e_2=e_2$&&\\ 

\hline 

$\mathfrak{WA}_{2}^2(\alpha,\beta)$ &  
${\bf E}_1(\alpha,\beta,\alpha,\beta)$ & 
$e_1e_1=e_1$ & $ e_1e_2=\alpha e_1+ \beta  e_2$  &&\\&& $e_2e_1=\alpha  e_1 + \beta  e_2$ & $ e_2e_2=e_2$&&\\ \hline
\end{longtable}
}

\subsubsection{3-dimensional nilpotent weakly associative   algebras}
The list of $3$-dimensional    nilpotent weakly associative algebras can be found in~\cite{ak21}.
Employing the graph of degenerations of~\cite{fkkv}, we obtain that the variety $\mathfrak{NWA}^3$ has two  irreducible  components and is defined by the following family of algebras

{\small
\vspace{-5mm}
\begin{longtable}{|l|l|lll|}
\caption*{ }   \\

\hline \multicolumn{1}{|c|}{${\mathcal A}$} & \multicolumn{1}{c|}{ } & \multicolumn{3}{c|}{Multiplication table} \\ \hline 
\endfirsthead

 \multicolumn{5}{l}%
{{\bfseries  continued from previous page}} \\
\hline \multicolumn{1}{|c|}{${\mathcal A}$} & \multicolumn{1}{c|}{ } & \multicolumn{3}{c|}{Multiplication table} \\ \hline 
\endhead

\hline \multicolumn{5}{|r|}{{Continued on next page}} \\ \hline
\endfoot

\hline 
\endlastfoot

$\mathfrak{NWA}^3_{1}$ &  ${\mathcal C}_{02}$ &  $e_1 e_1 = e_2$  & $e_2 e_2=e_3$  & \\\hline

$\mathfrak{NWA}^3_{2}(\alpha)$ &  ${\mathcal N}^{\alpha}_{02}$  &  $e_1 e_1 = \alpha e_3$   & $e_2 e_1=e_3$  & $e_2 e_2=e_3$
\end{longtable}
}

\subsubsection{4-dimensional nilpotent weakly associative algebras}
The algebraic and geometric classification of 
$4$-dimensional weakly associative  algebras are given in a paper by 
   Alvarez and       Kaygorodov  \cite{ak21}.
The variety $\mathfrak{NWA}^4$ has  five irreducible components:

\[{\rm Irr}(\mathfrak{NWA}^4)= \left\{\overline{{\mathcal O}\left(\mathfrak{NWA}_{1}^4\right)}\right\}\cup
\left\{\overline{\bigcup {\mathcal O}\left(\mathfrak{NWA}_{i}^4(\alpha)\right)} \right\}_{i=2}^3,\]
where
{\small
\vspace{-5mm}
\begin{longtable}{|l|l|lllll|}
\caption*{ }   \\

\hline \multicolumn{1}{|c|}{${\mathcal A}$} & \multicolumn{1}{c|}{ } & \multicolumn{5}{c|}{Multiplication table} \\ \hline 
\endfirsthead

 \multicolumn{7}{l}%
{{\bfseries  continued from previous page}} \\
\hline \multicolumn{1}{|c|}{${\mathcal A}$} & \multicolumn{1}{c|}{ } & \multicolumn{5}{c|}{Multiplication table} \\ \hline 
\endhead

\hline \multicolumn{7}{|r|}{{Continued on next page}} \\ \hline
\endfoot

\hline 
\endlastfoot

$\mathfrak{NWA}^4_{1}$ & 
${\mathcal S}_{01}$  &   $e_1e_1 = e_4$ & $e_1e_2 = e_3$ & $e_2e_1 = -e_3$ &&\\&& $e_2e_2 = e_4$ & $e_2e_3 = e_4$ & $e_3e_2 = -e_4$ &&\\   \hline
    
$\mathfrak{NWA}^4_{2}(\alpha)$ &$\mathfrak{W}_{06}^{\alpha}$ &  $e_1e_1=\alpha e_4$ & $e_1 e_2 =e_3+e_4$ & $e_2 e_1 =e_3$ & $e_2e_2=e_4$ & $e_3e_3 =e_4$\\
 \hline
   
$\mathfrak{NWA}^4_{3}(\alpha)$ & $ \mathfrak{C}_{19}(\alpha) $ & $  e_1 e_1 = e_2$ &$  e_1e_3=\alpha e_4$  &$  e_2 e_2=e_3$  &$ e_2e_3= e_4$  &$  e_3e_3=e_4$  \\ \hline

\end{longtable}
}
}

 \subsection{Terminal algebras} 

An algebra $\mathfrak{T}$ is called  {\it terminal} if it satisfies the identity 
{\footnotesize\begin{align*}
    &b(a(xy)-(ax)y-x(ay)) - a((bx)y) + (a(bx))y+(bx)(ay) -a(x(by))+(ax)(by)+x(a(by)) \\
&\quad= - \left(\frac 23ab+\frac 13 ba\right)(xy)+\left(\left(\frac 23ab+\frac 13 ba\right)x\right)y + x\left(\left(\frac 23ab+\frac 13 ba\right)y\right).
\end{align*}}
Note that there exists a simpler definition in terms of the product of bilinear maps. We will denote the variety of terminal algebras by $\mathfrak{Ter}$.

\subsubsection{2-dimensional terminal algebras}

The complete list of algebras and the graph of degenerations of the variety $\mathfrak{Ter}^2$ were constructed in~\cite{cfk19}. Basing on the general classification of~\cite{kv16}, 
it was proven that $\mathfrak{Ter}^2$ has four irreducible components, namely:

\[{\rm Irr}(\mathfrak{Ter}^2)=
\left\{ \overline{{\mathcal O}\big(\mathfrak{T}_{1}^2\big)} \right\}
\cup
\left\{\overline{\bigcup {\mathcal O}\big(\mathfrak{T}_{i}^2(\alpha)\big)}\right\}_{i=2}^{4},\]
where

{\small
\vspace{-5mm}
\begin{longtable}{|l|l|llll|}
\caption*{ }   \\

\hline \multicolumn{1}{|c|}{${\mathcal A}$} & \multicolumn{1}{c|}{ } & \multicolumn{4}{c|}{Multiplication table} \\ \hline 
\endfirsthead

 \multicolumn{6}{l}%
{{\bfseries  continued from previous page}} \\
\hline \multicolumn{1}{|c|}{${\mathcal A}$} & \multicolumn{1}{c|}{ } & \multicolumn{4}{c|}{Multiplication table} \\ \hline 
\endhead

\hline \multicolumn{6}{|r|}{{Continued on next page}} \\ \hline
\endfoot

\hline 
\endlastfoot

$\mathfrak{T}_{1}^2$ & $\mathfrak{T}_{09}$ & $e_1e_1=e_1$ &  $e_2e_2=e_2$ && \\ \hline
$\mathfrak{T}_{2}^2(\alpha)$ & $\mathfrak{T}_{07}(\alpha)$ & $e_1e_1=e_1$ &  $e_1e_2=\alpha e_2$ && \\ \hline
$\mathfrak{T}_{3}^2(\alpha)$ & $\mathfrak{T}_{08}(\alpha)$ & $e_1e_1=e_1$ &  $e_1e_2=\alpha e_2$ & $e_2e_1=(3-2\alpha)e_2$ &  \\ \hline
$\mathfrak{T}_{4}^2(\alpha)$ & $\mathfrak{T}_{10}(\alpha)$ & $e_1e_1=e_1$ &  $e_1e_2=(1-\alpha) e_1+ \alpha  e_2$ & $e_2e_1=\alpha  e_1 + (1-\alpha)e_2$ &  $e_2e_2=e_2$  \\ \hline
\end{longtable}
}

\subsubsection{3-dimensional nilpotent terminal algebras}
The list of $3$-dimensional nilpotent terminal algebras appeared in~\cite{kkp19geo}.
Their degenerations can be extracted from the classification of all the  nilpotent algebras of dimension $3$ (see~\cite{fkkv}). In the variety $\mathfrak{NTer}^3$, there are two irreducible components:
\[{\rm Irr}(\mathfrak{NTer}^3)= 
\left\{\overline{{\mathcal O}\left(\mathfrak{T}_{1}^3\right)}\right\}\cup
\left\{\overline{\bigcup {\mathcal O}\left(\mathfrak{T}_{2}^3(\alpha)\right)} \right\},\]
where

{\small
\vspace{-5mm}
\begin{longtable}{|l|l|lll|}
\caption*{ }   \\

\hline \multicolumn{1}{|c|}{${\mathcal A}$} & \multicolumn{1}{c|}{ } & \multicolumn{3}{c|}{Multiplication table} \\ \hline 
\endfirsthead

 \multicolumn{5}{l}%
{{\bfseries  continued from previous page}} \\
\hline \multicolumn{1}{|c|}{${\mathcal A}$} & \multicolumn{1}{c|}{ } & \multicolumn{3}{c|}{Multiplication table} \\ \hline 
\endhead

\hline \multicolumn{5}{|r|}{{Continued on next page}} \\ \hline
\endfoot

\hline 
\endlastfoot

$\mathfrak{NT}^3_{1}$ & $\mathfrak{N}_3$ & $e_{1}e_{1}=e_{2}$ & $e_{2}e_{1}=e_{3}$ & \\ \hline
$\mathfrak{NT}^3_{2}(\alpha)$ & $\mathfrak{N}_4(\alpha)$ & $e_{1}e_{1}=e_{2}$ & $e_{1}e_{2}=e_{3}$ & $e_{2}e_{1}=\alpha e_{3}$  \\ \hline

\end{longtable}
}





\subsubsection{4-dimensional nilpotent terminal algebras}
The algebraic and geometric classifications of $4$-dimensional nilpotent terminal algebras were obtained in~\cite{kkp19geo}. 
It was shown that the variety $\mathfrak{NTer}^4$ has three irreducible components, namely
\[{\rm Irr}\left(\mathfrak{NTer}^4\right)=\left\{
\overline{\bigcup {\mathcal O}\big(\mathfrak{NT}_{1}^4(\alpha)\big)}\right\}
\cup\left\{\overline{\bigcup {\mathcal O}\big(\mathfrak{NT}_{2}^4(\alpha,\beta)\big)}\right\}
\cup\left\{\overline{\bigcup {\mathcal O}\big(\mathfrak{NT}_{3}^4(\alpha,\beta,\gamma)\big)}
\right\},\]
where

{\small
\vspace{-5mm}
\begin{longtable}{|l|l|llll|}
\caption*{ }   \\

\hline \multicolumn{1}{|c|}{${\mathcal A}$} & \multicolumn{1}{c|}{ } & \multicolumn{4}{c|}{Multiplication table} \\ \hline 
\endfirsthead

 \multicolumn{6}{l}%
{{\bfseries  continued from previous page}} \\
\hline \multicolumn{1}{|c|}{${\mathcal A}$} & \multicolumn{1}{c|}{ } & \multicolumn{4}{c|}{Multiplication table} \\ \hline 
\endhead

\hline \multicolumn{6}{|r|}{{Continued on next page}} \\ \hline
\endfoot

\hline 
\endlastfoot

$\mathfrak{NT}_{1}^4(\alpha)$ & $\mathbf{T}^4_{41}(\alpha)$ & $e_1 e_1 = e_2$ & $e_1 e_2=\alpha e_4$ & $e_1e_3=-e_4$&\\&& $e_2 e_1=e_3$ & $e_2e_3 = e_4$ & $e_3e_1=3e_4$& \\ \hline
$\mathfrak{NT}_{2}^4(\alpha,\beta)$ & $\mathbf{T}^4_{43}(\alpha,\beta)$ & $e_1 e_1 = e_2$ & $e_1 e_2=\alpha e_3$ & $e_1e_3 = (\beta(\alpha-1)+1)e_4$ &\\&& $e_2 e_1=e_3$  & $e_2e_2=e_4$& $e_3e_1 = 3\beta e_4$& \\ \hline

$\mathfrak{NT}_{3}^4(\alpha,\beta,\gamma)$ & $\mathbf{D}^4_{01}(\alpha,\beta,\gamma)$ & $e_1 e_1 = \alpha e_3 + e_4$ & $e_1 e_3 = \beta e_4$ & $e_2 e_1=e_3$ &\\&& $e_2 e_2 = e_3$  & $e_2 e_3 = \gamma e_4$ & $e_3e_1 = e_4$&  \\ \hline
\end{longtable}
}

\subsection{Lie algebras}

An anticommutative algebra $\mathfrak{L}$ is called a {\it Lie} algebra  if it satisfies the identity 
\[J(x,y,z)=0.\] Let  $\mathfrak{Lie}$ be the variety of Lie algebras.

\subsubsection{2-dimensional Lie algebras}\label{2dlie}
For the variety $\mathfrak{Lie}^2$, we rely on the classification from~\cite{kv16}. It is determined by the rigid algebra:

{\small

\begin{longtable}{|l|l|ll|}
\caption*{ }   \\

\hline \multicolumn{1}{|c|}{${\mathcal A}$} & \multicolumn{1}{c|}{ } & \multicolumn{2}{c|}{Multiplication table} \\ \hline 
\endfirsthead

 \multicolumn{4}{l}%
{{\bfseries  continued from previous page}} \\
\hline \multicolumn{1}{|c|}{${\mathcal A}$} & \multicolumn{1}{c|}{ } & \multicolumn{2}{c|}{Multiplication table} \\ \hline 
\endhead

\hline \multicolumn{4}{|r|}{{Continued on next page}} \\ \hline
\endfoot

\hline 
\endlastfoot

$\mathfrak{L}_{1}^2$ & $\mathfrak{B}_3$ & $e_1e_2=e_2$ & \\ \hline

\end{longtable}
}

\subsubsection{3-dimensional nilpotent Lie algebras}
For the variety $\mathfrak{NLie}^3$, we  rely on the classification from~\cite{fkkv}. It is determined by the rigid algebra:

{\small

\begin{longtable}{|l|l|ll|}
\caption*{ }   \\

\hline \multicolumn{1}{|c|}{${\mathcal A}$} & \multicolumn{1}{c|}{ } & \multicolumn{2}{c|}{Multiplication table} \\ \hline 
\endfirsthead

 \multicolumn{4}{l}%
{{\bfseries  continued from previous page}} \\
\hline \multicolumn{1}{|c|}{${\mathcal A}$} & \multicolumn{1}{c|}{ } & \multicolumn{2}{c|}{Multiplication table} \\ \hline 
\endhead

\hline \multicolumn{4}{|r|}{{Continued on next page}} \\ \hline
\endfoot

\hline 
\endlastfoot

$\mathfrak{NL}_{1}^3$ & $\mathcal{N}_7$ & $e_1e_2=e_3$ & \\ \hline

\end{longtable}
}

\subsubsection{3-dimensional Lie algebras}
The classification of the Lie algebras of dimension $3$ is well-known and can be found in numerous books; for example, {\it Lie algebras} of Jacobson (1962). The four irreducible components of $\mathfrak{Lie}^3$ were found in~\cite{BC99}:
\[{\rm Irr}\left(\mathfrak{Lie}^3\right)=\left\{
\overline{{\mathcal O}\big(\mathfrak{L}_{i}^3\big)}\right\}_{i=1}^{3} \cup 
\left\{\overline{\bigcup {\mathcal O}\big(\mathfrak{L}_{4}^3(\alpha)\big)}\right\},\]
where

{\small

\begin{longtable}{|l|l|lll|}
\caption*{ }   \\

\hline \multicolumn{1}{|c|}{${\mathcal A}$} & \multicolumn{1}{c|}{ } & \multicolumn{3}{c|}{Multiplication table} \\ \hline 
\endfirsthead

 \multicolumn{5}{l}%
{{\bfseries  continued from previous page}} \\
\hline \multicolumn{1}{|c|}{${\mathcal A}$} & \multicolumn{1}{c|}{ } & \multicolumn{3}{c|}{Multiplication table} \\ \hline 
\endhead

\hline \multicolumn{5}{|r|}{{Continued on next page}} \\ \hline
\endfoot

\hline 
\endlastfoot

$\mathfrak{L}_{1}^3$ & $\mathfrak{r}_2(\mathbb{C})\oplus \mathbb{C}$ & $e_1e_2=e_1$ &&\\ \hline
$\mathfrak{L}_{2}^3$ & $\mathfrak{r}_3(\mathbb{C})$ & $e_1e_2=e_2$ & $e_1e_3=e_2+e_3$ & \\ \hline
$\mathfrak{L}_{3}^3$ & $\mathfrak{sl}_2(\mathbb{C})$ & $e_1e_2=e_3$ & $e_1e_3=-2e_1$ & $e_2e_3= 2e_2$  \\ \hline
$\mathfrak{L}_{4}^3(\alpha)$ & $\mathfrak{r}_{3,\alpha}(\mathbb{C})$ & $e_1e_2=e_2$ & $e_1e_3=\alpha e_3$ & \\ \hline
\end{longtable}
}

\subsubsection{4-dimensional nilpotent Lie algebras}
For the variety $\mathfrak{NLie}^4$, we  rely on the classification from~\cite{fkkv}. It is determined by the rigid algebra:

{\small

\begin{longtable}{|l|l|ll|}
\caption*{ }   \\

\hline \multicolumn{1}{|c|}{${\mathcal A}$} & \multicolumn{1}{c|}{ } & \multicolumn{2}{c|}{Multiplication table} \\ \hline 
\endfirsthead

 \multicolumn{4}{l}%
{{\bfseries  continued from previous page}} \\
\hline \multicolumn{1}{|c|}{${\mathcal A}$} & \multicolumn{1}{c|}{ } & \multicolumn{2}{c|}{Multiplication table} \\ \hline 
\endhead

\hline \multicolumn{4}{|r|}{{Continued on next page}} \\ \hline
\endfoot

\hline 
\endlastfoot

$\mathfrak{NL}_{1}^4$ & $\mathcal{A}_2$ & $e_1e_2=e_3$ & $e_1e_3=e_4$  \\ \hline

\end{longtable}
}

\subsubsection{4-dimensional Lie algebras}
Lie algebras of dimension $4$ were classified up to isomorphism by Steinhoff in 1997, and their degenerations were studied in~\cite{BC99}. The variety $\mathfrak{Lie}^4$ has seven irreducible components:
\[{\rm Irr}\left(\mathfrak{Lie}^4\right)=\left\{
\overline{{\mathcal O}\left(\mathfrak{L}_{i}^4\right)}\right\}_{i=1}^{3}\cup
\left\{\overline{\bigcup {\mathcal O}\left(\mathfrak{L}_{i}^4(\alpha)\right)}\right\}_{i=4}^{6}\cup
\left\{\overline{\bigcup {\mathcal O}\left(\mathfrak{L}_{7}^4(\alpha,\beta)\right)} \right\},\]
where

{\small

\begin{longtable}{|l|l|llll|}
\caption*{ }   \\

\hline \multicolumn{1}{|c|}{${\mathcal A}$} & \multicolumn{1}{c|}{ } & \multicolumn{4}{c|}{Multiplication table} \\ \hline 
\endfirsthead

 \multicolumn{6}{l}%
{{\bfseries  continued from previous page}} \\
\hline \multicolumn{1}{|c|}{${\mathcal A}$} & \multicolumn{1}{c|}{ } & \multicolumn{4}{c|}{Multiplication table} \\ \hline 
\endhead

\hline \multicolumn{6}{|r|}{{Continued on next page}} \\ \hline
\endfoot

\hline 
\endlastfoot

$\mathfrak{L}_{1}^4$ & $\mathfrak{r}_2(\mathbb{C})\oplus \mathfrak{r}_2(\mathbb{C})$ & $e_1e_2=e_1$ & $e_3e_4=e_3$ && \\ \hline
$\mathfrak{L}_{2}^4$ & $\mathfrak{r}_2(\mathbb{C})\oplus\mathbb{C}$ & $e_1e_2=e_3$ & $e_1e_3=-2e_1$ & $e_2e_3= 2e_2$ & \\ \hline
$\mathfrak{L}_{3}^4$ & $\mathfrak{g}_4$ & $e_1e_2=e_3$ & $e_1e_3=e_4$ & $e_1e_4= e_2$ & \\ \hline
$\mathfrak{L}_{4}^4(\alpha)$ & $\mathfrak{r}_{3,\alpha}(\mathbb{C})\oplus \mathbb{C}$ & $e_1e_2=e_2$ & $e_1e_3=\alpha e_3$ && \\ \hline
$\mathfrak{L}_{5}^4(\alpha)$ & $\mathfrak{g}_3(\alpha)$ & $e_1e_2=e_3$ & $e_1e_3= e_4$ & $e_1e_4=\alpha(e_2+e_3)$ & \\ \hline
$\mathfrak{L}_{6}^4(\alpha)$ & $\mathfrak{g}_8(\alpha)$ & $e_1e_2=e_3$ & $e_1e_3=-\alpha e_2+ e_3$& $e_1e_4=e_4$ & $e_2e_3=e_4$  \\ \hline
$\mathfrak{L}_{7}^4(\alpha,\beta)$ & $\mathfrak{g}_2(\alpha,\beta)$ & $e_1e_2=e_3$ & $e_1e_3=e_4$ & $e_1e_4=\alpha e_2 - \beta e_3 + e_4$ & \\ \hline
\end{longtable}
}


\subsubsection{5-dimensional nilpotent Lie algebras}
In her thesis (1966), Vergne obtained the list of nilpotent Lie algebras up to dimension $6$.
The degenerations of the variety $\mathfrak{NLie}^5$ were studied in~\cite{GRH}. The authors found out that the variety is irreducible: it consists of the orbit closure of the Lie algebra 

{\small

\begin{longtable}{|l|l|llll|}
\caption*{ }   \\

\hline \multicolumn{1}{|c|}{${\mathcal A}$} & \multicolumn{1}{c|}{ } & \multicolumn{4}{c|}{Multiplication table} \\ \hline 
\endfirsthead

 \multicolumn{6}{l}%
{{\bfseries  continued from previous page}} \\
\hline \multicolumn{1}{|c|}{${\mathcal A}$} & \multicolumn{1}{c|}{ } & \multicolumn{4}{c|}{Multiplication table} \\ \hline 
\endhead

\hline \multicolumn{6}{|r|}{{Continued on next page}} \\ \hline
\endfoot

\hline 
\endlastfoot

$\mathfrak{NL}_{1}^5$ & $\mathfrak{g}_5^6$ & $e_1e_2=e_3$ & $e_1e_3=e_4$ & $e_1e_4=e_5$ & $e_2e_3=e_5$  \\ \hline
\end{longtable}
}

The complete graph of degenerations can be found in~\cite{GRH}.










\subsubsection{6-dimensional nilpotent Lie algebras}
The study of the degenerations of the variety $\mathfrak{NLie}^6$ was initiated in~\cite{GRH2}, and corrected and completed in~\cite{S90}. However, it was known since~\cite{vergne} that  this variety is irreducible, defined by the Lie algebra $\mathfrak{NL}_{1}^6$ with product

{\small

\begin{longtable}{|l|l|llllll|}
\caption*{ }   \\

\hline \multicolumn{1}{|c|}{${\mathcal A}$} & \multicolumn{1}{c|}{ } & \multicolumn{6}{c|}{Multiplication table} \\ \hline 
\endfirsthead

 \multicolumn{8}{l}%
{{\bfseries  continued from previous page}} \\
\hline \multicolumn{1}{|c|}{${\mathcal A}$} & \multicolumn{1}{c|}{ } & \multicolumn{6}{c|}{Multiplication table} \\ \hline 
\endhead

\hline \multicolumn{8}{|r|}{{Continued on next page}} \\ \hline
\endfoot

\hline 
\endlastfoot

$\mathfrak{NL}_{1}^6$ & $g_{6,6}$ & $e_1e_2=e_3$ & $e_1e_3=e_4$ & $e_1e_4=e_5$ & $e_2e_3=e_5$ & $e_2e_5=e_6$ & $e_3e_4=-e_6$ \\ \hline

\end{longtable}
}

The complete graph of degenerations can be found in~\cite{S90}.



 \subsection{Malcev algebras} 
An anticommutative algebra $\mathfrak{M}$ is called a {\it Malcev} algebra  if it satisfies the identity 
\[J(x,y,xz)=J(x,y,z)x.\]
We will denote by $\mathfrak{Mal}$ the variety of Malcev algebras.

Every Lie algebra is a Malcev algebra. All Malcev algebras of dimension $\le 3$ are Lie algebras.

\subsubsection{4-dimensional Malcev algebras}
In 1970, Kuzmin proved that there exists only one non-Lie Malcev algebra of dimension $4$.
The graph of degenerations of all  $4$-dimensional Malcev algebras can be extracted from~\cite{kpv}. In particular, the variety $\mathfrak{Mal}^4$ has five irreducible components:
\[{\rm Irr}(\mathfrak{Mal}^4)=\left\{
\overline{{\mathcal O}\left(\mathfrak{M}_{i}^4\right)}\right\}_{i=1}^{3}\cup
\left\{\overline{\bigcup  {\mathcal O}\left(\mathfrak{M}_{4}^4(\alpha)\right)}\right\}
\cup\left\{\overline{\bigcup  {\mathcal O}\left(\mathfrak{M}_{5}^4(\alpha,\beta)\right)}\right\},\]
where

{\small

\begin{longtable}{|l|l|llll|}
\caption*{ }   \\

\hline \multicolumn{1}{|c|}{${\mathcal A}$} & \multicolumn{1}{c|}{ } & \multicolumn{4}{c|}{Multiplication table} \\ \hline 
\endfirsthead

 \multicolumn{6}{l}%
{{\bfseries  continued from previous page}} \\
\hline \multicolumn{1}{|c|}{${\mathcal A}$} & \multicolumn{1}{c|}{ } & \multicolumn{4}{c|}{Multiplication table} \\ \hline 
\endhead

\hline \multicolumn{6}{|r|}{{Continued on next page}} \\ \hline
\endfoot

\hline 
\endlastfoot

$\mathfrak{M}_{1}^4$ & $sl_2\oplus\mathbb{C}$ & $e_1e_2=e_2$ & $e_1e_3=-e_3$ & $e_2e_3=e_1$ & \\ \hline
$\mathfrak{M}_{2}^4$ & $r_2\oplus r_2$ & $e_1e_2=e_2$ & $e_3e_4=e_4$ &&\\ \hline
$\mathfrak{M}_{3}^4$ & $g_3(-1)$ & $e_1e_2=e_2$ & $e_1e_3=e_3$ & $e_1e_4=-e_4$ & $e_2e_3= e_4$  \\ \hline
$\mathfrak{M}_{4}^4(\alpha)$ & $g_5(\alpha)$ & $e_1e_2=e_2$ & $e_1e_3=e_2+\alpha e_3$ & $e_1e_4=(\alpha+1)e_4$  & $e_2e_3=e_4$ \\ \hline
$\mathfrak{M}_{5}^4(\alpha,\beta)$ & $g_4(\alpha,\beta)$ & $e_1e_2=e_2$ & $e_1e_3=e_2+\alpha e_3$ & $e_1e_4=e_3+\beta e_4$ & \\ \hline
\end{longtable}
}

\subsubsection{5-dimensional nilpotent Malcev algebras}
Also in 1970, Kuzmin classified the Malcev algebras of dimension $5$ up to isomorphism. Combining his results with~\cite{GRH}, the list of $5$-dimensional nilpotent Malcev algebras is easily obtained (see~\cite{kpv}).
Also in~\cite{kpv}, the authors constructed the graph of degenerations of  $\mathfrak{NMal}^5$, variety which has two irreducible components:
\[{\rm Irr}\left(\mathfrak{NMal}^5\right)=\left\{\overline{{\mathcal O}\big(\mathfrak{M}_{i}^5\big)}\right\}_{i=1}^{2},\]
where

{\small

\begin{longtable}{|l|l|llll|}
\caption*{ }   \\

\hline \multicolumn{1}{|c|}{${\mathcal A}$} & \multicolumn{1}{c|}{ } & \multicolumn{4}{c|}{Multiplication table} \\ \hline 
\endfirsthead

 \multicolumn{6}{l}%
{{\bfseries  continued from previous page}} \\
\hline \multicolumn{1}{|c|}{${\mathcal A}$} & \multicolumn{1}{c|}{ } & \multicolumn{4}{c|}{Multiplication table} \\ \hline 
\endhead

\hline \multicolumn{6}{|r|}{{Continued on next page}} \\ \hline
\endfoot

\hline 
\endlastfoot

$\mathfrak{M}_{1}^5$ & $g_{5,6}$ & $e_1e_2=e_3$ & $e_1e_3=e_4$ & $e_1e_4=e_5$ & $e_2e_3=e_5$  \\ \hline
$\mathfrak{M}_{2}^5$ & $M_5$ & $e_1e_2=e_4$ & $e_3e_4=e_5$ && \\ \hline
\end{longtable}
}

\subsubsection{6-dimensional nilpotent Malcev algebras}
The $6$-dimensional nilpotent Malcev algebras were also studied in~\cite{kpv}, employing the algebraic classification obtained by Kuzmin in 1970. 

The variety $\mathfrak{NMal}^6$ has two irreducible components:
\[{\rm Irr}\left(\mathfrak{NMal}^6\right)=\left\{\overline{{\mathcal O}\big(\mathfrak{M}_{i}^6\big)}\right\}_{i=1}^{2},\]
where

{\small

\begin{longtable}{|l|l|llllll|}
\caption*{ }   \\

\hline \multicolumn{1}{|c|}{${\mathcal A}$} & \multicolumn{1}{c|}{ } & \multicolumn{6}{c|}{Multiplication table} \\ \hline 
\endfirsthead

 \multicolumn{8}{l}%
{{\bfseries  continued from previous page}} \\
\hline \multicolumn{1}{|c|}{${\mathcal A}$} & \multicolumn{1}{c|}{ } & \multicolumn{6}{c|}{Multiplication table} \\ \hline 
\endhead

\hline \multicolumn{8}{|r|}{{Continued on next page}} \\ \hline
\endfoot

\hline 
\endlastfoot

$\mathfrak{M}_{1}^6$ & $g_6$ & $e_1e_2=e_3$ & $e_1e_3=e_4$ & $e_1e_4=e_5$
& $e_2e_3=e_5$ & $e_2e_5=e_6$ & $e_3e_4=-e_6$ \\ \hline
$\mathfrak{M}_{2}^6(\alpha)$ & $M_6^{\alpha}$ & $e_1e_2=e_3$  &  $e_1e_3=  e_5$ & $e_1e_5=e_6$ & $e_2e_4=\alpha e_5$  & $e_3e_4=e_6$ & \\ \hline
\end{longtable}
}

The graph of degenerations can be seen in~\cite{kpv}.

 \subsection{Binary Lie algebras}

Recall that an algebra $\mathfrak{A}$ is said to be {\it binary Lie} if all its $2$-generated subalgebras are Lie algebras. Let us denote this variety by $\mathfrak{BLie}$.
It was shown by Gainov in 1957 that $\mathfrak{A}$ is a binary Lie algebra if and only if it is anticommutative and satisfies the identity
\[ J(x,y,xy)=0.\]
In particular, all Lie and Malcev algebras are binary Lie.

It is straightforward that every $2$-dimensional binary Lie algebra is a Lie algebra. In 1963, Gainov 
proved that the same holds in dimension $3$. 

\subsubsection{4-dimensional binary Lie algebras}
The algebraic classification of  $4$-dimensional binary Lie algebras was obtained in the works of Gainov (1963) and Kuzmin (1998), and in~\cite{kpv}, the authors constructed the graph of degenerations. In particular, this variety 
$\mathfrak{BLie}^4$ has five irreducible components:
\[{\rm Irr}\left(\mathfrak{BLie}^4\right)=\left\{
\overline{{\mathcal O}\big(\mathfrak{BL}_{i}^4\big)}\right\}_{i=1}^{3}\cup
\left\{\overline{\bigcup {\mathcal O}\big(\mathfrak{BL}_{4}^4(\alpha)\big)}\right\}
\cup\left\{\overline{\bigcup {\mathcal O}\big(\mathfrak{BL}_{5}^4(\alpha,\beta)\big)} \right\}, \]
where

{\small

\begin{longtable}{|l|l|llll|}
\caption*{ }   \\

\hline \multicolumn{1}{|c|}{${\mathcal A}$} & \multicolumn{1}{c|}{ } & \multicolumn{4}{c|}{Multiplication table} \\ \hline 
\endfirsthead

 \multicolumn{6}{l}%
{{\bfseries  continued from previous page}} \\
\hline \multicolumn{1}{|c|}{${\mathcal A}$} & \multicolumn{1}{c|}{ } & \multicolumn{4}{c|}{Multiplication table} \\ \hline 
\endhead

\hline \multicolumn{6}{|r|}{{Continued on next page}} \\ \hline
\endfoot

\hline 
\endlastfoot
$\mathfrak{BL}_{1}^4$ & $sl_2\oplus\mathbb{C}$ & $e_1e_2=e_2$ &  $e_1e_3=-e_3$ & $e_2e_3=e_1$ & \\ \hline
$\mathfrak{BL}_{2}^4$ & $r_2\oplus r_2$ & $e_1e_2=e_2$ & $e_3e_4=e_4$ && \\ \hline
$\mathfrak{BL}_{3}^4(\alpha)$ & $g_3(\alpha)$ &  $e_1e_2=e_2$ & $e_1e_3=e_3$ & $e_1e_4=\alpha e_4$ & $e_2e_3= e_4$  \\ \hline
$\mathfrak{BL}_{4}^4(\alpha)$ & $g_5(\alpha)$ & $e_1e_2=e_2$ & $e_1e_3=e_2+\alpha e_3$ & $e_1e_4=(\alpha+1)e_4$ & $e_2e_3=e_4$  \\ \hline
$\mathfrak{BL}_{5}^4(\alpha,\beta)$ & $g_4(\alpha,\beta)$ & $e_1e_2=e_2$ & $e_1e_3=e_2+\alpha e_3$ & $e_1e_4=e_3+\beta e_4$ & \\ \hline
\end{longtable}
}

\subsubsection{5-dimensional nilpotent binary Lie algebras}
Every $5$-dimensional nilpotent binary Lie algebra is a nilpotent Malcev algebra.

\subsubsection{6-dimensional nilpotent binary Lie algebras}
The algebraic and geometric classification of the variety $\mathfrak{NBLie}^6$  can be found in~\cite{ack}. In particular, $\mathfrak{NBLie}^6$ has two irreducible components:
\[{\rm Irr}\left(\mathfrak{NBLie}^6\right)=\left\{\overline{{\mathcal O}\big(\mathfrak{NBL}_{i}^6\big)}\right\}_{i=1}^{2},\]

where
{\small
 
\begin{longtable}{|l|l|llllll|}
\caption*{ }   \\

\hline \multicolumn{1}{|c|}{${\mathcal A}$} & \multicolumn{1}{c|}{ } & \multicolumn{6}{c|}{Multiplication table} \\ \hline 
\endfirsthead

 \multicolumn{8}{l}%
{{\bfseries  continued from previous page}} \\
\hline \multicolumn{1}{|c|}{${\mathcal A}$} & \multicolumn{1}{c|}{ } & \multicolumn{6}{c|}{Multiplication table} \\ \hline 
\endhead

\hline \multicolumn{8}{|r|}{{Continued on next page}} \\ \hline
\endfoot

\hline 
\endlastfoot
$\mathfrak{NBL}_{1}^6$ & $\mathbf{B}_{6,3}$ & $e_{1}e_{2}=e_{3}$ & $e_{3}e_{4}=e_{5}$ &  $e_{1}e_{3}=e_{6}$ & $e_{4}e_{5}=e_{6}$ & & \\ \hline
$\mathfrak{NBL}_{2}^6$ & $g_6$ & $e_1e_2=e_3$ & $e_1e_3=e_4$ & $e_1e_4=e_5$  &
$e_2e_3=e_5$ & $e_2e_5=e_6$ & $e_3e_4=-e_6$ \\ \hline 
\end{longtable}
}

\subsection{Tortkara algebras} 

An anticommutative algebra $\mathfrak{A}$ is called a {\it Tortkara} algebra if it satisfies the identity 
\[(xy)(zy)=J(x,y,z)y.\]
We will denote this variety by $\mathfrak{Tor}$.

\subsubsection{2-dimensional Tortkara algebras}
Checking the classification of $2$-dimensional algebras~\cite{kv16}, we have that the variety $\mathfrak{Tor}^2$ has only one non-trivial algebra

{\small

\begin{longtable}{|l|l|l|}
\caption*{ }   \\

\hline \multicolumn{1}{|c|}{${\mathcal A}$} & \multicolumn{1}{c|}{ } & \multicolumn{1}{c|}{Multiplication table} \\ \hline 
\endfirsthead

 \multicolumn{3}{l}%
{{\bfseries  continued from previous page}} \\
\hline \multicolumn{1}{|c|}{${\mathcal A}$} & \multicolumn{1}{c|}{ } & \multicolumn{1}{c|}{Multiplication table} \\ \hline 
\endhead

\hline \multicolumn{3}{|r|}{{Continued on next page}} \\ \hline
\endfoot

\hline 
\endlastfoot
$\mathfrak{T}_{1}^2$ & $\mathbf{B}_3 $ & $ e_1e_2=e_2$ \\ \hline
\end{longtable}
}\noindent so it is irreducible.

\subsubsection{3-dimensional nilpotent Tortkara algebras}
For the variety $\mathfrak{NTor}^3$, we rely on the classification from~\cite{fkkv}. The rigid algebra determines it:

{\small

\begin{longtable}{|l|l|ll|}
\caption*{ }   \\

\hline \multicolumn{1}{|c|}{${\mathcal A}$} & \multicolumn{1}{c|}{ } & \multicolumn{2}{c|}{Multiplication table} \\ \hline 
\endfirsthead

 \multicolumn{4}{l}%
{{\bfseries  continued from previous page}} \\
\hline \multicolumn{1}{|c|}{${\mathcal A}$} & \multicolumn{1}{c|}{ } & \multicolumn{2}{c|}{Multiplication table} \\ \hline 
\endhead

\hline \multicolumn{4}{|r|}{{Continued on next page}} \\ \hline
\endfoot

\hline 
\endlastfoot

$\mathfrak{NT}_{1}^3$ & $\mathcal{N}_7$ & $e_1e_2=e_3$ & \\ \hline

\end{longtable}
}

\subsubsection{3-dimensional Tortkara algebras}
In~\cite{gkk19alg}, $3$-dimensional Tortkara algebras were selected from the list of $3$-dimensional anticommutative algebras of~\cite{ikv18}. A consequence of the graph of degenerations of the anticommutative algebras (see~\cite{ikv18}) is that $\mathfrak{Tor}^3$ is irreducible~\cite{gkk19alg}, defined by the rigid algebra 

{\small

\begin{longtable}{|l|l|ll|}
\caption*{ }   \\

\hline \multicolumn{1}{|c|}{${\mathcal A}$} & \multicolumn{1}{c|}{ } & \multicolumn{2}{c|}{Multiplication table} \\ \hline 
\endfirsthead

 \multicolumn{4}{l}%
{{\bfseries  continued from previous page}} \\
\hline \multicolumn{1}{|c|}{${\mathcal A}$} & \multicolumn{1}{c|}{ } & \multicolumn{2}{c|}{Multiplication table} \\ \hline 
\endhead

\hline \multicolumn{4}{|r|}{{Continued on next page}} \\ \hline
\endfoot

\hline 
\endlastfoot
$\mathfrak{T}_{1}^3$ & $\mathfrak{A}_{1}^{0}$ & $e_1e_2=e_3$ & $e_1e_3=e_1+e_3$  \\ \hline
\end{longtable}
}

%



\subsubsection{4-dimensional nilpotent Tortkara algebras}

In the variety $\mathfrak{NTor}^4$ there are only two non-trivial algebras and one irreducible component~\cite{gkks} determines by the rigid algebra:

{\small

\begin{longtable}{|l|l|ll|}
\caption*{ }   \\

\hline \multicolumn{1}{|c|}{${\mathcal A}$} & \multicolumn{1}{c|}{ } & \multicolumn{2}{c|}{Multiplication table} \\ \hline 
\endfirsthead

 \multicolumn{4}{l}%
{{\bfseries  continued from previous page}} \\
\hline \multicolumn{1}{|c|}{${\mathcal A}$} & \multicolumn{1}{c|}{ } & \multicolumn{2}{c|}{Multiplication table} \\ \hline 
\endhead

\hline \multicolumn{4}{|r|}{{Continued on next page}} \\ \hline
\endfoot

\hline 
\endlastfoot

$\mathfrak{NT}_{1}^4$ & $\mathbb{T}_{02}^5$ & $e_1e_2=e_3$ & $e_1e_3=e_4$ \\ \hline

\end{longtable}
}

\subsubsection{5-dimensional nilpotent Tortkara algebras}

The algebraic and geometric classifications of the variety $\mathfrak{NTor}^5$ were given in~\cite{gkks}. Again, there is a unique irreducible component defined by

{\small

\begin{longtable}{|l|l|lll|}
\caption*{ }   \\

\hline \multicolumn{1}{|c|}{${\mathcal A}$} & \multicolumn{1}{c|}{ } & \multicolumn{3}{c|}{Multiplication table} \\ \hline 
\endfirsthead

 \multicolumn{5}{l}%
{{\bfseries  continued from previous page}} \\
\hline \multicolumn{1}{|c|}{${\mathcal A}$} & \multicolumn{1}{c|}{ } & \multicolumn{3}{c|}{Multiplication table} \\ \hline 
\endhead

\hline \multicolumn{5}{|r|}{{Continued on next page}} \\ \hline
\endfoot

\hline 
\endlastfoot
$\mathfrak{NT}_{1}^5$ & $\mathbb{T}_{10}^5$ & $e_1e_2=e_3$ & $e_1e_3=e_4$ & $e_2e_4=e_5$  \\ \hline
\end{longtable}
}

\subsubsection{6-dimensional nilpotent Tortkara algebras}

More recently, in~\cite{gkk19}, it was provided the geometric classification of $6$-dimensional nilpotent Tortkara algebras, which is  based on the description of all $6$-dimensional nilpotent Tortkara algebras by Gorshkov, Kaygorodov and Khrypchenko (2019) and on the description of all degenerations of $6$-dimensional nilpotent Malcev algebras~\cite{kpv}.
In particular, there exist three irreducible components in the variety $\mathfrak{NTor}^6$: 
\[{\rm Irr}\left(\mathfrak{NTor}^6\right)=\left\{\overline{{\mathcal O}\big(\mathfrak{NT}_{i}^6\big)}\right\}_{i=1}^{3},\]
where

{\small

\begin{longtable}{|l|l|lllll|}
\caption*{ }   \\

\hline \multicolumn{1}{|c|}{${\mathcal A}$} & \multicolumn{1}{c|}{ } & \multicolumn{5}{c|}{Multiplication table} \\ \hline 
\endfirsthead

 \multicolumn{7}{l}%
{{\bfseries  continued from previous page}} \\
\hline \multicolumn{1}{|c|}{${\mathcal A}$} & \multicolumn{1}{c|}{ } & \multicolumn{5}{c|}{Multiplication table} \\ \hline 
\endhead

\hline \multicolumn{7}{|r|}{{Continued on next page}} \\ \hline
\endfoot

\hline 
\endlastfoot
$\mathfrak{NT}_{1}^6$ & $\mathbb{T}_{10}^6$ & $e_1e_2=e_3$ & $e_1e_3=e_6$ & $e_1e_4=e_5$ & $e_2e_3=e_5$ & $e_4e_5=e_6$ \\ \hline
$\mathfrak{NT}_{2}^6$ & $\mathbb{T}_{17}^6$ & $e_1e_2=e_3$ & $e_1e_3=e_4$ & $e_1e_4=e_5$ & $e_2e_3=e_5$ & $e_2e_5=e_6$ \\ \hline
$\mathfrak{NT}_{3}^6$ & $\mathbb{T}_{19}^6$ & $e_1e_2=e_3$ & $e_1e_3=e_4$ & $e_1e_5=e_6$ & $e_2e_4=e_5$ & $e_3e_4=e_6$  \\ \hline
\end{longtable}
}


 \subsection{Dual mock-Lie algebras}

An anticommutative algebra $\mathfrak{D}$ is called a {\it dual mock-Lie}  algebra if it satisfies the identity 
\[(xy)z =-x(yz).\]
Let $\mathfrak{DML}$ denote the variety of dual mock-Lie algebras.
The main subclass of dual mock-Lie algebras is $2$-step nilpotent Lie algebras.
The first example of non-Lie dual mock-Lie algebra appears in dimension $7$.

\subsubsection{7-dimensional dual mock-Lie algebras}

In~\cite{ckls}, the authors determined all the  $7$-dimensional  dual mock-Lie algebras up to isomorphism and found the degenerations between them. 
This variety $\mathfrak{DML}^7$ has three irreducible components:
\[{\rm Irr}\left(\mathfrak{DML}^7\right)=
\left\{\overline{{\mathcal O}\left(\mathfrak{D}_{i}^7\right)}\right\}_{i=1}^{3},\]
where

{\small

\begin{longtable}{|l|l|llllll|}
\caption*{ }   \\

\hline \multicolumn{1}{|c|}{${\mathcal A}$} & \multicolumn{1}{c|}{ } & \multicolumn{6}{c|}{Multiplication table} \\ \hline 
\endfirsthead

 \multicolumn{8}{l}%
{{\bfseries  continued from previous page}} \\
\hline \multicolumn{1}{|c|}{${\mathcal A}$} & \multicolumn{1}{c|}{ } & \multicolumn{6}{c|}{Multiplication table} \\ \hline 
\endhead

\hline \multicolumn{8}{|r|}{{Continued on next page}} \\ \hline
\endfoot

\hline 
\endlastfoot

$\mathfrak{D}_{1}^7$ & $\mathfrak{D}^7_{09}$ & $e_1e_2=e_6$ & $e_1e_5=e_7$ & $e_2e_3=e_7$ & $e_3e_4=e_6$ && \\ \hline
$\mathfrak{D}_{2}^7$ & $\mathfrak{D}^7_{13}$ & $e_1e_2=e_5$ & $e_1e_3=e_6$ & $e_2e_4=e_7$ & $e_3e_4=e_5$ &&  \\ \hline
$\mathfrak{D}_{3}^7$ & $\mathfrak{D}^7_{14}$ & $e_1e_2=e_4$ & $e_1e_3=e_5$ & $e_1e_6=e_7$ & $e_2e_3=e_6$ & $e_2e_5=-e_7$ & $e_3e_4=e_7$
\end{longtable}
}

\subsubsection{\texorpdfstring{$8$-dimensional dual mock-Lie algebras}{8-dimensional dual mock-Lie algebras}}

Also in~\cite{ckls}, the authors obtained the algebraic and geometric classifications of the $8$-dimensional  dual mock-Lie algebras. This variety $\mathfrak{DML}^8$ has four irreducible components:
\[{\rm Irr}(\mathfrak{DML}^8)=
\left\{\overline{{\mathcal O}\big(\mathfrak{D}_{i}^8\big)}\right\}_{i=1}^{4},\]
where

{\small

\begin{longtable}{|l|l|llllll|}
\caption*{ }   \\

\hline \multicolumn{1}{|c|}{${\mathcal A}$} & \multicolumn{1}{c|}{ } & \multicolumn{6}{c|}{Multiplication table} \\ \hline 
\endfirsthead

 \multicolumn{8}{l}%
{{\bfseries  continued from previous page}} \\
\hline \multicolumn{1}{|c|}{${\mathcal A}$} & \multicolumn{1}{c|}{ } & \multicolumn{6}{c|}{Multiplication table} \\ \hline 
\endhead

\hline \multicolumn{8}{|r|}{{Continued on next page}} \\ \hline
\endfoot

\hline 
\endlastfoot

$\mathfrak{D}_{1}^8$ & $\mathfrak{D}_{17}^8$ & $e_1e_2=e_7$ &  $e_3e_4=e_8$ & $e_5e_6=e_7+e_8$ &&&  \\ \hline
$\mathfrak{D}_{2}^8$ & $\mathfrak{D}_{30}^8$ & $e_1e_2=e_6$  & $e_1e_3=e_7$ & $e_1e_4=e_8$ & $e_2e_3=e_8$  & $e_2e_5=e_7$ & $e_4e_5=e_6$ \\ \hline 
$\mathfrak{D}_{3}^8$ & $\mathfrak{D}_{33}^8$ & $e_1e_2=e_5$ & $e_2e_3=e_6$ & $e_3e_4=e_7$ & $e_4e_1=e_8$ && \\ \hline
$\mathfrak{D}_{4}^8$ & $\mathfrak{D}_{36}^8$ & $e_1e_2 = e_4$ & $e_1e_3 = e_5$ & $e_1e_6 = e_8$ & $e_2e_3 = e_6$& & \\ 
&& $e_2e_5=-e_8$  & $e_3e_4 = e_8$  & $e_3e_7 = e_8$ &&& \\ \hline
\end{longtable}
}

\subsection{\texorpdfstring{$\mathfrak{CD}$-algebras}{CD-algebras}}

The class of  $\mathfrak{CD}$-algebras is defined by the 
property that the commutator of any pair of multiplication operators is a derivation; 
namely, an algebra $\mathfrak{A}$ is a $\mathfrak{CD}$-algebra if and only if  
\[ [T_x,T_y]   \in \mathfrak{Der} (\mathfrak{A}),\]
for all $x,y \in \mathfrak{A}$, where  $T_z \in \{ R_z,L_z\}$. Here we use the notation $R_z$ (resp. $L_z$) for the operator of right (resp. left) multiplication in $\mathfrak{A}$. We will denote the variety of $\mathfrak{CD}$-algebras by $\mathfrak{CD}$.
In terms of identities, the class of $\mathfrak{CD}$-algebras is defined by the following three:
\[((xy)a)b-((xy)b)a=((xa)b-(xb)a)y+x((ya)b-(yb)a),\]
\[(a(xy))b-a((xy)b)=((ax)b-a(xb))y+x((ay)b-a(yb)),\]
\[a(b(xy))-b(a(xy))=(a(bx)-b(ax))y+x(a(by)-b(ay)).\]
In the commutative and anticommutative cases, they are reduced to the first identity. All Lie and Jordan algebras are $\mathfrak{CD}$-algebras. On the other hand, each anticommutative $\mathfrak{CD}$-algebra is a binary Lie algebra.

\subsubsection{\texorpdfstring{2-dimensional $\mathfrak{CD}$-algebras}{2-dimensional CD-algebras}}

Analyzing the table of all $2$-dimensional algebras from~\cite{kv16}, we obtain the classification of all $2$-dimensional $\mathfrak{CD}$-algebras:
\vspace{-7mm}
{\small

\begin{longtable}{|l|l|lll|}
\caption*{ }   \\

\hline \multicolumn{1}{|c|}{${\mathcal A}$} & \multicolumn{1}{c|}{ } & \multicolumn{3}{c|}{Multiplication table} \\ \hline 
\endfirsthead

 \multicolumn{5}{l}%
{{\bfseries  continued from previous page}} \\
\hline \multicolumn{1}{|c|}{${\mathcal A}$} & \multicolumn{1}{c|}{ } & \multicolumn{3}{c|}{Multiplication table} \\ \hline 
\endhead

\hline \multicolumn{5}{|r|}{{Continued on next page}} \\ \hline
\endfoot

\hline 
\endlastfoot
$\mathfrak{CD}_{1}^2$ & $\mathbf{A}_2$ & $e_1e_1=e_2$ & $e_1e_2=e_2$ & $e_2e_1=-e_2$  \\ \hline
$\mathfrak{CD}_{2}^2$ & $\mathbf{A}_3$ & $e_1e_1=e_2$ && \\ \hline
$\mathfrak{CD}_{3}^2$ & $\mathbf{B}_3$ & $e_1e_2=e_2$ & $e_2e_1=-e_2$ & \\ \hline
$\mathfrak{CD}_{4}^2$ & $\mathbf{D}_2(0,0)$ & $e_1e_1=e_1$ && \\ \hline
$\mathfrak{CD}_{5}^2$ & $\mathbf{D}_2(1,1)$ & $e_1e_1=e_1$ & $e_1e_2=e_2$ & $e_2e_1=e_2$ \\ \hline
$\mathfrak{CD}_{6}^2$ & $\mathbf{E}_1(0,0,0,0)$ & $e_1e_1=e_1$ & $e_2e_2=e_2$ & \\ \hline
$\mathfrak{CD}_{7}^2(\alpha)$ & $\mathbf{A}_1(\alpha)$ & $e_1e_1=e_1+e_2$ & $e_1e_2=\alpha e_2$ & $e_2e_1=(1-\alpha)e_2$  \\ \hline
$\mathfrak{CD}_{8}^2(\alpha)$ & $\mathbf{E}_5(\alpha)$ & $e_1e_1=e_1$ & $e_1e_2=(1-\alpha)e_1+\alpha e_2$ & \\
&& $e_2e_1=\alpha e_1+(1-\alpha)e_2$  & $e_2e_2=e_2$ & \\ \hline

\end{longtable}
}

Basing on the full description of degenerations of $2$-dimensional algebras~\cite{kv16}, we conclude that
\[{\rm Irr}\left(\mathfrak{CD}^2\right)=\left\{\overline{{\mathcal O}\left(\mathfrak{CD}_{1}^2\right)}\right\}\cup
\left\{\overline{{\mathcal O}\left(\mathfrak{CD}_{6}^2\right)}\right\}
\cup\left\{\overline{\bigcup {\mathcal O}\left(\mathfrak{CD}_{7}^2(\alpha)\right)}\right\}.\]

\subsubsection{\texorpdfstring{$3$-dimensional nilpotent $\mathfrak{CD}$-algebras}{3-dimensional nilpotent CD-algebras}}

It is easy to see that each $3$-dimensional nilpotent algebra is a $\mathfrak{CD}$-algebra. The degenerations of $3$-dimensional nilpotent algebras were fully described in~\cite{fkkv}, where it was proved that this variety is irreducible and determined by the rigid algebra

{\small

\begin{longtable}{|l|l|lll|}
\caption*{ }   \\

\hline \multicolumn{1}{|c|}{${\mathcal A}$} & \multicolumn{1}{c|}{ } & \multicolumn{3}{c|}{Multiplication table} \\ \hline 
\endfirsthead

 \multicolumn{5}{l}%
{{\bfseries  continued from previous page}} \\
\hline \multicolumn{1}{|c|}{${\mathcal A}$} & \multicolumn{1}{c|}{ } & \multicolumn{3}{c|}{Multiplication table} \\ \hline 
\endhead

\hline \multicolumn{5}{|r|}{{Continued on next page}} \\ \hline
\endfoot

\hline 
\endlastfoot
$\mathfrak{CD}_{1}^3$ & $\mathfrak{N}_2$ & $e_1 e_1 = e_2$ & $e_2 e_1= e_3$ & $e_2 e_2=e_3$  \\ \hline
\end{longtable}
}

{

\subsubsection{\texorpdfstring{$4$-dimensional nilpotent $\mathfrak{CD}$-algebras}{4-dimensional nilpotent CD-algebras}}

The classification of the nilpotent   $\mathfrak{CD}$-algebras of dimension $4$ is a  result of the work by 
   Kaygorodov and Khrypchenko  (2022). Later, their geometric classification appeared in~\cite{kkcd21}. In particular, the variety $\mathfrak{NCD}^4$  has two irreducible components and is determined by the following algebras

{\small

\begin{longtable}{|l|l|lll|}
\caption*{ }   \\

\hline \multicolumn{1}{|c|}{${\mathcal A}$} & \multicolumn{1}{c|}{ } & \multicolumn{3}{c|}{Multiplication table} \\ \hline 
\endfirsthead

 \multicolumn{5}{l}%
{{\bfseries  continued from previous page}} \\
\hline \multicolumn{1}{|c|}{${\mathcal A}$} & \multicolumn{1}{c|}{ } & \multicolumn{3}{c|}{Multiplication table} \\ \hline 
\endhead

\hline \multicolumn{5}{|r|}{{Continued on next page}} \\ \hline
\endfoot

\hline 
\endlastfoot

$\mathfrak{CD}_{1}^4(\alpha, \beta)$ &$\mathfrak{CD}^4_{12}(\alpha, \beta)$ & 
$e_1 e_1 = e_2$ & $e_1 e_2=e_3$ & $e_1e_3=(\beta-2)e_4$ \\&& $e_2 e_1=\beta e_3$ & $e_2e_2=\alpha e_4$ & $e_3e_1=(1-2\beta) e_4$ \\
\hline

$\mathfrak{CD}_{2}^4(\alpha,\beta,\gamma,\delta)$ &$\mathfrak{CD}^{4}_{112}(\delta, \alpha, \beta, \gamma)$&  
$e_1 e_1 = \delta e_3+e_4$ & $e_1e_3=\alpha e_4$ & $e_2 e_1=e_3+\beta e_4$  \\
&&$e_2 e_2=e_3$& $e_2e_3=\gamma e_4$&$e_3e_3=e_4$\\
\hline 
\end{longtable}
}

}

 \subsection{\texorpdfstring{Commutative $\mathfrak{CD}$-algebras}{Commutative CD-algebras}}

\subsubsection{\texorpdfstring{$2$-dimensional commutative $\mathfrak{CD}$-algebras}{2-dimensional commutative CD-algebras}}

In dimension 2 we have the following commutative $\mathfrak{CD}$-algebras:

\

{\small

\begin{longtable}{|l|l|lll|}
\caption*{ }   \\

\hline \multicolumn{1}{|c|}{${\mathcal A}$} & \multicolumn{1}{c|}{ } & \multicolumn{3}{c|}{Multiplication table} \\ \hline 
\endfirsthead

 \multicolumn{5}{l}%
{{\bfseries  continued from previous page}} \\
\hline \multicolumn{1}{|c|}{${\mathcal A}$} & \multicolumn{1}{c|}{ } & \multicolumn{3}{c|}{Multiplication table} \\ \hline 
\endhead

\hline \multicolumn{5}{|r|}{{Continued on next page}} \\ \hline
\endfoot

\hline 
\endlastfoot

$\mathfrak{CCD}_{1}^2$ & $\mathfrak{CD}_{2}^2$ & $e_1e_1=e_2$ && \\ \hline
$\mathfrak{CCD}_{2}^2$ & $\mathfrak{CD}_{4}^2$ & $e_1e_1=e_1$ && \\ \hline
$\mathfrak{CCD}_{3}^2$ & $\mathfrak{CD}_{5}^2$ & $e_1e_1=e_1$ & $e_1e_2=e_2$ & \\ \hline
$\mathfrak{CCD}_{4}^2$ & $\mathfrak{CD}_{6}^2$ & $e_1e_1=e_1$ & $e_2e_2=e_2$ &\\ \hline
$\mathfrak{CCD}_{5}^2$ & $\mathfrak{CD}_{7}^2(\frac 12)$ & $e_1e_1=e_1+e_2$ & $e_1e_2=\frac 12 e_2$ & \\ \hline
$\mathfrak{CCD}_{6}^2$ & $\mathfrak{CD}_{8}^2(\frac 12)$ & $e_1e_1=e_1$ & $e_1e_2=\frac 12 e_1+\frac 12 e_2$ &  $e_2e_2=e_2$  \\ \hline
\end{longtable}
}

Hence
\[{\rm Irr}\left(\mathfrak{CCD}^2\right)=\left\{\overline{{\mathcal O}\left(\mathfrak{CCD}_{i}^2\right)}\right\}_{i=4}^{5}.\]

\subsubsection{\texorpdfstring{$3$-dimensional nilpotent commutative $\mathfrak{CD}$-algebras}{3-dimensional nilpotent commutative CD-algebras}}

Choosing the commutative algebras from the list of all nilpotent $3$-dimensional algebras~\cite{fkkv}, we obtain the next classification:

{\small

\begin{longtable}{|l|l|ll|}
\caption*{ }   \\

\hline \multicolumn{1}{|c|}{${\mathcal A}$} & \multicolumn{1}{c|}{ } & \multicolumn{2}{c|}{Multiplication table} \\ \hline 
\endfirsthead

 \multicolumn{4}{l}%
{{\bfseries  continued from previous page}} \\
\hline \multicolumn{1}{|c|}{${\mathcal A}$} & \multicolumn{1}{c|}{ } & \multicolumn{2}{c|}{Multiplication table} \\ \hline 
\endhead

\hline \multicolumn{4}{|r|}{{Continued on next page}} \\ \hline
\endfoot

\hline 
\endlastfoot

$\mathfrak{NCCD}_{1}^3$ & $\mathfrak{N}_1$ & $e_1 e_1 = e_2$ & $e_2 e_2=e_3$ \\ \hline
$\mathfrak{NCCD}_{2}^3$ & $\mathfrak{N}_4(1)$ &  $e_1 e_1 = e_2$ & $e_1 e_2=e_3$  \\ \hline
$\mathfrak{NCCD}_{3}^3$ & $\mathfrak{N}_5$ & $e_1 e_1 = e_2$ & \\ \hline
$\mathfrak{NCCD}_{4}^3$ & $\mathfrak{N}_6$ & $e_1 e_1 = e_3$ & $e_2 e_2=e_3$  \\ \hline
\end{longtable}
}

It follows from the graph of degenerations of $3$-dimensional nilpotent algebras~\cite{fkkv} that the variety $\mathfrak{NCCD}^3$ is irreducible and it is determined by the rigid algebra $\mathfrak{NCCD}_{1}^3$.

\subsubsection{\texorpdfstring{$4$-dimensional nilpotent commutative $\mathfrak{CD}$-algebras}{4-dimensional nilpotent commutative CD-algebras}}

By direct verification we see that only the algebras $\mathfrak{C}_{13}-\mathfrak{C}_{19}$ and $\mathfrak{C}_{21}-\mathfrak{C}_{24}$ from the list of all $4$-dimensional nilpotent commutative algebras~\cite{fkkv} are not $\mathfrak{CD}$. 
Hence, the description of all the degenerations of $4$-dimensional nilpotent commutative algebras~\cite{fkkv} implies that the subvariety formed by $\mathfrak{CD}$-algebras is irreducible and determined by the rigid algebra

{\small

\begin{longtable}{|l|l|llll|}
\caption*{ }   \\

\hline \multicolumn{1}{|c|}{${\mathcal A}$} & \multicolumn{1}{c|}{ } & \multicolumn{4}{c|}{Multiplication table} \\ \hline 
\endfirsthead

 \multicolumn{6}{l}%
{{\bfseries  continued from previous page}} \\
\hline \multicolumn{1}{|c|}{${\mathcal A}$} & \multicolumn{1}{c|}{ } & \multicolumn{4}{c|}{Multiplication table} \\ \hline 
\endhead

\hline \multicolumn{6}{|r|}{{Continued on next page}} \\ \hline
\endfoot

\hline 
\endlastfoot

$\mathfrak{NCCD}_{1}^4$ & $\mathfrak{C}_{29}$ & $e_1 e_1 = e_4$ & $e_1 e_2=e_3$ & $e_2e_2=e_4$ & $e_3e_3=e_4$  \\ \hline
\end{longtable}
}

{

\subsubsection{\texorpdfstring{$5$-dimensional nilpotent commutative $\mathfrak{CD}$-algebras}{5-dimensional nilpotent commutative CD-algebras}}

The classification of the nilpotent commutative $\mathfrak{CD}$-algebras of dimension $4$ is a joint result of the works by 
 Jumaniyozov,  Kaygorodov and Khudoyberdiyev (2021)and Abdelwahab and Hegazi (2016). Later, their geometric classification appeared in~\cite{jkk22}. In particular, the variety $\mathfrak{NCCD}^5$  has ten irreducible components and is determined by the following algebras

{\small

\begin{longtable}{|l|l|llll|}
\caption*{ }   \\

\hline \multicolumn{1}{|c|}{${\mathcal A}$} & \multicolumn{1}{c|}{ } & \multicolumn{4}{c|}{Multiplication table} \\ \hline 
\endfirsthead

 \multicolumn{6}{l}%
{{\bfseries  continued from previous page}} \\
\hline \multicolumn{1}{|c|}{${\mathcal A}$} & \multicolumn{1}{c|}{ } & \multicolumn{4}{c|}{Multiplication table} \\ \hline 
\endhead

\hline \multicolumn{6}{|r|}{{Continued on next page}} \\ \hline
\endfoot

\hline 
\endlastfoot

$\mathfrak{NCCD}_{1}^5$ &${\mathcal J}_{21}$ &     $e_1e_1=e_5$ & $e_1e_2=e_4$ & $e_2e_2=e_5$ &\\
&& $e_3e_3=e_4$ & $e_3e_4=e_5$ && \\
\hline

$\mathfrak{NCCD}_{2}^5(\alpha)$ &$\mathfrak{C}_{16}^5(\alpha)$&   $e_1e_1=e_2$ & $e_1e_2=e_4$ &$e_1e_4= (\alpha+1) e_5$&\\
&&$e_2e_2=\alpha e_5$ &$e_3e_3=e_4$ &&\\
\hline

$\mathfrak{NCCD}_{3}^5(\alpha,\beta)$ & $\mathfrak{C}^5_{26}(\alpha,\beta)$ &  $e_1e_1=\alpha e_5$ & $e_1e_2=e_3$ & $e_2e_2=\beta e_5$ &\\ 
&&$e_1e_3=e_4+e_5$ & $e_2e_3=e_4$ & $e_3e_3=e_5 $ &\\
\hline

 $\mathfrak{NCCD}_{4}^5(\alpha)$ & $\mathfrak{C}_{49}^5(\alpha)$& 
 $e_1e_1=e_3$ & $e_1e_2=e_5$ &$e_2e_2=e_4$ &\\
 &&$e_3e_3=\alpha e_5$ & $e_3e_4=e_5$& $e_4e_4=e_5$ &\\
\hline

$\mathfrak{NCCD}_{5}^5(\alpha)$ & ${\mathfrak{C}}_{69}^{5}(\alpha)$ &    $e_1 e_1 = e_4$&$e_1e_2=\alpha e_5$ &$e_1e_3=e_5$ &\\
&& $e_2e_2=e_5$
&$ e_2 e_3=e_4$ &$e_4e_4=e_5$&\\
\hline 

$\mathfrak{NCCD}_{6}^5$ &${\mathfrak{C}}_{72}^{5}$ &    $e_1 e_1 = e_4$&$e_2e_2=e_5$&$ e_2 e_3=e_4+e_5$ &$e_4e_4=e_5$ \\
\hline 

$\mathfrak{NCCD}_{7}^5$ & $\mathfrak{C}_{76}^5$&  $e_1e_1=e_2$ & $e_1e_2=e_4$  &$e_1e_4= e_5$ &\\
&&$e_2e_2= - 2 e_5$ &$e_3e_3=e_4+3e_5$&&\\
\hline 

$\mathfrak{NCCD}_{8}^5$ & $\mathfrak{C}_{77}^5$&   $e_1e_1=e_2$ & $e_1e_2=e_4$ &$e_1e_4= e_5$ &\\
&&$e_2e_3=  e_5$ &$e_3e_3=e_4$ &&\\
\hline

$\mathfrak{NCCD}_{9}^5(\alpha)$ & $\mathfrak{C}_{80}^5(\alpha)$&   
$e_1 e_1 = e_2$& 
$e_1 e_2=e_3$&
$e_1e_3= \alpha e_4$&\\&& 
$e_1 e_4=e_5$&
$e_2e_2= (\alpha +1)e_4$&
$e_2 e_3=(\alpha+3)e_5$ &\\
\hline 

$\mathfrak{NCCD}_{10}^5$ & $\mathfrak{C}_{81}^5 $&  
$e_1 e_1 = e_2$& 
$e_1 e_2=e_3$&
$e_1e_3=  e_4$& \\
&&$e_2e_2= 2 e_4$&
$e_2e_4=   e_5$ && \\
\hline 

\end{longtable}
}}

\subsection{\texorpdfstring{Anticommutative $\mathfrak{CD}$-algebras}{Anticommutative CD-algebras}}
\subsubsection{\texorpdfstring{$2$-dimensional anticommutative $\mathfrak{CD}$-algebras}{2-dimensional anticommutative CD-algebras}}

There is only one $2$-dimensional anticommutative $\mathfrak{CD}$-algebra:

{\small

\begin{longtable}{|l|l|ll|}
\caption*{ }   \\

\hline \multicolumn{1}{|c|}{${\mathcal A}$} & \multicolumn{1}{c|}{ } & \multicolumn{2}{c|}{Multiplication table} \\ \hline 
\endfirsthead

 \multicolumn{4}{l}%
{{\bfseries  continued from previous page}} \\
\hline \multicolumn{1}{|c|}{${\mathcal A}$} & \multicolumn{1}{c|}{ } & \multicolumn{2}{c|}{Multiplication table} \\ \hline 
\endhead

\hline \multicolumn{4}{|r|}{{Continued on next page}} \\ \hline
\endfoot

\hline 
\endlastfoot

$\mathfrak{ACD}_{1}^2$ & $\mathfrak{CD}_{3}^2$ & $e_1e_2=e_2$ & $e_2e_1=-e_2$ \\ \hline
\end{longtable}
}

\subsubsection{\texorpdfstring{$3$-dimensional anticommutative $\mathfrak{CD}$-algebras}{3-dimensional anticommutative CD-algebras}}

It was proved that each $3$-dimensional binary Lie algebra is a Lie algebra. Thus, the variety of $3$-dimensional anticommutative $\mathfrak{CD}$-algebras coincides with the variety of $3$-dimensional Lie algebras. 

\subsubsection{\texorpdfstring{$4$-dimensional anticommutative $\mathfrak{CD}$-algebras}{4-dimensional anticommutative CD-algebras}}

The full description of degenerations of all $4$-dimensional binary Lie algebras was made in~\cite{kpv}. Observe that almost all of such algebras are anticommutative $\mathfrak{CD}$-algebras, except $g_3(\beta)$ for $\beta\not\in\{0,2\}$ and $g_6$. It follows that 
\[
{\rm Irr}\left(\mathfrak{ACD}^4\right)=\left\{
\overline{{\mathcal O}\left(\mathfrak{ACD}_{i}^4\right)}\right\}_{i=1}^3
\cup  
\left\{\overline{\bigcup {\mathcal O}\left(\mathfrak{ACD}_{4}^4(\alpha)\right)}\right\}
\cup\left\{\overline{\bigcup {\mathcal O}\left(\mathfrak{ACD}_{5}^4(\alpha,\beta)\right)}
\right\},
\]
where

{\small
\begin{longtable}{|l|l|llll|}
\caption*{ }   \\

\hline \multicolumn{1}{|c|}{${\mathcal A}$} & \multicolumn{1}{c|}{ } & \multicolumn{4}{c|}{Multiplication table} \\ \hline 
\endfirsthead

 \multicolumn{6}{l}%
{{\bfseries  continued from previous page}} \\
\hline \multicolumn{1}{|c|}{${\mathcal A}$} & \multicolumn{1}{c|}{ } & \multicolumn{4}{c|}{Multiplication table} \\ \hline 
\endhead

\hline \multicolumn{6}{|r|}{{Continued on next page}} \\ \hline
\endfoot

\hline 
\endlastfoot

$\mathfrak{ACD}_{1}^4$ & $sl_2\oplus\mathbb{C}$ & $e_1e_2=e_2$ & $e_1e_3=-e_3$ & $e_2e_3=e_1$ & \\ \hline
$\mathfrak{ACD}_{2}^4$ & $r_2\oplus r_2$ & $e_1e_2=e_2$ & $e_3e_4=e_4$ && \\ \hline
$\mathfrak{ACD}_{3}^4$ & $g_3(0)$ & $e_1e_2=e_2$ & $e_1e_3=e_3$ & $e_2e_3=e_4$ & \\ \hline
$\mathfrak{ACD}_{4}^4(\alpha)$ & $g_5(\alpha)$ & $e_1e_2=e_2$ & $e_1e_3=e_2+\alpha e_3$ & $e_1e_4=(\alpha+1)e_4$ & $e_2e_3=e_4$  \\ \hline
$\mathfrak{ACD}_{5}^4(\alpha,\beta)$ & $g_4(\alpha,\beta)$ & $e_1e_2=e_2$ & $e_1e_3=e_2+\alpha e_3$ & $e_1e_4=e_3+\beta e_4$ & \\ \hline
\end{longtable}
}

\subsubsection{\texorpdfstring{$5$-dimensional nilpotent anticommutative $\mathfrak{CD}$-algebras}{5-dimensional nilpotent anticommutative CD-algebras}}

Using the algebraic classification of $5$-dimensional nilpotent binary Lie algebras~\cite{ack}, we see that all such algebras are anticommutative $\mathfrak{CD}$-algebras. Moreover, they are exactly all the $5$-dimensional nilpotent Malcev algebras, so their geometric classification can be deduced from the full description of degenerations of $5$-dimensional nilpotent Malcev algebras~\cite{kpv}:
\[
{\rm Irr}\left(\mathfrak{NACD}^5\right)=\left\{
\overline{{\mathcal O}\left(\mathfrak{NACD}_{i}^5\right)} \right\}_{i=1}^2,
\]
where

{\small

\begin{longtable}{|l|l|llll|}
\caption*{ }   \\

\hline \multicolumn{1}{|c|}{${\mathcal A}$} & \multicolumn{1}{c|}{ } & \multicolumn{4}{c|}{Multiplication table} \\ \hline 
\endfirsthead

 \multicolumn{6}{l}%
{{\bfseries  continued from previous page}} \\
\hline \multicolumn{1}{|c|}{${\mathcal A}$} & \multicolumn{1}{c|}{ } & \multicolumn{4}{c|}{Multiplication table} \\ \hline 
\endhead

\hline \multicolumn{6}{|r|}{{Continued on next page}} \\ \hline
\endfoot

\hline 
\endlastfoot

$\mathfrak{NACD}_{1}^5$ & $g_{5,6}$ & $e_1e_2=e_3$ & $e_1e_3=e_4$ & $e_1e_4=e_5$ & $e_2e_3=e_5$  \\ \hline
$\mathfrak{NACD}_{2}^5$ & $M_5$ & $e_1e_2=e_4$ & $e_3e_4=e_5$ && \\ \hline
\end{longtable}
}

\subsubsection{\texorpdfstring{$6$-dimensional nilpotent anticommutative $\mathfrak{CD}$-algebras}{6-dimensional nilpotent anticommutative CD-algebras}}

As it is explained in~\cite{ack}, checking the list of $6$-dimensional nilpotent binary Lie algebras yields that the $6$-dimensional nilpotent anticommutative $\mathfrak{CD}$-algebras are exactly the $6$-dimensional nilpotent Malcev algebras and $\mathbf{B}_{6,1}^{\alpha}$.
Then, the irreducible components of the variety $\mathfrak{NACD}^6$ are deduced as a corollary~\cite{ack} from those of $\mathfrak{NBL}^6$:
\[{\rm Irr}(\mathfrak{NACD}^6)=\left\{
\overline{{\mathcal O}\big(\mathfrak{NACD}_{i}^6\big)}\right\}_{i=1}^2 \cup \left\{ \overline{\bigcup {\mathcal O}\big(\mathfrak{NACD}_{3}^6(\alpha)\big)}\right\},\]
where

{\small

\begin{longtable}{|l|l|lllll|}
\caption*{ }   \\

\hline \multicolumn{1}{|c|}{${\mathcal A}$} & \multicolumn{1}{c|}{ } & \multicolumn{5}{c|}{Multiplication table} \\ \hline 
\endfirsthead

 \multicolumn{7}{l}%
{{\bfseries  continued from previous page}} \\
\hline \multicolumn{1}{|c|}{${\mathcal A}$} & \multicolumn{1}{c|}{ } & \multicolumn{5}{c|}{Multiplication table} \\ \hline 
\endhead

\hline \multicolumn{7}{|r|}{{Continued on next page}} \\ \hline
\endfoot

\hline 
\endlastfoot

$\mathfrak{NACD}_{1}^6$ & $g_6$ & $e_1e_2=e_3$ & $e_1e_3=e_4$ & $e_1e_4=e_5$ & $e_2e_3=e_5$ & $e_2e_5=e_6$  \\&& $e_3e_4=-e_6$ &&&&\\ \hline
$\mathfrak{NACD}_{2}^6$ & $\mathbf{B}_{6,1}^1$ & $e_{1}e_{2}=e_{4}$ & $e_{1}e_{3}=e_{5}$ & $e_{2}e_{3}=e_{6}$ & $e_{4}e_{5}=e_{6}$ & \\ \hline 
$\mathfrak{NACD}_{3}^6(\alpha)$ & $\mathbf{M}_6^{\alpha}$ & $e_1e_2=e_3$  &  $e_1e_3=  e_5$ & $e_1e_5=e_6$ & $e_2e_4=\alpha e_5$  & $e_3e_4=e_6$  \\ \hline
\end{longtable}
}

{

check

 \subsection{Symmetric Leibniz algebras} 
An   algebra $\mathfrak{L}$ is called  {\it symmetric Leibniz} if it satisfies the identities 
\[(xy)z=(xz)y+x(yz), \ x(yz)=(xy)z+y(xz).\] 
Let $\mathfrak{SLeib}$ denote the variety of symmetric Leibniz algebras.

\subsubsection{2-dimensional symmetric Leibniz algebras}
The algebraic classification of $2$-dimensional symmetric Leibniz algebras can be found in the work of Mohd Atan and Rakhimov (2012).
Analyzing the the graph of degenerations  of all $2$-dimensional algebras from~\cite{kv16}, we obtain the geometric classification of all $2$-dimensional symmetric Leibniz algebras:

\[{\rm Irr}\left(\mathfrak{SLeib}^2\right)=\left\{
\overline{{\mathcal O}\big(\mathfrak{SL}_{i}^2\big)} \right\}_{i=1}^2,\]
where

{\small
\vspace{-5mm}
\begin{longtable}{|l|l|ll|}
\caption*{ }   \\

\hline \multicolumn{1}{|c|}{${\mathcal A}$} & \multicolumn{1}{c|}{ } & \multicolumn{2}{c|}{Multiplication table} \\ \hline 
\endfirsthead

 \multicolumn{4}{l}%
{{\bfseries  continued from previous page}} \\
\hline \multicolumn{1}{|c|}{${\mathcal A}$} & \multicolumn{1}{c|}{ } & \multicolumn{2}{c|}{Multiplication table} \\ \hline 
\endhead

\hline \multicolumn{4}{|r|}{{Continued on next page}} \\ \hline
\endfoot

\hline 
\endlastfoot
$\mathfrak{SL}_{1}^2$ & ${\bf A}_3$ & $e_1e_1 = e_2$ &\\ \hline
$\mathfrak{SL}_{2}^2$ & ${\bf B}_3$ & $e_1e_2=  e_2$ & $e_2e_1=-e_2$ \\ \hline
\end{longtable}
}

\subsubsection{3-dimensional nilpotent symmetric Leibniz algebras}
The full graph of degenerations of Leibniz algebras in dimension $3$ was studied in~\cite{ikv18}. 
The restriction to the nilpotent symmetric Leibniz case gives the geometric classification of nilpotent symmetric Leibniz algebras.
The variety $\mathfrak{NSLeib}^3$  is irreducible:
\[{\rm Irr}\left(\mathfrak{NLeib}^3\right)=\left\{
\overline{{\mathcal O}\big(\mathfrak{NSL}_{1}^3\big)}\right\},\]
where

{\small
\vspace{-5mm}
\begin{longtable}{|l|l|lllll|}
\caption*{ }   \\

\hline \multicolumn{1}{|c|}{${\mathcal A}$} & \multicolumn{1}{c|}{ } & \multicolumn{5}{c|}{Multiplication table} \\ \hline 
\endfirsthead

 \multicolumn{7}{l}%
{{\bfseries  continued from previous page}} \\
\hline \multicolumn{1}{|c|}{${\mathcal A}$} & \multicolumn{1}{c|}{ } & \multicolumn{5}{c|}{Multiplication table} \\ \hline 
\endhead

\hline \multicolumn{7}{|r|}{{Continued on next page}} \\ \hline
\endfoot

\hline 
\endlastfoot

$\mathfrak{NSL}_{1}^3(\alpha)$ & $\mathfrak{L}^{\alpha}_1 $& 
$e_2e_2 =\alpha e_1$ & $e_3e_2=e_1$ &$e_3e_3 =e_1$ && \\ 
\end{longtable}

}

\subsubsection{3-dimensional symmetric Leibniz algebras}
The full graph of degenerations of Leibniz algebras in dimension $3$ was studied in~\cite{ikv18}. 
The restriction to the symmetric Leibniz case gives the geometric classification of symmetric Leibniz algebras.
The variety $\mathfrak{SLeib}^3$  has four irreducible components:
\[{\rm Irr}\left(\mathfrak{Leib}^3\right)=\left\{
\overline{{\mathcal O}\big(\mathfrak{SL}_{i}^3\big)}\right\}_{i=1}^{2}\cup
\left\{\overline{\bigcup {\mathcal O}\big(\mathfrak{SL}_{i}^3(\alpha)\big)} \right\}_{i=3}^{4},\]
where

{\small
\vspace{-5mm}
\begin{longtable}{|l|l|lllll|}
\caption*{ }   \\

\hline \multicolumn{1}{|c|}{${\mathcal A}$} & \multicolumn{1}{c|}{ } & \multicolumn{5}{c|}{Multiplication table} \\ \hline 
\endfirsthead

 \multicolumn{7}{l}%
{{\bfseries  continued from previous page}} \\
\hline \multicolumn{1}{|c|}{${\mathcal A}$} & \multicolumn{1}{c|}{ } & \multicolumn{5}{c|}{Multiplication table} \\ \hline 
\endhead

\hline \multicolumn{7}{|r|}{{Continued on next page}} \\ \hline
\endfoot

\hline 
\endlastfoot
$\mathfrak{SL}_{1}^3$ & $\mathfrak{g}_4$&
$e_1e_2 =e_3$ &$e_1e_3=-e_2$ & $e_2e_1=-e_3$ &&\\&&
$e_2e_3=e_1$ & $e_3e_1=e_2$ & $e_3e_2=-e_1$&&\\
 \hline

$\mathfrak{SL}_{2}^3$ &$\mathfrak{L}^{0}_4$ & 
$e_2e_3=e_2$ &$e_3e_2=-e_2$ & $e_3e_3=e_1$  && \\ \hline

$\mathfrak{SL}_{3}^3(\alpha)$ & 
$\mathfrak{g}^{\alpha}_3$&  $e_1e_3 =e_1+e_2$ & $e_2e_3=\alpha e_2$&
$e_3e_1 =-e_1-e_2$ & $e_3e_2=-\alpha e_2$ &\\
  \hline

$\mathfrak{SL}_{4}^3(\alpha)$ & $\mathfrak{L}^{\alpha}_1 $& 
$e_2e_2 =\alpha e_1$ & $e_3e_2=e_1$ &$e_3e_3 =e_1$ && \\ 
\end{longtable}

}

\subsubsection{4-dimensional nilpotent symmetric Leibniz algebras}
The classification of the nilpotent symmetric Leibniz algebras of dimension $4$ can be found in a paper by Alvarez and Kaygorodov \cite{ak21}.
Their geometric classification appeared in~\cite{ak21}. In particular, the variety $\mathfrak{NSLeib}^4$  has three irreducible components:

\[{\rm Irr}\left(\mathfrak{NSLeib}^4\right)= \left\{\overline{{\mathcal O}\left(\mathfrak{NSL}_{1}^4\right)}\right\} \cup
\left\{ \overline{\bigcup {\mathcal O}\left(\mathfrak{NSL}_{i}^4(\alpha)\right)}\right\}_{i=2}^{3},\]
where

{\small
\vspace{-5mm}
\begin{longtable}{|l|l|lllll|}
\caption*{ }   \\

\hline \multicolumn{1}{|c|}{${\mathcal A}$} & \multicolumn{1}{c|}{ } & \multicolumn{5}{c|}{Multiplication table} \\ \hline 
\endfirsthead

 \multicolumn{7}{l}%
{{\bfseries  continued from previous page}} \\
\hline \multicolumn{1}{|c|}{${\mathcal A}$} & \multicolumn{1}{c|}{ } & \multicolumn{5}{c|}{Multiplication table} \\ \hline 
\endhead

\hline \multicolumn{7}{|r|}{{Continued on next page}} \\ \hline
\endfoot

\hline 
\endlastfoot

$\mathfrak{NSL}_{1}^4(\alpha)$ & $\mathfrak{S}_{01}$ & $ e_1e_1=e_4$ & $e_1 e_2 =e_3$ & $e_2 e_1 =-e_3$ &&\\ && $e_2e_2=e_4$ & $e_2e_3 =e_4$ & $e_3e_2=-e_4$&&\\
\hline

$\mathfrak{NSL}_{2}^4(\alpha)$ &  ${\mathfrak N}_2(\alpha)$  &  $e_1e_1 = e_3$ & $e_1e_2 = e_4$ & $e_2e_1 = -\alpha e_3$ & $e_2e_2 = -e_4$ & \\
\hline

$\mathfrak{NSL}_{3}^4(\alpha)$ & ${\mathfrak N}_3(\alpha)$ &   $e_1e_1 = e_4$ & $e_1e_2 = \alpha e_4$ & $e_2e_1 = -\alpha e_4$ & $e_2e_2 = e_4$ & $e_3e_3=e_4$
\end{longtable}
}

\subsubsection{4-dimensional symmetric Leibniz algebras}
The classification of the symmetric Leibniz algebras of dimension $4$ can be found in a paper by Alvarez and Kaygorodov \cite{ak21}.
Their geometric classification appeared in~\cite{ak21}. In particular, the variety $\mathfrak{SLeib}^4$  has five irreducible components:

\[{\rm Irr}\left(\mathfrak{SLeib}^4\right)= \left\{\overline{{\mathcal O}\left(\mathfrak{SL}_{1}^4\right)}\right\} \cup
\left\{ \overline{\bigcup {\mathcal O}\left(\mathfrak{SL}_{i}^4(\alpha)\right)}\right\}_{i=2}^{5},\]
where

{\small

\begin{longtable}{|l|l|lllll|}
\caption*{ }   \\

\hline \multicolumn{1}{|c|}{${\mathcal A}$} & \multicolumn{1}{c|}{ } & \multicolumn{5}{c|}{Multiplication table} \\ \hline 
\endfirsthead

 \multicolumn{7}{l}%
{{\bfseries  continued from previous page}} \\
\hline \multicolumn{1}{|c|}{${\mathcal A}$} & \multicolumn{1}{c|}{ } & \multicolumn{5}{c|}{Multiplication table} \\ \hline 
\endhead

\hline \multicolumn{7}{|r|}{{Continued on next page}} \\ \hline
\endfoot

\hline 
\endlastfoot

$\mathfrak{SL}_{1}^4(\alpha)$ & $\mathfrak{L}_{02}$ & $
  e_1e_1=e_4$ & $e_1e_2=-e_2$ & $e_1e_3=e_3$ & $e_2e_1=e_2$ & \\
  && $e_2e_3=e_4$ &  $e_3e_1=-e_3$ &  $e_3e_2=-e_4$ &&\\ 
\hline

$\mathfrak{SL}_{2}^4(\alpha)$ & $\mathfrak{L}_{15}^{\alpha}$ &  $  
e_1e_1=\alpha e_4$ & $e_1e_2=e_4$ & $e_1e_3=-e_3$ &  $e_2e_2=e_4$ &   $e_3e_1=e_3$\\ 
\hline

$\mathfrak{SL}_{3}^4(\alpha)$ & $\mathfrak{L}_{24}^{\alpha}$ & 
$e_1e_1=e_4$ & $e_1e_2= - e_2$ &  $ e_1e_3= -\alpha e_3$ & $e_2e_1=e_2$ & $e_3e_1= \alpha e_3$\\ 
\hline

$\mathfrak{SL}_{4}^4(\alpha)$ &  ${\mathfrak N}_2(\alpha)$  &  $e_1e_1 = e_3$ & $e_1e_2 = e_4$ & $e_2e_1 = -\alpha e_3$ & $e_2e_2 = -e_4$ & \\
\hline

$\mathfrak{SL}_{5}^4(\alpha)$ & ${\mathfrak N}_3(\alpha)$ &   $e_1e_1 = e_4$ & $e_1e_2 = \alpha e_4$ & $e_2e_1 = -\alpha e_4$ & $e_2e_2 = e_4$ & $e_3e_3=e_4$
\end{longtable}
}

\subsubsection{5-dimensional nilpotent symmetric Leibniz algebras}
The classification of the nilpotent symmetric Leibniz algebras of dimension $5$ can be found in a paper by Alvarez and Kaygorodov \cite{ak21}.
Their geometric classification appeared in~\cite{ak21}. In particular, the variety $\mathfrak{NSLeib}^4$  has six irreducible components:

\[{\rm Irr}\left(\mathfrak{NSLeib}^5\right)= \left\{\overline{{\mathcal O}\left(\mathfrak{NSL}_{i}^5(\alpha)\right)}\right\}_{i=1}^{3} \cup
\left\{ \overline{\bigcup {\mathcal O}\left(\mathfrak{NSL}_{i}^5(\alpha,\beta)\right)}\right\}_{i=4}^{5}\cup
\left\{ \overline{\bigcup {\mathcal O}\left(\mathfrak{NSL}_{6}^5(\overline{\mu})\right)}\right\},\]
where

{\small

\begin{longtable}{|l|l|llll|}
\caption*{ }   \\

\hline \multicolumn{1}{|c|}{${\mathcal A}$} & \multicolumn{1}{c|}{ } & \multicolumn{4}{c|}{Multiplication table} \\ \hline 
\endfirsthead

 \multicolumn{6}{l}%
{{\bfseries  continued from previous page}} \\
\hline \multicolumn{1}{|c|}{${\mathcal A}$} & \multicolumn{1}{c|}{ } & \multicolumn{4}{c|}{Multiplication table} \\ \hline 
\endhead

\hline \multicolumn{6}{|r|}{{Continued on next page}} \\ \hline
\endfoot

\hline 
\endlastfoot

$\mathfrak{NSL}_{1}^5(\alpha)$ & $\mathbb{S}_{22}^{\alpha}$ & $e_1e_1=e_5$ &  $e_1 e_2 =e_3$  & $e_1e_3=e_5$ & $e_2 e_1 =-e_3$ \\&& $e_2e_2=\alpha e_5$ &  $e_2e_4=e_5$ & $e_3e_1=-e_5$ & $e_4e_4=e_5$\\
\hline

$\mathfrak{NSL}_{2}^5(\alpha)$ & $\mathbb{S}_{41}^{\alpha}$ & $e_1e_1=e_5$ & $e_1 e_2 =e_3$  & $e_1e_3=e_5$ & $e_2 e_1 =-e_3$ \\&& 
$e_2e_2=\alpha e_5$  & $e_2e_3=e_4$ & $e_2e_4=e_5$ &\\&& 
$e_3e_1=-e_5$ & $e_3e_2=-e_4$ & $e_4e_2=-e_5$&\\
\hline

 $\mathfrak{NSL}_{3}^5(\alpha)$ &${\mathfrak V}_{2+3}$ & 
$e_1e_1 = e_3 + \alpha e_5$& $e_1e_2 = e_3$ & $e_2e_1 = e_4$& $e_2e_2 = e_5$\\
\hline 

$\mathfrak{NSL}_{4}^5(\alpha,\beta)$ & $\mathbb{S}_{21}^{\alpha,\beta}$ & $e_1e_1=\alpha e_5$ &  $e_1 e_2 =e_3+e_4+\beta e_5$  & $e_1e_3=e_5$ & $e_2 e_1 =-e_3$ \\&& 
$e_2e_2=e_5$  & $e_2e_3 =e_4$ & $e_3e_1=-e_5$ & $e_3e_2=-e_4$\\
\hline 

$\mathfrak{NSL}_{5}^5(\alpha,\beta)$ & ${\mathfrak V}_{4+1}$ &   
$e_1e_2=e_5$& $e_2e_1=\alpha e_5$ &$e_3e_4=e_5$&$e_4e_3=\beta e_5$\\
\hline

$\mathfrak{NSL}_{6}^5(\overline{\mu})$ &  ${\mathfrak V}_{3+2}$ & 
$e_1e_1 =  e_4$& $e_1e_2 = \mu_1 e_5$ & $e_1e_3 =\mu_2 e_5$ &\\&& 
$e_2e_1 = \mu_3 e_5$  & $e_2e_2 = \mu_4 e_5$  
 & $e_2e_3 = \mu_5 e_5$  & \\&& 
 $e_3e_1 = \mu_6 e_5$  &   $e_3e_2 = \mu_0 e_4+ \mu_7 e_5$   & $e_3e_3 =  e_5$ &  \\

\hline

\end{longtable}
}

}

 \subsection{Leibniz algebras} 
An   algebra $\mathfrak{L}$ is called  {\it Leibniz} if it satisfies the identity 
\[(xy)z=(xz)y+x(yz).\] Let $\mathfrak{Leib}$ denote the variety of Leibniz algebras.

\subsubsection{2-dimensional Leibniz algebras}
The algebraic classification of $2$-dimensional Leibniz algebras can be found in the work of Mohd Atan and Rakhimov (2012).
Analyzing the the graph of degenerations  of all $2$-dimensional algebras from~\cite{kv16}, we obtain the geometric classification of all $2$-dimensional Leibniz algebras:

\[{\rm Irr}\left(\mathfrak{Leib}^2\right)=\left\{
\overline{{\mathcal O}\big(\mathfrak{L}_{i}^2\big)} \right\}_{i=1}^2,\]
where

{\small

\begin{longtable}{|l|l|ll|}
\caption*{ }   \\

\hline \multicolumn{1}{|c|}{${\mathcal A}$} & \multicolumn{1}{c|}{ } & \multicolumn{2}{c|}{Multiplication table} \\ \hline 
\endfirsthead

 \multicolumn{4}{l}%
{{\bfseries  continued from previous page}} \\
\hline \multicolumn{1}{|c|}{${\mathcal A}$} & \multicolumn{1}{c|}{ } & \multicolumn{2}{c|}{Multiplication table} \\ \hline 
\endhead

\hline \multicolumn{4}{|r|}{{Continued on next page}} \\ \hline
\endfoot

\hline 
\endlastfoot
$\mathfrak{L}_{1}^2$ & ${\bf B}_2(0)$ & $e_1e_2 = e_1$ &\\ \hline
$\mathfrak{L}_{2}^2$ & ${\bf B}_3$ & $e_1e_2=  e_2$ & $e_2e_1=-e_2$ \\ \hline
\end{longtable}
}

\subsubsection{3-dimensional nilpotent Leibniz algebras}
The algebraic and geometric classification of $3$-dimensional nilpotent Leibniz algebras can be obtained from the classification and description of degenerations of $3$-dimensional nilpotent algebras given in \cite{fkkv}.
Hence, we have  that the variety   $\mathfrak{NLeib}^3$ has two irreducible components:
\[{\rm Irr}(\mathfrak{NAss}^3)=
\left\{\overline{{\mathcal O}\left(\mathfrak{NL}_{i}^{3}\right)}\right\}_{i=1}^2\cup\left\{
\overline{\bigcup {\mathcal O}\left(\mathfrak{NL}_{3}^3(\alpha)\right)}\right\},\]
where

{\small
\vspace{-4mm}
\begin{longtable}{|l|l|llllll|}
\caption*{ }   \\

\hline \multicolumn{1}{|c|}{${\mathcal A}$} & \multicolumn{1}{c|}{ } & \multicolumn{6}{c|}{Multiplication table} \\ \hline 
\endfirsthead

 \multicolumn{8}{l}%
{{\bfseries  continued from previous page}} \\
\hline \multicolumn{1}{|c|}{${\mathcal A}$} & \multicolumn{1}{c|}{ } & \multicolumn{6}{c|}{Multiplication table} \\ \hline 
\endhead

\hline \multicolumn{8}{|r|}{{Continued on next page}} \\ \hline
\endfoot

\hline 
\endlastfoot
$\mathfrak{NL}_{1}^3$ & $\mathcal{N}_3$ &      $e_1 e_1 = e_2$ & $e_2 e_1=e_3$ & && &  \\\hline

$\mathfrak{NL}_{2}^3$ & $\mathcal{N}_6$ & $e_1 e_1 = e_3$ & $e_2 e_2=e_3$  &&& &\\ \hline

$\mathfrak{NL}_{3}^3(\alpha)$ &  $\mathcal{N}_{8}(\alpha)$  & $e_1e_1=\alpha e_3$ & $e_2e_1=e_3$ & $e_2e_2=e_3$ && &\\ \hline
\end{longtable}
}

\subsubsection{3-dimensional Leibniz algebras}
Dimension $3$ was studied in~\cite{ikv18}, employing the algebraic classification by Mohd Atan and  Rakhimov (2012).
The variety $\mathfrak{Leib}^3$  has five irreducible components:
\[{\rm Irr}\left(\mathfrak{Leib}^3\right)=\left\{
\overline{{\mathcal O}\big(\mathfrak{L}_{i}^3\big)}\right\}_{i=1}^{2}\cup
\left\{\overline{\bigcup {\mathcal O}\big(\mathfrak{L}_{i}^3(\alpha)\big)} \right\}_{i=3}^{5},\]
where

{\small
\vspace{-4mm}
\begin{longtable}{|l|l|lllll|}
\caption*{ }   \\

\hline \multicolumn{1}{|c|}{${\mathcal A}$} & \multicolumn{1}{c|}{ } & \multicolumn{5}{c|}{Multiplication table} \\ \hline 
\endfirsthead

 \multicolumn{7}{l}%
{{\bfseries  continued from previous page}} \\
\hline \multicolumn{1}{|c|}{${\mathcal A}$} & \multicolumn{1}{c|}{ } & \multicolumn{5}{c|}{Multiplication table} \\ \hline 
\endhead

\hline \multicolumn{7}{|r|}{{Continued on next page}} \\ \hline
\endfoot

\hline 
\endlastfoot
$\mathfrak{L}_{1}^3$ & $\mathfrak{g}_4$ & $e_1e_2 =e_3$ & $e_1e_3=-e_2$ & $e_2e_1=-e_3$ &&\\&& $e_2e_3=e_1$  & $e_3e_1=e_2$  & $e_3e_2=-e_1$ && \\ \hline
$\mathfrak{L}_{2}^3$ & $\mathfrak{L}_5$ &  $e_1e_3=2e_1$ &  $e_2e_2=e_1$&  $e_2e_3=e_2$ &  $e_3e_2=-e_2$ &  $e_3e_3=e_1$  \\ \hline
$\mathfrak{L}_{3}^3(\alpha)$ & $\mathfrak{g}^{\alpha}_3$ & $e_1e_3 =e_1+e_2$ & $e_2e_3=\alpha e_2$ & $e_3e_1=-e_1-e_2$ & $e_3e_2=-\alpha e_2$ & \\ \hline
$\mathfrak{L}_{4}^3(\alpha)$ & $\mathfrak{L}^{\alpha}_4$ & $e_1e_3=\alpha e_1$ & $e_2e_3=e_2$ & $e_3e_2=-e_2$ & $e_3e_3=e_1$ & \\ \hline
$\mathfrak{L}_{5}^3(\alpha)$ & $\mathfrak{L}^{\alpha}_6$ & $e_1e_3=\alpha e_1$ & $e_2e_3=e_2$ &&& \\ \hline
\end{longtable}
}

\subsubsection{4-dimensional Leibniz algebras}
The classification of the Leibniz algebras of dimension $4$ is a joint result of the works by Albeverio, Omirov, and Rakhimov (2006), Ca\~nete and Khudoyberdiyev (2013), and Omirov, Rakhimov and Turdibaev (2013). Later, their geometric classification appeared in~\cite{ikv17}. In particular, the variety $\mathfrak{Leib}^4$  has seventeen irreducible components:

\[{\rm Irr}\left(\mathfrak{Leib}^4\right)= \left\{\overline{{\mathcal O}\left(\mathfrak{L}_{i}^4\right)}\right\}_{i=1}^{6}\cup
\left\{ \overline{\bigcup {\mathcal O}\left(\mathfrak{L}_{i}^4(\alpha)\right)}\right\}_{i=7}^{13}\cup
\left\{\overline{\bigcup {\mathcal O}\left(\mathfrak{L}_{i}^4(\alpha,\beta)\right)} \right\}_{i=14}^{17},\]
where

{\small

\vspace{-4mm}
\begin{longtable}{|l|l|lll|}
\caption*{ }   \\

\hline \multicolumn{1}{|c|}{${\mathcal A}$} & \multicolumn{1}{c|}{ } & \multicolumn{3}{c|}{Multiplication table} \\ \hline 
\endfirsthead

 \multicolumn{5}{l}%
{{\bfseries  continued from previous page}} \\
\hline \multicolumn{1}{|c|}{${\mathcal A}$} & \multicolumn{1}{c|}{ } & \multicolumn{3}{c|}{Multiplication table} \\ \hline 
\endhead

\hline \multicolumn{5}{|r|}{{Continued on next page}} \\ \hline
\endfoot

\hline 
\endlastfoot
$\mathfrak{L}_{1}^4$ & $\mathfrak{sl}_2\oplus\mathbb{C}$ & $e_1e_2=e_2$ & $e_1e_3=-e_3$ & $e_2e_1=-e_2$ \\&& $e_2e_3=e_1$  & $e_3e_1=e_3$& $e_3e_2=-e_1$ \\ \hline
$\mathfrak{L}_{2}^4$ & $\mathfrak{R}_1$ & $e_3e_1 = e_3$ & $e_4e_2=e_4$ & \\ \hline
$\mathfrak{L}_{3}^4$ & $\mathfrak{R}_2$ & $e_1e_2=e_2$ & $e_2e_1=-e_2$ & $e_3e_4=e_4$ \\&& $e_4e_3=-e_4$&&  \\ \hline
$\mathfrak{L}_{4}^4$ & $\mathfrak{R}_3$ & $e_2e_4=-e_4$ & $e_3e_1 = e_3$ & $e_4e_2=e_4$   \\ \hline
$\mathfrak{L}_{5}^4$ & $\mathfrak{L}_2$ & $e_1e_1=e_4$ & $e_1e_2=-e_2$ & $e_1e_3=e_3$ \\&& $e_2e_1=e_2$ & $e_2e_3=e_4$ & $e_3e_1=-e_3$ \\&&  $e_3e_2=-e_4$ && \\ \hline
$\mathfrak{L}_{6}^4$ & $\mathfrak{L}_{44}$ & $e_1e_2 =- e_2$ & $e_2e_1=  e_2$ & $e_2e_2=e_3$ \\&& $e_3e_1= 2e_3$ & $e_3e_2= e_4$ & $e_4e_1=3e_4$ \\ \hline
$\mathfrak{L}_{7}^4(\alpha)$ & $g_5(\alpha)$ & $e_1e_2=e_2$ & $e_1e_3=e_2+\alpha e_3$ & $e_1e_4=(\alpha+1)e_4$ \\&& $e_2e_1=-e_2$ & $e_2e_3=e_4$ & $e_3e_1=-e_2-\alpha e_3$ \\&& $e_3e_2=-e_4$ & $e_4e_1=-(\alpha+1)e_4$ & \\ \hline
$\mathfrak{L}_{8}^4(\alpha)$ & $\mathfrak{L}_4^{\alpha}$ & $e_1e_2=-e_2$ & $e_2e_1=e_2$ & $e_3e_1= \alpha  e_3$ \\&& $e_3e_2=e_4$ & $e_4e_1=(1+\alpha)e_4$ & \\ \hline
$\mathfrak{L}_{9}^4(\alpha)$ & $\mathfrak{L}_8^{\alpha}$ & $e_1e_2=-e_2$ &  $e_1e_3=-\alpha e_3$ & $e_2e_1=e_2$ \\&& $e_2e_3=\alpha e_4$ & $e_3e_1=\alpha  e_3$ & $e_3e_2=e_4$ \\&&  $e_4e_1=(\alpha+1)e_4$ & & \\ \hline
$\mathfrak{L}_{10}^4(\alpha)$ & $\mathfrak{L}_9^{\alpha}$ & $e_1e_2=-e_2$ & $e_1e_3=-\alpha e_3$ &  $e_2e_1=e_2$ \\&& $e_2e_2=e_4$ &  $e_3e_1=\alpha e_3$ & $e_4e_1=2e_4$ \\ \hline
$\mathfrak{L}_{11}^4(\alpha)$ & $\mathfrak{L}_{10}^{\alpha}$ & $e_1e_2=-e_2$ & $e_2e_1=e_2$ & $e_2e_2=e_4$ \\&& $e_3e_1=\alpha e_3$ & $e_4e_1=2e_4$ & \\ \hline
$\mathfrak{L}_{12}^4(\alpha)$ & $\mathfrak{L}_{15}^{\alpha}$ & $e_1e_1=\alpha e_4$ & $e_1e_2=e_4$ & $e_1e_3=-e_3$ \\&& $e_2e_2=e_4$ & $e_3e_1=e_3$ &  \\ \hline
$\mathfrak{L}_{13}^4(\alpha)$ & $\mathfrak{L}_{18}^{\alpha}$ & $e_1e_1=\alpha  e_4$ & $e_1e_2=e_4$ & $e_2e_2=e_4$ \\&& $e_3e_1=e_3$ && \\ \hline
$\mathfrak{L}_{14}^4(\alpha,\beta)$ & $g_4(\alpha,\beta)$ & $e_1e_2=e_2$ & $e_1e_3=e_2+\alpha e_3$ & $e_1e_4=e_3+\beta e_4$ \\&& $e_2e_1=-e_2$  & $e_3e_1=-e_2-\alpha e_3$ & $e_4e_1=-e_3-\beta e_4$ \\ \hline
$\mathfrak{L}_{15}^4(\alpha,\beta)$ & $\mathfrak{L}_{21}^{\alpha,\beta}$ & $e_1e_2=-e_2$ & $e_1e_3=- \alpha  e_3$ & $e_2e_1=e_2$ \\&& $e_3e_1= \alpha  e_3$ &  $e_4e_1= \beta  e_4$ & \\ \hline
$\mathfrak{L}_{16}^4(\alpha,\beta)$ & $\mathfrak{L}_{22}^{\alpha,\beta}$ & $e_1e_2=-e_2$ & $e_2e_1=e_2$ & $e_3e_1= \alpha e_3$ \\&& $e_4e_1=\beta e_4$ && \\ \hline 
$\mathfrak{L}_{17}^4(\alpha,\beta)$ & $\mathfrak{L}_{23}^{\alpha,\beta}$ & $e_2e_1=e_2$ & $e_3e_1= \alpha e_3$ &  $e_4e_1= \beta e_4$  \\ \hline
\end{longtable}
}

Focusing on $4$-dimensional nilpotent Leibniz algebras, whose algebraic classification was ultimately given by Albeverio, Omirov and Rakhimov (2006), we find in~\cite{kppv} the following geometric classification:
\[{\rm Irr}\left(\mathfrak{NLeib}^4\right)= \left\{\overline{{\mathcal O}\left(\mathfrak{NL}_{i}^4\right)}\right\}_{i=1}^{3}\cup
\left\{ \overline{\bigcup {\mathcal O}\left(\mathfrak{NL}_{4}^4(\alpha)\right)}\right\},\]
where

{\small

\begin{longtable}{|l|l|lllll|}
\caption*{ }   \\

\hline \multicolumn{1}{|c|}{${\mathcal A}$} & \multicolumn{1}{c|}{ } & \multicolumn{5}{c|}{Multiplication table} \\ \hline 
\endfirsthead

 \multicolumn{7}{l}%
{{\bfseries  continued from previous page}} \\
\hline \multicolumn{1}{|c|}{${\mathcal A}$} & \multicolumn{1}{c|}{ } & \multicolumn{5}{c|}{Multiplication table} \\ \hline 
\endhead

\hline \multicolumn{7}{|r|}{{Continued on next page}} \\ \hline
\endfoot

\hline 
\endlastfoot
$\mathfrak{NL}_{1}^4$ & $\mathfrak{L}_2$ & $e_1e_1 = e_2$ & $e_2e_1 = e_3$  & $e_3e_1 = e_4$ && \\ \hline
$\mathfrak{NL}_{2}^4$ & $\mathfrak{L}_5$ & $e_1e_1 = e_3$ & $e_2e_1 = e_3$ & $e_2e_2 = e_4$ & $e_3e_1 = e_4$ & \\ \hline
$\mathfrak{NL}_{3}^4$ & $\mathfrak{L}_{11}$ & $e_1e_1=e_4$ & $e_1e_2 = -e_3$ & $e_1e_3 = -e_4$ &&  \\&&$e_2e_1 = e_3$ & $e_2e_2 = e_4$ & $e_3e_1=e_4$ &&  \\ \hline
$\mathfrak{NL}_{4}^4(\alpha)$ & $\mathfrak{N}_3(\alpha)$  & $e_1e_1 = e_4$ & $e_1e_2 = \alpha e_4$ & $e_2e_1 = -\alpha e_4$ & $e_2e_2 = e_4$ & $e_3e_3 = e_4$  \\ \hline
\end{longtable}
}

{

\subsubsection{5-dimensional nilpotent Leibniz algebras}
The classification of nilpotent Leibniz algebras of dimension $5$ is a  result of the work by Abdurasulov,   Kaygorodov, and Khudoyberdiyev  (2023). Their geometric classification appeared in~\cite{akk23}. In particular, the variety $\mathfrak{NLeib}^5$  has ten irreducible components:
\begin{align*}
&{\rm Irr}\left(\mathfrak{NLeib}^5\right) = \\
&\left\{\overline{{\mathcal O}\left(\mathfrak{NL}_{1}^5\right)}\right\}\cup
\left\{ \overline{\bigcup {\mathcal O}\left(\mathfrak{NL}_{i}^5(\alpha)\right)}\right\}_{i=2}^{5} \cup \left\{\overline{\bigcup {\mathcal O}\left(\mathfrak{NL}_{i}^5(\alpha,\beta)\right)} \right\}_{i=6}^{9} \cup
\left\{\overline{\bigcup {\mathcal O}\left(\mathfrak{NL}_{i}^5(\overline{\mu})\right)} \right\}.
\end{align*}

where

{\small

\begin{longtable}{|l|l|llll|}
\caption*{ }   \\

\hline \multicolumn{1}{|c|}{${\mathcal A}$} & \multicolumn{1}{c|}{ } & \multicolumn{4}{c|}{Multiplication table} \\ \hline 
\endfirsthead

 \multicolumn{6}{l}%
{{\bfseries  continued from previous page}} \\
\hline \multicolumn{1}{|c|}{${\mathcal A}$} & \multicolumn{1}{c|}{ } & \multicolumn{4}{c|}{Multiplication table} \\ \hline 
\endhead

\hline \multicolumn{6}{|r|}{{Continued on next page}} \\ \hline
\endfoot

\hline 
\endlastfoot

$\mathfrak{NL}_{1}^5$ &$\mathbb{L}_{82}$ &   $e_1e_1=e_2$ &  $e_2e_1=e_3$ & $e_3e_1=e_4$ & $e_4e_1=e_5$\\
\hline

$\mathfrak{NL}_{2}^5(\alpha)$ & $\mathbb{L}_{28}^{\alpha}$ &   $e_1e_1=e_3$  & $e_1e_2=e_3$ & $e_1e_4=\alpha e_5$ &  $e_2e_2=e_5$ \\

 \hline

$\mathfrak{NL}_{3}^5(\alpha)$ &$\mathbb{L}_{79}^{\alpha}$ &  $e_1e_1=e_3$ & $e_1e_2=e_4$ & $e_2e_1=e_3$   & $e_2e_2=e_4+e_5$ \\
&&  $e_3e_1=e_4+\alpha e_5$ & $e_3e_2=e_5$ & $e_4e_1=e_5$ & \\
\hline

$\mathfrak{NL}_{4}^5(\alpha)$ & $\mathbb{S}_{22}^{\alpha}$ & $e_1e_1=e_5$ &  $e_1 e_2 =e_3$  & $e_1e_3=e_5$ & $e_2 e_1 =-e_3$ \\&& $e_2e_2=\alpha e_5$  &  $e_2e_4=e_5$ & $e_3e_1=-e_5$ & $e_4e_4=e_5$\\
\hline

$\mathfrak{NL}_{5}^5(\alpha)$ & $\mathbb{S}_{41}^{\alpha}$ &  $e_1e_1=e_5$ & $e_1 e_2 =e_3$  & $e_1e_3=e_5$ & $e_2 e_1 =-e_3$ \\
&& $e_2e_2=\alpha e_5$  & $e_2e_3=e_4$ & $e_2e_4=e_5$&\\
& &  $e_3e_1=-e_5$ & $e_3e_2=-e_4$ & $e_4e_2=-e_5$&\\
\hline

$\mathfrak{NL}_{6}^5(\alpha,\beta)$ &$\mathbb{L}_{47}^{\alpha, \beta}$ &   $e_1e_1=e_3$  & $e_1e_2=e_4$ & $e_2e_1=-\alpha e_3$ &\\
&&  $e_2e_2=-e_4$ & $e_3e_1=e_5$ & $e_4e_2=\beta e_5$ &  \\
\hline

$\mathfrak{NL}_{7}^5(\alpha,\beta)$ &$\mathbb{L}_{52}^{\alpha, \beta}$ &  $e_1e_2=e_3$ & $e_1e_3=-e_5$ & $e_1e_4=e_5$ & $e_2e_1=e_4$ \\
&&  $e_2e_3=\beta e_5$  & $e_2e_4=-\beta e_5$ & $e_3e_1=e_5$ & \\
&& $e_3e_2=e_5$ &  $e_4e_1=\alpha e_5$&  $e_4e_2=\beta e_5$  & \\
\hline

$\mathfrak{NL}_{8}^5(\alpha,\beta)$ & $\mathbb{S}_{21}^{\alpha,\beta}$  & 
$e_1e_1=\alpha e_5$ &  $e_1 e_2 =e_3+e_4+\beta e_5$  & $e_1e_3=e_5$ & $e_2 e_1 =-e_3$ \\&& $e_2e_2=e_5$  & $e_2e_3 =e_4$ & $e_3e_1=-e_5$ & $e_3e_2=-e_4$\\
\hline 

$\mathfrak{NL}_{9}^5(\alpha, \beta)$ &${\mathfrak V}_{4+1}$ &    
$e_1e_2=e_5$& $e_2e_1=\alpha e_5$ &$e_3e_4=e_5$&$e_4e_3=\beta e_5$\\

 \hline

$\mathfrak{NL}_{10}^5(\overline{\mu})$ &${\mathfrak V}_{3+2}$ & 
$e_1e_1 =  e_4$& $e_1e_2 = \mu_1 e_5$ & $e_1e_3 =\mu_2 e_5$&\\
&&
$e_2e_1 = \mu_3 e_5$  & $e_2e_2 = \mu_4 e_5$   & $e_2e_3 = \mu_5 e_5$  &\\
&&$e_3e_1 = \mu_6 e_5$  &$e_3e_2 = \mu_0 e_4+ \mu_7 e_5$  & $e_3e_3 =  e_5$  &\\
 \hline
\end{longtable}
}
 }

 
 \subsection{Zinbiel algebras} 
An   algebra $\mathfrak{Z}$ is called  {\it Zinbiel} if it satisfies the identity 
\[(xy)z=x(yz+zy).\] We will denote this variety by $\mathfrak{Zin}$.

\subsubsection{2-dimensional Zinbiel algebras}
Dzhumadildaev and Tulenbaev proved in 2005 that every finite-dimensional   Zinbiel algebra is nilpotent.
Also, the lists of Zinbiel algebras of dimension $2$ and $3$ were given in that paper. In fact, there is just one $2$-dimensional  Zinbiel algebra, namely $\mathfrak{Z}_{1}^2$ with $e_1e_1=e_2$.

\subsubsection{3-dimensional Zinbiel algebras}

Regarding dimension $3$, the geometric classification can be extracted from~\cite{kppv}. We find three irreducible components in $\mathfrak{Zin}^3$:
\[{\rm Irr}(\mathfrak{Zin}^3)=
\left\{\overline{{\mathcal O}\left(\mathfrak{Z}_{i}^3\right)}\right\}_{i=1}^2
\cup
\left\{\overline{\bigcup {\mathcal O}\left(\mathfrak{Z}_{3}^3(\alpha)\right)}\right\},\]
where

{\small

\begin{longtable}{|l|l|lll|}
\caption*{ }   \\

\hline \multicolumn{1}{|c|}{${\mathcal A}$} & \multicolumn{1}{c|}{ } & \multicolumn{3}{c|}{Multiplication table} \\ \hline 
\endfirsthead

 \multicolumn{5}{l}%
{{\bfseries  continued from previous page}} \\
\hline \multicolumn{1}{|c|}{${\mathcal A}$} & \multicolumn{1}{c|}{ } & \multicolumn{3}{c|}{Multiplication table} \\ \hline 
\endhead

\hline \multicolumn{5}{|r|}{{Continued on next page}} \\ \hline
\endfoot

\hline 
\endlastfoot

$\mathfrak{Z}_{1}^3$ & $\mathfrak{Z}_1^{\mathbb{C}}$ & $e_1e_1 = e_2$ &  e$_1e_2 =\frac{1}{2} e_3$ & $e_2e_1=e_3$ \\ \hline
$\mathfrak{Z}_{2}^3$ & $\mathfrak{N}_3^{\mathbb{C}}$ & $e_1e_1=  e_3$ & $e_1e_2=e_3$ & $e_2e_1=e_3$  \\ \hline
$\mathfrak{Z}_{3}^3(\alpha)$ & $\mathfrak{N}_2^{\mathbb{C}}(\alpha)$ & $e_1e_1=  e_3$ & $e_1e_2=e_3$ & $e_2e_2=\alpha e_3$  \\ \hline
\end{longtable}
}

\subsubsection{4-dimensional Zinbiel algebras}
The algebraic classification of $4$-dimensional Zinbiel algebras is given in a paper by Adashev, Khudoyberdiyev and Omirov (2010).
After that, it was constructed the graph of degenerations of this variety $\mathfrak{Zin}^4$  in~\cite{kppv}.
In particular, there exist five irreducible components in $\mathfrak{Zin}^4$: 
\[{\rm Irr}(\mathfrak{Zin}^4)=
\left\{\overline{{\mathcal O}\big(\mathfrak{Z}_{i}^4\big)}\right\}_{i=1}^3 
\cup
\left\{\overline{\bigcup {\mathcal O}\big(\mathfrak{Z}_{i}^4(\alpha)\big)}\right\}_{i=4}^5,\]
where

{\small

\begin{longtable}{|l|l|lllll|}
\caption*{ }   \\

\hline \multicolumn{1}{|c|}{${\mathcal A}$} & \multicolumn{1}{c|}{ } & \multicolumn{5}{c|}{Multiplication table} \\ \hline 
\endfirsthead

 \multicolumn{7}{l}%
{{\bfseries  continued from previous page}} \\
\hline \multicolumn{1}{|c|}{${\mathcal A}$} & \multicolumn{1}{c|}{ } & \multicolumn{5}{c|}{Multiplication table} \\ \hline 
\endhead

\hline \multicolumn{7}{|r|}{{Continued on next page}} \\ \hline
\endfoot

\hline 
\endlastfoot
$\mathfrak{Z}_{1}^4$ & $\mathfrak{Z}_1$ & $e_1e_1 = e_2$ & $e_1e_2 = e_3$ & $e_1e_3 = e_4$ &&\\&& $e_2e_1 = 2e_3$ & $e_2e_2 = 3e_4$   & $e_3e_1 = 3e_4$&&  \\ \hline
$\mathfrak{Z}_{2}^4$ & $\mathfrak{Z}_3$ & $e_1e_1 = e_3$ & $e_1e_3 = e_4$ & $e_2e_2 = e_4$ & $e_3e_1 = 2e_4$ &  \\ \hline
$\mathfrak{Z}_{3}^4$ & $\mathfrak{Z}_5$ & $e_1e_2 = e_3$ & $e_1e_3 = e_4$ & $e_2e_1 = -e_3$ & $e_2e_2 = e_4$ &  \\ \hline
$\mathfrak{Z}_{4}^4(\alpha)$ & $\mathfrak{N}_2(\alpha)$ & $e_1e_1 = e_3$ & $e_1e_2 = e_4$ &  $e_2e_1 = -\alpha e_3$ & $e_2e_2 = -e_4$ &  \\ \hline
$\mathfrak{Z}_{5}^4(\alpha)$ & $\mathfrak{N}_3(\alpha)$ & $e_1e_1 = e_4$ & $e_1e_2 = \alpha e_4$ & $e_2e_1 = -\alpha e_4$ & $e_2e_2 = e_4$ & $e_3e_3 = e_4$  \\ \hline
\end{longtable}
}

{ 

\subsubsection{5-dimensional Zinbiel algebras}
The algebraic and geometric classification of 
$5$-dimensional Zinbiel algebras is given in a paper by 
 Alvarez,   Fehlberg J\'{u}nior and  Kaygorodov \cite{afk22}.
The variety $\mathfrak{Zin}^5$ has  sixteen irreducible components: 
\[{\rm Irr}(\mathfrak{Zin}^5)=
\left\{\overline{{\mathcal O}\big(\mathfrak{Z}_{i}^5\big)}\right\}_{i=1}^{11} 
\cup
\left\{\overline{\bigcup {\mathcal O}\big(\mathfrak{Z}_{i}^5(\alpha)\big)}\right\}_{i=12}^{14}
\cup
\left\{\overline{\bigcup {\mathcal O}\big(\mathfrak{Z}_{15}^5(\alpha,\beta)\big)}\right\} 
\cup
\left\{\overline{\bigcup {\mathcal O}\big(\mathfrak{Z}_{16}^5(\overline{\mu})\big)}\right\},\]
where

{\small

\begin{longtable}{|l|l|lllll|}
\caption*{ }   \\

\hline \multicolumn{1}{|c|}{${\mathcal A}$} & \multicolumn{1}{c|}{ } & \multicolumn{5}{c|}{Multiplication table} \\ \hline 
\endfirsthead

 \multicolumn{7}{l}%
{{\bfseries  continued from previous page}} \\
\hline \multicolumn{1}{|c|}{${\mathcal A}$} & \multicolumn{1}{c|}{ } & \multicolumn{5}{c|}{Multiplication table} \\ \hline 
\endhead

\hline \multicolumn{7}{|r|}{{Continued on next page}} \\ \hline
\endfoot

\hline 
\endlastfoot
$\mathfrak{Z}_{1}^5$ & $[\mathfrak{N}_1]^2_{08}$  &$e_1e_1=e_4$& $e_1e_2=  e_3$ &$e_1e_3=e_5$ &&\\
&&$ e_2e_1= -e_3$ &$e_2e_2=e_4$ &$e_2e_3=e_4$ &&\\
 \hline

$\mathfrak{Z}_{2}^5$ & 
$[\mathfrak{N}_1^{\mathbb{C}}]^2_{06}$ &  
$ e_1e_1=e_2$ &$e_1 e_2 =e_4$ &$e_1 e_3 =e_4+e_5$ &&\\
&&$e_2 e_1 =2e_4$ &$e_3 e_3 =e_5$ &&&\\ \hline

$\mathfrak{Z}_{3}^5$ &$\mathfrak{Z}_{05}$ & 
$ e_1e_1=e_3$ & $e_1 e_3 =e_5$ & $e_2e_2 =e_4$ && \\
&& $e_2 e_4 =e_5$ & $e_3e_1=2e_5$ & $e_4e_2=2e_5$&& \\ \hline

$\mathfrak{Z}_{4}^5$ & $\mathfrak{Z}_{22}$ & 
$e_1e_1=e_5$ & $e_1e_2=e_3$ & $e_2e_1 =-e_3$ &&\\
&& $e_2e_2=e_5$ & $e_2e_4=e_5$  & $e_4e_3=e_5$ && \\ 
 \hline

$\mathfrak{Z}_{5}^5$ &$\mathfrak{Z}_{23}$ & 
$ e_1e_2=e_3$ & $e_1 e_3 =e_5$ & $e_1e_4 =-e_5$ & $e_2e_1 =e_4$ &\\
&& $e_2 e_2 =-e_3$  & $e_2e_3=-e_5$ & $e_2e_4=e_5$ & $e_3e_2=-2e_5$& 
 \\ \hline

$\mathfrak{Z}_{6}^5$ & $\mathfrak{Z}_{24}$ &
$e_1e_1=e_3$ & $e_1 e_2 =e_4$ & $e_1e_4=-e_5$ && \\
&& $e_2e_1 =-e_3$ &$e_2e_2 =-e_4$ & $e_2 e_4 =e_5$ &&\\
&&$e_3 e_2 =-e_5$ & $e_4e_1=-e_5$ & $e_4e_2=2 e_5$ &&\\ 
 \hline

$\mathfrak{Z}_{7}^5$ &$\mathfrak{Z}_{27}$ &
$ e_1e_2=e_3$ & $e_1 e_3 =-e_5$ & $e_1 e_4 =e_5$&& \\
&& $e_2e_1=e_4$ & $e_2e_3 =-e_5$ & $e_2e_4=e_5$&&

 \\ \hline
$\mathfrak{Z}_{8}^5$ &$\mathfrak{Z}_{34}$ & 
$e_1e_1=e_5$ & $e_1e_2=e_4$ & $e_1 e_3 =\frac{1}{2}e_5$ & $e_2 e_1 =-\frac{1}{2}e_4$ &\\
&&$e_2e_2=e_3$ & $e_2e_3=e_5$ & $e_2e_4 =-e_5$ && \\
&& $e_3e_1=-\frac{1}{2}e_5$ & $e_3e_2=2e_5$ & $e_4e_2=e_5$&&
\\ \hline
 
$\mathfrak{Z}_{9}^5$ &$\mathfrak{Z}_{35}$ & 
$ e_1e_2=e_4$ & $e_1 e_4 =e_5$ & $e_2 e_1 =-e_4$ && \\
&& $e_2e_2=e_3$ &  $e_2e_3=e_5$ & $e_3 e_2 =2e_5$&&

 \\ \hline
$\mathfrak{Z}_{10}^5$ &$\mathfrak{Z}_{38}$ & 
$ e_1e_1=e_2$ & $e_1 e_2 =e_3$ & $e_1e_3=e_5$ & $e_2e_1=2e_3$& \\
&& $e_2e_2=3e_5$  & $e_3e_1=3e_5$ & $e_4e_4=e_5$ &&\\
 \hline

$\mathfrak{Z}_{11}^5$ &$\mathfrak{Z}_{40}$ & 
$ e_1e_1=e_2$ & $e_1 e_2 =\frac{1}{2}e_3$ & $e_1e_3=2e_4$ & $e_1e_4=e_5$ &\\ 
&& $e_2e_1=e_3$ & $e_2e_2=3e_4$ & $e_2e_3=8e_5$ &&\\ 
&&$e_3e_1=6e_4$ & $e_3e_2=12e_5$ & $e_4e_1=4e_5$ &&\\
 \hline

$\mathfrak{Z}_{12}^5(\alpha)$ &$\mathfrak{Z}_{02}^\alpha$ &  
$ e_1e_1=e_2$ & $e_1 e_2 =e_5$ & $e_2e_1=2e_5$ &&\\
&&$e_3e_4=e_5$ & $e_4e_3 =\alpha e_5$ &&&\\
 \hline

$\mathfrak{Z}_{13}^5(\alpha)$ &$\mathfrak{Z}_{14}^\alpha$ & 
$e_1e_1=\alpha e_5$ & $e_1e_2=e_3$ & $e_1e_4=e_5$ &&\\
&& $e_2e_1 =-e_3$ & $e_2 e_3 =e_5$  & $e_4e_4=e_5$ &&\\ 

 \hline

$\mathfrak{Z}_{14}^5(\alpha)$ &$\mathfrak{Z}_{30}^\alpha $ & 
$ e_1e_2=e_4$ & $e_1 e_3 =(\alpha+1)e_5$ & $e_2 e_1 =\alpha e_4$ & $e_2e_2=e_3$ &\\
&& $e_2e_4 =2\alpha e_5$  
 & 
$e_3e_1=2\alpha(\alpha+1)e_5$ & 
$e_4e_2=2(\alpha+1)e_5$&&\\
 \hline

$\mathfrak{Z}_{15}^5(\alpha, \beta)$ &${\mathfrak V}_{4+1}$ &    
$e_1e_2=e_5$& $e_2e_1=\alpha e_5$ &$e_3e_4=e_5$&$e_4e_3=\beta e_5$&\\

 \hline

$\mathfrak{Z}_{16}^5(\overline{\mu})$ &${\mathfrak V}_{3+2}$ & 
$e_1e_1 =  e_4$& $e_1e_2 = \mu_1 e_5$ & $e_1e_3 =\mu_2 e_5$&&\\
&&
$e_2e_1 = \mu_3 e_5$  & $e_2e_2 = \mu_4 e_5$   & $e_2e_3 = \mu_5 e_5$  &&\\
&&$e_3e_1 = \mu_6 e_5$  &$e_3e_2 = \mu_0 e_4+ \mu_7 e_5$  & $e_3e_3 =  e_5$  &&\\
 \hline
\end{longtable}
}
}

\subsection{Novikov algebras} 

An   algebra $\mathfrak{N}$ is called  {\it Novikov}  if it satisfies the identities
\[\begin{array}{rcl}
(xy)z = (xz)y, & (x,y,z)= (y,x,z).
\end{array} \] The variety of Novikov algebras will be denoted by $\mathfrak{Nov}$.

\subsubsection{2-dimensional Novikov algebras}
The algebraic classification of  $2$-dimensional Novikov algebras can be found in Burde (1992), and their graph of degenerations was given in~\cite{bb09}. The variety $\mathfrak{Nov}^2$ has three irreducible components:
\[{\rm Irr}(\mathfrak{Nov}^2)=\left\{\overline{{\mathcal O}\left(\mathfrak{N}_{i}^2\right)}\right\}_{i=1}^{2}\cup\left\{\overline{\bigcup {\mathcal O}\left(\mathfrak{N}_{3}^2(\alpha)\right)}\right\},\]
where

\hspace*{2cm}

{\small

\begin{longtable}{|l|l|lll|}
\caption*{ }  \\

\hline \multicolumn{1}{|c|}{${\mathcal A}$} & \multicolumn{1}{c|}{ } & \multicolumn{3}{c|}{Multiplication table} \\ \hline 
\endfirsthead

 \multicolumn{5}{l}%
{{\bfseries  continued from previous page}} \\
\hline \multicolumn{1}{|c|}{${\mathcal A}$} & \multicolumn{1}{c|}{ } & \multicolumn{3}{c|}{Multiplication table} \\ \hline 
\endhead

\hline \multicolumn{5}{|r|}{{Continued on next page}} \\ \hline
\endfoot

\hline 
\endlastfoot

$\mathfrak{N}_{1}^2 $ & $A_3$ & $e_1e_1 = e_1$ &  $e_2e_2 =e_2$ & \\ \hline
$\mathfrak{N}_{2}^2 $ & $B_5$ & $e_1e_2=  e_1$ & $e_2e_2=e_1+e_2$ & \\ \hline
$\mathfrak{N}_{3}^2(\alpha) $ & $B_2(\alpha)$ & $e_1e_2= \alpha e_1$ & $e_2e_1=(\alpha-1)e_1$ & $e_2e_2=\alpha e_2$  \\ \hline
\end{longtable}
}

\subsubsection{3-dimensional nilpotent Novikov algebras}
The algebraic and geometric classification of $3$-dimensional nilpotent Novikov algebras can be obtained from the classification and description of degenerations of $3$-dimensional nilpotent algebras given in \cite{fkkv}.
Hence, we have  that the variety   $\mathfrak{NNov}^3$ has two irreducible components:
\[{\rm Irr}(\mathfrak{NNov}^3)=
\left\{\overline{{\mathcal O}\left(\mathfrak{NN}_{1}^{3}\right)}\right\}\cup\left\{
\overline{\bigcup {\mathcal O}\left(\mathfrak{NN}_{2}^3(\alpha)\right)}\right\},\]
where

{\small

\begin{longtable}{|l|l|llllll|}
\caption*{ }   \\

\hline \multicolumn{1}{|c|}{${\mathcal A}$} & \multicolumn{1}{c|}{ } & \multicolumn{6}{c|}{Multiplication table} \\ \hline 
\endfirsthead

 \multicolumn{8}{l}%
{{\bfseries  continued from previous page}} \\
\hline \multicolumn{1}{|c|}{${\mathcal A}$} & \multicolumn{1}{c|}{ } & \multicolumn{6}{c|}{Multiplication table} \\ \hline 
\endhead

\hline \multicolumn{8}{|r|}{{Continued on next page}} \\ \hline
\endfoot

\hline 
\endlastfoot
$\mathfrak{NN}_{1}^3$ & $\mathcal{N}_3$ &              $e_1 e_1 = e_2$ & $e_2 e_1=e_3$ &&&&  \\\hline

$\mathfrak{NN}_{2}^3(\alpha)$ &  $\mathcal{N}_4(\alpha)$&     $e_1 e_1 = e_2$ & $e_1 e_2=e_3$ & $e_2 e_1=\alpha e_3$ &&& \\ \hline
\end{longtable}
}

\subsubsection{3-dimensional Novikov algebras}
The geometric classification in dimension $3$ was obtained some years later in~\cite{bb14}, although the algebraic classification was known since the work of Bai and Meng (2001). The variety $\mathfrak{Nov}^3$ has eleven irreducible components, namely:
\[{\rm Irr}(\mathfrak{Nov}^3)=
\left\{\overline{{\mathcal O}\left(\mathfrak{N}_{i}^3\right)}\right\}_{i=1}^{6}\cup
\left\{\overline{\bigcup {\mathcal O}\left(\mathfrak{N}_{i}^3(\alpha)\right)}\right\}_{7}^{10}\cup
\left\{\overline{\bigcup {\mathcal O}\left(\mathfrak{N}_{11}^3(\alpha,\beta)\right)}\right\},\]
where

{\small

\begin{longtable}{|l|l|lll|}
\caption*{ }   \\

\hline \multicolumn{1}{|c|}{${\mathcal A}$} & \multicolumn{1}{c|}{ } & \multicolumn{3}{c|}{Multiplication table} \\ \hline 
\endfirsthead

 \multicolumn{5}{l}%
{{\bfseries  continued from previous page}} \\
\hline \multicolumn{1}{|c|}{${\mathcal A}$} & \multicolumn{1}{c|}{ } & \multicolumn{3}{c|}{Multiplication table} \\ \hline 
\endhead

\hline \multicolumn{5}{|r|}{{Continued on next page}} \\ \hline
\endfoot

\hline 
\endlastfoot

$\mathfrak{N}_{1}^3$ & $ A_4$ & $e_1e_1 = e_1$ &  $e_2e_2 =e_2$ & $e_3e_3=e_3$ \\ \hline
$\mathfrak{N}_{2}^3$ & $B_2$ & $e_1e_1=  e_1$ & $e_1e_2=e_2+e_3$ & $e_1e_3=e_3$ \\&& $e_2e_1=e_2$ & $e_2e_2=e_3$ & $e_3e_1=e_3$ \\ \hline
$\mathfrak{N}_{3}^3$ & $C_1$ & $e_1e_1=-e_1+e_2$ & $e_2e_1=-e_2$ & $e_3e_3= e_3$  \\ \hline
$\mathfrak{N}_{4}^3$ & $D_1$ & $e_1e_1 =- e_1+e_3$ & $e_1e_3=e_2$ &  $e_2e_1 =-e_2$  \\&& $e_3e_1=-e_3$&&  \\ \hline
$\mathfrak{N}_{5}^3$ & $E_3$ & $e_1e_1=  -\frac{1}{2}e_1 + e_3$ & $e_1e_2=\frac{1}{2}e_2$ & $e_1e_3=e_2$ \\&& $e_2e_1=-\frac{1}{2}e_2$ & $e_3e_1=e_2-\frac{1}{2}e_3$ & \\ \hline
$\mathfrak{N}_{6}^3$ & $E_4$ & $e_1e_1= - e_1+e_2$ & $e_1e_3=-\frac{1}{2}e_3$ & $e_2e_1=-e_2$ \\&& $e_3e_1= -e_3$ & $e_3e_3=e_2$ & \\ \hline 
$\mathfrak{N}_{7}^3(\alpha)$ & $C_6(\alpha)$ & $e_1e_1 = \alpha e_1$ &  $e_1e_2 =(\alpha +1)e_2$ & $e_2e_1=\alpha e_2$  \\&& $e_3e_3=e_3$ && \\ \hline
$\mathfrak{N}_{8}^3(\alpha)$ & $D_2(\alpha)$ & $e_1e_1 =\alpha e_1$ &  $e_1e_2 =(\alpha +1)e_2$ & $e_1e_3=e_2 +(\alpha +1)e_3$ \\&& $e_2e_1=\alpha e_2$ & $e_3e_1=\alpha e_3$ & \\ \hline
$\mathfrak{N}_{9}^3(\alpha)$ & $E_{2,\alpha}$ & $e_1e_1= - e_1+e_2$ & $e_1e_3= (\alpha-1) e_3$ & $e_2e_1=-e_2$ \\& & $e_3e_1=-e_3$&&  \\ \hline
$\mathfrak{N}_{10}^3(\alpha)$ & $E_5(\alpha)$ & $e_1e_1 =\alpha e_1$ & $e_1e_2 =(\alpha +1)e_2$ & $e_1e_3=(\alpha+\frac{1}{2})e_3$ \\&& $e_2e_1=\alpha e_2$ & $e_3e_1=\alpha e_3$ & $e_3e_3=e_2$ \\ \hline
$\mathfrak{N}_{11}^3(\alpha,\beta) $ & $E_{1,\beta}(\alpha)$ & $e_1e_1=\beta e_1$ & $e_1e_2= (\beta +1) e_2$ & $e_1e_3=(\alpha + \beta)e_3$ \\&& $e_2e_1=\beta e_2$ & $e_3e_1=\beta e_3$ &  \\ \hline
\end{longtable}
}

The entire degeneration system is detailed in~\cite{bb14}.

\subsubsection{4-dimensional nilpotent Novikov algebras}
The algebraic classification of  $4$-dimensional Novikov algebras has not been obtained yet. However, 
in~\cite{kkk18}, the authors determined all the $4$-dimensional nilpotent Novikov algebras up to isomorphism and also studied the geometric decomposition. They proved that  $\mathfrak{NNov}^4$ has two irreducible components, defined by two families of algebras:
\[{\rm Irr}(\mathfrak{NNov}^4)=
\left\{\overline{\bigcup {\mathcal O}\left(\mathfrak{NN}_{i}^4(\alpha)\right) }\right\}_{i=1}^2,\]
where

{\small
\vspace{-5mm}
\begin{longtable}{|l|l|lllll|}
\caption*{ }   \\

\hline \multicolumn{1}{|c|}{${\mathcal A}$} & \multicolumn{1}{c|}{ } & \multicolumn{5}{c|}{Multiplication table} \\ \hline 
\endfirsthead

 \multicolumn{7}{l}%
{{\bfseries  continued from previous page}} \\
\hline \multicolumn{1}{|c|}{${\mathcal A}$} & \multicolumn{1}{c|}{ } & \multicolumn{5}{c|}{Multiplication table} \\ \hline 
\endhead

\hline \multicolumn{7}{|r|}{{Continued on next page}} \\ \hline
\endfoot

\hline 
\endlastfoot

$\mathfrak{NN}_{1}^4(\alpha)$ & $\mathfrak{N}_{20}(\alpha)$ & $e_1 e_1=\alpha e_4$ & $e_1e_2=e_3$  &  $e_1 e_3=e_4$& $e_2 e_2=e_4$ & $e_2e_3=e_4$  \\&& $e_3 e_2=-e_4$ &&&& \\ \hline
$\mathfrak{NN}_{2}^4(\alpha)$ & $\mathfrak{N}_{22}(\alpha)$ & $e_1 e_1 = e_2$ &  $e_1 e_2=e_3$ & $e_1 e_3=(2-\alpha)e_4$ &  $e_2 e_1= \alpha e_3$  &  $e_2 e_2=\alpha e_4$ \\&&  $e_3 e_1=\alpha e_4$ &&&&\\ \hline
\end{longtable}
}

 \subsection{Bicommutative algebras} 
An   algebra $\mathfrak{B}$ is called {\it bicommutative} if it satisfies the identities 
\[
\begin{array}{rcl}
(xy)z = (xz)y,  &  x(yz) = y(xz).
\end{array} \]
We will denote this variety by $\mathfrak{Bic}$. Note that bicommutative algebras are also known as {\it LR} algebras.

\subsubsection{2-dimensional bicommutative algebras}
The algebraic and geometric classifications of the  $2$-dimensional bicommutative algebras can be found in~\cite{kv16}.
This variety has two irreducible components:
\[{\rm Irr}(\mathfrak{Bic}^2)= \left\{
\overline{{\mathcal O}\left(\mathfrak{B}_{i}^2\right)} \right\}_{i=1}^{2},\]
where

{\small
\vspace{-5mm}
\begin{longtable}{|l|l|ll|}
\caption*{ }   \\

\hline \multicolumn{1}{|c|}{${\mathcal A}$} & \multicolumn{1}{c|}{ } & \multicolumn{2}{c|}{Multiplication table} \\ \hline 
\endfirsthead

 \multicolumn{4}{l}%
{{\bfseries  continued from previous page}} \\
\hline \multicolumn{1}{|c|}{${\mathcal A}$} & \multicolumn{1}{c|}{ } & \multicolumn{2}{c|}{Multiplication table} \\ \hline 
\endhead

\hline \multicolumn{4}{|r|}{{Continued on next page}} \\ \hline
\endfoot

\hline 
\endlastfoot
$\mathfrak{B}^2_{1}$ & $\mathbf{D}_1(0,0)$ & $e_{1}e_{1}=e_{1}$ & $e_{1}e_{2}=e_{1}$ \\ \hline
 $\mathfrak{B}^2_{2}$ & $\mathbf{E}_1(0,0,0,0)$ & $e_{1}e_{1}=e_{1}$ & $e_{2}e_{2}=e_{2}$ \\ \hline
\end{longtable}
}

\subsubsection{3-dimensional nilpotent bicommutative algebras}
For the algebraic classification of the  $3$-dimensional nilpotent bicommutative algebras, consult~\cite{kpv19}. The geometric one can be extracted from the graph of degenerations of all the  nilpotent algebras of dimension $3$ \cite{fkkv}. In $\mathfrak{NBic}^3$, there are two irreducible components:
\[{\rm Irr}(\mathfrak{NBic}^3)= 
\left\{\overline{{\mathcal O}\left(\mathfrak{NB}_{1}^3\right)}\right\}\cup
\left\{\overline{\bigcup {\mathcal O}\left(\mathfrak{NB}_{2}^3(\alpha)\right)} \right\},\]
where

{\small
\begin{longtable}{|l|l|lll|}
\caption*{ }   \\

\hline \multicolumn{1}{|c|}{${\mathcal A}$} & \multicolumn{1}{c|}{ } & \multicolumn{3}{c|}{Multiplication table} \\ \hline 
\endfirsthead

 \multicolumn{5}{l}%
{{\bfseries  continued from previous page}} \\
\hline \multicolumn{1}{|c|}{${\mathcal A}$} & \multicolumn{1}{c|}{ } & \multicolumn{3}{c|}{Multiplication table} \\ \hline 
\endhead

\hline \multicolumn{5}{|r|}{{Continued on next page}} \\ \hline
\endfoot

\hline 
\endlastfoot

$\mathfrak{NB}^3_{1}$ & $\mathfrak{N}_3$  & $e_{1}e_{1}=e_{2}$ & $e_{2}e_{1}=e_{3}$ & \\ \hline
$\mathfrak{NB}^3_{2}(\alpha)$ & $\mathfrak{N}_4(\alpha)$ & $e_{1}e_{1}=e_{2}$ & $e_{1}e_{2}=e_{3}$ & $e_{2}e_{1}=\alpha e_{3}$ \\ \hline
\end{longtable}
}

\subsubsection{4-dimensional nilpotent bicommutative algebras}
The variety $\mathfrak{NBic}^4$ was classified in~\cite{kpv19}, both algebraically and geometrically. It has two irreducible components:

\[{\rm Irr}(\mathfrak{NBic}^4)= \left\{\overline{{\mathcal O}\left(\mathfrak{NB}_{1}^4\right)}\right\}\cup
\left\{\overline{\bigcup {\mathcal O}\left(\mathfrak{NB}_{2}^4(\alpha)\right)} \right\},\]
where

{\small
\begin{longtable}{|l|l|lllll|}
\caption*{ }   \\

\hline \multicolumn{1}{|c|}{${\mathcal A}$} & \multicolumn{1}{c|}{ } & \multicolumn{5}{c|}{Multiplication table} \\ \hline 
\endfirsthead

 \multicolumn{7}{l}%
{{\bfseries  continued from previous page}} \\
\hline \multicolumn{1}{|c|}{${\mathcal A}$} & \multicolumn{1}{c|}{ } & \multicolumn{5}{c|}{Multiplication table} \\ \hline 
\endhead

\hline \multicolumn{7}{|r|}{{Continued on next page}} \\ \hline
\endfoot

\hline 
\endlastfoot

$\mathfrak{NB}^4_{1}$ & $\mathfrak{B}_{10}^4 $ & $e_1 e_2=e_3$  & $e_1 e_3=e_4$ & $e_2e_1=e_4$ & $e_3e_2=e_4$ & \\ \hline
$\mathfrak{NB}^4_{2}(\alpha)$ & $\mathfrak{B}_{24}^4(\alpha)$ & $e_1 e_1 = e_2$  & $e_1 e_2=e_3$  & $e_1e_3=e_4$ & $e_2 e_1 = \alpha e_3$ & $e_2e_2=\alpha e_4$ \\&& $e_3e_1=\alpha e_4$ &&&&\\ \hline
\end{longtable}
}

 \subsection{Assosymmetric algebras} 
An algebra $\mathfrak{A}$ is called  {\it assosymmetric} if it satisfies the  
 identities:
\[
\begin{array}{rcl} 
(x,y,z) = (x,z,y), & \ &  (x,y,z) = (y,x,z).
\end{array} \]
Let  $\mathfrak{Asso}$ denote the variety of assosymmetric algebras.

\subsubsection{2-dimensional   assosymmetric algebras}

All associative algebras are assosymmetric. Also, every assosymmetric algebra of dimension $2$ is associative. The variety of $2$-dimensional assosymmetric  algebras has three irreducible components:

\[{\rm Irr}(\mathfrak{Asso}^2)=\left\{\overline{{\mathcal O}\left(\mathfrak{AS}_{i}^2\right)} \right\}_{i=1}^{3},\]
where  

{\small
\begin{longtable}{|l|ll|}
\caption*{ }   \\

\hline \multicolumn{1}{|c|}{${\mathcal A}$} & \multicolumn{2}{c|}{Multiplication table} \\ \hline 
\endfirsthead

\multicolumn{3}{l}%
{{\bfseries continued from previous page}} \\
\hline \multicolumn{1}{|c|}{${\mathcal A}$} & \multicolumn{2}{c|}{Multiplication table} \\ \hline 
\endhead

\hline \multicolumn{3}{|r|}{{Continued on next page}} \\ \hline
\endfoot

\hline 
\endlastfoot

$\mathfrak{AS}_{1}^2$ & $e_1e_1 = e_1$ & $e_2e_2 = e_2$  \\ \hline
$\mathfrak{AS}_{2}^2$ & $e_1e_1 = e_1$ & $e_1e_2 =e_2$ \\ \hline
$\mathfrak{AS}_{3}^2$ & $e_1e_1 = e_1$ & $e_2e_1 = e_2$ \\ \hline

\end{longtable}
}

\subsubsection{3-dimensional nilpotent assosymmetric algebras}
The list of $3$-dimensional nilpotent assosymmetric algebras can be found in~\cite{ikm19}.
Employing the graph of degenerations of~\cite{fkkv}, we see that in the variety $\mathfrak{NAsso}^3$ there are two irreducible components:
\[{\rm Irr}(\mathfrak{NAsso}^3)= 
\left\{\overline{{\mathcal O}\left(\mathfrak{NAS}_{1}^3\right)}\right\}\cup
\left\{\overline{\bigcup {\mathcal O}\left(\mathfrak{NAS}_{2}^3(\alpha)\right)} \right\},\]
where

{\small
\begin{longtable}{|l|l|lll|}
\caption*{ }   \\

\hline \multicolumn{1}{|c|}{${\mathcal A}$} & \multicolumn{1}{c|}{ } & \multicolumn{3}{c|}{Multiplication table} \\ \hline 
\endfirsthead

 \multicolumn{5}{l}%
{{\bfseries  continued from previous page}} \\
\hline \multicolumn{1}{|c|}{${\mathcal A}$} & \multicolumn{1}{c|}{ } & \multicolumn{3}{c|}{Multiplication table} \\ \hline 
\endhead

\hline \multicolumn{5}{|r|}{{Continued on next page}} \\ \hline
\endfoot

\hline 
\endlastfoot

$ \mathfrak{NAS}^3_{1}$ & $\mathfrak{N}_3$ & $e_{1}e_{1}=e_{2}$ & $e_{2}e_{1}=e_{3}$ &\\ \hline 
$\mathfrak{NAS}^3_{2}(\alpha)$ & $\mathfrak{N}_4(\alpha)$ & $e_{1}e_{1}=e_{2}$ & $e_{1}e_{2}=e_{3}$ & $e_{2}e_{1}=\alpha e_{3}$ \\ \hline
\end{longtable}
}

\subsubsection{4-dimensional nilpotent assosymmetric algebras}

In~\cite{ikm19}, the authors determined all the $4$-dimensional   nilpotent assosymmetric algebras up to isomorphism and found all the degenerations between them. This variety $\mathfrak{NAsso}^4$ has four irreducible components:
\[{\rm Irr}(\mathfrak{NAsso}^4)=
\left\{\overline{{\mathcal O}\big(\mathfrak{NAS}_{1}^4\big)}\right\}_{i=1}^3
\cup 
\left\{\overline{\bigcup {\mathcal O}\left(\mathfrak{NAS}_{4}^4(\alpha)\right)}\right\}\]
where

{\small
\begin{longtable}{|l|l|lll|}
\caption*{ }   \\

\hline \multicolumn{1}{|c|}{${\mathcal A}$} & \multicolumn{1}{c|}{ } & \multicolumn{3}{c|}{Multiplication table} \\ \hline 
\endfirsthead

 \multicolumn{5}{l}%
{{\bfseries  continued from previous page}} \\
\hline \multicolumn{1}{|c|}{${\mathcal A}$} & \multicolumn{1}{c|}{ } & \multicolumn{3}{c|}{Multiplication table} \\ \hline 
\endhead

\hline \multicolumn{5}{|r|}{{Continued on next page}} \\ \hline
\endfoot

\hline 
\endlastfoot

$\mathfrak{NAS}_{1}^4$ & $\mathfrak{A}_{11}^4$ & $e_1 e_1 = e_3$ & $e_1e_3=e_4$ & $e_2e_1=e_3+e_4$ \\&& $e_2e_2=e_3$ & $e_2e_3=e_4$ & $e_3e_2=-e_4$ \\ \hline
$\mathfrak{NAS}_{2}^4$ & $\mathfrak{A}_{12}^4$ & $e_1 e_1 = e_3$ & $e_1e_3=e_4$ &   $e_2e_1=e_3$ \\&&  $e_2e_2=e_3+e_4$ & $e_2e_3=\frac{1+\sqrt{3}i}{2}e_4$ & $e_3e_1=\frac{-1+\sqrt{3}i}{2}e_4$ \\&& $e_3e_2=-e_4$ && \\ \hline
$\mathfrak{NAS}_{3}^4$ & $\mathfrak{A}_{14}^4$ & $e_1 e_1 = e_3$ &  $e_1e_3=e_4$ &   $e_2e_1=e_3$ \\&&  $e_2e_2=e_3+e_4$ & $e_2e_3=\frac{1-\sqrt{3}i}{2}e_4$ & $e_3e_1=\frac{-1-\sqrt{3}i}{2}e_4$ \\&& $e_3e_2=-e_4$ && \\ \hline
$\mathfrak{NAS}_{4}^4(\alpha)$ & $\mathfrak{A}_{18}^4(\alpha)$ & $e_1e_1=e_2$ & $e_1e_2=e_3$ & $e_1e_3=(2-\alpha)e_4$ \\&& $e_2e_1=\alpha e_3$ & $e_2e_2=(\alpha^2-\alpha+1)e_4$ & $e_3e_1=(2\alpha-1)e_4$ \\ \hline
\end{longtable}
}
{

 \subsection{{Antiassociative algebras}} 
An algebra $\mathfrak{A}$ is called  {\it antiassociative} if it satisfies the  
 identity:
\[
\begin{array}{rcl} 
(xy)z =- x(zy),  
\end{array} \]
Let  $\mathfrak{AA}$ denote the variety of assosymmetric algebras.

All antiassociative algebras are nilpotent.

\subsubsection{2-dimensional   antiassociative algebras}
For the variety $\mathfrak{AA}^2$, we  rely on the classification from~\cite{kv16}. The rigid algebra determines it:

{\small

\begin{longtable}{|l|l|ll|}
\caption*{ }   \\

\hline \multicolumn{1}{|c|}{${\mathcal A}$} & \multicolumn{1}{c|}{ } & \multicolumn{2}{c|}{Multiplication table} \\ \hline 
\endfirsthead

 \multicolumn{4}{l}%
{{\bfseries  continued from previous page}} \\
\hline \multicolumn{1}{|c|}{${\mathcal A}$} & \multicolumn{1}{c|}{ } & \multicolumn{2}{c|}{Multiplication table} \\ \hline 
\endhead

\hline \multicolumn{4}{|r|}{{Continued on next page}} \\ \hline
\endfoot

\hline 
\endlastfoot

$\mathfrak{AA}_{1}^2$ & $\mathfrak{A}_3$ & $e_1e_1=e_2$ & \\ \hline

\end{longtable}
}

\subsubsection{3-dimensional   antiassociative algebras}
The list of $3$-dimensional   antiassociative algebras can be found in~\cite{fkk22}.
Employing the graph of degenerations of~\cite{fkkv}, we see that in the variety $\mathfrak{AA}^3$ there are two irreducible components:
\[{\rm Irr}(\mathfrak{AA}^3)= 
\left\{\overline{{\mathcal O}\left(\mathfrak{AA}_{1}^3\right)}\right\}\cup
\left\{\overline{\bigcup {\mathcal O}\left(\mathfrak{AA}_{2}^3(\alpha)\right)} \right\},\]
where

\vspace{+4mm}
{\small

\begin{longtable}{|l|l|lll|}
\caption*{ }   \\

\hline \multicolumn{1}{|c|}{${\mathcal A}$} & \multicolumn{1}{c|}{ } & \multicolumn{3}{c|}{Multiplication table} \\ \hline 
\endfirsthead

 \multicolumn{5}{l}%
{{\bfseries  continued from previous page}} \\
\hline \multicolumn{1}{|c|}{${\mathcal A}$} & \multicolumn{1}{c|}{ } & \multicolumn{3}{c|}{Multiplication table} \\ \hline 
\endhead

\hline \multicolumn{5}{|r|}{{Continued on next page}} \\ \hline
\endfoot

\hline 
\endlastfoot

$ \mathfrak{AA}^3_{1}$ &$\mathfrak{N}_4(-1)$&      $e_1 e_1 = e_2$ & $e_1 e_2=e_3$ & $e_2 e_1=- e_3$  \\ \hline
 
$\mathfrak{AA}^3_{2}(\alpha)$ &$\mathfrak{N}_8(\alpha)$&     $e_1 e_1 = \alpha e_3$ & $e_2 e_1=e_3$  & $e_2 e_2=e_3$ \\ \hline
\end{longtable}
}

\subsubsection{4-dimensional   antiassociative algebras}

In~\cite{fkk22}, the authors determined all the $4$-dimensional   antiassociative algebras up to isomorphism and found   the geometric classification of them. 
This variety $\mathfrak{AA}^4$ has three irreducible components:
\[{\rm Irr}(\mathfrak{AA}^4)=
\left\{\overline{{\mathcal O}\big(\mathfrak{AA}_{1}^4\big)}\right\} 
\cup 
\left\{\overline{\bigcup {\mathcal O}\left(\mathfrak{AA}_{i}^4(\alpha)\right)}\right\}_{i=2}^3\]
where

{\small

\begin{longtable}{|l|l|lllll|}
\caption*{ }   \\

\hline \multicolumn{1}{|c|}{${\mathcal A}$} & \multicolumn{1}{c|}{ } & \multicolumn{5}{c|}{Multiplication table} \\ \hline 
\endfirsthead

 \multicolumn{5}{l}%
{{\bfseries  continued from previous page}} \\
\hline \multicolumn{1}{|c|}{${\mathcal A}$} & \multicolumn{1}{c|}{ } & \multicolumn{5}{c|}{Multiplication table} \\ \hline 
\endhead

\hline \multicolumn{7}{|r|}{{Continued on next page}} \\ \hline
\endfoot

\hline 
\endlastfoot

$\mathfrak{AA}_{1}^4$ & ${\mathbb{A}}_{4,3}$ &  $ e_1e_1=e_2$ & $e_1 e_2 =e_4$ &  $e_2 e_1 =-e_4$ & $e_3 e_3 = e_4$ &\\\hline

$\mathfrak{AA}_{2}^4(\alpha)$ & $\mathcal A_{4,8}^{\alpha}$ &  $e_1 e_1 = e_3$ & $  e_1 e_2 = e_4$ & $  e_2 e_1 =-\alpha e_3$ & $ e_2 e_2 =-e_4$ &  \\ \hline

$\mathfrak{AA}_{3}^4(\alpha)$ &   $\mathcal  A_{4,9}^{\alpha}$ & $e_1 e_1 = e_4$ & $ e_1 e_2 = \alpha e_4$ & $  e_2 e_1 =-\alpha e_4$ & $ e_2 e_2 = e_4$ & $e_3 e_3 = e_4$    \\

\end{longtable}
}

\subsubsection{5-dimensional   antiassociative algebras}

In~\cite{fkk22}, the authors determined all the $5$-dimensional   antiassociative algebras up to isomorphism and found the geometric classification of them. 
This variety $\mathfrak{AA}^5$ has three irreducible components:
\begin{align*}
&{\rm Irr}(\mathfrak{AA}^5) = \\
&\left\{\overline{{\mathcal O}\big(\mathfrak{AA}_{i}^5\big)}\right\}_{i=1}^4 
\cup 
\left\{\overline{\bigcup {\mathcal O}\left(\mathfrak{AA}_{i}^5(\alpha)\right)}\right\}_{i=5}^6 
\cup 
\left\{\overline{\bigcup {\mathcal O}\big(\mathfrak{AA}_{7}^5(\lambda,\mu)\big)}\right\}
\cup 
\left\{\overline{\bigcup {\mathcal O}\big(\mathfrak{AA}_{8}^5(\overline{\mu})\big)}\right\}.
\end{align*}

where

{\small

\begin{longtable}{|l|l|lllll|}
\caption*{ }   \\

\hline \multicolumn{1}{|c|}{${\mathcal A}$} & \multicolumn{1}{c|}{ } & \multicolumn{5}{c|}{Multiplication table} \\ \hline 
\endfirsthead

 \multicolumn{7}{l}%
{{\bfseries  continued from previous page}} \\
\hline \multicolumn{1}{|c|}{${\mathcal A}$} & \multicolumn{1}{c|}{ } & \multicolumn{5}{c|}{Multiplication table} \\ \hline 
\endhead

\hline \multicolumn{7}{|r|}{{Continued on next page}} \\ \hline
\endfoot

\hline 
\endlastfoot

$\mathfrak{AA}_{1}^5$ & $\mathbb A_{5,10}$&$ e_1  e_1=e_2$&$ e_1e_2=e_4$&$  e_1e_3=e_4$ & &\\
&&$  e_2e_1=-e_4$&$  e_3e_1=e_5$ &$  e_3e_3=e_5$ && \\
\hline

$\mathfrak{AA}_{2}^5$ & $\mathbb A_{5,19}$&$  e_1e_1=e_2$&$e_1e_2=e_5$&$ e_1e_3=  e_5$&$e_2e_1=-e_5$& \\
&&$  e_3 e_3=e_4$ &$ e_3e_4=e_5$&$e_4e_3=-e_5$  && \\
\hline

$\mathfrak{AA}_{3}^5$ & $\mathbb A_{5,21}$ & $  e_1 e_2 = e_3+e_5$&$  e_1e_4= e_5 $&$ e_2 e_1 = e_4$&&\\
&&$  e_2 e_2 =-e_3$&$  e_2 e_4=-e_5$ &$  e_3 e_1= -e_5$ &&\\
\hline

$\mathfrak{AA}_{4}^5$ &  $\mathbb A_{5,23}$&$   e_1 e_2 = e_3$&$ e_1e_4=e_5$&$ e_2 e_1 = e_4$&&\\
&&$ e_2e_3=e_5$&$e_3e_1=-e_5$ &$e_4e_2=-e_5$&& \\
 \hline

 $\mathfrak{AA}_{5}^5(\alpha)$ & $\mathbb A_{5,14}^{\alpha}$&   
 $e_1e_1=e_2$ &   $ e_1e_2=e_5$& $e_2e_1=-e_5$ &&\\
 &&$ e_3e_4=e_5$ & $ e_4e_3=\alpha e_5 $&&&\\
\hline
 
 $\mathfrak{AA}_{6}^5(\alpha)$ & $\mathbb A_{5,26}^{\alpha}$&$ e_1 e_2 = e_4$&$ e_1e_3=e_5$&$ e_2 e_1 = \alpha e_4$&$ e_2 e_2 = e_3$&
 \\&&$ e_2e_4=\alpha e_5 $ & $ e_3e_1=-\alpha^2e_5$&$ e_4e_2=-e_5$ &&\\
\hline
 
 $\mathfrak{AA}_{7}^5(\lambda, \mu)$ &${\mathfrak V}_{4+1}$ &    
$e_1e_2=e_5$& $e_2e_1=\lambda e_5$ &$e_3e_4=e_5$&$e_4e_3=\mu e_5$&\\

 \hline

$\mathfrak{AA}_{8}^5(\overline{\mu})$ &${\mathfrak V}_{3+2}$ & 
$e_1e_1 =  e_4$& $e_1e_2 = \mu_1 e_5$ & $e_1e_3 =\mu_2 e_5$&&\\
&&
$e_2e_1 = \mu_3 e_5$  & $e_2e_2 = \mu_4 e_5$   & $e_2e_3 = \mu_5 e_5$  &&\\
&&$e_3e_1 = \mu_6 e_5$  &$e_3e_2 = \mu_0 e_4+ \mu_7 e_5$  & $e_3e_3 =  e_5$  &&\\
 \hline
 
\end{longtable}
}
}

 \subsection{Left-symmetric algebras} 

An algebra $\mathfrak{LS}$ is called {\it left-symmetric} (or {\it pre-Lie}) when
it satisfies the  identity:
\[
(x,y,z) = (y,x,z).  \]
The variety of left-symmetric algebras will be denoted by $\mathfrak{LS}$.

\subsubsection{2-dimensional left-symmetric algebras}
The algebraic classification of $2$-dimensional  left-symmetric algebras can be found in Burde (1992). In~\cite{bb09}, it was established that the variety $\mathfrak{LS}^2$ has six irreducible components
\[{\rm Irr}\left(\mathfrak{LS}^2\right)=\left\{
\overline{{\mathcal O}\big(\mathfrak{LS}_{i}^2\big)}\right\}_{i=1}^{4}\cup
\left\{\overline{\bigcup {\mathcal O}\big(\mathfrak{LS}_{i}^2(\alpha)\big)} \right\}_{i=5}^{6},\]
where

{\small

\begin{longtable}{|l|l|lll|}
\caption*{ }   \\

\hline \multicolumn{1}{|c|}{${\mathcal A}$} & \multicolumn{1}{c|}{ } & \multicolumn{3}{c|}{Multiplication table} \\ \hline 
\endfirsthead

 \multicolumn{5}{l}%
{{\bfseries  continued from previous page}} \\
\hline \multicolumn{1}{|c|}{${\mathcal A}$} & \multicolumn{1}{c|}{ } & \multicolumn{3}{c|}{Multiplication table} \\ \hline 
\endhead

\hline \multicolumn{5}{|r|}{{Continued on next page}} \\ \hline
\endfoot

\hline 
\endlastfoot

$\mathfrak{LS}_{1}^2$ & $A_3$ & $e_1e_1 = e_1$ &  $e_2e_2 =e_2$ &  \\ \hline
$\mathfrak{LS}_{2}^2$ & $B_3$ & $e_2e_1 =- e_1$ &  $e_2e_2 =e_1-e_2$ &  \\ \hline
$\mathfrak{LS}_{3}^2$ & $B_4$ & $e_1e_1 = e_2$ & $e_2e_1 =- e_1$ &  $e_2e_2 =-2e_2$  \\ \hline
$\mathfrak{PL}_{4}^2$ & $B_5$ & $e_1e_2 = e_1$ &  $e_2e_2 =e_1+e_2$ & \\ \hline
$\mathfrak{LS}_{5}^2(\alpha)$ & $B_1(\alpha)$ & $e_2e_1 =- e_1$ &  $e_2e_2 =\alpha e_2$ &  \\ \hline
$\mathfrak{LS}_{6}^2(\alpha)$ & $B_2(\alpha)$ & $e_1e_2=\alpha e_1$ & $e_2e_1 =(\alpha -1) e_1$ &  $e_2e_2 =\alpha e_2$  \\ \hline
\end{longtable}
}

{

\subsubsection{3-dimensional nilpotent left-symmetric algebras}
The list of $3$-dimensional    nilpotent left-symmetric algebras can be found in~\cite{akks21}.
Employing the graph of degenerations of~\cite{fkkv}, we obtain that the variety $\mathfrak{NLS}^3$ is  irreducible and defined by the following family of algebras

{\small
\vspace{-5mm}
\begin{longtable}{|l|l|lll|}
\caption*{ }   \\

\hline \multicolumn{1}{|c|}{${\mathcal A}$} & \multicolumn{1}{c|}{ } & \multicolumn{3}{c|}{Multiplication table} \\ \hline 
\endfirsthead

 \multicolumn{5}{l}%
{{\bfseries  continued from previous page}} \\
\hline \multicolumn{1}{|c|}{${\mathcal A}$} & \multicolumn{1}{c|}{ } & \multicolumn{3}{c|}{Multiplication table} \\ \hline 
\endhead

\hline \multicolumn{5}{|r|}{{Continued on next page}} \\ \hline
\endfoot

\hline 
\endlastfoot

$\mathfrak{NLSC}^3_{1}(\alpha)$ & ${\bf L}^{3*}_{06}(\alpha)$ &  $e_1 e_1 = e_2$ & $e_1 e_2=e_3$ & $e_2 e_1=\alpha e_3$ \\ \hline
\end{longtable}
}

\subsubsection{4-dimensional nilpotent left-symmetric  algebras}
The algebraic and geometric classifications of 
$4$-dimensional left-symmetric algebras are
 given in a paper by 
   Adashev,     Kaygorodov,     Khudoyberdiyev, and   Sattarov \cite{akks22}.
The variety $\mathfrak{NLS}^4$ has  three irreducible components:

\[{\rm Irr}(\mathfrak{NLS}^4)= \left\{\overline{{\mathcal O}\left(\mathfrak{NLS}_{i}^4(\alpha)\right)}\right\}_{i=1}^2\cup
\left\{\overline{\bigcup {\mathcal O}\left(\mathfrak{NLS}_{3}^4(\alpha,\beta)\right)} \right\},\]
where
{\small
\vspace{-5mm}
\begin{longtable}{|l|l|llll|}
\caption*{ }   \\

\hline \multicolumn{1}{|c|}{${\mathcal A}$} & \multicolumn{1}{c|}{ } & \multicolumn{4}{c|}{Multiplication table} \\ \hline 
\endfirsthead

 \multicolumn{6}{l}%
{{\bfseries  continued from previous page}} \\
\hline \multicolumn{1}{|c|}{${\mathcal A}$} & \multicolumn{1}{c|}{ } & \multicolumn{4}{c|}{Multiplication table} \\ \hline 
\endhead

\hline \multicolumn{6}{|r|}{{Continued on next page}} \\ \hline
\endfoot

\hline 
\endlastfoot

$\mathfrak{NLS}^4_{1}(\alpha)$ & 
${\bf L}^4_{12}(\alpha) $ & 
$e_1 e_1 = \alpha e_3$& $e_1 e_2=e_4$&  $e_2 e_1=e_3$&\\
&&$e_2e_2=e_3 $ & $e_2e_3=e_4 $ & $e_3e_1=\alpha e_4 $ &\\  \hline
 
$\mathfrak{NLS}^4_{2}(\alpha)$ & 
${\bf L}^4_{21}(\alpha)$ & 
$e_1 e_1 = e_2$& $e_1 e_2=e_4$&$e_1 e_3=\alpha e_4$&\\
&&$e_2 e_1=e_3$&$e_2 e_3=e_4$& $e_3 e_1=-\alpha e_4$&\\   \hline
    
$\mathfrak{NLS}^4_{3}(\alpha,\beta)$ &${\bf L}^4_{23}(\beta,\alpha)$ & 
$e_1 e_1 = e_2$&$e_1 e_2=e_3$& \multicolumn{2}{l|}{$e_1 e_3=\big((2-\beta)\alpha+1\big) e_4$} \\
&& $e_2 e_1=\beta e_3$& $e_2 e_2=(\beta\alpha+1)e_4$&
\multicolumn{2}{l|}{$e_3 e_1=(\beta\alpha-1)e_4$} \\ \hline

\end{longtable}
}

}

\subsection{{Right alternative algebras}}
Recall that an algebra is said to be {\it right alternative} if it satisfies the identity \[(xy)y=x(yy).\] We will denote this variety by $\mathfrak{RAlt}$.

\subsubsection{2-dimensional   right alternative algebras}

It is easy to see that every $2$-dimensional right alternative algebra is associative.
The variety of $2$-dimensional right alternative   algebras has three irreducible components:

\[{\rm Irr}(\mathfrak{RAlt}^2)=\left\{\overline{{\mathcal O}\left(\mathfrak{RA}_{i}^2\right)} \right\}_{i=1}^{3},\]
where  

{\small
\vspace{-5mm}
\begin{longtable}{|l|ll|}
\caption*{ }   \\

\hline \multicolumn{1}{|c|}{${\mathcal A}$} & \multicolumn{2}{c|}{Multiplication table} \\ \hline 
\endfirsthead

\multicolumn{3}{l}%
{{\bfseries continued from previous page}} \\
\hline \multicolumn{1}{|c|}{${\mathcal A}$} & \multicolumn{2}{c|}{Multiplication table} \\ \hline 
\endhead

\hline \multicolumn{3}{|r|}{{Continued on next page}} \\ \hline
\endfoot

\hline 
\endlastfoot

$\mathfrak{RA}_{1}^2$ & $e_1e_1 = e_1$ & $e_2e_2 = e_2$  \\ \hline
$\mathfrak{RA}_{2}^2$ & $e_1e_1 = e_1$ & $e_1e_2 =e_2$ \\ \hline
$\mathfrak{RA}_{3}^2$ & $e_1e_1 = e_1$ & $e_2e_1 = e_2$ \\ \hline

\end{longtable}
}

{

\subsubsection{3-dimensional nilpotent right alternative algebras}
The list of $3$-dimensional    nilpotent right alternative algebras can be found in~\cite{ikm22}.
Employing the graph of degenerations of~\cite{fkkv}, we obtain that the variety $\mathfrak{NRAlt}^3$ has two irreducible components:

\[{\rm Irr}(\mathfrak{NRAlt}^3)= \left\{\overline{{\mathcal O}\left(\mathfrak{NRA}_{1}^3\right)}\right\}\cup\left\{\overline{\bigcup {\mathcal O}\left(\mathfrak{NRA}_{2}^3(\alpha)\right)} \right\},\]
where

{\small
\vspace{-5mm}
\begin{longtable}{|l|l|lll|}
\caption*{ }   \\

\hline \multicolumn{1}{|c|}{${\mathcal A}$} & \multicolumn{1}{c|}{ } & \multicolumn{3}{c|}{Multiplication table} \\ \hline 
\endfirsthead

 \multicolumn{5}{l}%
{{\bfseries  continued from previous page}} \\
\hline \multicolumn{1}{|c|}{${\mathcal A}$} & \multicolumn{1}{c|}{ } & \multicolumn{3}{c|}{Multiplication table} \\ \hline 
\endhead

\hline \multicolumn{5}{|r|}{{Continued on next page}} \\ \hline
\endfoot

\hline 
\endlastfoot

$\mathfrak{NRA}^3_{1}$ & ${\mathcal N}_{4}(1)$ & $e_1 e_1 = e_2$  & $e_1 e_2=e_3$  & $e_2 e_1 = e_3$  \\ \hline
$\mathfrak{NRA}^{3}_{2}(\alpha)$ & ${\mathcal N}_{8}(\alpha)$ & $e_1 e_1 = \alpha e_3$  & $e_2 e_1=e_3$  & $e_2 e_2=e_3$  \\ \hline
\end{longtable}
}

\subsubsection{4-dimensional nilpotent right alternative algebras}
The algebraic and geometric classifications of 
$4$-dimensional right alternative algebras are given in a paper by 
 Ismailov,     Kaygorodov and Mustafa \cite{ikm22}.
The variety $\mathfrak{NRA}^4$ has  five irreducible components:

\[{\rm Irr}(\mathfrak{NRAlt}^4)= \left\{\overline{{\mathcal O}\left(\mathfrak{NRA}_{i}^4\right)}\right\}_{i=1}^4\cup
\left\{\overline{\bigcup {\mathcal O}\left(\mathfrak{NRA}_{5}^4(\alpha)\right)} \right\},\]
where
{\small
\vspace{-5mm}
\begin{longtable}{|l|l|lllll|}
\caption*{ }   \\

\hline \multicolumn{1}{|c|}{${\mathcal A}$} & \multicolumn{1}{c|}{ } & \multicolumn{5}{c|}{Multiplication table} \\ \hline 
\endfirsthead

 \multicolumn{7}{l}%
{{\bfseries  continued from previous page}} \\
\hline \multicolumn{1}{|c|}{${\mathcal A}$} & \multicolumn{1}{c|}{ } & \multicolumn{5}{c|}{Multiplication table} \\ \hline 
\endhead

\hline \multicolumn{7}{|r|}{{Continued on next page}} \\ \hline
\endfoot

\hline 
\endlastfoot

$\mathfrak{NRA}^4_{1}$ & 
${\mathcal R}^4_{5}$  &  
 $e_1 e_2 = e_3$& $e_1 e_3=e_4$  &$e_2e_1=-e_3$  &&\\ \hline
 
$\mathfrak{NRA}^4_{2}$ & ${\mathcal R}^4_{6}$  & 
$e_1 e_2 = e_3$  & $e_1 e_3=e_4$ & $e_2e_1=-e_3$ & $e_2 e_2=e_4$& \\ \hline
    
$\mathfrak{NRA}^4_{3}$ & ${\mathcal R}^4_{8} $ &  
$e_1 e_1=e_4$ &  $e_2 e_1 = e_3 $& $e_2e_2=e_3$ &  $e_2 e_3=e_4$ & $e_3e_2=e_4$ \\ \hline
   
$\mathfrak{NRA}^4_{4}$ & ${\mathcal R}^4_{9}$  &  
$e_1 e_1 = e_2$ &  $e_1e_2=e_3$&  $e_1 e_3=e_4$ &&\\
&&$e_2 e_1=e_3 $ & $e_2e_2=e_4$ & $e_3e_1=e_4$ && \\ \hline

$\mathfrak{NRA}^4_{5}(\alpha)$ & $\mathfrak{N}_3(\alpha)$  & $e_1e_1 = e_4$ &$e_1e_2 = \alpha e_4$  &$e_2e_1 = -\alpha e_4$ &$e_2e_2 = e_4$  &$e_3e_3 = e_4$
\\ \hline
\end{longtable}
}

}

\subsection{{Right commutative algebras}}
Recall that an algebra is said to be {\it right commutative} if it satisfies the identity \[(xy)z=(xz)y.\] We will denote this variety by $\mathfrak{RC}$.


{

\subsubsection{2-dimensional  right commutative algebras}
The variety of $2$-dimensional right commutative algebras    is   irreducible and defined by the following family of algebras:

{\small
\vspace{-5mm}
\begin{longtable}{|l|l|llll|}
\caption*{ }   \\

\hline \multicolumn{1}{|c|}{${\mathcal A}$} & \multicolumn{1}{c|}{ } & \multicolumn{4}{c|}{Multiplication table} \\ \hline 
\endfirsthead

 \multicolumn{6}{l}%
{{\bfseries  continued from previous page}} \\
\hline \multicolumn{1}{|c|}{${\mathcal A}$} & \multicolumn{1}{c|}{ } & \multicolumn{4}{c|}{Multiplication table} \\ \hline 
\endhead

\hline \multicolumn{6}{|r|}{{Continued on next page}} \\ \hline
\endfoot

\hline 
\endlastfoot

$\mathfrak{RC}_{1}^2(\alpha,\beta)$ & ${\bf E}_1(\alpha,0,0,\beta)$ & $e_1e_1 = e_1$ & $e_2e_1 = \alpha e_1$
&$e_2e_1 = \beta e_2$ &
$e_2e_2 = e_2$ \\ \hline

\end{longtable}
}

\subsubsection{3-dimensional nilpotent right commutative algebras}
The list of $3$-dimensional    nilpotent right commutative algebras can be found in~\cite{akks22}.
Employing the graph of degenerations of~\cite{fkkv}, we obtain that the variety $\mathfrak{NRC}^3$ is  irreducible and defined by the following family of algebras:

{\small
\vspace{-5mm}
\begin{longtable}{|l|l|lll|}
\caption*{ }   \\

\hline \multicolumn{1}{|c|}{${\mathcal A}$} & \multicolumn{1}{c|}{ } & \multicolumn{3}{c|}{Multiplication table} \\ \hline 
\endfirsthead

 \multicolumn{5}{l}%
{{\bfseries  continued from previous page}} \\
\hline \multicolumn{1}{|c|}{${\mathcal A}$} & \multicolumn{1}{c|}{ } & \multicolumn{3}{c|}{Multiplication table} \\ \hline 
\endhead

\hline \multicolumn{5}{|r|}{{Continued on next page}} \\ \hline
\endfoot

\hline 
\endlastfoot

$\mathfrak{NRC}^3_{1}(\alpha)$ & ${\bf R}^{3*}_{06}(\alpha)$ &  $e_1 e_1 = e_2$ & $e_1 e_2=e_3$ & $e_2 e_1=\alpha e_3$ \\ \hline
\end{longtable}
}

\subsubsection{4-dimensional nilpotent right commutative algebras}
The algebraic and geometric classifications of 
$4$-dimensional right commutative  algebras are given in a paper by 
   Adashev,     Kaygorodov,     Khudoyberdiyev, and   Sattarov \cite{akks22}.
The variety $\mathfrak{NRC}^4$ has  five irreducible components:

\[{\rm Irr}(\mathfrak{NRC}^4)= \left\{\overline{{\mathcal O}\left(\mathfrak{NRC}_{i}^4(\alpha)\right)}\right\}_{i=1}^4\cup
\left\{\overline{\bigcup {\mathcal O}\left(\mathfrak{NRC}_{5}^4(\alpha, \beta)\right)} \right\},\]
where
{\small
\vspace{-5mm}
\begin{longtable}{|l|l|llll|}
\caption*{ }   \\

\hline \multicolumn{1}{|c|}{${\mathcal A}$} & \multicolumn{1}{c|}{ } & \multicolumn{4}{c|}{Multiplication table} \\ \hline 
\endfirsthead

 \multicolumn{6}{l}%
{{\bfseries  continued from previous page}} \\
\hline \multicolumn{1}{|c|}{${\mathcal A}$} & \multicolumn{1}{c|}{ } & \multicolumn{4}{c|}{Multiplication table} \\ \hline 
\endhead

\hline \multicolumn{6}{|r|}{{Continued on next page}} \\ \hline
\endfoot

\hline 
\endlastfoot

$\mathfrak{NRC}^4_{1}(\alpha)$ & 
${\bf R}^4_{12}(\alpha)$ & 
$e_1 e_1 = \alpha e_3$& $e_1 e_3=e_4$&  $e_2 e_1=e_3+e_4$&  $e_2e_2=e_3 $  \\  \hline
 
$\mathfrak{NRC}^4_{2}(\alpha)$ & 
${\bf R}^4_{18}(\alpha)$ &  
$e_1e_3=\alpha e_4$& $e_2 e_1=e_3+e_4$ &   $e_2e_2=e_3$  &\\ && $e_3e_1=(1-\alpha) e_4$ & $e_3e_2=(1-\alpha) e_4$ &&\\   \hline

$\mathfrak{NRC}^4_{3}(\alpha)$ & ${\bf R}^4_{29}(\alpha)$  &
$e_1 e_1 = e_2$ & $e_1 e_2=e_3$& $e_1e_3=2e_4$ &\\ &&$e_2e_1=\alpha e_4$ & $e_3e_2=e_4$&&\\ \hline

$\mathfrak{NRC}^4_{4}(\alpha)$ & ${\mathcal N}^4_{20}(\alpha)$  &  $e_1e_2 = e_3$ & $e_1e_1 = \alpha e_4$ & $e_1e_3 = e_4$ &\\
&&$e_2e_2 = e_4$ & $e_2e_3 = e_4$ & $e_3e_2 = -e_4$ &\\\hline

$\mathfrak{NRC}^4_{5}(\alpha,\beta)$ &${\bf R}^4_{27}(\beta, \alpha)$ & 
$e_1 e_1 = e_2$ & $e_1 e_2=e_3$ & \multicolumn{2}{l|}{$e_1e_3=  (2\alpha-\alpha\beta-1) e_4$}\\ 
&& $e_2 e_1=\beta e_3$ & $e_2e_2=(\alpha \beta+1) e_4$ & \multicolumn{2}{l|}{$e_3e_1= (\alpha\beta+1) e_4$}\\
 \hline

\end{longtable}
}

}

 \subsection{Filippov algebras} 
An  algebra $\mathfrak{F}$ with an anticommutative $n$-ary  multiplication is called a {\it Filippov} algebra if it satisfies the identity
\[
\begin{array}{rcl} 
[[x_1,\ldots, x_n],y_2,\ldots, y_n] &=& \sum_{i=1}^{n}[x_1,\ldots, x_{i-1}, [x_i, y_2, \ldots, y_n], x_{i+1}, \ldots, x_n].
\end{array} \] Let us denote this variety by $\mathfrak{Fil}$.

\subsubsection{\texorpdfstring{$(n+1)$-dimensional Filippov ($n$-Lie) algebras}{(n+1)-dimensional Filippov (n-Lie) algebras}}
 
In~\cite{kv19}, the authors gave the geometric classification of the $n$-ary Filippov algebras of dimension $n+1$, namely the variety $\mathfrak{Fil}_n^{n+1}$. For that purpose, they based on the algebraic classification given in a paper by Filippov in 1985. The variety $\mathfrak{Fil}_n^{n+1}$ has two irreducible components:
\[{\rm Irr}(\mathfrak{Fil}_n^{n+1})=\left\{
\overline{{\mathcal O}\big(\mathfrak{F}_{n,1}^{n+1}\big)}\right\}\cup
\left\{\overline{\bigcup {\mathcal O}\big(\mathfrak{F}_{n,2}^{n+1}(\alpha)\big)}\right\}, \]
where

{\small

\begin{longtable}{|l|l|ll|}
\caption*{ }   \\

\hline \multicolumn{1}{|c|}{${\mathcal A}$} & \multicolumn{1}{c|}{ } & \multicolumn{2}{c|}{Multiplication table} \\ \hline 
\endfirsthead

 \multicolumn{4}{l}%
{{\bfseries  continued from previous page}} \\
\hline \multicolumn{1}{|c|}{${\mathcal A}$} & \multicolumn{1}{c|}{ } & \multicolumn{2}{c|}{Multiplication table} \\ \hline 
\endhead

\hline \multicolumn{4}{|r|}{{Continued on next page}} \\ \hline
\endfoot

\hline 
\endlastfoot

$\mathfrak{F}_{n,1}^{n+1}$ & $D_{n+1}$ & $[e_1,\ldots, e_{i-1}, e_{i+1},\ldots, e_{n+1}]=e_i$ & \\ \hline
$\mathfrak{F}_{n,2}^{n+1}(\alpha)$ & $C_2(\alpha)$ & $[e_2,\ldots, e_{n+1}]=\alpha e_1 + e_2$ & 
$[e_1,e_3,\ldots, e_{n+1}]=e_2,$  \\ \hline
\end{longtable}
}\noindent for $1\leq i \leq n+1$.

 \subsection{{Lie triple systems}}
 An  algebra $\mathfrak{A}$ with a $3$-ary  multiplication is called a {\it Lie triple system}  if it satisfies the identities
\[
 [x, y,z] =- [y, x,z], \ \
 [x, y,z] + [y,z, x] + [z, x, y] = 0,\]
\[ [u, v,[x, y,z]] = [[u, v, x], y,z] + [x,[u, v, y] ,z] + [x, y,[u, v,z]]. \] 
Let us denote this variety by $\mathfrak{LTS}$.

\subsubsection{3-dimensional nilpotent Lie triple systems}
The list of $3$-dimensional    nilpotent Lie triple systems can be found in~\cite{abcf232}.
The variety $\mathfrak{NLTS}^3$ is  irreducible:

\[{\rm Irr}(\mathfrak{NLTS}^3)= \left\{\overline{{\mathcal O}\left(\mathfrak{NLTS}_{1}^3\right)}\right\},\]
where
{\small

\begin{longtable}{|l|l|lll|}
\caption*{ }   \\

\hline \multicolumn{1}{|c|}{${\mathcal A}$} & \multicolumn{1}{c|}{ } & \multicolumn{3}{c|}{Multiplication table} \\ \hline 
\endfirsthead

 \multicolumn{5}{l}%
{{\bfseries  continued from previous page}} \\
\hline \multicolumn{1}{|c|}{${\mathcal A}$} & \multicolumn{1}{c|}{ } & \multicolumn{3}{c|}{Multiplication table} \\ \hline 
\endhead

\hline \multicolumn{5}{|r|}{{Continued on next page}} \\ \hline
\endfoot

\hline 
\endlastfoot

$\mathfrak{NLTS}^3_{1}$ & $\mathfrak{T}_{3,2}$ &  $[e_1, e_2, e_1] = e_3$  && \\ \hline
\end{longtable}
}

\subsubsection{4-dimensional nilpotent Lie triple systems}
The list of $4$-dimensional    nilpotent Lie triple systems can be found in~\cite{abcf232}.
The variety $\mathfrak{NLTS}^4$ has two irreducible components:

\[{\rm Irr}(\mathfrak{NLTS}^4)= \left\{\overline{{\mathcal O}\left(\mathfrak{NLTS}_{1}^4\right)}\right\} \cup
\left\{\overline{{\mathcal O}\left(\mathfrak{NLTS}_{2}^4(\alpha)\right)}\right\} ,\]
where
{\small

\begin{longtable}{|l|l|lll|}
\caption*{ }   \\

\hline \multicolumn{1}{|c|}{${\mathcal A}$} & \multicolumn{1}{c|}{ } & \multicolumn{3}{c|}{Multiplication table} \\ \hline 
\endfirsthead

 \multicolumn{5}{l}%
{{\bfseries  continued from previous page}} \\
\hline \multicolumn{1}{|c|}{${\mathcal A}$} & \multicolumn{1}{c|}{ } & \multicolumn{3}{c|}{Multiplication table} \\ \hline 
\endhead

\hline \multicolumn{5}{|r|}{{Continued on next page}} \\ \hline
\endfoot

\hline 
\endlastfoot

$\mathfrak{NLTS}^4_{1}$ & $\mathfrak{T}_{4,7}$ &  
$[e_1, e_2, e_1] = e_3$ & $[e_1, e_2, e_3] = e_4$ &$[e_1, e_3, e_2] = e_4$   \\ \hline

$\mathfrak{NLTS}^4_{2}(\alpha)$ & $\mathfrak{T}_{4,6}(\alpha)$ &  
$[e_1, e_2, e_3] = - (\alpha  + 1) e_4$ &
$[e_2, e_3, e_1] = \alpha e_4$ & $[e_3, e_1, e_2] = e_4$
  \\ \hline
\end{longtable}
}

 \subsection{{Anticommutative ternary algebras}}
 An  algebra $\mathfrak{A}$ with a $3$-ary  multiplication is called a {\it anticommutative}  if it satisfies the identities
\[
 [x, y,z] =- [y, x,z] = [y, z, x]. \] 
Let us denote this variety by $\mathfrak{AT}$.

\subsubsection{3-dimensional anticommutative ternary algebras}
The list of $3$-dimensional   anticommutative ternary algebras can be found in~\cite{kv23}.
The variety $\mathfrak{AT}^3$ is  irreducible:

\[{\rm Irr}(\mathfrak{AT}^3)= \left\{\overline{{\mathcal O}\left(\mathfrak{AT}_{1}^3\right)}\right\},\]
where
{\small

\begin{longtable}{|l|l|llll|}
\caption*{ }   \\

\hline \multicolumn{1}{|c|}{${\mathcal A}$} & \multicolumn{1}{c|}{ } & \multicolumn{4}{c|}{Multiplication table} \\ \hline 
\endfirsthead

 \multicolumn{6}{l}%
{{\bfseries  continued from previous page}} \\
\hline \multicolumn{1}{|c|}{${\mathcal A}$} & \multicolumn{1}{c|}{ } & \multicolumn{4}{c|}{Multiplication table} \\ \hline 
\endhead

\hline \multicolumn{6}{|r|}{{Continued on next page}} \\ \hline
\endfoot

\hline 
\endlastfoot

$\mathfrak{AT}^3_{1}$ &  & 
 
$[e_1,e_2,e_3] =   e_3$& &&
 
\end{longtable}
}

\subsubsection{4-dimensional anticommutative ternary algebras}
The list of $4$-dimensional   anticommutative ternary algebras can be found in~\cite{kv23}.
The variety $\mathfrak{AT}^4$ is  irreducible:

\[{\rm Irr}(\mathfrak{AT}^4)= \left\{\overline{{\mathcal O}\left(\mathfrak{AT}_{1}^4(\alpha,\beta)\right)}\right\},\]
where
{\small

\begin{longtable}{|l|l|llll|}
\caption*{ }   \\

\hline \multicolumn{1}{|c|}{${\mathcal A}$} & \multicolumn{1}{c|}{ } & \multicolumn{4}{c|}{Multiplication table} \\ \hline 
\endfirsthead

 \multicolumn{6}{l}%
{{\bfseries  continued from previous page}} \\
\hline \multicolumn{1}{|c|}{${\mathcal A}$} & \multicolumn{1}{c|}{ } & \multicolumn{4}{c|}{Multiplication table} \\ \hline 
\endhead

\hline \multicolumn{6}{|r|}{{Continued on next page}} \\ \hline
\endfoot

\hline 
\endlastfoot

$\mathfrak{AT}^4_{1}(\alpha,\beta)$ & $\mathfrak{H}_{13}^{\alpha;\beta}$ & 
 
$[e_1,e_2,e_3] = -\alpha e_3$&
$[e_1,e_2,e_4] = e_4$  & 
$[e_1,e_3,e_4] = -\beta e_1$ & 
$[e_2,e_3,e_4] = e_2$   \\ \hline
\end{longtable}
}

 \subsection{Superalgebras}
 A {\it superalgebra} is a ${\mathbb Z}_2$-graded algebra $\mathcal A=\mathcal A_0\oplus \mathcal A_1$. If $\dim\mathcal A_0=m$ and $\dim\mathcal A_1=n$, we say that $\mathcal A$ has dimension $(m,n)$. We define $|a|=i$ for $a\in \mathcal A_i$. We will denote by $\{e_1,\dots,e_m\}$ a fixed basis of $\mathfrak{A}_0$ and by $\{f_1,\dots, f_n\}$ a fixed basis of $\mathfrak{A}_1$.

\subsubsection{Lie  superalgebras}
A superalgebra $\mathcal A$ is a {\it Lie} superalgebra if it satisfies $ab=-(-1)^{|a||b|}ba$ and
\begin{align*}
    (-1)^{|a||c|}(ab)c + (-1)^{|a||b|}(bc)a + (-1)^{|b||c|}(ca)b = 0.
\end{align*}
The variety of Lie superalgebras of dimension $(m,n)$ will be denoted by $\mathfrak{SLie}^{m,n}$. 

The algebraic and geometric classification of $\mathfrak{SLie}^{2,2}$ have been obtained in~\cite{aleis}:
\[
{\rm Irr}\left(\mathfrak{SLie}^{2,2}\right)=\left\{\overline{{\mathcal O}\big(\mathfrak{SL}^{2,2}_i\big)}\right\}_{i=1}^3\cup\left\{\overline{\bigcup {\mathcal O}\big(\mathfrak{SL}^{2,2}_i(\alpha)\big)}\right\}_{i=4}^6\cup\left\{\overline{\bigcup {\mathcal O}\big(\mathfrak{SL}^{2,2}_7(\alpha,\beta)\big)}\right\},
\]
where

{\small

\begin{longtable}{|l|l|llll|}
\caption*{ }   \\

\hline \multicolumn{1}{|c|}{${\mathcal A}$} & \multicolumn{1}{c|}{ } & \multicolumn{4}{c|}{Multiplication table} \\ \hline 
\endfirsthead

 \multicolumn{6}{l}%
{{\bfseries  continued from previous page}} \\
\hline \multicolumn{1}{|c|}{${\mathcal A}$} & \multicolumn{1}{c|}{ } & \multicolumn{4}{c|}{Multiplication table} \\ \hline 
\endhead

\hline \multicolumn{6}{|r|}{{Continued on next page}} \\ \hline
\endfoot

\hline 
\endlastfoot

$\mathfrak{SL}^{2,2}_{1}$ & $\mathfrak{LS}_{19}$ & $e_1e_2 = e_1$ & $e_1f_2 = f_1$ & $e_2f_1 = -f_1$ & $f_1f_2 = e_1$  \\&& $f_2f_2 = 2e_2$ &&& \\ \hline
$\mathfrak{SL}^{2,2}_{2}$ & $\mathfrak{LS}_{1}$ & $f_1f_1 = e_1$ & $f_2f_2 = e_2$ && \\ \hline
$\mathfrak{SL}^{2,2}_{3}$ & $\mathfrak{LS}_{5}$ & $e_1f_1 = f_1$ & $e_2f_2 = f_2$ && \\ \hline
$\mathfrak{SL}^{2,2}_{4}(\alpha)$ & $\mathfrak{LS}^\alpha_{14}$ & $e_1e_2 = e_1$ & $e_2f_1 = \alpha f_1$ & $e_2f_2 = -(\alpha+1)f_2$ & $f_1f_2 = e_1$  \\ \hline
$\mathfrak{SL}^{2,2}_{5}(\alpha)$ & $\mathfrak{LS}^\alpha_{15}$ & $e_1e_2 = e_1$ & $e_2f_1 = \alpha f_1$ & $e_2f_2 = -\frac 12 f_2$ & $f_2f_2 = e_1$   \\ \hline
$\mathfrak{SL}^{2,2}_{6}(\alpha)$ & $\mathfrak{LS}^\alpha_{18}$ & $e_1e_2 = e_1$ & $e_1f_2 = f_1$ & $e_2f_1 = \alpha f_1$ & $e_2f_2 = (\alpha+1)f_2$ \\ \hline
$\mathfrak{SL}^{2,2}_{7}(\alpha,\beta)$ & $\mathfrak{LS}^{\alpha,\beta}_{13}$ & $e_1e_2 = e_1$ & $e_2f_1 = \alpha f_1$ & $e_2f_2 = \beta f_2$ & \\ \hline
\end{longtable}
}

The varieties $\mathfrak{NSLie}^{m,n}$ with $m+n=5$ have been classified algebraically and geometrically in~\cite{aleis2}. We will not consider the cases $(5,0)$ and $(0,5)$ since the first one gives usual Lie algebras, and the second one gives an algebra with zero multiplication. 

The variety $\mathfrak{NSLie}^{4,1}$ is irreducible, determined by the rigid superalgebra

{\small

\begin{longtable}{|l|l|lll|}
\caption*{ }   \\

\hline \multicolumn{1}{|c|}{${\mathcal A}$} & \multicolumn{1}{c|}{ } & \multicolumn{3}{c|}{Multiplication table} \\ \hline 
\endfirsthead

 \multicolumn{5}{l}%
{{\bfseries  continued from previous page}} \\
\hline \multicolumn{1}{|c|}{${\mathcal A}$} & \multicolumn{1}{c|}{ } & \multicolumn{3}{c|}{Multiplication table} \\ \hline 
\endhead

\hline \multicolumn{5}{|r|}{{Continued on next page}} \\ \hline
\endfoot

\hline 
\endlastfoot

$\mathfrak{NSL}^{4,1}_{1}$ & $(4|1)_6$ & $e_1e_2 = e_3$ & $e_1e_3 = e_4$ & $f_1f_1 = e_4$  \\ \hline
\end{longtable}
}

The variety $\mathfrak{NSLie}^{1,4}$ has two irreducible components whose rigid superalgebras are

{\small
\vspace{-5mm}
\begin{longtable}{|l|l|llll|}
\caption*{ }   \\

\hline \multicolumn{1}{|c|}{${\mathcal A}$} & \multicolumn{1}{c|}{ } & \multicolumn{4}{c|}{Multiplication table} \\ \hline 
\endfirsthead

 \multicolumn{6}{l}%
{{\bfseries  continued from previous page}} \\
\hline \multicolumn{1}{|c|}{${\mathcal A}$} & \multicolumn{1}{c|}{ } & \multicolumn{4}{c|}{Multiplication table} \\ \hline 
\endhead

\hline \multicolumn{6}{|r|}{{Continued on next page}} \\ \hline
\endfoot

\hline 
\endlastfoot

$\mathfrak{NSL}^{1,4}_{1}$ & $(1|4)_4$ & $f_1f_1 = e_1$ & $f_2f_2 = e_1$ & $f_3f_3 = e_1$ & $f_4f_4=e_1$\\ \hline
$\mathfrak{NSL}^{1,4}_{2}$ & $(1|4)_7$  & $e_1f_2 = f_1$ & $e_1f_3 = f_2$ & $e_1f_4 = f_3$ & \\ \hline
\end{longtable}
}

The number of irreducible components of $\mathfrak{NSLie}^{3,2}$ is also two. They are determined by the rigid superalgebras

{\small

\begin{longtable}{|l|l|llll|}
\caption*{ }   \\

\hline \multicolumn{1}{|c|}{${\mathcal A}$} & \multicolumn{1}{c|}{ } & \multicolumn{4}{c|}{Multiplication table} \\ \hline 
\endfirsthead

 \multicolumn{6}{l}%
{{\bfseries  continued from previous page}} \\
\hline \multicolumn{1}{|c|}{${\mathcal A}$} & \multicolumn{1}{c|}{ } & \multicolumn{4}{c|}{Multiplication table} \\ \hline 
\endhead

\hline \multicolumn{6}{|r|}{{Continued on next page}} \\ \hline
\endfoot

\hline 
\endlastfoot
$\mathfrak{NSL}^{3,2}_{1}$ & $(3|2)_5$ & $f_1f_1 = e_2$ & $f_1f_2 = e_1$ & $f_2f_2 = e_3$ &\\ \hline
$\mathfrak{NSL}^{3,2}_{2}$ & $(3|2)_{13}$ & $e_1e_2 = e_3$ & $e_1f_2 = f_1$ & $f_1f_2 = e_3$ & $f_2f_2 = 2e_2$  \\ \hline
\end{longtable}
}

Finally,
\[
{\rm Irr}\left(\mathfrak{NSLie}^{2,3}\right)=\left\{\overline{{\mathcal O}\big(\mathfrak{NSL}^{2,3}_i\big)}\right\}_{i=1}^5,
\]
where

{\small
\vspace{-5mm}
\begin{longtable}{|l|l|llll|}
\caption*{ }   \\

\hline \multicolumn{1}{|c|}{${\mathcal A}$} & \multicolumn{1}{c|}{ } & \multicolumn{4}{c|}{Multiplication table} \\ \hline 
\endfirsthead

 \multicolumn{6}{l}%
{{\bfseries  continued from previous page}} \\
\hline \multicolumn{1}{|c|}{${\mathcal A}$} & \multicolumn{1}{c|}{ } & \multicolumn{4}{c|}{Multiplication table} \\ \hline 
\endhead

\hline \multicolumn{6}{|r|}{{Continued on next page}} \\ \hline
\endfoot

\hline 
\endlastfoot

$\mathfrak{NSL}^{2,3}_{1}$ & $(2|3)_6$  & $f_1f_1 = e_1$ & $f_2f_2 = e_2$ & $f_3f_3 = e_1+e_2$ &\\ \hline
$\mathfrak{NSL}^{2,3}_{2}$ & $(2|3)_{18}$ & $e_1f_3 = f_1$ & $e_2f_2 = f_1$ & $f_2f_2 = 2e_1$ & $f_2f_3 = -e_2$\\ \hline
$\mathfrak{NSL}^{2,3}_{3}$ & $(2|3)_{19}$ & $e_1f_3 = f_1$ & $e_2f_2 = f_1$ & $f_2f_3 = -e_1$ & $f_3f_3 = 2e_2$\\ \hline
$\mathfrak{NSL}^{2,3}_{4}$ & $(2|3)_{23}$ & $e_1f_2 = f_1$ & $e_1f_3 = f_2$ & $f_1f_3 = -e_2$ & $f_2f_2 = e_2$\\ \hline
$\mathfrak{NSL}^{2,3}_{5}$ & $(2|3)_{24}$ & $e_1f_2 = f_1$ & $e_1f_3 = f_2$ & $e_2f_3 = f_1$ & \\ \hline
\end{longtable}
}

\subsubsection{Jordan superalgebras}
A superalgebra $\mathcal A$ is a {\it Jordan} superalgebra if it satisfies $ab=(-1)^{|a||b|}ba$ and
\begin{align*}
    &(ab)(cd) + (-1)^{|b||c|}(ac)(bd) + (-1)^{|b||d|+|c||d|}(ad)(bc)\\
    &\quad = ((ab)c)d + (-1)^{|c||d|+|b||c|}((ad)c)b + (-1)^{|a||b|+|a||c|+|a||d|+|c||d|}((bd)c)a.
\end{align*}
Denote by $\mathfrak{SJord}^{m,n}$ the variety of Jordan superalgebras of dimension $(m,n)$. 

Using the algebraic classification from Martin (2017), it was proven in~\cite{maria} that
\[
{\rm Irr}\left(\mathfrak{SJord}^{1,2}\right)=\left\{\overline{{\mathcal O}\big(\mathfrak{SJ}^{1,2}_i\big)}\right\}_{i=1}^7,
\]
where

{\small
\vspace{-9mm}
\begin{longtable}{|l|l|llll|}
\caption*{ }   \\

\hline \multicolumn{1}{|c|}{${\mathcal A}$} & \multicolumn{1}{c|}{ } & \multicolumn{4}{c|}{Multiplication table} \\ \hline 
\endfirsthead

 \multicolumn{6}{l}%
{{\bfseries  continued from previous page}} \\
\hline \multicolumn{1}{|c|}{${\mathcal A}$} & \multicolumn{1}{c|}{ } & \multicolumn{4}{c|}{Multiplication table} \\ \hline 
\endhead

\hline \multicolumn{6}{|r|}{{Continued on next page}} \\ \hline
\endfoot

\hline 
\endlastfoot

$\mathfrak{SJ}^{1,2}_{1}$ & $U_1^s$ & $e_1e_1 = e_1$ &&& \\ \hline
$\mathfrak{SJ}^{1,2}_{2}$ & $S_1^2$ & $e_1e_1 = e_1$ & $e_1f_1 = \frac 12 f_1$ &&\\ \hline
$\mathfrak{SJ}^{1,2}_{3}$ & $S_1^3$ & $e_1f_1 = f_2$ & $f_1f_2 = e_1$  &&\\ \hline
$\mathfrak{SJ}^{1,2}_{4}$ & $S_2^2$ & $e_1f_1 = e_1$ & $e_1f_1 = f_1$  &&\\ \hline
$\mathfrak{SJ}^{1,2}_{5}$ & $S_4^3$ & $e_1e_1 = e_1$ & $e_1f_1 = f_1$ & $e_1f_2 = \frac 12 f_2$ & \\ \hline
$\mathfrak{SJ}^{1,2}_{6}$ & $S_7^3$ & $e_1e_1 = e_1$ & $e_1f_1 = \frac 12 f_1$ & $e_1f_2 = \frac 12 f_2$ & $f_1f_2 = e_1$  \\ \hline
$\mathfrak{SJ}^{1,2}_{7}$ & $S_8^3$ & $e_1e_1 = e_1$ & $e_1f_1 = f_1$ & $e_1f_2 = f_2$ & $f_1f_2 = e_1$ \\ \hline
\end{longtable}
}

Again employing the classification from Martin (2017), it was proven in~\cite{maria} that the variety $\mathfrak{SJord}^{2,1}$ also has seven irreducible components:
\[
{\rm Irr}\left(\mathfrak{SJord}^{2,1}\right)=\left\{\overline{{\mathcal O}\big(\mathfrak{SJ}^{2,1}_i\big)}\right\}_{i=1}^7,
\]
where

{\small

\begin{longtable}{|l|l|llll|}
\caption*{ }   \\

\hline \multicolumn{1}{|c|}{${\mathcal A}$} & \multicolumn{1}{c|}{ } & \multicolumn{4}{c|}{Multiplication table} \\ \hline 
\endfirsthead

 \multicolumn{6}{l}%
{{\bfseries  continued from previous page}} \\
\hline \multicolumn{1}{|c|}{${\mathcal A}$} & \multicolumn{1}{c|}{ } & \multicolumn{4}{c|}{Multiplication table} \\ \hline 
\endhead

\hline \multicolumn{6}{|r|}{{Continued on next page}} \\ \hline
\endfoot

\hline 
\endlastfoot

$\mathfrak{SJ}^{2,1}_{1}$ & $2U_1^s$ & $e_1e_1 = e_1$ & $e_2e_2 = e_2$ &&\\ \hline
$\mathfrak{SJ}^{2,1}_{2}$ & $B_2^s$ & $e_1e_1 = e_1$ & $e_1e_2 = \frac 12 e_2$ &&\\ \hline
$\mathfrak{SJ}^{2,1}_{3}$ & $S_1^2\oplus U_1^s$ & $e_1e_1 = e_1$ & $e_2e_2 = e_2$ & $e_1f_1 = \frac 12 f_1$ & \\ \hline
$\mathfrak{SJ}^{2,1}_{4}$ & $S_2^2\oplus U_1^s$ & $e_1e_1 = e_1$ & $e_2e_2 = e_2$ & $e_1f_1 = f_1$ & \\ \hline
$\mathfrak{SJ}^{2,1}_{5}$ & $S_{11}^3$ & $e_1e_1 = e_1$ & $e_1e_2 = \frac 12 e_2$ & $e_1f_1 = \frac 12 f_1$ & \\ \hline
$\mathfrak{SJ}^{2,1}_{6}$ & $S_{12}^3$& $e_1e_1 = e_1$ & $e_1e_2 = \frac 12 e_2$ & $e_1f_1 = f_1$ & \\ \hline
$\mathfrak{SJ}^{2,1}_{7}$ & $S_{13}^3$ & $e_1e_1 = e_1$ & $e_1f_1 = \frac 12 f_1$ & $e_2e_2 = e_2$  & $e_2f_1 = \frac 12 f_1$  \\ \hline
\end{longtable}
}

As the Jordan superalgebras of dimension $(3,0)$ are nothing but the $3$-dimensional Jordan algebras, and the unique Jordan superalgebra of dimension $(0,3)$ is trivial, the classification of Jordan superalgebras of dimension $(m,n)$ with $m+n=3$ is complete.

\subsection{Poisson algebras}
An   algebra $(\mathfrak{P}, \cdot, \{ \cdot,\cdot \})$  is called {\it Poisson} if it satisfies the identities 
\[
(x,y,z) = 0, \ \ xy = yx,   \ \  \{x,y\}=-\{y,x\},\]
\[\{x,\{y,z\}\}=\{\{x,y\},z\}+\{y,\{x,z\}\},  \
\{xy,z\}=\{x,z\}y+x\{y,z\}. \]
We will denote this variety by $\mathfrak{Poiss}$.

\subsubsection{2-dimensional Poisson algebras}
There are no non-trivial Poisson algebras in dimension $2$ (i.e. there are no Poisson algebras with both non-zero multiplications) \cite{afm22}. Hence, based on subsections \ref{2dasc} and \ref{2dlie}, we have
  
\[
{\rm Irr}\left(\mathfrak{Poiss}^{2}\right)=\left\{\overline{\bigcup {\mathcal O}\big(\mathfrak{P}^2_{i}\big)}\right\}_{i=1}^2,
\]
where

{\small

\begin{longtable}{|l|l|ll|}
\caption*{ }   \\

\hline \multicolumn{1}{|c|}{${\mathcal A}$} & \multicolumn{1}{c|}{ } & \multicolumn{2}{c|}{Multiplication table} \\ \hline 
\endfirsthead

 \multicolumn{4}{l}%
{{\bfseries  continued from previous page}} \\
\hline \multicolumn{1}{|c|}{${\mathcal A}$} & \multicolumn{1}{c|}{ } & \multicolumn{2}{c|}{Multiplication table} \\ \hline 
\endhead

\hline \multicolumn{4}{|r|}{{Continued on next page}} \\ \hline
\endfoot

\hline 
\endlastfoot

$\mathfrak{P}^2_{1}$ & $\mathfrak{CA}_{1}^2$   & $e_1e_1=e_1$ & $e_2e_2=e_2$  \\ \hline

$\mathfrak{P}^2_{2} $ & $\mathfrak{L}_{1}^2$ &  $\{e_1,e_2\}=e_2$ & \\ \hline

\end{longtable}
}

\subsubsection{3-dimensional nilpotent Poisson algebras}

The full graph of degenerations of nilpotent Poisson algebras in dimension $3$ was studied in~\cite{abcf23}. 
The variety $\mathfrak{NPoiss}^{3}$  has two  irreducible components \cite{abcf23}, corresponding to the algebras:

\[
{\rm Irr}\left(\mathfrak{NPoiss}^{3}\right)=
\left\{\overline{ {\mathcal O}\big(\mathfrak{NP}^3_{1}\big)}\right\} 
\cup
\left\{\overline{\bigcup {\mathcal O}\big(\mathfrak{NP}^3_{2}(\alpha)\big)}\right\},
\]
where
{\small

\begin{longtable}{|l|l|llll|}
\caption*{ }   \\

\hline \multicolumn{1}{|c|}{${\mathcal A}$} & \multicolumn{1}{c|}{ } & \multicolumn{4}{c|}{Multiplication table} \\ \hline 
\endfirsthead

$ \mathfrak{NP}^3_{1} $ & $ \mathfrak{P}_{3,6}$ & 
$e_1  e_1 = e_2 $& $ e_1  e_2 = e_3$ &&
 \\  \hline 

$ \mathfrak{NP}^3_{2}(\alpha) $ & $ \mathfrak{P}_{3,3}^{\alpha}$ & 
$e_1  e_2 = \alpha e_3$ & $\{e_1, e_2\} = e_3 $  & & \\  \hline 
 
\end{longtable}
}

\subsubsection{3-dimensional Poisson algebras}
The full graph of degenerations of Poisson algebras in dimension $3$ was studied in~\cite{afm22}. 
The variety $\mathfrak{Poiss}^{3}$  has six irreducible components \cite{afm22}, corresponding to the algebras:

\[
{\rm Irr}\left(\mathfrak{Poiss}^{3}\right)=
\left\{\overline{ {\mathcal O}\big(\mathfrak{P}^3_{i}\big)}\right\}_{i=1}^4
\cup
\left\{\overline{\bigcup {\mathcal O}\big(\mathfrak{P}^3_{i}(\alpha)\big)}\right\}_{i=5}^6,
\]
where
{\small

\begin{longtable}{|l|l|llll|}
\caption*{ }   \\

\hline \multicolumn{1}{|c|}{${\mathcal A}$} & \multicolumn{1}{c|}{ } & \multicolumn{4}{c|}{Multiplication table} \\ \hline 
\endfirsthead

$ \mathfrak{P}^3_{1}$ & $ \mathfrak{P}_{3,5} $ & 
$\left\{ e_{1},e_{2}\right\} =e_{3}$ &$\left\{
e_{1},e_{3}\right\} =-2e_{1}$&$\left\{ e_{2},e_{3}\right\} =2e_{2}$ & \\ \hline 

$ \mathfrak{P}^3_{2}$ & $ \mathfrak{P}_{3,7} $ & 
$e_{1} e_{1}=e_{1}$ & 
$e_{2} e_{2}=e_{2} $& 
$e_{3} e_{3}=e_{3}$  & \\ \hline 

$ \mathfrak{P}^3_{3}$ & $ \mathfrak{P}_{3,18} $ & 
$e_{1} e_{1}=e_{1}$ & $e_{1} e_{2}=e_{2}$ & $e_{1} e_{3}=e_{3}$  & $\left\{ e_{2},e_{3}\right\} =e_{2}$  
 \\ \hline 

$ \mathfrak{P}^3_{4}$ & $ \mathfrak{P}_{3,20} $ &  
$e_{1} e_{1}=e_{1}$ &  
$\left\{ e_{2},e_{3}\right\} =e_{2}$
&&   \\ \hline 

$ \mathfrak{P}^3_{5}(\alpha)$ & $ \mathfrak{P}_{3,4}^{\alpha} $ & $\left\{ e_{1},e_{2}\right\}=e_{2} $& 
$\left\{ e_{1},e_{3}\right\} =\alpha e_{3}$ &&\\  \hline 

$ \mathfrak{P}^3_{6}(\alpha)$ & $ \mathfrak{P}_{3,16}^{\alpha } $ & $e_{1} e_{2}=e_{3}$ & 
$\left\{ e_{1},e_{2}\right\} =\alpha e_{3}$ && \\ \hline 

\end{longtable}
}

\subsubsection{4-dimensional nilpotent Poisson algebras}

The full graph of degenerations of nilpotent Poisson algebras in dimension $3$ was studied in~\cite{abcf23}. 
The variety $\mathfrak{NPoiss}^{4}$  has five irreducible components \cite{abcf23}, corresponding to the algebras:

\[
{\rm Irr}\left(\mathfrak{NPoiss}^{4}\right)=
\left\{\overline{ {\mathcal O}\big(\mathfrak{NP}^4_{i}\big)}\right\}_{i=1}^3
\cup
\left\{\overline{\bigcup {\mathcal O}\big(\mathfrak{NP}^4_{i}(\alpha)\big)}\right\}_{i=4}^5,
\]
where
{\small

\begin{longtable}{|l|l|llll|}
\caption*{ }   \\

\hline \multicolumn{1}{|c|}{${\mathcal A}$} & \multicolumn{1}{c|}{ } & \multicolumn{4}{c|}{Multiplication table} \\ \hline 
\endfirsthead

$ \mathfrak{NP}^4_{1} $ & $ \mathfrak{P}_{4,20}$ & 
$e_1  e_1 = e_2$&$e_1  e_2 = e_3$ &$  e_1  e_3 = e_4$ & $e_2 
e_2 = e_4$ \\  \hline 

$ \mathfrak{NP}^4_{2} $ & $ \mathfrak{P}_{4,12}$ & 
$e_1  e_1 = e_2$  &$e_1  e_2 = e_4 $& $ e_3  e_3 = e_4$ &  $\{e_1, e_3\} = e_4$   \\  \hline 

$ \mathfrak{NP}^4_{3} $ & $ \mathfrak{P}_{4,15}$ & 
$e_1  e_1 = e_4$ & $e_2  e_2 = e_4$ &  $\{e_1, e_2\} = e_3$ & $ \{e_1, e_3\} = e_4$   \\  \hline 

 $ \mathfrak{NP}^4_{4}(\alpha) $ & $ \mathfrak{P}_{4,10}^{\alpha}$ & 
$e_1  e_2 = e_4$ & $e_3  e_3 = e_4$ &  $\{e_1, e_3\} = e_4$ & $ \{e_2, e_3\} = \alpha e_4$   \\  

\hline 
 $ \mathfrak{NP}^4_{5}(\alpha) $ & $ \mathfrak{P}_{4,26}^{\alpha}$ & 
$e_1  e_1 = e_3$ & $e_2  e_2 = \alpha e_3$ & $e_1  e_2 = e_4$  
&  $\{e_1, e_2\} = e_3$   \\  

\hline 
\end{longtable}
}

\subsection{Transposed Poisson algebras}

An   algebra $(\mathfrak{P}, \cdot, \{ , \})$  is called {\it transposed Poisson} if it satisfies the identities 
\[ (x,y,z) = 0, \ \ xy = yx,   \ \  \{x,y \}=-\{y,x\},\]
\[ \{x,\{y,z\}\} =\{\{x,y\},z\}+\{y,\{ x,z \}\},  \
2 x  \{y, z\}  = \{x y, z\} + \{y, x z\}. \]
We will denote this variety by $\mathfrak{TPoiss}$.

\subsubsection{2-dimensional transposed  Poisson algebras}
The full graph of degenerations of transposed  Poisson algebras in dimension $2$ was studied in~\cite{afk24}. 
The variety $\mathfrak{TPoiss}^{2}$  has two  irreducible components \cite{afk24}, corresponding to the algebras:

\[
{\rm Irr}\left(\mathfrak{TPoiss}^{2}\right)=
\left\{\overline{ {\mathcal O}\big(\mathfrak{TP}^2_{1}\big)}\right\}  
\cup
\left\{\overline{\bigcup {\mathcal O}\big(\mathfrak{TP}^2_{2}(\alpha)\big)}\right\},
\]
where
{\small

\begin{longtable}{|l|l|llll|}
\caption*{ }   \\

\hline \multicolumn{1}{|c|}{${\mathcal A}$} & \multicolumn{1}{c|}{ } & \multicolumn{4}{c|}{Multiplication table} \\ \hline 
\endfirsthead

$ \mathfrak{TLP}^2_{1} $ &${\rm T}_{3}$& 
$e_1  e_1  = e_1$ & $e_2  e_2  =e_2$ &&\\ \hline

$ \mathfrak{TLP}^2_{2}(\alpha) $ & ${\rm T}_{5}^{\alpha}$& 
$e_1  e_1   =  e_1$&$ e_1  e_2   = e_2$ & 
$ \{e_{1},e_{2} \}   = \alpha e_{2}$& \\
\hline
 
\hline
\end{longtable}
}

\subsubsection{3-dimensional transposed Poisson algebras}

The variety $\mathfrak{TPoiss}^{3}$ of $3$-dimensional transposed-Poisson algebras has five irreducible components \cite{dfk22}.

\[
{\rm Irr}\left(\mathfrak{TPoiss}^{3}\right)=
\left\{\overline{ {\mathcal O}\big(\mathfrak{TP}^3_{i}\big)}\right\}_{i=1}^2
\cup 
\left\{\overline{\bigcup {\mathcal O}\big(\mathfrak{TP}^3_{i}(\alpha)\big)}\right\}_{i=3}^4\cup
\left\{\overline{\bigcup {\mathcal O}\big(\mathfrak{TP}^3_{1}(\alpha, \beta)\big)}\right\},
\]
where
{\small

\begin{longtable}{|l|l|llll|}
\caption*{ }   \\

\hline \multicolumn{1}{|c|}{${\mathcal A}$} & \multicolumn{1}{c|}{ } & \multicolumn{4}{c|}{Multiplication table} \\ \hline 
\endfirsthead

$ \mathfrak{TP}^3_{1}$ & $ {\rm T}_{01} $ & 
$ \{e_1, e_{2} \}=e_3$ & $ \{e_{1}, e_3 \}=-e_2$ &  $ \{e_2, e_{3} \}=e_1$  &\\ \hline 

$ \mathfrak{TP}^3_{2}$ & $ {\rm T}_{20} $ & $e_1 e_1 = e_1$ & $e_2 e_2 = e_2$ & $e_3 e_3 = e_3$ & \\ \hline

$ \mathfrak{TP}^3_{3}(\alpha)$ & $ {\rm T}_{17}^{\alpha} $ & 

$e_1  e_1=e_2$ & $e_1  e_2=- e_2$ & $e_1  e_3=\alpha  e_1$ &\\
&& $e_2  e_2=   e_2$ &  $e_2  e_3= \alpha  e_2$ &
$e_3  e_3= \alpha  e_3$ & \\ 
&&  $ \{e_1, e_{3} \}=e_1+e_2$ &&& \\ \hline

$ \mathfrak{TP}^3_{4}(\alpha)$ & $ {\rm T}_{12}^{\alpha} $ & 
$ e_1  e_1 =  e_2$&$e_1  e_3 =  \alpha e_1 $&$e_2  e_3 =   \alpha e_2$ & $ e_3  e_3 =   \alpha  e_3$ \\ 
&& $ \{e_1, e_{3} \}=e_1+e_2$ & $ \{e_2, e_3 \}= 2 e_2$ && \\ \hline

$ \mathfrak{TP}^3_{5}(\alpha, \beta)$ & ${\rm T}_{09}^{\alpha, \beta} $ & 
$ e_1  e_3= \beta  e_1$ & $ e_2  e_3 = \beta  e_2$ & $  e_3  e_3 = \beta e_3$ & \\ 
&& $ \{e_1, e_{3} \}=e_1+e_2$ & $\{e_2, e_3\}= \alpha e_2$ &&\\ \hline 
\end{longtable}
}


\subsection{{Generic Poisson   algebras}}
An   algebra $(\mathfrak{P}, \cdot, \{\cdot,\cdot\})$  is called {\it generic Poisson} if it satisfies the identities 
\[
(x,y,z) = 0, \ \ xy = yx,  \  \ \{x,y \}=- \{y,x \},  \]
\[ 
\{x y, z\}  =   \{x, z \} y + x \{y, z \}. \]
We will denote this variety by $\mathfrak{GPoiss}$.

\subsubsection{2-dimensional generic Poisson algebras}
Each anticommutative $2$-dimensional algebra is Lie.
Hence, each $2$-dimensional generic Poisson algebra is Poisson.
There are no non-trivial Poisson algebras in dimension $2$ (i.e. there are no Poisson algebras with both non-zero multiplications) \cite{afm22}. Hence, based on subsections \ref{2dasc} and \ref{2dlie}, we have
  
\[
{\rm Irr}\left(\mathfrak{GPoiss}^{2}\right)=\left\{\overline{\bigcup {\mathcal O}\big(\mathfrak{GP}^2_{i}\big)}\right\}_{i=1}^2,
\]
where

{\small

\begin{longtable}{|l|l|ll|}
\caption*{ }   \\

\hline \multicolumn{1}{|c|}{${\mathcal A}$} & \multicolumn{1}{c|}{ } & \multicolumn{2}{c|}{Multiplication table} \\ \hline 
\endfirsthead

 \multicolumn{4}{l}%
{{\bfseries  continued from previous page}} \\
\hline \multicolumn{1}{|c|}{${\mathcal A}$} & \multicolumn{1}{c|}{ } & \multicolumn{2}{c|}{Multiplication table} \\ \hline 
\endhead

\hline \multicolumn{4}{|r|}{{Continued on next page}} \\ \hline
\endfoot

\hline 
\endlastfoot

$\mathfrak{GP}^2_{1}$ & $\mathfrak{CA}_{1}^2$   & $e_1e_1=e_1$ & $e_2e_2=e_2$  \\ \hline

$\mathfrak{GP}^2_{2} $ & $\mathfrak{L}_{1}^2$ &  $ \{e_1,e_2 \}=e_2$ & \\ \hline

\end{longtable}
}

\subsubsection{3-dimensional generic Poisson algebras}
There is a one-to-one correspondence between generic Poisson algebras and Kokoris algebras \cite{aak24}.
Hence, the algebraic and geometric classifications of generic Poisson algebras  follow from the algebraic and geometric classifications of Kokoris algebras.
The algebraic and geometric classifications of $3$-dimensional   Kokoris can be found in \cite{aak24}.  In particular, it follows that the variety  
$\mathfrak{GPoiss}^3$ has five irreducible components:
\[{\rm Irr}(\mathfrak{GPoiss}^3)=\left\{
\overline{{\mathcal O}\big(\mathfrak{GP}_{i}^3\big)}\right\}_{i=1}^3 \cup \left\{
\overline{\bigcup {\mathcal O}\big(\mathfrak{GP}_{i}^3(\alpha)\big)} \right\}_{i=4}^5 ,\]
where  

{\small

\begin{longtable}{|l|l|llll|}
\caption*{ }   \\

\hline \multicolumn{1}{|c|}{${\mathcal A}$} & \multicolumn{1}{c|}{ } & \multicolumn{4}{c|}{Multiplication table} \\ \hline 
\endfirsthead

 \multicolumn{6}{l}%
{{\bfseries  continued from previous page}} \\
\hline \multicolumn{1}{|c|}{${\mathcal A}$} & \multicolumn{1}{c|}{ } & \multicolumn{4}{c|}{Multiplication table} \\ \hline 
\endhead

\hline \multicolumn{6}{|r|}{{Continued on next page}} \\ \hline
\endfoot

\hline 
\endlastfoot
$\mathfrak{GP}_{1}^3$& ${\bf A}_{04} $ &   $e_{1} e_{1}=e_{1}$ & $ e_{2} e_{2}=e_{2}$ && \\ 
\hline

$\mathfrak{GP}_{2}^3$&
${\bf A}_{29}$ &   $e_{1} e_{1}=e_{1}$ & $e_{1} e_{2}=e_{2}$ &   $e_{1} e_{3}=e_{3}$ 
&   $\{e_{2}, e_{3} \}=e_{3}$   \\
\hline

$\mathfrak{GP}_{3}^3$&${\bf A}_{30}$ &   
$e_{1} e_{1}=e_{1}$ & $\{ e_{2}, e_{3} \}=e_{3}$ &  &\\
\hline

$\mathfrak{GP}_{4}^3(\alpha)$&  ${\bf A}_{02}$ &  
$e_1 e_2= e_3$ & $ \{e_1, e_2\}= \alpha  e_3$ && \\

\hline 
$\mathfrak{GP}_{5}^3(\alpha)$&${\bf A}_{24}^{\alpha }$ &   
$  \{e_{1},e_{2}\}  =e_{3}$ & 
$ \{ e_{1},e_{3} \}  =e_{1}+e_{3}$ &$ \{e_{2},e_{3} \} =\alpha e_{2}$ & \\
\hline 

\end{longtable}

\subsection{{Generic Poisson-Jordan  algebras}}
An   algebra $(\mathfrak{P}, \cdot, \{\cdot,\cdot\})$  is called {\it generic Poisson-Jordan} if it satisfies the identities 
\[
(x^2,y,x) = 0, \ \ xy = yx,  \  \ \{x,y \}=- \{y,x \},  \]
\[ 
\{x y, z\}  =   \{x, z \} y + x \{y, z \}. \]
We will denote this variety by $\mathfrak{GPJ}$. 

\subsubsection{2-dimensional generic Poisson-Jordan algebras}
There is a one-to-one correspondence between generic Poisson-Jordan algebras and noncommutative Jordan algebras \cite{aak24}.
Hence, the algebraic and geometric classifications of generic Poisson-Jordan algebras 
follow from the algebraic and geometric classifications of noncommutative Jordan algebras.
The algebraic and geometric classification of $2$-dimensional noncommutative Jordan algebras can be found in \cite{jkk19}.  In particular, it is proven that the variety  $\mathfrak{GPJ}^2$ has two irreducible components:
\[{\rm Irr}(\mathfrak{GPJ}^2)=\left\{
\overline{{\mathcal O}\big(\mathfrak{GPJ}_{1}^2\big)} \right\} \cup 
\left\{\overline{\bigcup {\mathcal O}\big(\mathfrak{GPJ}_{2}^2(\alpha)\big)}\right\} ,\]
where the algebras $\mathfrak{GPJ}_{1}^2$ and $\mathfrak{GPJ}_{2}^2(\alpha)$ are defined as follows:

{\small

\begin{longtable}{|l|l|lll|}
\caption*{ }   \\

\hline \multicolumn{1}{|c|}{${\mathcal A}$} & \multicolumn{1}{c|}{ } & \multicolumn{3}{c|}{Multiplication table} \\ \hline 
\endfirsthead

 \multicolumn{5}{l}%
{{\bfseries  continued from previous page}} \\
\hline \multicolumn{1}{|c|}{${\mathcal A}$} & \multicolumn{1}{c|}{ } & \multicolumn{3}{c|}{Multiplication table} \\ \hline 
\endhead

\hline \multicolumn{5}{|r|}{{Continued on next page}} \\ \hline
\endfoot

\hline 
\endlastfoot

$\mathfrak{GPJ}_{1}^2$ & ${\bf E}_{1}(0,0,0,0)$ & $e_1 e_1 = e_1$ & $e_2 e_2=e_2$ & \\ \hline
$\mathfrak{GPJ}_{2}^2(\alpha)$ & ${\bf E}_{5}(\alpha)$ & $e_1 e_1 = e_1$ & 
$e_1 e_2=e_1+ e_2$ & \\&&  $e_2 e_2 = e_2$ &   $\{e_1, e_2\}= \alpha e_1 -\alpha e_2$ &  \\ \hline 
\end{longtable}
}

\subsubsection{3-dimensional generic Poisson-Jordan algebras}
There is a one-to-one correspondence between generic Poisson-Jordan algebras and noncommutative Jordan algebras \cite{aak24}.
Hence, the algebraic and geometric classifications of generic Poisson-Jordan algebras 

  follow from the algebraic and geometric classifications of noncommutative Jordan algebras.
The algebraic and geometric classification of $3$-dimensional   noncommutative Jordan algebras can be found in \cite{aak24}.  In particular, it is proven that the variety  $\mathfrak{GPJ}^3$ has eight irreducible components:
\[{\rm Irr}(\mathfrak{GPJ}^3)=\left\{
\overline{{\mathcal O}\big(\mathfrak{GPJ}_{i}^3\big)}\right\}_{i=1}^5 \cup \left\{
\overline{\bigcup {\mathcal O}\big(\mathfrak{GPJ}_{i}^3(\alpha)\big)} \right\}_{i=6}^8 ,\]
where  

{\small

\begin{longtable}{|l|l|llll|}
\caption*{ }   \\

\hline \multicolumn{1}{|c|}{${\mathcal A}$} & \multicolumn{1}{c|}{ } & \multicolumn{4}{c|}{Multiplication table} \\ \hline 
\endfirsthead

 \multicolumn{6}{l}%
{{\bfseries  continued from previous page}} \\
\hline \multicolumn{1}{|c|}{${\mathcal A}$} & \multicolumn{1}{c|}{ } & \multicolumn{4}{c|}{Multiplication table} \\ \hline 
\endhead

\hline \multicolumn{6}{|r|}{{Continued on next page}} \\ \hline
\endfoot

\hline 
\endlastfoot
$\mathfrak{GPJ}_{1}^3$& ${\bf A}_{04} $ &   $e_{1} e_{1}=e_{1}$ & $ e_{2} e_{2}=e_{2}$ && \\ 
\hline

$\mathfrak{GPJ}_{2}^3$& ${\bf A}_{12}$ &  
$e_{1} e_{1}=e_{1}$ & $e_{1} e_{3}= \frac{1}{2}e_{3}$ & $e_{2} e_{2}=e_{2}$ & \\
&&  $e_{2} e_{3}=\frac{1}{2}e_{3}$ &
  $ e_{3} e_{3}=e_{1}+e_{2}$ & &
 \\
\hline

$\mathfrak{GPJ}_{3}^3$& ${\bf A}_{16}$ &   
$e_{1} e_{1}=e_{1}$ & $e_{1} e_{2}= \frac{1}{2}e_{2}$ & $e_{1} e_{3}=e_{3}$   &   $e_{2}  e_{2}=e_{3}$  \\
\hline

$\mathfrak{GPJ}_{4}^3$&${\bf A}_{30}$ &   
$e_{1} e_{1}=e_{1}$ & $ \{e_{2}, e_{3}\}=e_{3}$ &  &\\
\hline 

$\mathfrak{GPJ}_{5}^3$& ${\bf A}_{32}$ &  
$e_{1} e_{1}=e_{1}$ & $e_{1}  e_{2}=\frac{1}{2}e_{2} $ & $e_{1} e_{3}=\frac{1}{2}e_{3}$   &\\
&& $\{e_{1}, e_{2} \}=e_{3}$   & $\{e_{2},  e_{3}\}=e_{2}$ &  &  \\

\hline 

$\mathfrak{GPJ}_{6}^3(\alpha)$& ${\bf A}_{17}^{\alpha}$ &  
$e_{1} e_{1}=e_{1}$  & $e_{1} e_{3}=\frac{1}{2} e_{3}$ & $e_{2} e_{2}=e_{2}$ & \\
&& $e_{2} e_{3}=\frac{1}{2}e_{3}$& 
$\{ e_{1} ,e_{3} \}= \alpha e_{3}$   &  $\{e_{2}, e_{3}\}=-\alpha e_{3}$ & \\
\hline 

$\mathfrak{GPJ}_{7}^3(\alpha)$&${\bf A}_{19}^{\alpha }$ &   
$e_{1} e_{1}=e_{1}$ & $e_{1} e_{3}=\frac{1}{2}e_{3}$ & $e_{2} e_{2}=e_{2}$  & 
$\{ e_{1}, e_{3} \}= \alpha e_{3}$   \\
\hline 
$\mathfrak{GPJ}_{8}^3(\alpha)$&${\bf A}_{24}^{\alpha }$ &   
$\{e_{1}, e_{2} \}  =e_{3}$ & 
$\{e_{1},e_{3} \}=e_{1}+e_{3}$ &
$ \{e_{2},e_{3}\} =\alpha e_{2}$ &\\
\hline 

\end{longtable}
}

\subsection{{Poisson-type  algebras}}

\subsubsection{2-dimensional Leibniz–Poisson algebras}
An   algebra $(\mathfrak{P}, \cdot, \{\cdot,\cdot\})$  is called {\it Leibniz-Poisson} if it satisfies the identities 
\[
(x,y,z) = 0, \ \ xy = yx,    \]
\[ 
\{ \{  x,y \} ,z\} =\{ \{ x,z \},y \}+ \{ x,\{ y,z\}\},  \
\{x y, z\}  =   \{x, z\} y + x \{y, z\}. \]
We will denote this variety by $\mathfrak{LPoiss}$.

The full graph of degenerations of   Leibniz–Poisson algebras in dimension $2$ was studied in~\cite{afk24}. 
The variety $\mathfrak{LPoiss}^{2}$  has four  irreducible components \cite{afk24}, corresponding to the algebras:

\[
{\rm Irr}\left(\mathfrak{LPoiss}^{2}\right)=
\left\{\overline{ {\mathcal O}\big(\mathfrak{LP}^2_{i}\big)}\right\}_{i=1}^2 
\cup
\left\{\overline{\bigcup {\mathcal O}\big(\mathfrak{LP}^2_{i}(\alpha)\big)}\right\}_{i=3}^4,
\]
where

{\small

\begin{longtable}{|l|l|llll|}
\caption*{ }   \\

\hline \multicolumn{1}{|c|}{${\mathcal A}$} & \multicolumn{1}{c|}{ } & \multicolumn{4}{c|}{Multiplication table} \\ \hline 
\endfirsthead

$ \mathfrak{LP}^2_{1} $ & ${\rm L}_{3}$& $\{ e_{1},e_{2} \} = e_{2},$ &$ \{e_{2},e_{1}\} = -e_{2}$ &&
 \\  \hline 

$ \mathfrak{LP}^2_{2} $ & ${\rm L}_{4} $& $e_1  e_1 = e_1$ &$ e_2  e_2 = e_2$ &&\\
\hline 

$ \mathfrak{LP}^2_{3}(\alpha) $ & ${\rm L}_{5}^{\alpha}$&   $e_1  e_1  =  e_1$&$ e_1  e_2   = e_2$  &
$\{ e_{2},e_{1} \}  =  \alpha e_{2}$ & \\
\hline

 $ \mathfrak{LP}^2_{4}(\alpha) $ &${\rm L}_{6}^{\alpha}$& 
 $e_1  e_1  =  e_1$  &
$\{ e_{2},e_{1} \}  =  \alpha e_{2}$ &&
\\
\hline
\end{longtable}
}

\subsubsection{2-dimensional  transposed Leibniz–Poisson   algebras}

An   algebra $(\mathfrak{P}, \cdot, [\cdot,\cdot])$  is called {\it transposed Leibniz-Poisson} if it satisfies the identities 
\[
(x,y,z) = 0, \ \ xy = yx,  \ \  2 x \{y, z\}  =   \{ xy, z \} + \{y, xz\}, \]
\[ 
\{\{  x,y\} ,z\} =\{\{ x,z \},y \} +\{x, \{ y,z\}\},  \ \
\{x, \{y,z \}\}  = \{\{  x,y \},z\} +\{ y, \{ x,z \}\}.     \]
We will denote this variety by $\mathfrak{TLPoiss}$.

The full graph of degenerations of transposed  Leibniz–Poisson algebras in dimension $2$ was studied in~\cite{afk24}. 
The variety $\mathfrak{TLPoiss}^{2}$  has three  irreducible components \cite{afk24}, corresponding to the algebras:

\[
{\rm Irr}\left(\mathfrak{TLPoiss}^{2}\right)=
\left\{\overline{ {\mathcal O}\big(\mathfrak{TLP}^2_{i}\big)}\right\}_{i=1}^2 
\cup
\left\{\overline{\bigcup {\mathcal O}\big(\mathfrak{TLP}^2_{3}(\alpha)\big)}\right\},
\]
where \relax
{\small

\begin{longtable}{|l|l|llll|}
\caption*{ }   \\

\hline \multicolumn{1}{|c|}{${\mathcal A}$} & \multicolumn{1}{c|}{ } & \multicolumn{4}{c|}{Multiplication table} \\ \hline 
\endfirsthead

$ \mathfrak{TLP}^2_{1} $ &${\rm T}_{3}$& 
$e_1  e_1  = e_1$ & $e_2  e_2  =e_2$ &&\\ \hline 

$ \mathfrak{TLP}^2_{2} $ & ${\rm T}_{4}$& 
$e_1  e_1   = e_1$&$ e_1  e_2   = e_2$ &
$ \{ e_{1},e_{1} \}   =   e_{2}$ & \\
\hline 

$ \mathfrak{TLP}^2_{3}(\alpha) $ & ${\rm T}_{5}^{\alpha}$& 
$e_1  e_1   =  e_1$&$ e_1  e_2   = e_2$ & 
$\{e_{1},e_{2} \}   = \alpha e_{2}$&$ \{ e_{2},e_{1} \}  =  -\alpha e_{2}$ \\

\hline
\end{longtable}
}

\subsubsection{2-dimensional Novikov–Poisson algebras}
An   algebra $(\mathfrak{P}, \cdot, \circ)$  is called {\it  Novikov-Poisson} if it satisfies the identities 
\[
(x,y,z) = 0, \ \ xy = yx,  \ \
(x, y, z)_{\circ}=(y,x, z)_{\circ}, \ \
 ( x  \circ  y)\circ  z = ( x \circ z)\circ y, \]
\[(x  y) \circ z   =  x ( y \circ z), \ \
(x \circ y) z - (y \circ x) z  =  x\circ (y  z)- y \circ (x  z).  \]
We will denote this variety by $\mathfrak{NPoiss}$.

The full graph of degenerations of    Novikov–Poisson algebras in dimension $2$ was studied in~\cite{afk24}. 
The variety $\mathfrak{NPoiss}^{2}$  has two  irreducible components \cite{afk24}, corresponding to the algebras:

\[
{\rm Irr}\left(\mathfrak{NPoiss}^{2}\right)=
\left\{\overline{ {\mathcal O}\big(\mathfrak{NP}^2_{i}(\alpha,\beta)\big)}\right\}_{i=1}^2 
,
\]
where \relax
{\small

\begin{longtable}{|l|l|llll|}
\caption*{ }   \\

\hline \multicolumn{1}{|c|}{${\mathcal A}$} & \multicolumn{1}{c|}{ } & \multicolumn{4}{c|}{Multiplication table} \\ \hline 
\endfirsthead

$ \mathfrak{NP}^2_{1}(\alpha,\beta)$ & ${\rm N}_{07}^{\alpha, \beta}$ & 
$e_1  e_1 = e_1$&$ e_2  e_2 =  e_2$ & 
$ e_{1}\circ e_{1}   =  \alpha e_{1}  $&$  e_{2}\circ e_{2}   =  \beta e_{2}$ \\  \hline 

$ \mathfrak{NP}^2_{2}(\alpha,\beta) $ & ${\rm N}_{08}^{\alpha, \beta} $& 
$e_1  e_1   =  e_1$&$ e_1  e_2  = e_2$ &&\\ &&
$ e_{1}\circ e_{1}   =  \alpha e_{1} + e_{2}$&$  e_{1}\circ e_{2}   =  \beta e_{2}$&$e_{2}\circ e_{1} =  \alpha e_{2}$ &\\
\hline 
 
\end{longtable}
}

\subsubsection{2-dimensional pre-Lie Poisson algebras}
An   algebra $(\mathfrak{P}, \cdot, \circ)$  is called {\it  pre-Lie Poisson} if it satisfies the identities 
\[
(x,y,z) = 0, \ \ xy = yx,  \ \
(x, y, z)_{\circ}=(y,x, z)_{\circ},  \]
\[(x  y) \circ z   =  x ( y \circ z), \ \
(x \circ y) z - (y \circ x) z  =  x\circ (y  z)- y \circ (x  z).  \]
We will denote this variety by $\mathfrak{pLPoiss}$.

The variety $\mathfrak{pLPoiss}^{2}$  has two  irreducible components \cite{afk24}, corresponding to the algebras:

\[
{\rm Irr}\left(\mathfrak{pLPoiss}^{2}\right)=
\left\{\overline{ {\mathcal O}\big(\mathfrak{pLP}^2_{1}\big)}\right\} \cup 
\left\{\overline{ {\mathcal O}\big(\mathfrak{pLP}^2_{2}(\alpha)\big)}\right\}
\cup \left\{\overline{ {\mathcal O}\big(\mathfrak{pLP}^2_{i}(\alpha,\beta)\big)}\right\}_{i=3}^4,
\]
where \relax
{\small

\begin{longtable}{|l|l|llll|}
\caption*{ }   \\

\hline \multicolumn{1}{|c|}{${\mathcal A}$} & \multicolumn{1}{c|}{ } & \multicolumn{4}{c|}{Multiplication table} \\ \hline 
\endfirsthead

$\mathfrak{pLP}^2_{1} $ & ${\rm P}_{02}$&  
$e_1 \circ e_1  =  e_2$ &$ e_2 \circ e_1  =  -e_1 $&$e_2 \circ e_2  =  -2e_2$& \\
 \hline

$ \mathfrak{pLP}^2_{2}(\alpha)$ & ${\rm P}_{01}^{\alpha}$&
$e_2 \circ e_1  =  -e_1$&$ e_2 \circ e_2  =  \alpha e_2$ && \\ 
\hline

$ \mathfrak{pLP}^2_{3}(\alpha,\beta)$ & ${\rm N}_{07}^{\alpha, \beta}$ & 
$e_1  e_1 = e_1$&$ e_2  e_2 =  e_2$ & 
$ e_{1}\circ e_{1}   =  \alpha e_{1}  $&$  e_{2}\circ e_{2}   =  \beta e_{2}$ \\  \hline 

$ \mathfrak{pLP}^2_{4}(\alpha,\beta) $ & ${\rm N}_{08}^{\alpha, \beta} $& 
$e_1  e_1   =  e_1$&$ e_1  e_2  = e_2$ &&\\ &&
$ e_{1}\circ e_{1}   =  \alpha e_{1} + e_{2}$&$  e_{1}\circ e_{2}   =  \beta e_{2}$&$e_{2}\circ e_{1} =  \alpha e_{2}$ &\\
\hline 
 
\end{longtable}
}

\subsubsection{2-dimensional  commutative pre-Lie algebras}
An   algebra $(\mathfrak{P}, \cdot, \circ)$  is called {\it  commutative pre-Lie} if it satisfies the identities 
\[
(x,y,z) = 0, \ \ xy = yx,  \ \
(x, y, z)_{\circ}=(y,x, z)_{\circ},  
   \]
\[x\circ   \left( y  z\right) =
\left( x\circ   y\right) z+y  \left( x\circ   z\right).  \]
We will denote this variety by $\mathfrak{CPL}$.

The full graph of degenerations of commutative pre-Lie in dimension $2$ was studied in~\cite{afk24}. 
The variety $\mathfrak{CPL}^{2}$  has eleven  irreducible components \cite{afk24}, corresponding to the algebras:

\[
{\rm Irr}\left(\mathfrak{CPL}^{2}\right)=
\left\{\overline{ {\mathcal O}\big(\mathfrak{CPL}^2_{i} \big)}\right\}_{i=1}^5 \cup
\left\{\overline{ {\mathcal O}\big(\mathfrak{CPL}^2_{i}(\alpha)\big)}\right\}_{i=6}^{11},
\]
where \relax
{\small

\begin{longtable}{|l|l|llll|}
\caption*{ }   \\

\hline \multicolumn{1}{|c|}{${\mathcal A}$} & \multicolumn{1}{c|}{ } & \multicolumn{4}{c|}{Multiplication table} \\ \hline 
\endfirsthead

$ \mathfrak{CPL}^2_{1} $ & ${\rm C}_{07}$&$  e_{1}\circ   e_{1}  = e_{1}$&$  e_{2}\circ   e_{2} = e_{2}$ &&\\ \hline

$ \mathfrak{CPL}^2_{2} $ &  ${\rm C}_{09}$& $ e_1  e_1 = e_1$ &$e_2  e_2 = e_2$ && \\ 
\hline 

 $ \mathfrak{CPL}^2_{3} $ &  ${\rm C}_{11}$&  
$e_1  e_1   = e_1 $&$ e_1  e_2   = e_2$ & 
$ e_{2} \circ   e_{2}  =  e_{2}$%
 &\\ \hline

$ \mathfrak{CPL}^2_{4} $ &  ${\rm C}_{13}$& 
$e_1  e_1   =  e_1$ &
$e_{2}\circ   e_{2}  =  e_{2}$ && \\
\hline

$ \mathfrak{CPL}^2_{5} $ &  ${\rm C}_{14}$& 
$e_1  e_1  =  e_2$ & 
$e_{1}\circ   e_{1}  =  e_{1}$&$  e_{1}\circ   e_{2}  =  2 e_{2}$ & \\ \hline

$ \mathfrak{CPL}^2_{6}(\alpha) $ &  ${\rm C}_{05}^{\alpha}$& 
$e_{1}\circ   e_{1}   = e_{1}$&$  e_{1}\circ   e_{2} =  \alpha e_{2}$ &&\\
 \hline

$ \mathfrak{CPL}^2_{7}(\alpha) $ & ${\rm C}_{06}^{\alpha}$& 
$e_{1} \circ   e_{1}   = e_{1}$&$  e_{1}  \circ   e_{2}   =  \alpha e_{2} $&$ e_{2} \circ   e_{1}  = e_{2}$ & \\
 \hline

$ \mathfrak{CPL}^2_{8}(\alpha) $ &  ${\rm C}_{10}^{\alpha} $ &
$e_1  e_1  =  e_1 $&$ e_1  e_2   = e_2$ &
$ e_{1} \circ   e_{2}  =  \alpha e_{2}$
 & \\
 \hline

$ \mathfrak{CPL}^2_{9}(\alpha) $ &  ${\rm C}_{12}^{\alpha} $& 
$e_1  e_1   =  e_1$ & 
$ e_{1}\circ   e_{2}  =  \alpha e_{2}$
 && \\
 \hline

$ \mathfrak{CPL}^2_{10}(\alpha) $ &  ${\rm C}_{15}^{\alpha}$&
$ e_1  e_1  =  e_2$ & &&\\ && 
$e_{1}\circ   e_{1}  =  \alpha e_{1}$&$  e_{1}\circ   e_{2}  =  2 \alpha e_{2}$ &$e_{2}\circ   e_{1}  =  e_1 + \alpha e_{2}$&$ e_{2}\circ   e_{2}  =  2 e_{2}$%
 \\
\hline

 $ \mathfrak{CPL}^2_{11}(\alpha) $ & ${\rm C}_{17}^{\alpha} $& 
$e_1  e_1   =  e_2$ &  $e_{1}\circ   e_{1} =\alpha e_{2}$ 
 &&\\
 \hline
 
\end{longtable}
}
\subsubsection{2-dimensional anti-Pre-Lie  Poisson  algebras}

An   algebra $(\mathfrak{P}, \cdot, \circ)$  is called {\it anti-pre-Lie Poisson} if it satisfies the identities 
\[
(x,y,z) = 0, \ \ xy = yx,  \ \
x\circ \left( y\circ z\right) -y\circ \left( x\circ z\right) =
(y \circ x - x \circ y)\circ z, \ \ 
J( x, y,z)_{\circ}=J(y,x,z)_{\circ},\]
\[ 
2\left( x\circ y- y\circ x\right) z = y\cdot
\left( x\circ z\right) -x\left( y\circ z\right),  \ \
2   x\circ\left( yz\right) = \left( zx\right) \circ y+z\cdot
\left( x\circ y\right).
\]
We will denote this variety by $\mathfrak{apLPoiss}$.

The variety $\mathfrak{apLPoiss}^{2}$  has three  irreducible components \cite{afk24}, corresponding to the algebras:

\[
{\rm Irr}\left(\mathfrak{apLPoiss}^{2}\right)=
\left\{\overline{ {\mathcal O}\big(\mathfrak{apLP}^2_{1}\big)}\right\} \cup 
\left\{\overline{ {\mathcal O}\big(\mathfrak{apLP}^2_{2}(\alpha,\beta)\big)}\right\}_{i=2}^3,
\]
where \relax
{\small

\begin{longtable}{|l|l|llll|}
\caption*{ }   \\

\hline \multicolumn{1}{|c|}{${\mathcal A}$} & \multicolumn{1}{c|}{ } & \multicolumn{4}{c|}{Multiplication table} \\ \hline 
\endfirsthead

$\mathfrak{apLP}^2_{1} $ &  $ {\rm A}_{05}$& 
$e_{1}\circ e_{1}  =  -e_{2}$&$ e_{2}\circ e_{1}  =  -e_{1}$&& \\
 \hline

$ \mathfrak{apLP}^2_{2}(\alpha,\beta)$ & ${\rm A}_{10}^{\alpha, \beta} $& 
$e_{1}e_{1} = e_{1}$ & $e_{2}e_{2} = e_{2}$ & 
$e_{1}\circ e_{1} = \alpha e_{1}$ & $e_{2}\circ e_{2} = \beta e_{2}$\\
\hline

$ \mathfrak{apLP}^2_{3}(\alpha,\beta)$ & 
 ${\rm A}_{11}^{\alpha, \beta}$& 
$e_{1}e_{1} = e_{1}$ & $e_{1}e_{2} = e_{2}$ & &  \\ &&
$e_{1}\circ e_{1} = \left( 2\alpha -\beta \right) e_{1}+e_{2}$ & $e_{1}\circ
e_{2} = \alpha e_{2}$ & $e_{2}\circ e_{1} = \beta e_{2}$ &
\\ \hline 
 
\end{longtable}
}

\subsubsection{2-dimensional  pre-Poisson algebras}
An   algebra $(\mathfrak{P}, \cdot, \circ)$  is called {\it pre-Poisson} if it satisfies the identities 
\[
x \left(    y    z\right)= 
\left( y   x + x   y\right)  z, \ \
(x\circ y - y \circ x ) z  =  x \circ (y z) -y (x \circ z),\] \[
(xy + y x ) \circ z = x (y \circ z) +y (x \circ z).\]
We will denote this variety by $\mathfrak{pPoiss}$.

The variety $\mathfrak{pPoiss}^{2}$  has five  irreducible components \cite{afk24}, corresponding to the algebras:

\[
{\rm Irr}\left(\mathfrak{pPoiss}^{2}\right)=
\left\{\overline{ {\mathcal O}\big(\mathfrak{pP}^2_{i}\big)}\right\}_{i=1}^3 \cup 
\left\{\overline{ {\mathcal O}\big(\mathfrak{pP}^2_{2}(\alpha)\big)}\right\}_{i=4}^5,
\]
where \relax
{\small

\begin{longtable}{|l|l|llll|}
\caption*{ }   \\

\hline \multicolumn{1}{|c|}{${\mathcal A}$} & \multicolumn{1}{c|}{ } & \multicolumn{4}{c|}{Multiplication table} \\ \hline 
\endfirsthead

$\mathfrak{pP}^2_{1} $ & 
${\rm C}_{07}$ &
$ e_{1}\circ   e_{1}  = e_{1}$&$  e_{2}\circ   e_{2} =e_{2}$&& \\
\hline 

$\mathfrak{pP}^2_{2} $ & ${\rm C}_{08}$&$   e_{1}\circ   e_{1}   = e_{1}$&$  e_{1}\circ   e_{2}  =  2 e_{2}$&$ e_{2}\circ   e_{1}   = \frac{1}{2} e_{1} + e_2 $&$  e_{2}\circ   e_{2} = e_{2}$ \\
\hline

$ \mathfrak{pP}^2_{3}$ & 
 $\rm{P}_{10}$& 
$e_{1}e_{1} = e_{2}$   &
$e_{1}\circ e_{1}= e_{1}$ & $e_{1}\circ e_{2}  = e_{2}$ & $e_{2}\circ e_{1} = e_{2}$ \\
\hline

$\mathfrak{pP}^2_{4}(\alpha) $ & ${\rm C}_{05}^{\alpha}$ &$  e_{1}\circ   e_{1}  = e_{1}$&$  e_{1}\circ   e_{2}  = \alpha e_{2}$&&\\
\hline 

$\mathfrak{pP}^2_{5}(\alpha) $ & ${\rm C}_{06}^{\alpha}$& $e_{1} \circ   e_{1}   = e_{1}$&$  e_{1}  \circ   e_{2}   =  \alpha e_{2}$&$ e_{2} \circ   e_{1} = e_{2}$ & \\
\hline
\end{longtable}
}

\subsection{Compatible   algebras}

\subsubsection{2-dimensional compatible commutative associative  algebras}
An   algebra $(\mathfrak{C}, \cdot, \ast)$  is called {\it compatible commutative associative} if it satisfies the identities 
\[
x   y = y  x,  \  x \ast y= y \ast x,\ 
( x,   y,   z) =0,   \ 
(x,y,   z)_{\ast} = 0, \]
\[(x\ast y) z+(x y)\ast z =
x\ast (y z)+x (y\ast z). 
    \]
We will denote this variety by $\mathfrak{CCA}$.

The variety $\mathfrak{CCA}^{2}$  has two  irreducible components \cite{akm24}, corresponding to the algebras:

\[
{\rm Irr}\left(\mathfrak{CCA}^{2}\right)=
\left\{\overline{ \bigcup {\mathcal O}\big(\mathfrak{CCA}^2_{1}(\alpha,\beta)\big)}\right\} 
\cup
\left\{\overline{\bigcup {\mathcal O}\big(\mathfrak{CCA}^2_{2}(\alpha,\beta,\gamma)\big)}\right\},
\]
where
{\small
\setlength{\tabcolsep}{2pt} 

\begin{longtable}{|l|l|llll|}
\caption*{ }   \\

\hline \multicolumn{1}{|c|}{${\mathcal A}$} & \multicolumn{1}{c|}{ } & \multicolumn{4}{c|}{Multiplication table} \\ \hline 
\endfirsthead

$ \mathfrak{CCA}^2_{1}(\alpha,\beta)$ & $\mathcal{C}_{39}^{\alpha ,\beta }$& 

$e_{1} e_{1}=e_{1}$ & $e_{2} e_{2}=e_{2}$ & & \\ &&
$e_{1}\ast e_{1}=\alpha e_{1}$ & $e_{2}\ast e_{2}=\beta e_{2}$ && \\
\hline

$ \mathfrak{CCA}^2_{2}(\alpha,\beta,\gamma)$ &  
$\mathcal{C}_{38}^{\alpha ,\beta ,\gamma }$ &$
e_{1} e_{1}=e_{1} $ & $e_{2} e_{2}=e_{2}$   &   
$e_{1}\ast e_{1}=\left( \gamma +\beta -\alpha \right) e_{1}-\beta e_{2}$ & \\ && 
$e_{1}\ast e_{2}=\alpha e_{1}+\beta e_{2}$  
& $e_{2}\ast e_{1}=\alpha e_{1}+\beta e_{2}$ & $e_{2}\ast e_{2}=-\alpha e_{1}+\gamma e_{2}$  &\\

 \hline

\end{longtable}
}

\subsubsection{2-dimensional compatible   associative  algebras}
An   algebra $(\mathfrak{C}, \cdot, \ast)$  is called {\it compatible   associative} if it satisfies the identities 
\[
( x,   y,   z) =0,   \ 
(x,y,   z)_{\ast} = 0, \]
\[(x\ast y) z+(x y)\ast z =
x\ast (y z)+x (y\ast z). 
    \]
We will denote this variety by $\mathfrak{CA}$.

The variety $\mathfrak{CA}^{2}$  has four  irreducible components \cite{akm24}, 
corresponding to the algebras:

\[
{\rm Irr}\left(\mathfrak{CA}^{2}\right)=
\left\{\overline{ {\mathcal O}\big(\mathfrak{CA}^2_i\big)}\right\}_{i=1}^2 
\cup
\left\{\overline{ \bigcup {\mathcal O}\big(\mathfrak{CA}^2_{3}(\alpha,\beta)\big)}\right\} 
\cup
\left\{\overline{\bigcup {\mathcal O}\big(\mathfrak{CA}^2_{4}(\alpha,\beta,\gamma)\big)}\right\},
\]
where
{\small
\vspace{-5mm}
\setlength{\tabcolsep}{2pt} 

\begin{longtable}{|l|l|llll|}
\caption*{ }   \\

\hline \multicolumn{1}{|c|}{${\mathcal A}$} & \multicolumn{1}{c|}{ } & \multicolumn{4}{c|}{Multiplication table} \\ \hline 
\endfirsthead

$ \mathfrak{CA}^2_{1}$ & $\mathcal{C}_{28}$ &
$e_{1} e_{1}=e_{1}$ & $e_{1} e_{2}=e_{2}$ &&\\ && 
$e_{2}\ast e_{1}=e_{1}$ & $e_{2}\ast e_{2}=e_{2}$ &&\\ 
 \hline

$ \mathfrak{CA}^2_{2}$ &$\mathcal{C}_{33}$&$e_{1} e_{1}=e_{1}$ & $e_{2} e_{1}=e_{2}$ &&\\&& 
$e_{1}\ast e_{2}=e_{1}$ & $e_{2}\ast e_{2}=e_{2}$ &&\\
\hline

$ \mathfrak{CA}^2_{3}(\alpha,\beta)$ & $\mathcal{C}_{39}^{\alpha ,\beta }$& 

$e_{1} e_{1}=e_{1}$ & $e_{2} e_{2}=e_{2}$ & & \\ &&
$e_{1}\ast e_{1}=\alpha e_{1}$ & $e_{2}\ast e_{2}=\beta e_{2}$ && \\
\hline

$ \mathfrak{CA}^2_{4}(\alpha,\beta,\gamma)$ &  
$\mathcal{C}_{38}^{\alpha ,\beta ,\gamma }$ &$
e_{1} e_{1}=e_{1} $ & $e_{2} e_{2}=e_{2}$   &   
$e_{1}\ast e_{1}=\left( \gamma +\beta -\alpha \right) e_{1}-\beta e_{2}$ & \\ && 
$e_{1}\ast e_{2}=\alpha e_{1}+\beta e_{2}$  
& $e_{2}\ast e_{1}=\alpha e_{1}+\beta e_{2}$ & $e_{2}\ast e_{2}=-\alpha e_{1}+\gamma e_{2}$  &\\

 \hline

\end{longtable}
}

\subsubsection{2-dimensional compatible   Novikov  algebras}
An   algebra $(\mathfrak{C}, \cdot, \ast)$  is called {\it compatible   Novikov} if it satisfies the identities 
\[
( x ,   y,    z) = ( y,   x,    z), \
(x  y)z = (x  z) y, \ 
( x,   y,    z)_{\ast}= ( y,    x,     z)_{\ast},\
(x \ast y) \ast z  = (x \ast z)\ast y, \]
\[(x\ast y) z-x\ast (y z)+(x y)\ast z-x (y\ast z) = 
(y\ast x) z-y\ast (x z)+(y x)\ast z-y  (x\ast z), \]
\[(x \ast y)z +(x  y) \ast z  = (x \ast z)y+(x  z)\ast y.\]

We will denote this variety by $\mathfrak{CN}$.

The variety $\mathfrak{CN}^{2}$  has six  irreducible components \cite{akm24}, 
corresponding to the algebras:

\[
{\rm Irr}\left(\mathfrak{CN}^{2}\right)=
\left\{\overline{ {\mathcal O}\big(\mathfrak{CN}^2_1\big)}\right\} 
\cup
\left\{\overline{ \bigcup {\mathcal O}\big(\mathfrak{CN}^2_{i}(\alpha,\beta)\big)}\right\} _{i=2}^4 
\cup
\left\{\overline{\bigcup {\mathcal O}\big(\mathfrak{CN}^2_{i}(\alpha,\beta,\gamma)\big)}\right\}_{i=5}^6,
\]
where
{\small
\vspace{-5mm}
\setlength{\tabcolsep}{4pt} 

\begin{longtable}{|l|l|llll|}
\caption*{ }   \\

\hline \multicolumn{1}{|c|}{${\mathcal A}$} & \multicolumn{1}{c|}{ } & \multicolumn{4}{c|}{Multiplication table} \\ \hline 
\endfirsthead

$ \mathfrak{CN}^2_{1}$ &$\mathcal{C}_{33}$&$e_{1} e_{1}=e_{1}$ & $e_{2} e_{1}=e_{2}$ &&\\&& 
$e_{1}\ast e_{2}=e_{1}$ & $e_{2}\ast e_{2}=e_{2}$ &&\\
\hline

$ \mathfrak{CN}^2_{2}(\alpha,\beta)$ & 
$\mathcal{C}_{09}^{\alpha ,\beta}$& 
$e_{1} e_{1}=e_{1}+e_{2}$  & $e_{2} e_{1}=e_{2}$ & &  \\ 
&& $e_{1}\ast e_{1}=\alpha e_{1}$ & $e_{1}\ast e_{2}=\beta e_{2}$ &  
$e_{2}\ast e_{1}=\alpha e_{2}$ &\\ 
\hline

$ \mathfrak{CN}^2_{3}(\alpha,\beta)$ & $\mathcal{C}_{22}^{\alpha,\beta }$&
$e_{2} e_{1}=e_{1}$ &  
$e_{1}\ast e_{2}=\alpha e_{1}$ & &\\ 
&& $e_{2}\ast e_{1}=\beta e_{1}$ & $e_{2}\ast e_{2}=e_{1}+\alpha e_{2}$ &&\\
\hline

$ \mathfrak{CN}^2_{4}(\alpha,\beta)$ & $\mathcal{C}_{39}^{\alpha ,\beta }$& 
$e_{1} e_{1}=e_{1}$ & $e_{2} e_{2}=e_{2}$ & & \\ &&
$e_{1}\ast e_{1}=\alpha e_{1}$ & $e_{2}\ast e_{2}=\beta e_{2}$ && \\
\hline

$ \mathfrak{CN}^2_{5}(\alpha,\beta,\gamma)$ &   $\mathcal{C}_{31}^{\alpha ,\beta ,\gamma }$ &
$e_{1} e_{1}=e_{1}$ & $e_{1} e_{2}=\alpha e_{2}$ & $e_{2}
e_{1}=e_{2}$ &\\
&& $e_{1}\ast e_{1}=\beta e_{1}+e_{2}$ &
 $e_{1}\ast e_{2}=\gamma e_{2}$ & $%
e_{2}\ast e_{1}=\beta e_{2}$  & \\ 
 \hline

$ \mathfrak{CN}^2_{6}(\alpha,\beta,\gamma)$ &  
$\mathcal{C}_{38}^{\alpha ,\beta ,\gamma }$ &$
e_{1} e_{1}=e_{1} $ & $e_{2} e_{2}=e_{2}$   &   
$e_{1}\ast e_{1}=\left( \gamma +\beta -\alpha \right) e_{1}-\beta e_{2}$ & \\ && 
$e_{1}\ast e_{2}=\alpha e_{1}+\beta e_{2}$  
& $e_{2}\ast e_{1}=\alpha e_{1}+\beta e_{2}$ & $e_{2}\ast e_{2}=-\alpha e_{1}+\gamma e_{2}$  &\\

 \hline

\end{longtable}
}

\subsubsection{2-dimensional compatible   pre-Lie  algebras}
An   algebra $(\mathfrak{C}, \cdot, \ast)$  is called {\it compatible   pre-Lie} if it satisfies the identities 
\[
( x ,   y,    z) = ( y,   x,    z), \
( x,   y,    z)_{\ast}= ( y,    x,     z)_{\ast},\]
\[(x\ast y) z-x\ast (y z)+(x y)\ast z-x (y\ast z) = 
(y\ast x) z-y\ast (x z)+(y x)\ast z-y  (x\ast z).\]

We will denote this variety by $\mathfrak{CpL}$.

The variety $\mathfrak{CpL}^{2}$  has fourteen  irreducible components \cite{akm24}, 
corresponding to the algebras:

\begin{align*}
&{\rm Irr}\left(\mathfrak{CpL}^{2}\right) =\\ 
&\left\{\overline{ {\mathcal O}\big(\mathfrak{CpL}^2_i\big)}\right\} _{i=1}^2 
\cup
\left\{\overline{ \bigcup {\mathcal O}\big(\mathfrak{CpL}^2_{3}(\alpha)\big)}\right\}  \cup
\left\{\overline{ \bigcup {\mathcal O}\big(\mathfrak{CpL}^2_{i}(\alpha,\beta)\big)}\right\} _{i=4}^{10} 
\cup
\left\{\overline{\bigcup {\mathcal O}\big(\mathfrak{CpL}^2_{i}(\alpha,\beta,\gamma)\big)}\right\}_{i=11}^{14}.
\end{align*}

where
{\small
\setlength{\tabcolsep}{3pt} 

\begin{longtable}{|l|l|llll|}
\caption*{ }   \\

\hline \multicolumn{1}{|c|}{${\mathcal A}$} & \multicolumn{1}{c|}{ } & \multicolumn{4}{c|}{Multiplication table} \\ \hline 
\endfirsthead

$ \mathfrak{CpL}^2_{1}$ &
$\mathcal{C}_{28}$ & 
$e_{1} e_{1}=e_{1}$ & $e_{1} e_{2}=e_{2}$  &&\\ 
&& $e_{2}\ast e_{1}=e_{1}$ & $e_{2}\ast e_{2}=e_{2}$
&& \\
\hline

$ \mathfrak{CpL}^2_{2}$ &$\mathcal{C}_{33}$&$e_{1} e_{1}=e_{1}$ & $e_{2} e_{1}=e_{2}$ &&\\&& 
$e_{1}\ast e_{2}=e_{1}$ & $e_{2}\ast e_{2}=e_{2}$ &&\\
\hline

$ \mathfrak{CpL}^2_{3}(\alpha)$ &
$\mathcal{C}_{37}^{\alpha }$&
$e_{1} e_{1}=e_{1}$ & $e_{1} e_{2}=2e_{2}$ & $e_{2}
e_{1}=e_{2}$ &  \\ 
&&  $e_{1}\ast e_{1}=\alpha e_{1}$ & \multicolumn{2}{l}{$e_{1}\ast e_{2}=e_{1}+2\alpha e_{2}$ }& \\ 
&&$%
e_{2}\ast e_{1}=2e_{1}+\alpha e_{2}$ & $e_{2}\ast e_{2}=e_{2}$ &&\\
\hline

$ \mathfrak{CpL}^2_{4}(\alpha,\beta)$ & 
$\mathcal{C}_{09}^{\alpha ,\beta}$& 
$e_{1} e_{1}=e_{1}+e_{2}$  & $e_{2} e_{1}=e_{2}$ & &  \\ 
&& $e_{1}\ast e_{1}=\alpha e_{1}$ & $e_{1}\ast e_{2}=\beta e_{2}$ &  
$e_{2}\ast e_{1}=\alpha e_{2}$ &\\ 
\hline

$\mathfrak{CpL}^2_{5}(\alpha,\beta)$ & 
$\mathcal{C}_{11}^{\alpha ,\beta  }$ &
$e_{1} e_{1}=e_{1}+e_{2}$ & $e_{1} e_{2}=e_{2}$ && \\ 
&& $e_{1}\ast e_{1}=\alpha e_{1}$ & $e_{1}\ast e_{2}=\beta e_{2}$ && \\
\hline

$ \mathfrak{CpL}^2_{6}(\alpha,\beta)$ & $\mathcal{C}_{22}^{\alpha,\beta }$&
$e_{2} e_{1}=e_{1}$ &  
$e_{1}\ast e_{2}=\alpha e_{1}$ & &\\ 
&& $e_{2}\ast e_{1}=\beta e_{1}$ & \multicolumn{2}{l}{$e_{2}\ast e_{2}=e_{1}+\alpha e_{2}$ }&\\
\hline

$ \mathfrak{CpL}^2_{7}(\alpha,\beta)$ &
$\mathcal{C}_{27}^{\alpha ,\beta }$ & 
$e_{1} e_{1}=e_{1}$ & $e_{1} e_{2}=\frac{1}{2}e_{2}$ &  $e_{1}\ast e_{1}=2\alpha e_{1}$ & \\
&& $e_{1}\ast e_{2}=e_{1}+\alpha e_{2}$ & $
e_{2}\ast e_{1}=2e_{1}$ & $e_{2}\ast e_{2}=\beta e_{1}+e_{2}$ &\\ 
\hline

$ \mathfrak{CpL}^2_{8}(\alpha,\beta)$ &
$\mathcal{C}_{36}^{\alpha ,\beta }$ &
$e_{1} e_{1}=e_{1}$ & $e_{1} e_{2}=2e_{2}$ & $e_{2}
e_{1}=e_{2}$ &  \\ 
&& $e_{1}\ast e_{1}=\alpha e_{1}+\beta e_{2}$ & $e_{1}\ast e_{2}=2\alpha e_{2}$ &&\\
&&    $e_{2}\ast e_{1}=e_{1}+\alpha e_{2}$ & $e_{2}\ast e_{2}=2e_{2}$ &&  \\
\hline

$ \mathfrak{CpL}^2_{9}(\alpha,\beta)$ & $\mathcal{C}_{39}^{\alpha ,\beta }$& 
$e_{1} e_{1}=e_{1}$ & $e_{2} e_{2}=e_{2}$ & & \\ &&
$e_{1}\ast e_{1}=\alpha e_{1}$ & $e_{2}\ast e_{2}=\beta e_{2}$ && \\
\hline

$ \mathfrak{CpL}^2_{10}(\alpha,\beta)$ &
$\mathcal{C}_{41}^{\alpha ,\beta }$& 
$e_{1} e_{1}=e_{1}$ & $e_{1} e_{2}=2e_{2}$ &$e_{2} e_{1}=\frac{1}{2}e_{1}+e_{2}$ & \\ 
&& $e_{2} e_{2}=e_{2}$ &    $e_{1}\ast e_{1}=\alpha e_{1}$ & 
$e_{1}\ast e_{2}=\beta e_{1}+2\alpha e_{2}$ &\\ 
&& $e_{2}\ast e_{1}=\left( 2\beta +\frac{1}{2}\alpha \right) e_{1}+\alpha
e_{2}$ & 
\multicolumn{2}{l}{$e_{2}\ast e_{2}=\frac{1}{2}\beta e_{1}+\left( \alpha +\beta
\right) e_{2}$} & \\
\hline

$ \mathfrak{CpL}^2_{11}(\alpha,\beta,\gamma)$ & 
$\mathcal{C}_{24}^{\alpha ,\beta ,\gamma }$ &
$e_{1} e_{1}=e_{1}$ & $e_{1} e_{2}=\alpha e_{2}$ && \\ 
&& $e_{1}\ast e_{1}=\beta e_{1}+e_{2}$ & $e_{1}\ast e_{2}=\gamma e_{2}$ && \\ 
\hline

$ \mathfrak{CpL}^2_{12}(\alpha,\beta,\gamma)$ &   $\mathcal{C}_{31}^{\alpha ,\beta ,\gamma }$ &
$e_{1} e_{1}=e_{1}$ & $e_{1} e_{2}=\alpha e_{2}$ & $e_{2}
e_{1}=e_{2}$ &\\
&& $e_{1}\ast e_{1}=\beta e_{1}+e_{2}$ &
 $e_{1}\ast e_{2}=\gamma e_{2}$ & $%
e_{2}\ast e_{1}=\beta e_{2}$  & \\ 
 \hline

$ \mathfrak{CpL}^2_{13}(\alpha,\beta,\gamma)$ &  
$\mathcal{C}_{38}^{\alpha ,\beta ,\gamma }$ &$
e_{1} e_{1}=e_{1} $ & $e_{2} e_{2}=e_{2}$   &  & \\ && 
$e_{1}\ast e_{1}=\left( \gamma +\beta -\alpha \right) e_{1}-\beta e_{2}$ &
\multicolumn{2}{l}{$e_{1}\ast e_{2}=\alpha e_{1}+\beta e_{2}$  }& \\ && 
 $e_{2}\ast e_{1}=\alpha e_{1}+\beta e_{2}$ & 
 \multicolumn{2}{l}{$e_{2}\ast e_{2}=-\alpha e_{1}+\gamma e_{2}$  }&\\

 \hline

$ \mathfrak{CpL}^2_{14}(\alpha,\beta,\gamma)$ &  
$\mathcal{C}_{40}^{\alpha ,\beta ,\gamma }$& 
$e_{1} e_{1}=e_{1}$ & $e_{1} e_{2}=2e_{2}$ & $e_{2} e_{1}=%
\frac{1}{2}e_{1}+e_{2}$ & \\ 
&& $e_{2} e_{2}=e_{2}$ & 
\multicolumn{2}{l}{$e_{1}\ast e_{1}=\alpha e_{1}+\beta e_{2}$} & \\
&& $e_{1}\ast e_{2}=2\alpha e_{2}$ & $e_{2}\ast e_{1}=\gamma e_{1}+\alpha e_{2}$ & $e_{2}\ast e_{2}=2\gamma
e_{2}$ &   \\
\hline

\end{longtable}
}

\section{The level classification of algebras}
Throughout this section, we summarize the level classification of different varieties of (not necessarily associative) algebras over the field $\mathbb{C}$. In what follows, we will not refer to the base field anymore.


Let us establish some notation that will be used throughout this section.
Let $m,m_1,\dots,m_k$ be positive integers such that $\sum_{i=1}^k m_i=m$, and let $\mathfrak{S}_q$ be the symmetric group on $q$ elements. There is a natural action of $\mathfrak{S}_{m_1}\times\dots\times\mathfrak{S}_{m_k}$ on $\mathbb{C}^m$: the symmetric group $\mathfrak{S}_{m_l}$ permutes the components located from the position number $\sum_{i=1}^{l-1}m_i+1$ to the position number $\sum_{i=1}^{l}m_i$ in $\mathbb{C}^m$. We will denote by $K_{m_1,\dots,m_k}$ the fixed set of representatives of this action. There is also an action of $\mathbb{C}^{\star}$ on $\mathbb{C}^m$ by multiplication. Both actions commute and stabilize the zero element. Then, $K_{m_1,\dots,m_k}^{\star}$ will denote a fixed set of representatives of the action of $\mathfrak{S}_{m_1}\times\dots\times\mathfrak{S}_{m_k}\times \mathbb{C}^{\star}$ on $\mathbb{C}^m\setminus (0,\dots,0)$.

\subsection{Algebras of level one} 
The complete classification of algebras of level one was given in~\cite{khud13}.
It was proven that every $2$-dimensional algebra of level one is isomorphic to one of the following algebras:

{\small
\vspace{-5mm}
\begin{longtable}{|l|l|lll|}
\caption*{ }   \\

\hline \multicolumn{1}{|c|}{${\mathcal A}$} & \multicolumn{1}{c|}{ } & \multicolumn{3}{c|}{Multiplication table} \\ \hline 
\endfirsthead

 \multicolumn{5}{l}%
{{\bfseries  continued from previous page}} \\
\hline \multicolumn{1}{|c|}{${\mathcal A}$} & \multicolumn{1}{c|}{ } & \multicolumn{3}{c|}{Multiplication table} \\ \hline 
\endhead

\hline \multicolumn{5}{|r|}{{Continued on next page}} \\ \hline
\endfoot

\hline 
\endlastfoot

$^1\mathbf{A}_{1}^{2}$ & $p_2^{-}$ & $e_1e_2=e_2$ & $e_2e_1=-e_2$ &\\ \hline
$^1\mathbf{A}_{2}^{2}$ & $\lambda_2$ & $e_1e_1=e_2$ && \\ \hline
$^1\mathbf{A}_{3}^{2}(\alpha)$ & $\nu_2(\alpha)$ & $e_1e_1=e_1$ & $e_1e_2=\alpha e_2$ & $e_2e_1=(1-\alpha)e_2$  \\ \hline
\end{longtable}
}

In dimension $n\geq 3$, the classification is the following:

{\small
\vspace{-5mm}
\begin{longtable}{|l|l|lll|}
\caption*{ }   \\

\hline \multicolumn{1}{|c|}{${\mathcal A}$} & \multicolumn{1}{c|}{ } & \multicolumn{3}{c|}{Multiplication table} \\ \hline 
\endfirsthead

 \multicolumn{5}{l}%
{{\bfseries  continued from previous page}} \\
\hline \multicolumn{1}{|c|}{${\mathcal A}$} & \multicolumn{1}{c|}{ } & \multicolumn{3}{c|}{Multiplication table} \\ \hline 
\endhead

\hline \multicolumn{5}{|r|}{{Continued on next page}} \\ \hline
\endfoot

\hline 
\endlastfoot

$^1\mathbf{A}_{1}^{n}$ & $p_n^{-}$ & $e_1e_i=e_i$ & $e_ie_1=-e_i$ & \\ \hline
$^1\mathbf{A}_{2}^{n}$ & $\lambda_2\oplus a_{n-2}$ & $e_1e_1=e_2$ && \\ \hline
$^1\mathbf{A}_{3}^{n}$ & $n_3^{-}\oplus a_{n-3}$ & $e_1e_2=e_3$ & $e_2e_1=-e_3$ & \\ \hline
$^1\mathbf{A}_{4}^{n}(\alpha)$ & $\nu_n(\alpha)$ & $e_1e_1=e_1$ & $e_1e_i=\alpha e_i$ & $e_ie_1=(1-\alpha)e_i,$  \\ \hline
\end{longtable}
}\noindent for $2\leq i \leq n$.

\subsection{Algebras of level two} 
The classification of all  algebras of level two is considerably more complex than that of level one, and was accomplished in~\cite{kv17}.

Each $2$-dimensional algebra of level two is isomorphic to one of the following algebras:

{\small

\begin{longtable}{|l|l|lll|}
\caption*{ }   \\

\hline \multicolumn{1}{|c|}{${\mathcal A}$} & \multicolumn{1}{c|}{ } & \multicolumn{3}{c|}{Multiplication table} \\ \hline 
\endfirsthead

 \multicolumn{5}{l}%
{{\bfseries  continued from previous page}} \\
\hline \multicolumn{1}{|c|}{${\mathcal A}$} & \multicolumn{1}{c|}{ } & \multicolumn{3}{c|}{Multiplication table} \\ \hline 
\endhead

\hline \multicolumn{5}{|r|}{{Continued on next page}} \\ \hline
\endfoot

\hline 
\endlastfoot

$^2\mathbf{A}_{1}^{2}$ & $\mathbf{A}_2$ & $e_1e_1=e_2$ & $e_1e_2=e_2$ & $e_2e_1=-e_2$ \\ \hline
$^2\mathbf{A}_{2}^{2}$ & $\mathbf{E}_4$ & $e_1e_1=e_1$ & $e_1e_2=e_1+e_2$ & $e_2e_2=e_2$  \\ \hline
$^2\mathbf{A}_{3}^{2}(\alpha)$ & $\mathbf{A}_1^{\alpha}$ & $e_1e_1=e_1+e_2$ & $e_1e_2=\alpha e_2$ & $e_2e_1=(1-\alpha)e_2$ \\ \hline
$^2\mathbf{A}_{4}^{2}(\alpha)$ & $\mathbf{B}_2^{\alpha}$ & $e_1e_2=\alpha e_2$ & $e_2e_1=(1-\alpha)e_2$ & \\ \hline
$^2\mathbf{A}_{5}^{2}(\alpha,\beta),\ \alpha+\beta\neq 1$ & $\mathbf{D}_2^{\alpha,\beta}$ & $e_1e_1=e_1$ & $e_1e_2=\alpha e_2$ & $e_2e_1=\beta e_2$ \\ \hline
\end{longtable}
}

Each $3$-dimensional algebra of level two is isomorphic to one of the following algebras:

{\small
\vspace{-5mm}
\begin{longtable}{|l|l|ll|}
\caption*{ }   \\

\hline \multicolumn{1}{|c|}{${\mathcal A}$} & \multicolumn{1}{c|}{ } & \multicolumn{2}{c|}{Multiplication table} \\ \hline 
\endfirsthead

 \multicolumn{4}{l}%
{{\bfseries  continued from previous page}} \\
\hline \multicolumn{1}{|c|}{${\mathcal A}$} & \multicolumn{1}{c|}{ } & \multicolumn{2}{c|}{Multiplication table} \\ \hline 
\endhead

\hline \multicolumn{4}{|r|}{{Continued on next page}} \\ \hline
\endfoot

\hline 
\endlastfoot
$^2\mathbf{A}_{1}^{3}$ & $\mathbb{C}\rtimes_1\mathbf{A}_2$ & $e_1e_1=e_2$ & $e_1e_2=e_2$ \\&& $e_1e_3=e_3$ & $e_2e_1=-e_2$ \\&& $e_3e_1=-e_3$ &\\ \hline
$^2\mathbf{A}_{2}^{3}$ & $\mathbb{C}\rtimes\mathbf{E}_4$ & $e_1e_1=e_1$ & $e_1e_2=e_1+e_2$ \\&& $e_1e_3=e_3$ & $e_2e_2=e_2$ \\&& $e_3e_2=e_3$ &\\ \hline
$^2\mathbf{A}_{3}^{3}(\alpha)$ & $\mathbb{C}\rtimes_{\alpha} \mathbf{A}_1^{\alpha}$ & $e_1e_1=e_1+e_2$ & $e_1e_2=\alpha e_2$ \\&& $e_1e_3=\alpha e_3$ & $e_2e_1=(1-\alpha)e_2$ \\&& $e_3e_1=(1-\alpha)e_3$ &\\ \hline
$^2\mathbf{A}_{4}^{3}(\alpha)$ & $\mathbb{C}\rtimes_{0}^t\mathbf{B}_2^{\alpha}$ & $e_1e_2=\alpha e_2$ & $e_1e_3=\alpha e_3$ \\&& $e_2e_1=(1-\alpha)e_2$ & $e_3e_1=(1-\alpha)e_3$  \\ \hline
$^2\mathbf{A}_{5}^{3}(\alpha,\beta),\ \alpha+\beta\neq 1$ & $\mathbb{C}\rtimes_{0}^t\mathbf{D}_2^{\alpha,\beta}$ & $e_1e_1=e_1$ & $e_1e_2=\alpha e_2$ \\&& $e_1e_3=\alpha e_3$ & $e_2e_1=\beta e_2$ \\&& $e_3e_1=\beta e_3$ &\\ \hline
$^2\mathbf{A}_{6}^{3}(\alpha,\beta),\ (\alpha,\beta)\in K_{2}^{\star}$ & $\mathbf{F}^{\alpha,\beta}$ & $e_1e_1=e_3$ & $e_1e_2=\alpha e_3$ \\&& $e_2e_1=\beta e_3$ & \\ \hline
$^2\mathbf{A}_{7}^{3}(\alpha,\beta),\ (\alpha,\beta)\in K_{2}^{\star}$ & $\mathbf{T}_0^{2,\overline{\alpha,\beta}}$ & $e_1e_2=\alpha e_2+e_3$ & $e_1e_3=\beta e_3$ \\&& $e_2e_1=-\alpha e_2-e_3$ & $e_3e_1=-\beta e_3$ \\ \hline
$^2\mathbf{A}_{8}^{3}(\alpha,\beta),\ (\alpha,\beta)\in K_{2}$ & $\mathbf{T}_1^{2,\overline{\alpha,\beta}}$ & $e_1e_1=e_1$ & $e_1e_2=\alpha e_2+e_3$ \\&& $e_1e_3=\beta e_3$ & $e_2e_1=(1-\alpha) e_2-e_3$ \\&& $e_3e_1=(1-\beta)e_3. $ &\\ \hline
\end{longtable}
}

Every $4$-dimensional algebra of level two is isomorphic to one of the following algebras:

{\small

\begin{longtable}{|l|l|ll|}
\caption*{ }   \\

\hline \multicolumn{1}{|c|}{${\mathcal A}$} & \multicolumn{1}{c|}{ } & \multicolumn{2}{c|}{Multiplication table} \\ \hline 
\endfirsthead

 \multicolumn{4}{l}%
{{\bfseries  continued from previous page}} \\
\hline \multicolumn{1}{|c|}{${\mathcal A}$} & \multicolumn{1}{c|}{ } & \multicolumn{2}{c|}{Multiplication table} \\ \hline 
\endhead

\hline \multicolumn{4}{|r|}{{Continued on next page}} \\ \hline
\endfoot

\hline 
\endlastfoot

$^2\mathbf{A}_{1}^{4}$ & $\mathbb{C}^2\rtimes_1\mathbf{A}_2$ & $e_1e_1=e_2$ & $e_1e_2=e_2$ \\&& $e_1e_3=e_3$ & $e_1e_4=e_4$ \\&& $e_2e_1=-e_2$  & $e_3e_1=-e_3$ \\ && $e_4e_1=-e_4$ &\\ \hline
$^2\mathbf{A}_{2}^{4}$ & $\mathbb{C}^2\rtimes\mathbf{E}_4$ & $e_1e_1=e_1$ & $e_1e_2=e_1+e_2$ \\&& $e_1e_3=e_3$ 
& $e_1e_4=e_4$ 
\\&&  $e_2e_2=e_2$  & $e_3e_2=e_3$ 
\\&& $e_4e_2=e_4$ &\\ \hline
$^2\mathbf{A}_{3}^{4}$ & $T_0^3$ & $e_1e_2=e_3$ & $e_1e_3=e_4$ \\&& $e_2e_1=-e_3$ & $e_3e_1=-e_4$ \\ \hline
$^2\mathbf{A}_{4}^{4}(\alpha)$ & $\mathbb{C}^2\rtimes_{\alpha} \mathbf{A}_1^{\alpha}$  
& $e_1e_1=e_1+e_2$ & $e_1e_2=\alpha e_2$ \\&& $e_1e_3=\alpha e_3$ & $e_1e_4=\alpha e_4$ \\&& $e_2e_1=(1-\alpha)e_2$ & $e_3e_1=(1-\alpha)e_3$ \\ &&  $e_4e_1=(1-\alpha)e_4$ &\\ \hline
$^2\mathbf{A}_{5}^{4}(\alpha)$ & $\mathbb{C}^2\rtimes_{0}^t\mathbf{B}_2^{\alpha}$ & $e_1e_2=\alpha e_2$ & $e_1e_3=\alpha e_3$ \\&& $e_1e_4=\alpha e_4$ & $e_2e_1=(1-\alpha)e_2$ \\&&  $e_3e_1=(1-\alpha)e_3$  & $e_4e_1=(1-\alpha) e_4$ \\ \hline
$^2\mathbf{A}_{6}^{4}(\alpha,\beta),\ \alpha+\beta\neq 1$ & $\mathbb{C}^2\rtimes_{0}^t\mathbf{D}_2^{\alpha,\beta}$ & $e_1e_1=e_1$ & $e_1e_2=\alpha e_2$  \\&& $e_1e_3=\alpha e_3$ & $e_1e_4=\alpha e_4$ \\&& $e_2e_1=\beta e_2$  & $e_3e_1=\beta e_3$ \\&& $e_4e_1=\beta e_4$ & \\ \hline
$^2\mathbf{A}_{7}^{4}(\alpha,\beta),\ (\alpha,\beta)\in K_{2}^{\star}$ & $\mathbf{F}^{\alpha,\beta}\oplus \mathbb{C}$ & $e_1e_1=e_3$ & $e_1e_2=\alpha e_3$ \\&& $e_2e_1=\beta e_3$ & \\ \hline
$^2\mathbf{A}_{8}^{4}(\alpha,\beta),\ (\alpha,\beta)\in K_{1,1}^{\star}$ & $T_0^{2,\overline{\alpha,\beta}}$ & $e_1e_2=\alpha e_2+e_3$ & $e_1e_3=\beta e_3$ \\&& $e_1e_4=\beta e_4$ & $e_2e_1=-\alpha e_2-e_3$ \\&& $e_3e_1=-\beta e_3$ & $e_4e_1=-\beta e_4$\\ \hline
$^2\mathbf{A}_{9}^{4}(\alpha,\beta),\ (\alpha,\beta)\in K_{1,1}$ & $T_1^{2,\overline{\alpha,\beta}}$  & $e_1e_1=e_1$ & $e_1e_2=\alpha e_2+e_3$ \\&& $e_1e_3=\beta e_3$ & $e_1e_4=\beta e_4$ \\&& $e_2e_1=(1-\alpha) e_2-e_3$  & $e_3e_1=(1-\beta)e_3$ \\ && $e_4e_1=(1-\beta)e_4$& \\ \hline
\end{longtable}
}
\vspace{-3mm}

Every $n$-dimensional  algebra of level two, $n\geq 5$, is isomorphic to one of the following algebras:

{\small
\vspace{-9mm}
\begin{longtable}{|l|l|ll|}
\caption*{ }   \\

\hline \multicolumn{1}{|c|}{${\mathcal A}$} & \multicolumn{1}{c|}{ } & \multicolumn{2}{c|}{Multiplication table} \\ \hline 
\endfirsthead

 \multicolumn{4}{l}%
{{\bfseries  continued from previous page}} \\
\hline \multicolumn{1}{|c|}{${\mathcal A}$} & \multicolumn{1}{c|}{ } & \multicolumn{2}{c|}{Multiplication table} \\ \hline 
\endhead

\hline \multicolumn{4}{|r|}{{Continued on next page}} \\ \hline
\endfoot

\hline 
\endlastfoot

$^2\mathbf{A}_{1}^{n}$ & $\mathbb{C}^{n-2}\rtimes_1\mathbf{A}_2$  & $e_1e_1=e_2$ & $e_1e_i=e_i$ \\&& $e_ie_1=-e_i$ & \\ \hline
$^2\mathbf{A}_{2}^{n}$ & $\mathbb{C}^{n-2}\rtimes\mathbf{E}_4$ & $e_1e_1=e_1$ & $e_1e_2=e_1+e_2$ \\&& $e_1e_j=e_j$ &  $e_2e_2=e_2$ \\&& $e_je_2=e_j$ &\\ \hline
$^2\mathbf{A}_{3}^{n}$ & $T_0^{2,2}$ & $e_1e_2=e_3$ & $e_1e_4=e_5$ \\&& $e_2e_1=-e_3$ & $e_4e_1=-e_5$ \\ \hline
$^2\mathbf{A}_{4}^{n}$ & $\eta_2\oplus \mathbb{C}^{n-5}$ & $e_1e_2=e_5$  & $e_2e_1=-e_5$ \\&& $e_3e_4=e_5$ & $e_4e_3=-e_5$ \\ \hline
$^2\mathbf{A}_{5}^{n}(\alpha)$ & $\mathbb{C}^{n-2}\rtimes_{\alpha} \mathbf{A}_1^{\alpha}$ & $e_1e_1=e_1+e_2$ & $e_1e_i=\alpha e_i$ \\&& $e_ie_1=(1-\alpha)e_i$ & \\ \hline
$^2\mathbf{A}_{6}^{n}(\alpha)$ & $\mathbb{C}^{n-2}\rtimes_{0}^t\mathbf{B}_2^{\alpha}$ & $e_1e_i=\alpha e_i$ &  $e_ie_1=(1-\alpha)e_i$ \\ \hline
$^2\mathbf{A}_{7}^{n}(\alpha,\beta),\ \alpha+\beta\neq 1$ & $\mathbb{C}^{n-2}\rtimes_{0}^t\mathbf{D}_2^{\alpha,\beta}$ & $e_1e_1=e_1$ & $e_1e_i=\alpha e_i$  \\&& $e_ie_1=\beta e_i$ & \\ \hline
$^2\mathbf{A}_{8}^{n}(\alpha,\beta),\ (\alpha,\beta)\in K_{2}^{\star}$ & $\mathbf{F}^{\alpha,\beta}\oplus \mathbb{C}^{n-3}$ & $e_1e_1=e_3$ & $e_1e_2=\alpha e_3$ \\&& $e_2e_1=\beta e_3$ & \\ \hline
$^2\mathbf{A}_{9}^{n}(\alpha,\beta),\ (\alpha,\beta)\in K_{1,1}^{\star}$ & $T_0^{2,\overline{\alpha,\beta}}$ & $e_1e_2=\alpha e_2+e_3$ & $e_1e_j=\beta e_j$ \\&& $e_2e_1=-\alpha e_2-e_3$ & $e_je_1=-\beta e_j$\\ \hline
$^2\mathbf{A}_{10}^{n}(\alpha,\beta),\ (\alpha,\beta)\in K_{1,1}$ & $T_1^{2,\overline{\alpha,\beta}}$ & $e_1e_1=e_1$ & $e_1e_2=\alpha e_2+e_3$ \\&& $e_1e_j=\beta e_j$ & $e_2e_1=(1-\alpha) e_2-e_3$ \\&& $e_je_1=(1-\beta)e_j,$ & \\ \hline
\end{longtable}
}\noindent for $2\leq i\leq n$ and $3\leq j\leq n$.

\subsection{Nilpotent algebras} 
The nilpotent algebras of level one can be selected from the general results of~\cite{khud13}, and those of level two, from~\cite{kv17}.

\subsubsection{Nilpotent algebras of level one}
There exists only one  $2$-dimensional nilpotent algebras of level one, up to isomorphism:

{\small

\begin{longtable}{|l|l|l|}
\caption*{ }   \\

\hline \multicolumn{1}{|c|}{${\mathcal A}$} & \multicolumn{1}{c|}{ } & \multicolumn{1}{c|}{Multiplication table} \\ \hline 
\endfirsthead

 \multicolumn{3}{l}%
{{\bfseries  continued from previous page}} \\
\hline \multicolumn{1}{|c|}{${\mathcal A}$} & \multicolumn{1}{c|}{ } & \multicolumn{1}{c|}{Multiplication table} \\ \hline 
\endhead

\hline \multicolumn{3}{|r|}{{Continued on next page}} \\ \hline
\endfoot

\hline 
\endlastfoot
$^1\mathbf{N}_{1}^{2}$ & $\lambda_2$ & $e_1e_1=e_2$  \\ \hline
\end{longtable}
}

In dimension $n\geq 3$, there exist two:

{\small
\vspace{-5mm}
\begin{longtable}{|l|l|ll|}
\caption*{ }   \\

\hline \multicolumn{1}{|c|}{${\mathcal A}$} & \multicolumn{1}{c|}{ } & \multicolumn{2}{c|}{Multiplication table} \\ \hline 
\endfirsthead

 \multicolumn{4}{l}%
{{\bfseries  continued from previous page}} \\
\hline \multicolumn{1}{|c|}{${\mathcal A}$} & \multicolumn{1}{c|}{ } & \multicolumn{2}{c|}{Multiplication table} \\ \hline 
\endhead

\hline \multicolumn{4}{|r|}{{Continued on next page}} \\ \hline
\endfoot

\hline 
\endlastfoot

$^1\mathbf{N}_{1}^{n}$ & $\lambda_2\oplus a_{n-2}$ & $e_1e_1=e_2$ & \\ \hline
$^1\mathbf{N}_{2}^{n}$ & $n_3^{-}\oplus a_{n-3}$ & $e_1e_2=e_3$ & $e_2e_1=-e_3$ \\ \hline 
\end{longtable}
}

\subsubsection{Nilpotent algebras of level two}
In this section, we correct some inaccuracies of~\cite{fkrv}.

There are no nilpotent algebras of level two in dimension $2$; in dimension $3$, there exists only one family:

{\small

\begin{longtable}{|l|l|lll|}
\caption*{ }   \\

\hline \multicolumn{1}{|c|}{${\mathcal A}$} & \multicolumn{1}{c|}{ } & \multicolumn{3}{c|}{Multiplication table} \\ \hline 
\endfirsthead

 \multicolumn{5}{l}%
{{\bfseries  continued from previous page}} \\
\hline \multicolumn{1}{|c|}{${\mathcal A}$} & \multicolumn{1}{c|}{ } & \multicolumn{3}{c|}{Multiplication table} \\ \hline 
\endhead

\hline \multicolumn{5}{|r|}{{Continued on next page}} \\ \hline
\endfoot

\hline 
\endlastfoot

$^2\mathbf{N}_{1}^{3}(\alpha,\beta),\ (\alpha,\beta)\in K_{2}^{\star}$ & $\mathbf{F}^{\alpha,\beta}$  & $e_1e_1=e_3$ & $e_1e_2=\alpha e_3$ & $e_2e_1=\beta e_3$  \\ \hline 
\end{longtable}
}

In dimension $4$, we find one algebra and one family, namely:

{\small
\vspace{-5mm}
\begin{longtable}{|l|l|llll|}
\caption*{ }   \\

\hline \multicolumn{1}{|c|}{${\mathcal A}$} & \multicolumn{1}{c|}{ } & \multicolumn{4}{c|}{Multiplication table} \\ \hline 
\endfirsthead

 \multicolumn{6}{l}%
{{\bfseries  continued from previous page}} \\
\hline \multicolumn{1}{|c|}{${\mathcal A}$} & \multicolumn{1}{c|}{ } & \multicolumn{4}{c|}{Multiplication table} \\ \hline 
\endhead

\hline \multicolumn{6}{|r|}{{Continued on next page}} \\ \hline
\endfoot

\hline 
\endlastfoot

$^2\mathbf{N}_{1}^{4}$ & $T_0^3$ & $e_1e_2=e_3$ & $e_1e_3=e_4$ & $e_2e_1=-e_3$ & $e_3e_1=-e_4$ \\ \hline
$^2\mathbf{N}_{2}^{4}(\alpha,\beta),\ (\alpha,\beta)\in K_{2}^{\star}$ & $\mathbf{F}^{\alpha,\beta}\oplus \mathbb{C}$ & $e_1e_1=e_3$ & $e_1e_2=\alpha e_3$ & $e_2e_1=\beta e_3$ & \\ \hline
\end{longtable}
}

In dimension $n\geq 5$, the classification of nilpotent algebras of level two is as follows:

{\small

\begin{longtable}{|l|l|llll|}
\caption*{ }   \\

\hline \multicolumn{1}{|c|}{${\mathcal A}$} & \multicolumn{1}{c|}{ } & \multicolumn{4}{c|}{Multiplication table} \\ \hline 
\endfirsthead

 \multicolumn{6}{l}%
{{\bfseries  continued from previous page}} \\
\hline \multicolumn{1}{|c|}{${\mathcal A}$} & \multicolumn{1}{c|}{ } & \multicolumn{4}{c|}{Multiplication table} \\ \hline 
\endhead

\hline \multicolumn{6}{|r|}{{Continued on next page}} \\ \hline
\endfoot

\hline 
\endlastfoot
$^2\mathbf{N}_{1}^{n}$ & $ T_0^{2,2}$ & $e_1e_2=e_3$ & $e_1e_4=e_5$ & $e_2e_1=-e_3$ & $e_4e_1=-e_5$ \\ \hline
$^2\mathbf{N}_{2}^{n}$ & $\eta_2\oplus \mathbb{C}^{n-5}$ & $e_1e_2=e_5$  & $e_2e_1=-e_5$ & $e_3e_4=e_5$ & $e_4e_3=-e_5$  \\ \hline
$^2\mathbf{N}_{3}^{n}(\alpha,\beta),\ (\alpha,\beta)\in K_{2}^{\star}$ & $\mathbf{F}^{\alpha,\beta}\oplus \mathbb{C}^{n-3}$ & $e_1e_1=e_3$ & $e_1e_2=\alpha e_3$ & $e_2e_1=\beta e_3$ & \\ \hline
\end{longtable}
}

\subsection{Commutative algebras} 

Thanks to~\cite{khud13}, we know that every $n$-dimensional commutative algebra of level one, with $n\geq 2$,  is isomorphic to one of the following two algebras:

{\small

\begin{longtable}{|l|l|ll|}
\caption*{ }   \\

\hline \multicolumn{1}{|c|}{${\mathcal A}$} & \multicolumn{1}{c|}{ } & \multicolumn{2}{c|}{Multiplication table} \\ \hline 
\endfirsthead

 \multicolumn{4}{l}%
{{\bfseries  continued from previous page}} \\
\hline \multicolumn{1}{|c|}{${\mathcal A}$} & \multicolumn{1}{c|}{ } & \multicolumn{2}{c|}{Multiplication table} \\ \hline 
\endhead

\hline \multicolumn{4}{|r|}{{Continued on next page}} \\ \hline
\endfoot

\hline 
\endlastfoot
$^1\mathbf{C}_{1}^{2}$ & $\lambda_n \oplus a_{n-2}$ & $e_1e_1=e_2$  &\\ \hline
$^1\mathbf{C}_{2}^{2}$ & $\nu_n(\frac{1}{2})$ & $e_1e_1=e_1$ & $e_1e_i=\frac{1}{2} e_i,$ \\ \hline 
\end{longtable}
}\noindent for $2\leq i\leq n$.

In~\cite{kv17}, 
 the classification of the commutative algebras of level two was
also presented. 
Each commutative algebra of dimension $2$ and level two is isomorphic to one of the following algebras: 

{\small

\begin{longtable}{|l|l|ll|}
\caption*{ }   \\

\hline \multicolumn{1}{|c|}{${\mathcal A}$} & \multicolumn{1}{c|}{ } & \multicolumn{2}{c|}{Multiplication table} \\ \hline 
\endfirsthead

 \multicolumn{4}{l}%
{{\bfseries  continued from previous page}} \\
\hline \multicolumn{1}{|c|}{${\mathcal A}$} & \multicolumn{1}{c|}{ } & \multicolumn{2}{c|}{Multiplication table} \\ \hline 
\endhead

\hline \multicolumn{4}{|r|}{{Continued on next page}} \\ \hline
\endfoot

\hline 
\endlastfoot

$^2\mathbf{C}_{1}^{2}$ & $\mathbf{A}_1^{\frac{1}{2}}$ & $e_1e_1=e_1+e_2$ & $e_1e_2=\frac{1}{2} e_2$  \\ \hline
$^2\mathbf{C}_{2}^{2}$ & $\mathbf{B}_2^{\frac{1}{2}}$ & $e_1e_2=\frac{1}{2} e_2$ & \\ \hline
$^2\mathbf{C}_{3}^{2}(\alpha),\ \alpha\neq \frac{1}{2}$ & $\mathbf{D}_2^{\alpha,\alpha}$ & $e_1e_1=e_1$ & $e_1e_2=\alpha e_2$  \\ \hline
\end{longtable}
}

In dimension $n\geq 3$, the commutative algebras of level two are as follows: 

{\small

\begin{longtable}{|l|l|ll|}
\caption*{ }   \\

\hline \multicolumn{1}{|c|}{${\mathcal A}$} & \multicolumn{1}{c|}{ } & \multicolumn{2}{c|}{Multiplication table} \\ \hline 
\endfirsthead

 \multicolumn{4}{l}%
{{\bfseries  continued from previous page}} \\
\hline \multicolumn{1}{|c|}{${\mathcal A}$} & \multicolumn{1}{c|}{ } & \multicolumn{2}{c|}{Multiplication table} \\ \hline 
\endhead

\hline \multicolumn{4}{|r|}{{Continued on next page}} \\ \hline
\endfoot

\hline 
\endlastfoot

$^2\mathbf{C}_{1}^{n}$ & $\mathbb{C}^{n-2}\rtimes_{\frac{1}{2}} \mathbf{A}_1^{\frac{1}{2}}$ & $e_1e_1=e_1+e_2$ & $e_1e_i=\frac{1}{2} e_i$ \\ \hline
$^2\mathbf{C}_{2}^{n}$ & $\mathbb{C}^{n-2}\rtimes_{0}^t\mathbf{B}_2^{\frac{1}{2}}$ & $e_1e_i=\frac{1}{2} e_i$ & \\ \hline
$^2\mathbf{C}_{3}^{n}$ & $\mathbf{F}^{1,1}\oplus \mathbb{C}^{n-3}$ & $e_1e_1=e_3$ & $e_1e_2= e_3$  \\ \hline
$^2\mathbf{C}_{4}^{n}(\alpha),\ \alpha\neq \frac{1}{2}$ & $\mathbb{C}^{n-2}\rtimes_{0}^t\mathbf{D}_2^{\alpha,\alpha}$ & $e_1e_1=e_1$ & $e_1e_i=\alpha e_i$,  \\ \hline
\end{longtable}
}\noindent for $2\leq i\leq n$.

\subsection{Anticommutative algebras}\label{ss:la}
Anticommutative algebras of levels one and two are completely classified in~\cite{gorb91} and~\cite{kv17}, respectively. For higher levels, we will impose the condition of being Engel to classify the algebras (see~\cite{wolf1}).
The algebra ${\mathcal A}$ is called $m$-Engel if $(L_a)^m = 0$ for any $a \in {\mathcal A}.$
Here we use the notation $L_a$  for the operator of left multiplication in $\mathfrak{A}$.
We will call the algebra ${\mathcal A}$ Engel if it is $m$-Engel for some $m > 0$.

\subsubsection{Anticommutative algebras of level one}

The unique $2$-dimensional anticommutative algebra of level one is:

{\small

\begin{longtable}{|l|l|l|}
\caption*{ }   \\

\hline \multicolumn{1}{|c|}{${\mathcal A}$} & \multicolumn{1}{c|}{ } & \multicolumn{1}{c|}{Multiplication table} \\ \hline 
\endfirsthead

 \multicolumn{3}{l}%
{{\bfseries  continued from previous page}} \\
\hline \multicolumn{1}{|c|}{${\mathcal A}$} & \multicolumn{1}{c|}{ } & \multicolumn{1}{c|}{Multiplication table} \\ \hline 
\endhead

\hline \multicolumn{3}{|r|}{{Continued on next page}} \\ \hline
\endfoot

\hline 
\endlastfoot
$^{1}\mathbf{AC}_{1}^{2}$ & $p_2^{-}$ & $e_1e_2=e_2$  \\ \hline
\end{longtable}
}

In dimension $n\geq 3$, any anticommutative algebra of level one is isomorphic to one of the following two algebras: 

{\small

\begin{longtable}{|l|l|l|}
\caption*{ }   \\

\hline \multicolumn{1}{|c|}{${\mathcal A}$} & \multicolumn{1}{c|}{ } & \multicolumn{1}{c|}{Multiplication table} \\ \hline 
\endfirsthead

 \multicolumn{3}{l}%
{{\bfseries  continued from previous page}} \\
\hline \multicolumn{1}{|c|}{${\mathcal A}$} & \multicolumn{1}{c|}{ } & \multicolumn{1}{c|}{Multiplication table} \\ \hline 
\endhead

\hline \multicolumn{3}{|r|}{{Continued on next page}} \\ \hline
\endfoot

\hline 
\endlastfoot

$^1\mathbf{AC}_{1}^{n}$ & $p_n^{-}$  & $e_1e_i=e_i$ \\ \hline
$^1\mathbf{AC}_{2}^{n}$ & $n_3^{-}\oplus a_{n-3}$ & $e_1e_2=e_3,$  \\ \hline 
\end{longtable}
}\noindent for $2\leq i \leq n$.

Note that all these algebras are Lie. On the other hand, only $^1\mathbf{AC}_{2}^{n}$ is Engel, for $n\geq 3$.

\subsubsection{Anticommutative algebras of level two}
There are no anticommutative algebras of level two and dimension $2$, and in dimension $3$ there exists only one:

{\small

\begin{longtable}{|l|l|ll|}
\caption*{ }   \\

\hline \multicolumn{1}{|c|}{${\mathcal A}$} & \multicolumn{1}{c|}{ } & \multicolumn{2}{c|}{Multiplication table} \\ \hline 
\endfirsthead

 \multicolumn{4}{l}%
{{\bfseries  continued from previous page}} \\
\hline \multicolumn{1}{|c|}{${\mathcal A}$} & \multicolumn{1}{c|}{ } & \multicolumn{2}{c|}{Multiplication table} \\ \hline 
\endhead

\hline \multicolumn{4}{|r|}{{Continued on next page}} \\ \hline
\endfoot

\hline 
\endlastfoot
$^2\mathbf{AC}_{1}^{3}(\alpha,\beta),\ (\alpha,\beta)\in K_{2}^{\star}$ & $\mathbf{T}_0^{2,\overline{\alpha,\beta}}$  & $e_1e_2=\alpha e_2+e_3$ & $e_1e_3=\beta e_3$  \\ \hline
\end{longtable}
}

In dimension $4$, it turns out that there are, up to isomorphism, two anticommutative algebras:

{\small

\begin{longtable}{|l|l|lll|}
\caption*{ }   \\

\hline \multicolumn{1}{|c|}{${\mathcal A}$} & \multicolumn{1}{c|}{ } & \multicolumn{3}{c|}{Multiplication table} \\ \hline 
\endfirsthead

 \multicolumn{5}{l}%
{{\bfseries  continued from previous page}} \\
\hline \multicolumn{1}{|c|}{${\mathcal A}$} & \multicolumn{1}{c|}{ } & \multicolumn{3}{c|}{Multiplication table} \\ \hline 
\endhead

\hline \multicolumn{5}{|r|}{{Continued on next page}} \\ \hline
\endfoot

\hline 
\endlastfoot
$^2\mathbf{AC}_{1}^{4}$ & $T_0^3$ & $e_1e_2=e_3$ & $e_1e_3=e_4$ & \\ \hline
$^2\mathbf{AC}_{2}^{4}(\alpha,\beta),\ (\alpha,\beta)\in K_{1,1}^{\star}$ & $T_0^{2,\overline{\alpha,\beta}}$ & $e_1e_2=\alpha e_2+e_3$ & $e_1e_3=\beta e_3$ & $e_1e_4=\beta e_4$  \\\hline
\end{longtable}
}

Finally, for $n\geq 5$, every $n$-dimensional anticommutative algebra of level two is isomorphic to one of the following algebras:

{\small
\vspace{-5mm}
\begin{longtable}{|l|l|ll|}
\caption*{ }   \\

\hline \multicolumn{1}{|c|}{${\mathcal A}$} & \multicolumn{1}{c|}{ } & \multicolumn{2}{c|}{Multiplication table} \\ \hline 
\endfirsthead

 \multicolumn{4}{l}%
{{\bfseries  continued from previous page}} \\
\hline \multicolumn{1}{|c|}{${\mathcal A}$} & \multicolumn{1}{c|}{ } & \multicolumn{2}{c|}{Multiplication table} \\ \hline 
\endhead

\hline \multicolumn{4}{|r|}{{Continued on next page}} \\ \hline
\endfoot

\hline 
\endlastfoot

$^2\mathbf{AC}_{1}^{n}$ & $T_0^{2,2}$ & $e_1e_2=e_3$ & $e_1e_4=e_5$  \\\hline
$^2\mathbf{AC}_{2}^{n}$ & $\eta_2\oplus \mathbb{C}^{n-5}$ & $e_1e_2=e_5$  & $e_3e_4=e_5$ \\\hline
$^2\mathbf{AC}_{3}^{n}(\alpha,\beta),\ (\alpha,\beta)\in K_{1,1}^{\star}$ & $T_0^{2,\overline{\alpha,\beta}}$ & $e_1e_2=\alpha e_2+e_3$ & $e_1e_j=\beta e_j,$ \\\hline
\end{longtable}
}\noindent for $3\leq j\leq n$.

Note that all anticommutative algebras of the level two are Lie algebras.
Also, any Engel anticommutative algebra of level two is isomorphic to $^2\mathbf{AC}_{1}^{4}$, to $^2\mathbf{AC}_{1}^{n}$ or to $^2\mathbf{AC}_{2}^{n}$, for $n\geq 5$.

\subsubsection{Engel anticommutative algebras of level three}
There are no Engel anticommutative algebras of level three and dimension at most $4$, and there exists only one, up to isomorphism, of dimension~$5$:

{\small
\vspace{-5mm}
\begin{longtable}{|l|l|ll|}
\caption*{ }   \\

\hline \multicolumn{1}{|c|}{${\mathcal A}$} & \multicolumn{1}{c|}{ } & \multicolumn{2}{c|}{Multiplication table} \\ \hline 
\endfirsthead

 \multicolumn{4}{l}%
{{\bfseries  continued from previous page}} \\
\hline \multicolumn{1}{|c|}{${\mathcal A}$} & \multicolumn{1}{c|}{ } & \multicolumn{2}{c|}{Multiplication table} \\ \hline 
\endhead

\hline \multicolumn{4}{|r|}{{Continued on next page}} \\ \hline
\endfoot

\hline 
\endlastfoot
$^3\mathbf{EAC}_{1}^{5}$ & $T^{3}$ & $e_1e_2=e_3$ & $e_1e_3=e_5$ \\\hline
\end{longtable}
}

In dimension $6$, we find that every Engel anticommutative algebra of level three is isomorphic to one of the following algebras:

{\small
\vspace{-5mm}
\begin{longtable}{|l|l|lll|}
\caption*{ }   \\

\hline \multicolumn{1}{|c|}{${\mathcal A}$} & \multicolumn{1}{c|}{ } & \multicolumn{3}{c|}{Multiplication table} \\ \hline 
\endfirsthead

 \multicolumn{5}{l}%
{{\bfseries  continued from previous page}} \\
\hline \multicolumn{1}{|c|}{${\mathcal A}$} & \multicolumn{1}{c|}{ } & \multicolumn{3}{c|}{Multiplication table} \\ \hline 
\endhead

\hline \multicolumn{5}{|r|}{{Continued on next page}} \\ \hline
\endfoot

\hline 
\endlastfoot
$^3\mathbf{EAC}_{1}^{6}$ & $T^{3}$ & $e_1e_2=e_3$ & $e_1e_3=e_6$ & \\\hline
$^3\mathbf{EAC}_{2}^{6}$ & $T^{2,2}(\epsilon_{23}^4)$ & $e_1e_2=e_5$  & $e_1e_3=e_6$ & $e_2e_3=e_4$ \\\hline
$^3\mathbf{EAC}_{3}^{6}$ & $T^{2,2}(\epsilon_{24}^6)$ & $e_1e_2=e_5$ & $e_1e_3=e_6$ & $e_2e_4=e_6$  \\\hline
\end{longtable}
}

Finally, in dimension $n\geq 7$ there exist the following Engel anticommutative algebras, up to isomorphism:

{\small

\begin{longtable}{|l|l|lll|}
\caption*{ }   \\

\hline \multicolumn{1}{|c|}{${\mathcal A}$} & \multicolumn{1}{c|}{ } & \multicolumn{3}{c|}{Multiplication table} \\ \hline 
\endfirsthead

 \multicolumn{5}{l}%
{{\bfseries  continued from previous page}} \\
\hline \multicolumn{1}{|c|}{${\mathcal A}$} & \multicolumn{1}{c|}{ } & \multicolumn{3}{c|}{Multiplication table} \\ \hline 
\endhead

\hline \multicolumn{5}{|r|}{{Continued on next page}} \\ \hline
\endfoot

\hline 
\endlastfoot
$^3\mathbf{EAC}_{1}^{n}$ & $\eta_3$ & $e_1e_2=e_7$ & $e_3e_4=e_7$ & $e_5e_6=e_7$ \\\hline
$^3\mathbf{EAC}_{2}^{n}$ & $ T^{2,2,2}$  & $e_1e_2=e_{n-2}$ & $e_1e_3=e_{n-1}$ & $e_1e_4=e_n$ \\\hline
$^3\mathbf{EAC}_{3}^{n}$ & $T^{3}$ & $e_1e_2=e_3$ & $e_1e_3=e_{n}$ & \\ \hline
$^3\mathbf{EAC}_{4}^{n}$ & $T^{2,2}(\epsilon_{23}^{n-2})$ & $e_1e_2=e_{n-1}$  & $e_1e_3=e_n$ & $e_2e_3=e_{n-2}$  \\\hline
$^3\mathbf{EAC}_{5}^{n}$ & $T^{2,2}(\epsilon_{24}^n)$ & $e_1e_2=e_{n-1}$ & $e_1e_3=e_n$ & $e_2e_4=e_n$  \\ \hline
\end{longtable}
}

Note that every Engel anticommutative algebra of level three is a Lie algebra.

\subsubsection{Engel anticommutative algebras of level four}

The Engel anticommutative algebras of level four have dimension at least $5$. In dimension $5$ there exist three, up to isomorphism:

{\small

\begin{longtable}{|l|l|lll|}
\caption*{ }   \\

\hline \multicolumn{1}{|c|}{${\mathcal A}$} & \multicolumn{1}{c|}{ } & \multicolumn{3}{c|}{Multiplication table} \\ \hline 
\endfirsthead

 \multicolumn{5}{l}%
{{\bfseries  continued from previous page}} \\
\hline \multicolumn{1}{|c|}{${\mathcal A}$} & \multicolumn{1}{c|}{ } & \multicolumn{3}{c|}{Multiplication table} \\ \hline 
\endhead

\hline \multicolumn{5}{|r|}{{Continued on next page}} \\ \hline
\endfoot

\hline 
\endlastfoot
$^4\mathbf{EAC}_{1}^{5}$ & $T^{4}$ & $e_1e_2=e_3$ & $e_1e_3=e_4$ & $e_1e_4=e_5$  \\ \hline
$^4\mathbf{EAC}_{2}^{5}$ & $T^{3}(\epsilon_{23}^4)$ & $e_1e_2=e_3$  & $e_1e_3=e_5$ & $e_2e_3=e_4$ \\ \hline
$^4\mathbf{EAC}_{3}^{5}$ & $T^{3}(\epsilon_{24}^5)$ & $e_1e_2=e_3$ & $e_1e_3=e_5$ & $e_2e_4=e_5$\\\hline
\end{longtable}
}

In dimension $6$, there exist four:

{\small

\begin{longtable}{|l|l|lll|}
\caption*{ }   \\

\hline \multicolumn{1}{|c|}{${\mathcal A}$} & \multicolumn{1}{c|}{ } & \multicolumn{3}{c|}{Multiplication table} \\ \hline 
\endfirsthead

 \multicolumn{5}{l}%
{{\bfseries  continued from previous page}} \\
\hline \multicolumn{1}{|c|}{${\mathcal A}$} & \multicolumn{1}{c|}{ } & \multicolumn{3}{c|}{Multiplication table} \\ \hline 
\endhead

\hline \multicolumn{5}{|r|}{{Continued on next page}} \\ \hline
\endfoot

\hline 
\endlastfoot
$^4\mathbf{EAC}_{1}^{6}$ & $T^{3,2}$  & $e_1e_2=e_{5}$ & $e_1e_3=e_4$ & $e_1e_4=e_6$ \\\hline
$^4\mathbf{EAC}_{2}^{6}$ & $T^{3}(\epsilon_{23}^5)$ & $e_1e_2=e_3$  & $e_1e_3=e_6$ & $e_2e_3=e_5$ \\\hline
$^4\mathbf{EAC}_{3}^{6}$ & $T^{3}(\epsilon_{24}^6)4$ & $e_1e_2=e_3$ & $e_1e_3=e_6$ & $e_2e_4=e_6$\\\hline
$^4\mathbf{EAC}_{4}^{6}$ & $T^{2,2}(\epsilon_{34}^6)$ & $e_1e_2=e_5$ & $e_1e_3=e_6$ & $e_3e_4=e_6$\\\hline
\end{longtable}
}

In dimension $n=7,8$, we find the following list:

{\small

\begin{longtable}{|l|l|llll|}
\caption*{ }   \\

\hline \multicolumn{1}{|c|}{${\mathcal A}$} & \multicolumn{1}{c|}{ } & \multicolumn{4}{c|}{Multiplication table} \\ \hline 
\endfirsthead

 \multicolumn{6}{l}%
{{\bfseries  continued from previous page}} \\
\hline \multicolumn{1}{|c|}{${\mathcal A}$} & \multicolumn{1}{c|}{ } & \multicolumn{4}{c|}{Multiplication table} \\ \hline 
\endhead

\hline \multicolumn{6}{|r|}{{Continued on next page}} \\ \hline
\endfoot

\hline 
\endlastfoot
$^4\mathbf{EAC}_{1}^{n}$ & $T^{3,2}$ & $e_1e_2=e_{n-1}$ & $e_1e_3=e_4$ & $e_1e_4=e_n$ & \\\hline
$^4\mathbf{EAC}_{2}^{n}$ & $T^{3}(\epsilon_{23}^{n-1})$ & $e_1e_2=e_3$  & $e_1e_3=e_n$ & $e_2e_3=e_{n-1}$ &\\\hline
$^4\mathbf{EAC}_{3}^{n}$ & $T^{3}(\epsilon_{24}^{n})$ & $e_1e_2=e_3$ & $e_1e_3=e_n$ & $e_2e_4=e_n$ &\\ \hline
$^4\mathbf{EAC}_{4}^{n}$ & $T^{2,2}(\epsilon_{34}^{n})$ & $e_1e_2=e_{n-1}$ & $e_1e_3=e_{n}$ & $e_3e_4=e_n$ &\\ \hline
$^4\mathbf{EAC}_{5}^{n}$ & $T^{2,2,2}(\epsilon_{23}^{n})$ & $e_1e_2=e_{n-2}$ & $e_1e_3=e_{n-1}$ & $e_1e_4=e_{n}$ & $e_2e_3=e_n$ \\ \hline
\end{longtable}
}

Finally, in dimension $n\geq 9$, we find the following Engel anticommutative algebras of level four:

{\small

\begin{longtable}{|l|l|llll|}
\caption*{ }   \\

\hline \multicolumn{1}{|c|}{${\mathcal A}$} & \multicolumn{1}{c|}{ } & \multicolumn{4}{c|}{Multiplication table} \\ \hline 
\endfirsthead

 \multicolumn{6}{l}%
{{\bfseries  continued from previous page}} \\
\hline \multicolumn{1}{|c|}{${\mathcal A}$} & \multicolumn{1}{c|}{ } & \multicolumn{4}{c|}{Multiplication table} \\ \hline 
\endhead

\hline \multicolumn{6}{|r|}{{Continued on next page}} \\ \hline
\endfoot

\hline 
\endlastfoot
$^4\mathbf{EAC}_{1}^{n}$ & $\eta_4$ & $e_1e_2=e_9$ & $e_3e_4=e_9$ & $e_5e_6=e_9$ & $e_7e_8=e_9$ \\ \hline
$^4\mathbf{EAC}_{2}^{n}$ & $T^{3,2}$ & $e_1e_2=e_{n-1}$ & $e_1e_3=e_4$ & $e_1e_4=e_n$ & \\ \hline
$^4\mathbf{EAC}_{3}^{n}$ & $T^{2,2,2,2}$ & $e_1e_2=e_{n-3}$ & $e_1e_3=e_{n-2}$ & $e_1e_4=e_{n-1}$ & $e_1e_5=e_n$ \\ \hline
$^4\mathbf{EAC}_{4}^{n}$ & $T^{3}(\epsilon_{23}^{n-1})$  & $e_1e_2=e_3$  & $e_1e_3=e_n$ & $e_2e_3=e_{n-1}$ &\\ \hline
$^4\mathbf{EAC}_{5}^{n}$ & $T^{3}(\epsilon_{24}^{n})$ & $e_1e_2=e_3$ & $e_1e_3=e_n$ & $e_2e_4=e_n$ & \\ \hline
$^4\mathbf{EAC}_{6}^{n}$ & $T^{2,2}(\epsilon_{34}^{n})$ & $e_1e_2=e_{n-1}$ & $e_1e_3=e_{n}$ & $e_3e_4=e_n$ &\\ \hline
$^4\mathbf{EAC}_{7}^{n}$ & $ T^{2,2,2}(\epsilon_{23}^{n})$ & $e_1e_2=e_{n-2}$ & $e_1e_3=e_{n-1}$ & $e_1e_4=e_{n}$ & $e_2e_3=e_n$  \\ \hline
\end{longtable}
}

Note that every Engel anticommutative algebra of level four is a Lie algebra.

\subsubsection{Engel anticommutative algebras of level five}

There are no Engel anticommutative algebras of level five and dimension lower than $5$.
In dimension $5$, there exist two Engel anticommutative algebras of level five, up to isomorphism:

{\small
\vspace{-5mm}
\begin{longtable}{|l|l|llll|}
\caption*{ }   \\

\hline \multicolumn{1}{|c|}{${\mathcal A}$} & \multicolumn{1}{c|}{ } & \multicolumn{4}{c|}{Multiplication table} \\ \hline 
\endfirsthead

 \multicolumn{6}{l}%
{{\bfseries  continued from previous page}} \\
\hline \multicolumn{1}{|c|}{${\mathcal A}$} & \multicolumn{1}{c|}{ } & \multicolumn{4}{c|}{Multiplication table} \\ \hline 
\endhead

\hline \multicolumn{6}{|r|}{{Continued on next page}} \\ \hline
\endfoot

\hline 
\endlastfoot

$^5\mathbf{EAC}_{1}^{5}$ & $T^{3}(\epsilon_{34}^{5})$ & $e_1e_2=e_3$  & $e_1e_3=e_5$ & $e_3e_4=e_5$ & \\ \hline
$^5\mathbf{EAC}_{2}^{5}$ & $T^{4}(\epsilon_{23}^{5})$ & $e_1e_2=e_3$ & $e_1e_3=e_4$ & $e_1e_4=e_5$ & $e_2e_3=e_5$ \\\hline
\end{longtable}
}

In dimension $6$, every Engel anticommutative  algebra of level five is isomorphic to one of the following algebras:

{\small
\vspace{-7mm}
\begin{longtable}{|l|l|llll|}
\caption*{ }   \\

\hline \multicolumn{1}{|c|}{${\mathcal A}$} & \multicolumn{1}{c|}{ } & \multicolumn{4}{c|}{Multiplication table} \\ \hline 
\endfirsthead

 \multicolumn{6}{l}%
{{\bfseries  continued from previous page}} \\
\hline \multicolumn{1}{|c|}{${\mathcal A}$} & \multicolumn{1}{c|}{ } & \multicolumn{4}{c|}{Multiplication table} \\ \hline 
\endhead

\hline \multicolumn{6}{|r|}{{Continued on next page}} \\ \hline
\endfoot

\hline 
\endlastfoot

$^5\mathbf{EAC}_{1}^{6}$ & $T^{4}$ & $e_1e_2=e_{3}$ & $e_1e_3=e_4$ & $e_1e_4=e_6$ & \\ \hline
$^5\mathbf{EAC}_{2}^{6}$ & $T^{3}(\epsilon_{34}^{6})$ & $e_1e_2=e_3$  & $e_1e_3=e_6$ & $e_3e_4=e_{6}$ &\\ \hline
$^5\mathbf{EAC}_{3}^{6}$ & $T^{3}(\epsilon_{45}^{6})$ & $e_1e_2=e_3$ & $e_1e_3=e_6$ & $e_4e_5=e_6$&\\ \hline
$^5\mathbf{EAC}_{4}^{6}$ & $T^{3,2}(\epsilon_{23}^{6})$ & $e_1e_2=e_{5}$ & $e_1e_3=e_{4}$ & $e_1e_4=e_6$ & $e_2e_3=e_6$ \\ \hline
\end{longtable}
}

In dimension $7$, we find:

{\small
\vspace{-9mm}
\begin{longtable}{|l|l|llll|}
\caption*{ }   \\

\hline \multicolumn{1}{|c|}{${\mathcal A}$} & \multicolumn{1}{c|}{ } & \multicolumn{4}{c|}{Multiplication table} \\ \hline 
\endfirsthead

 \multicolumn{6}{l}%
{{\bfseries  continued from previous page}} \\
\hline \multicolumn{1}{|c|}{${\mathcal A}$} & \multicolumn{1}{c|}{ } & \multicolumn{4}{c|}{Multiplication table} \\ \hline 
\endhead

\hline \multicolumn{6}{|r|}{{Continued on next page}} \\ \hline
\endfoot

\hline 
\endlastfoot
$^5\mathbf{EAC}_{1}^{7}$ & $ T^{4}$ & $e_1e_2=e_{3}$ & $e_1e_3=e_4$ & $e_1e_4=e_7$ & \\\hline
$^5\mathbf{EAC}_{2}^{7}$ & $T^{3,3}$ & $e_1e_2=e_{3}$ & $e_1e_3=e_{4}$ & $e_1e_5=e_{6}$ & $e_1e_6=e_7$ \\ \hline
$^5\mathbf{EAC}_{3}^{7}$ & $T^{3}(\epsilon_{34}^{7})$ & $e_1e_2=e_3$  & $e_1e_3=e_7$ & $e_3e_4=e_{7}$ & \\ \hline
$^5\mathbf{EAC}_{4}^{7}$ & $T^{2,2}(\epsilon_{45}^{7})$ & $e_1e_2=e_6$ & $e_1e_3=e_7$ & $e_4e_5=e_7$ &\\ \hline
$^5\mathbf{EAC}_{5}^{7}$ & $T^{3,2}(\epsilon_{23}^{7})$ & $e_1e_2=e_{6}$ & $e_1e_3=e_{4}$ & $e_1e_4=e_7$ & $e_2e_3=e_7$ \\ \hline
$^5\mathbf{EAC}_{6}^{7}$ & $T^{2,2,2}(\epsilon_{24}^{7})$ & $e_1e_2=e_{5}$ & $e_1e_3=e_{6}$ & $e_1e_4=e_{7}$ & $e_2e_4=e_7$ \\ \hline
$^5\mathbf{EAC}_{7}^{7}$ & $T^{2,2,2}(\epsilon_{23}^{4}-\epsilon_{26}^{7}+\epsilon_{35}^{7})$  & $e_1e_2=e_{5}$ & $e_1e_3=e_{6}$ & $e_1e_4=e_{7}$ & $e_2e_3=e_4$ \\&& $e_2e_6=-e_7$ & $e_3e_5=e_7$ &&\\ \hline
\end{longtable}
}

The  Engel anticommutative algebras of level five and dimension $n$, with $n=8,9,10$, are the following, up to isomorphism:

{\small

\begin{longtable}{|l|l|llll|}
\caption*{ }   \\

\hline \multicolumn{1}{|c|}{${\mathcal A}$} & \multicolumn{1}{c|}{ } & \multicolumn{4}{c|}{Multiplication table} \\ \hline 
\endfirsthead

 \multicolumn{6}{l}%
{{\bfseries  continued from previous page}} \\
\hline \multicolumn{1}{|c|}{${\mathcal A}$} & \multicolumn{1}{c|}{ } & \multicolumn{4}{c|}{Multiplication table} \\ \hline 
\endhead

\hline \multicolumn{6}{|r|}{{Continued on next page}} \\ \hline
\endfoot

\hline 
\endlastfoot
$^5\mathbf{EAC}_{1}^{n}$ & $T^{4}$ & $e_1e_2=e_{3}$ & $e_1e_3=e_4$ & $e_1e_4=e_n$ & \\ \hline
$^5\mathbf{EAC}_{2}^{n}$ & $T^{3,2,2}$  & $e_1e_2=e_{n-2}$ & $e_1e_3=e_{n-1}$ & $e_1e_4=e_{5}$ & $e_1e_5=e_n$ \\ \hline
$^5\mathbf{EAC}_{3}^{n}$ & $T^{3}(\epsilon_{34}^{n})$ & $e_1e_2=e_3$  & $e_1e_3=e_n$ & $e_3e_4=e_{n}$ &\\ \hline
$^5\mathbf{EAC}_{4}^{n}$ & $T^{2,2}(\epsilon_{45}^{n})$ & $e_1e_2=e_{n-1}$ & $e_1e_3=e_n$ & $e_4e_5=e_n$&\\ \hline
$^5\mathbf{EAC}_{5}^{n}$ & $T^{3,2}(\epsilon_{23}^{n})$ & $e_1e_2=e_{n-1}$ & $e_1e_3=e_{4}$ & $e_1e_4=e_n$ & $e_2e_3=e_n$ \\ \hline
$^5\mathbf{EAC}_{6}^{n}$ & $T^{2,2,2}(\epsilon_{24}^{n})$ & $e_1e_2=e_{n-2}$ & $e_1e_3=e_{n-1}$ & $e_1e_4=e_{n}$ & $e_2e_4=e_n$ \\ \hline
\end{longtable}
}

Finally, in dimension $n\geq 11$, there exist, up to isomorphism, the following Engel anticommutative  algebras:

{\small

\begin{longtable}{|l|l|llll|}
\caption*{ }   \\

\hline \multicolumn{1}{|c|}{${\mathcal A}$} & \multicolumn{1}{c|}{ } & \multicolumn{4}{c|}{Multiplication table} \\ \hline 
\endfirsthead

 \multicolumn{6}{l}%
{{\bfseries  continued from previous page}} \\
\hline \multicolumn{1}{|c|}{${\mathcal A}$} & \multicolumn{1}{c|}{ } & \multicolumn{4}{c|}{Multiplication table} \\ \hline 
\endhead

\hline \multicolumn{6}{|r|}{{Continued on next page}} \\ \hline
\endfoot

\hline 
\endlastfoot

$^5\mathbf{EAC}_{1}^{n}$ & $\eta_5$ & $e_1e_2=e_{11}$ & $e_3e_4=e_{11}$ & $e_5e_6=e_{11}$ & $e_7e_8=e_{11}$ \\&& $e_9e_{10}=e_{11}$ &&&\\ \hline
$^5\mathbf{EAC}_{2}^{n}$ & $T^{4}$ & $e_1e_2=e_{3}$ & $e_1e_3=e_4$ & $e_1e_4=e_n$ & \\\hline
$^5\mathbf{EAC}_{3}^{n}$ & $T^{3,2,2}$ & $e_1e_2=e_{n-2}$ & $e_1e_3=e_{n-1}$ & $e_1e_4=e_{5}$ & $e_1e_5=e_n$  \\ \hline
$^5\mathbf{EAC}_{4}^{n}$ & $T^{2,2,2,2,2}$ & $e_1e_2=e_{n-4}$ & $e_1e_3=e_{n-3}$ & $e_1e_4=e_{n-2}$ & $e_1e_5=e_{n-1}$ \\&& $e_1e_6=e_{n}$ &&&  \\ \hline
$^5\mathbf{EAC}_{5}^{n}$ & $T^{3}(\epsilon_{34}^{n}) $  & $e_1e_2=e_3$  & $e_1e_3=e_n$ & $e_3e_4=e_{n}$ & \\ \hline
$^5\mathbf{EAC}_{6}^{n}$ & $T^{2,2}(\epsilon_{45}^{n})$ & $e_1e_2=e_{n-1}$ & $e_1e_3=e_n$ & $e_4e_5=e_n$ & \\ \hline
$^5\mathbf{EAC}_{7}^{n}$ & $T^{3,2}(\epsilon_{23}^{n})$ & $e_1e_2=e_{n-1}$ & $e_1e_3=e_{4}$ & $e_1e_4=e_n$ & $e_2e_3=e_n$ \\ \hline
$^5\mathbf{EAC}_{8}^{n}$ & $T^{2,2,2}(\epsilon_{24}^{n})$ & $e_1e_2=e_{n-2}$ & $e_1e_3=e_{n-1}$ & $e_1e_4=e_{n}$ & $e_2e_4=e_n$  \\ \hline
\end{longtable}
}

As it is pointed out in~\cite{wolf1}, any Engel anticommutative algebra of level at most five is a Lie algebra, except for $^5\mathbf{EAC}_{2}^{6}$, $^5\mathbf{EAC}_{3}^{7}$, $^5\mathbf{EAC}_{3}^{n}$ for $n\in\{8,9,10\}$, $^5\mathbf{EAC}_{5}^{n}$ for $n\geq 11$ (all of them grouped under the name $T^{3}(\epsilon_{34}^{n})$ in~\cite{wolf1}) and $^5\mathbf{EAC}_{7}^{7}$, which are Malcev. Also, any  Engel anticommutative algebra of level at most five is nilpotent.

\subsection{Lie algebras} 
For results about Lie and Engel Lie algebras, we refer the reader to Subsection~\ref{ss:la}. 

\subsection{Malcev algebras} 
For results about Malcev and Engel Malcev algebras, we refer the reader to Subsection~\ref{ss:la}.

\subsection{Jordan algebras} 
The classification of  Jordan algebras of level two was given in~\cite{khud15}. In the same article, the authors selected the Jordan algebras of level one from the classification of~\cite{khud13}. Also, the classification of Jordan algebras of level two can be seen as an easy corollary of the general results of~\cite{kv17}. In this survey, we will refer to the notation of the original work~\cite{khud15}.

\subsubsection{Jordan algebras of level one}
Every $n$-dimensional Jordan algebra of level one is isomorphic to one of the two following algebras:

{\small

\begin{longtable}{|l|l|lll|}
\caption*{ }   \\

\hline \multicolumn{1}{|c|}{${\mathcal A}$} & \multicolumn{1}{c|}{ } & \multicolumn{3}{c|}{Multiplication table} \\ \hline 
\endfirsthead

 \multicolumn{4}{l}%
{{\bfseries  continued from previous page}} \\
\hline \multicolumn{1}{|c|}{${\mathcal A}$} & \multicolumn{1}{c|}{ } & \multicolumn{3}{c|}{Multiplication table} \\ \hline 
\endhead

\hline \multicolumn{5}{|r|}{{Continued on next page}} \\ \hline
\endfoot

\hline 
\endlastfoot
$^1\mathbf{J}_{1}^{n}$ & $\lambda_2\oplus a_{n-2}$  & $e_1e_1=e_2$ & &\\ \hline
$^1\mathbf{J}_{2}^{n}(\frac{1}{2})$ & $\nu_n(\frac{1}{2})$ & $e_1e_1=e_1$ & $e_1e_i=\frac{1}{2} e_i,$ &\mbox{for }$2\leq i \leq n$  \\ \hline
\end{longtable}
} 
\subsubsection{Jordan algebras of level two}
Up to isomorphism, there exist two Jordan algebras of level and dimension $2$, namely:

{\small
\vspace{-5mm}
\begin{longtable}{|l|l|ll|}
\caption*{ }   \\

\hline \multicolumn{1}{|c|}{${\mathcal A}$} & \multicolumn{1}{c|}{ } & \multicolumn{2}{c|}{Multiplication table} \\ \hline 
\endfirsthead

 \multicolumn{4}{l}%
{{\bfseries  continued from previous page}} \\
\hline \multicolumn{1}{|c|}{${\mathcal A}$} & \multicolumn{1}{c|}{ } & \multicolumn{2}{c|}{Multiplication table} \\ \hline 
\endhead

\hline \multicolumn{4}{|r|}{{Continued on next page}} \\ \hline
\endfoot

\hline 
\endlastfoot
$^2\mathbf{J}_{1}^{2}$ & $J_1$ & $e_1e_1=e_1 $ & \\ \hline
$^2\mathbf{J}_{2}^{2}$ & $J_2$ & $e_1e_1=e_1$ & $e_1e_2= e_2$ \\ \hline
\end{longtable}
}

In dimension $n\geq 3$, we find the following list:

{\small

\begin{longtable}{|l|l|l|}
\caption*{ }   \\

\hline \multicolumn{1}{|c|}{${\mathcal A}$} & \multicolumn{1}{c|}{ } & \multicolumn{1}{c|}{Multiplication table} \\ \hline 
\endfirsthead

 \multicolumn{3}{l}%
{{\bfseries  continued from previous page}} \\
\hline \multicolumn{1}{|c|}{${\mathcal A}$} & \multicolumn{1}{c|}{ } & \multicolumn{1}{c|}{Multiplication table} \\ \hline 
\endhead

\hline \multicolumn{3}{|r|}{{Continued on next page}} \\ \hline
\endfoot

\hline 
\endlastfoot

$^2\mathbf{J}_{1}^{n}$ & $J_1$ & $e_1e_1=e_1$ \\ \hline
$^2\mathbf{J}_{2}^{n}$ & $J_2$ &  $e_1e_i= e_i$ \\ \hline
$^2\mathbf{J}_{3}^{n}$ & $J_3$ & $e_1e_2=e_3,$  \\ \hline 
\end{longtable}
}\noindent for $1\leq i \leq n$. Note that $^2\mathbf{J}_{3}^{n}$ and $\mathbf{F}^{1,1}\oplus\mathbb{C}^{n-3}$ (with the notation of~\cite{kv17}) are isomorphic, for $n\geq 3$.

\subsection{Left-alternative algebras}\label{ss:lla}

We select the left-alternative algebras of level one from the general classification of~\cite{khud13}. For level two, consult~\cite{kv17}.

\subsubsection{Left-alternative algebras of level one}
Up to isomorphism, the $2$-dimensional left-alternative algebras of level one are

{\small

\begin{longtable}{|l|l|ll|}
\caption*{ }   \\

\hline \multicolumn{1}{|c|}{${\mathcal A}$} & \multicolumn{1}{c|}{ } & \multicolumn{2}{c|}{Multiplication table} \\ \hline 
\endfirsthead

 \multicolumn{4}{l}%
{{\bfseries  continued from previous page}} \\
\hline \multicolumn{1}{|c|}{${\mathcal A}$} & \multicolumn{1}{c|}{ } & \multicolumn{2}{c|}{Multiplication table} \\ \hline 
\endhead

\hline \multicolumn{4}{|r|}{{Continued on next page}} \\ \hline
\endfoot

\hline 
\endlastfoot
$^1\mathbf{LA}_{1}^{2}$ & $\lambda_2$ & $e_1e_1=e_2$ & \\ \hline
$^1\mathbf{LA}_{2}^{2}4$ & $\nu_2(0)$ & $e_1e_1=e_1$ & $e_2e_1=e_2$  \\ \hline
$^1\mathbf{LA}_{3}^{2}$ & $\nu_2(1)$ & $e_1e_1=e_1$ & $e_1e_2=e_2$  \\ \hline
\end{longtable}
}

In dimension $n\geq 3$, the classification is the following:

{\small

\begin{longtable}{|l|l|ll|}
\caption*{ }   \\

\hline \multicolumn{1}{|c|}{${\mathcal A}$} & \multicolumn{1}{c|}{ } & \multicolumn{2}{c|}{Multiplication table} \\ \hline 
\endfirsthead

 \multicolumn{4}{l}%
{{\bfseries  continued from previous page}} \\
\hline \multicolumn{1}{|c|}{${\mathcal A}$} & \multicolumn{1}{c|}{ } & \multicolumn{2}{c|}{Multiplication table} \\ \hline 
\endhead

\hline \multicolumn{4}{|r|}{{Continued on next page}} \\ \hline
\endfoot

\hline 
\endlastfoot

$^1\mathbf{LA}_{1}^{n}$ & $\lambda_2\oplus a_{n-2}$ & $e_1e_1=e_2$ & \\ \hline
$^1\mathbf{LA}_{2}^{n}$ & $n_3^{-}\oplus a_{n-3}$ & $e_1e_2=e_3$ & $e_2e_1=-e_3$  \\ \hline
$^1\mathbf{LA}_{3}^{n}$ & $\nu_n(0)$ & $e_ie_1=e_i$ &\\ \hline
$^1\mathbf{LA}_{4}^{n}$ & $\nu_n(1)$ & $e_1e_i=e_i,$ & \\ \hline
\end{longtable}
}\noindent for $1\leq i \leq n$.

Note that all these algebras are associative.

\subsubsection{Left-alternative algebras of level two}
In dimension $2$, there exist only two  left-alternative algebras of level two, up to isomorphism:

{\small
\vspace{-5mm}
\begin{longtable}{|l|l|ll|}
\caption*{ }   \\

\hline \multicolumn{1}{|c|}{${\mathcal A}$} & \multicolumn{1}{c|}{ } & \multicolumn{2}{c|}{Multiplication table} \\ \hline 
\endfirsthead

 \multicolumn{4}{l}%
{{\bfseries  continued from previous page}} \\
\hline \multicolumn{1}{|c|}{${\mathcal A}$} & \multicolumn{1}{c|}{ } & \multicolumn{2}{c|}{Multiplication table} \\ \hline 
\endhead

\hline \multicolumn{4}{|r|}{{Continued on next page}} \\ \hline
\endfoot

\hline 
\endlastfoot
$^2\mathbf{LA}_{1}^{2}$ & $\mathbf{D}_2^{0,0}$  & $e_1e_1=e_1$ & \\ \hline
$^2\mathbf{LA}_{2}^{2}$ & $\mathbf{D}_2^{1,1}$ & $e_1e_2=e_2$ & $e_2e_1=e_2$ \\ \hline 
\end{longtable}
}

In dimension $3$, the classification is the following:

{\small
\vspace{-5mm}
\begin{longtable}{|l|l|lll|}
\caption*{ }   \\

\hline \multicolumn{1}{|c|}{${\mathcal A}$} & \multicolumn{1}{c|}{ } & \multicolumn{3}{c|}{Multiplication table} \\ \hline 
\endfirsthead

 \multicolumn{5}{l}%
{{\bfseries  continued from previous page}} \\
\hline \multicolumn{1}{|c|}{${\mathcal A}$} & \multicolumn{1}{c|}{ } & \multicolumn{3}{c|}{Multiplication table} \\ \hline 
\endhead

\hline \multicolumn{5}{|r|}{{Continued on next page}} \\ \hline
\endfoot

\hline 
\endlastfoot
$^2\mathbf{LA}_{1}^{3}$ & $\mathbb{C}\rtimes_{0}^t\mathbf{D}_2^{0,0}$ & $e_1e_1=e_1$ && \\ \hline
$^2\mathbf{LA}_{2}^{3}$  & $\mathbb{C}\rtimes_{0}^t\mathbf{D}_2^{1,1}$ & $e_1e_1=e_1$ & $e_1e_2= e_2$  & $e_1e_3= e_3$ \\&& $e_2e_1=e_2$  & $e_3e_1=e_3$ & \\ \hline
$^2\mathbf{LA}_{3}^{3}(\alpha,\beta),\ (\alpha,\beta)\in K_{2}^{\star}$ & $\mathbf{F}^{\alpha,\beta}$ & $e_1e_1=e_3$ & $e_1e_2=\alpha e_3$ & $e_2e_1=\beta e_3 $ \\ \hline
$^2\mathbf{LA}_{4}^{3}$ & $\mathbf{T}_1^{2,\overline{1,0}}$  & $e_1e_1=e_1$ & $e_1e_2=e_2+e_3$ & $e_2e_1=-e_3 $ \\&& $e_3e_1=e_3$ && \\ \hline
$^2\mathbf{LA}_{5}^{3}$ & $\mathbf{T}_1^{2,\overline{0,1}}$ & $e_1e_1=e_1$ & $e_1e_2=e_3$ & $e_1e_3=e_3$ \\&& $e_2e_1=e_2-e_3$ && \\ \hline
\end{longtable}
}

In dimension $4$, we find:

{\small
\vspace{-5mm}
\begin{longtable}{|l|l|lll|}
\caption*{ }   \\

\hline \multicolumn{1}{|c|}{${\mathcal A}$} & \multicolumn{1}{c|}{ } & \multicolumn{3}{c|}{Multiplication table} \\ \hline 
\endfirsthead

 \multicolumn{5}{l}%
{{\bfseries  continued from previous page}} \\
\hline \multicolumn{1}{|c|}{${\mathcal A}$} & \multicolumn{1}{c|}{ } & \multicolumn{3}{c|}{Multiplication table} \\ \hline 
\endhead

\hline \multicolumn{5}{|r|}{{Continued on next page}} \\ \hline
\endfoot

\hline 
\endlastfoot
$^2\mathbf{LA}_{1}^{4}$ & $\mathbb{C}\rtimes_{0}^t\mathbf{D}_2^{0,0}$  & $e_1e_1=e_1$ && \\ \hline
$^2\mathbf{LA}_{2}^{4}$ & $\mathbb{C}\rtimes_{0}^t\mathbf{D}_2^{1,1}$ & $e_1e_1=e_1$ & $e_1e_2= e_2$  & $e_1e_3= e_3$ \\&& $e_2e_1= e_2$  & $e_3e_1= e_3$ & \\ \hline
$^2\mathbf{LA}_{3}^{4}(\alpha,\beta),\ (\alpha,\beta)\in K_{2}^{\star}$ & $\mathbf{F}^{\alpha,\beta}\oplus \mathbb{C}$  & $e_1e_1=e_3$ & $e_1e_2=\alpha e_3$ & $e_2e_1=\beta e_3$  \\ \hline
$^2\mathbf{LA}_{4}^{4}$ & $\mathbf{T}_1^{2,\overline{1,0}}$ & $e_1e_1=e_1$ & $e_1e_2= e_2+e_3$ & $e_2e_1=-e_3$ \\&& $e_3e_1=e_3$&&  \\ \hline
$^2\mathbf{LA}_{5}^{4}$ & $\mathbf{T}_1^{2,\overline{0,1}}$ & $e_1e_1=e_1$ & $e_1e_2=e_3$ & $e_1e_3=e_3$ \\&& $e_2e_1=e_2-e_3$ && \\ \hline
\end{longtable}
}\clearpage

Finally, the classification of  $n$-dimensional left-alternative algebras of level two, for $n\geq 5$, is the following:

{\small
\begin{longtable}{|l|l|lll|}
\caption*{ }   \\

\hline \multicolumn{1}{|c|}{${\mathcal A}$} & \multicolumn{1}{c|}{ } & \multicolumn{3}{c|}{Multiplication table} \\ \hline 
\endfirsthead

 \multicolumn{5}{l}%
{{\bfseries  continued from previous page}} \\
\hline \multicolumn{1}{|c|}{${\mathcal A}$} & \multicolumn{1}{c|}{ } & \multicolumn{3}{c|}{Multiplication table} \\ \hline 
\endhead

\hline \multicolumn{5}{|r|}{{Continued on next page}} \\ \hline
\endfoot

\hline 
\endlastfoot
$^2\mathbf{LA}_{1}^{n}$ & $T_0^{2,2}$ & $e_1e_2=e_3$ & $e_1e_4=e_5$ & $e_2e_1=-e_3$ \\&& $e_4e_1=-e_5$ && \\ \hline
$^2\mathbf{LA}_{2}^{n}$ & $\eta_2\oplus \mathbb{C}^{n-5}$ & $e_1e_2=e_5$  & $e_2e_1=-e_5$ & $e_3e_4=e_5$ \\&& $e_4e_3=-e_5$ && \\ \hline
$^2\mathbf{LA}_{3}^{n}$ & $\mathbb{C}^{n-2}\rtimes_{0}^t\mathbf{D}_2^{0,0}$ & $e_1e_1=e_1$ && \\ \hline
$^2\mathbf{LA}_{4}^{n}$ & $\mathbb{C}^{n-2}\rtimes_{0}^t\mathbf{D}_2^{1,1}$ & $e_1e_1=e_1$ & $e_1e_i= e_i$  & $e_ie_1=e_i$ \\ \hline
$^2\mathbf{LA}_{5}^{n}(\alpha,\beta),\ (\alpha,\beta)\in K_{2}^{\star}$ & $\mathbf{F}^{\alpha,\beta}\oplus \mathbb{C}^{n-3}$ & $e_1e_1=e_3$ & $e_1e_2=\alpha e_3$ & $e_2e_1=\beta e_3$ \\ \hline
$^2\mathbf{LA}_{6}^{n}$ & $T_1^{2,\overline{1,0}}$ & $e_1e_1=e_1$ & $e_1e_2= e_2+e_3$  & $e_2e_1=-e_3$  \\&& $e_je_1=e_j$  && \\ \hline
$^2\mathbf{LA}_{7}^{n}$ & $T_1^{2,\overline{0,1}}$ & $e_1e_1=e_1$ & $e_1e_2=e_3$ & $e_1e_j=e_j$ \\&& $e_2e_1=e_2-e_3,$ && \\ \hline
\end{longtable}
}\noindent for $2\leq i\leq n$ and $3\leq j\leq n$.

Note that all these algebras are associative.

\subsection{Associative algebras} 
We refer the reader to Subsection~\ref{ss:lla}, as the associative and left-alternative algebras of levels one and two coincide.

\subsection{Leibniz algebras} 
In this section, we will deal with Leibniz algebras of levels one and two. The algebras are selected from the general classifications of~\cite{khud13} and~\cite{kv17}, respectively.

\subsubsection{Leibniz algebras of level one}
There exist two Leibniz algebras of level one and dimension $2$: one of them is Lie, and the other one is not.

{\small
\begin{longtable}{|l|l|ll|}
\caption*{ }   \\

\hline \multicolumn{1}{|c|}{${\mathcal A}$} & \multicolumn{1}{c|}{ } & \multicolumn{2}{c|}{Multiplication table} \\ \hline 
\endfirsthead

 \multicolumn{4}{l}%
{{\bfseries  continued from previous page}} \\
\hline \multicolumn{1}{|c|}{${\mathcal A}$} & \multicolumn{1}{c|}{ } & \multicolumn{2}{c|}{Multiplication table} \\ \hline 
\endhead

\hline \multicolumn{4}{|r|}{{Continued on next page}} \\ \hline
\endfoot

\hline 
\endlastfoot
$^1\mathbf{L}_{1}^{2}$ & $p_2^{-}$ & $e_1e_2=e_2$ & $e_2e_1=-e_2$  \\ \hline
$^1\mathbf{L}_{2}^{2}$ & $\lambda_2$ & $e_1e_1=e_2$ & \\ \hline
\end{longtable}
}

In dimension $n\geq 3$, the classification is the following:

{\small
\begin{longtable}{|l|l|ll|}
\caption*{ }   \\

\hline \multicolumn{1}{|c|}{${\mathcal A}$} & \multicolumn{1}{c|}{ } & \multicolumn{2}{c|}{Multiplication table} \\ \hline 
\endfirsthead

 \multicolumn{4}{l}%
{{\bfseries  continued from previous page}} \\
\hline \multicolumn{1}{|c|}{${\mathcal A}$} & \multicolumn{1}{c|}{ } & \multicolumn{2}{c|}{Multiplication table} \\ \hline 
\endhead

\hline \multicolumn{4}{|r|}{{Continued on next page}} \\ \hline
\endfoot

\hline 
\endlastfoot
$^1\mathbf{L}_{1}^{n}$ & $p_n^{-}$ & $e_1e_i=e_i$ & $e_ie_1=-e_i$ \\ \hline
$^1\mathbf{L}_{2}^{n}$ & $\lambda_2\oplus a_{n-2}$  & $e_1e_1=e_2$ & \\ \hline
$^1\mathbf{L}_{3}^{n}$ & $n_3^{-}\oplus a_{n-3}$ & $e_1e_2=e_3$ & $e_2e_1=-e_3,$  \\ \hline
\end{longtable}
}\noindent for $2\leq i \leq n$.

\subsubsection{Leibniz algebras of level two}
In this section, we correct some inaccuracies of~\cite{fkrv}.

Up to isomorphism, there exists   one Leibniz non-Lie algebra of level two and dimension $2$:

{\small
\begin{longtable}{|l|l|l|}
\caption*{ }   \\

\hline \multicolumn{1}{|c|}{${\mathcal A}$} & \multicolumn{1}{c|}{ } & \multicolumn{1}{c|}{Multiplication table} \\ \hline 
\endfirsthead

 \multicolumn{3}{l}%
{{\bfseries  continued from previous page}} \\
\hline \multicolumn{1}{|c|}{${\mathcal A}$} & \multicolumn{1}{c|}{ } & \multicolumn{1}{c|}{Multiplication table} \\ \hline 
\endhead

\hline \multicolumn{3}{|r|}{{Continued on next page}} \\ \hline
\endfoot

\hline 
\endlastfoot
$^2\mathbf{L}_{1}^{2}$ & $\mathbf{B}_2^{0}$ & $e_2e_1=e_2$  \\ \hline
\end{longtable}
}

In dimension $3$, the classification is:

{\small
\begin{longtable}{|l|l|lll|}
\caption*{ }   \\

\hline \multicolumn{1}{|c|}{${\mathcal A}$} & \multicolumn{1}{c|}{ } & \multicolumn{3}{c|}{Multiplication table} \\ \hline 
\endfirsthead

 \multicolumn{5}{l}%
{{\bfseries  continued from previous page}} \\
\hline \multicolumn{1}{|c|}{${\mathcal A}$} & \multicolumn{1}{c|}{ } & \multicolumn{3}{c|}{Multiplication table} \\ \hline 
\endhead

\hline \multicolumn{5}{|r|}{{Continued on next page}} \\ \hline
\endfoot

\hline 
\endlastfoot
$^2\mathbf{L}_{1}^{3}$ & $\mathbb{C}\rtimes_{0}^t\mathbf{B}_2^{0}$ & $e_2e_1=e_2$ &  $e_3e_1=e_3$ & \\ \hline
$^2\mathbf{L}_{2}^{3}(\alpha,\beta),\ (\alpha,\beta)\in K_{2}^{\star}$ & $\mathbf{F}^{\alpha,\beta}$ & $e_1e_1=e_3$ & $e_1e_2=\alpha e_3$ & $e_2e_1=\beta e_3$ \\ \hline
$^2\mathbf{L}_{3}^{3}(\alpha,\beta),\ (\alpha,\beta)\in K_{2}^{\star}$ & $\mathbf{T}_0^{2,\overline{\alpha,\beta}}$  & $e_1e_2=\alpha e_2+e_3$ & $e_1e_3=\beta e_3$ & $e_2e_1=-\alpha e_2-e_3$ \\&& $e_3e_1=-\beta e_3$ && \\ \hline
\end{longtable}
}

The $4$-dimensional Leibniz algebras of level two are, up to isomorphism, the following ones:

{\small
\begin{longtable}{|l|l|lll|}
\caption*{ }   \\

\hline \multicolumn{1}{|c|}{${\mathcal A}$} & \multicolumn{1}{c|}{ } & \multicolumn{3}{c|}{Multiplication table} \\ \hline 
\endfirsthead

 \multicolumn{5}{l}%
{{\bfseries  continued from previous page}} \\
\hline \multicolumn{1}{|c|}{${\mathcal A}$} & \multicolumn{1}{c|}{ } & \multicolumn{3}{c|}{Multiplication table} \\ \hline 
\endhead

\hline \multicolumn{5}{|r|}{{Continued on next page}} \\ \hline
\endfoot

\hline 
\endlastfoot
$^2\mathbf{L}_{1}^{4}$ & $T_0^3$ & $e_1e_2=e_3$ & $e_1e_3=e_4$ & $e_2e_1=-e_3$ \\&& $e_3e_1=-e_4$ && \\ \hline
$^2\mathbf{L}_{2}^{4}$ & $\mathbb{C}^2\rtimes_{0}^t\mathbf{B}_2^{0}$ & $e_2e_1=e_2$ &  $e_3e_1=e_3$ & $e_4e_1= e_4$ \\ \hline
$^2\mathbf{L}_{3}^{4}(\alpha,\beta),\ (\alpha,\beta)\in K_{2}^{\star}$ & $\mathbf{F}^{\alpha,\beta}\oplus \mathbb{C}$ & $e_1e_1=e_3$ & $e_1e_2=\alpha e_3$ & $e_2e_1=\beta e_3$  \\ \hline
$^2\mathbf{L}_{4}^{4}(\alpha,\beta),\ (\alpha,\beta)\in K_{1,1}^{\star}$ & $T_0^{2,\overline{\alpha,\beta}}$ & $e_1e_2=\alpha e_2+e_3$ & $e_1e_3=\beta e_3$ & $e_1e_4=\beta e_4$ \\&& $e_2e_1=-\alpha e_2-e_3$ & $e_3e_1=-\beta e_3$ & $e_4e_1=-\beta e_4$ \\ \hline
\end{longtable}
}

Finally, in dimension $n\geq 5$, the classification of Leibniz algebras of level two up to isomorphism is as follows:

{\small
\begin{longtable}{|l|l|lll|}
\caption*{ }   \\

\hline \multicolumn{1}{|c|}{${\mathcal A}$} & \multicolumn{1}{c|}{ } & \multicolumn{3}{c|}{Multiplication table} \\ \hline 
\endfirsthead

 \multicolumn{5}{l}%
{{\bfseries  continued from previous page}} \\
\hline \multicolumn{1}{|c|}{${\mathcal A}$} & \multicolumn{1}{c|}{ } & \multicolumn{3}{c|}{Multiplication table} \\ \hline 
\endhead

\hline \multicolumn{5}{|r|}{{Continued on next page}} \\ \hline
\endfoot

\hline 
\endlastfoot
$^2\mathbf{L}_{1}^{n}$ & $T_0^{2,2}$ & $e_1e_2=e_3$ & $e_1e_4=e_5$ & $e_2e_1=-e_3$ \\&& $e_4e_1=-e_5$ && \\ \hline
$^2\mathbf{L}_{2}^{n}$ & $\eta_2\oplus \mathbb{C}^{n-5}$ & $e_1e_2=e_5$  & $e_2e_1=-e_5$ & $e_3e_4=e_5$ \\&& $e_4e_3=-e_5$ && \\ \hline
$^2\mathbf{L}_{3}^{n}$ & $\mathbb{C}^{n-2}\rtimes_{0}^t\mathbf{B}_2^{0}$ & $e_ie_1=e_i$ && \\ \hline
$^2\mathbf{L}_{4}^{n}(\alpha,\beta),\ (\alpha,\beta)\in K_{2}^{\star}$ & $\mathbf{F}^{\alpha,\beta}\oplus \mathbb{C}^{n-3}$ & $e_1e_1=e_3$ & $e_1e_2=\alpha e_3$ & $e_2e_1=\beta e_3$  \\ \hline
$^2\mathbf{L}_{5}^{n}(\alpha,\beta),\ (\alpha,\beta)\in K_{1,1}^{\star}$ & $T_0^{2,\overline{\alpha,\beta}}$  & $e_1e_2=\alpha e_2+e_3$ & $e_1e_j=\beta e_j$ &  $e_2e_1=-\alpha e_2-e_3$ \\&& $e_je_1=-\beta e_j,$ &&\\ \hline
\end{longtable}
}\noindent for $2\leq i\leq n$ and $3\leq j\leq n$.

\subsection{\texorpdfstring{$n$-ary algebras}{n-ary algebras}}

In~\cite{wolf2}, the author described all the $n$-ary algebras of level one. In particular, he gave an explicit classification for $n=2$, which coincides with the one in~\cite{khud13}, and for $n=3$, which we present below.
\vspace{+4mm}

\subsubsection{Ternary algebras of level one}
Up to isomorphism, there exist the following  ternary algebras of level one and dimension $2$:

{\small
\begin{longtable}{|l|l|ll|}
\caption*{ }   \\

\hline \multicolumn{1}{|c|}{${\mathcal A}$} & \multicolumn{1}{c|}{ } & \multicolumn{2}{c|}{Multiplication table} \\ \hline 
\endfirsthead

 \multicolumn{4}{l}%
{{\bfseries  continued from previous page}} \\
\hline \multicolumn{1}{|c|}{${\mathcal A}$} & \multicolumn{1}{c|}{ } & \multicolumn{2}{c|}{Multiplication table} \\ \hline 
\endhead

\hline \multicolumn{4}{|r|}{{Continued on next page}} \\ \hline
\endfoot

\hline 
\endlastfoot
$^1\mathbf{T}_{1}^{2}$ & $p=(3,0,\dots)$ & $[e_1,e_1,e_1]=e_2$ & \\ \hline
$^1\mathbf{T}_{2}^{2}(\epsilon,\beta_1,\beta_2,\beta_3)$ & $p=(2,0,\dots)$  
& $[e_1,e_1,e_1]=\epsilon e_1$ 
& $[e_1,e_1,e_2]=\beta_3e_2$ \\&
& $[e_1,e_2,e_1]=\beta_2e_2$ 
& $[e_2,e_1,e_1]=\beta_1e_2$ \\ \hline
$^1\mathbf{T}_{3}^{2}(\alpha_1,\alpha_2,\alpha_3)$ & $p=(1,1,0,\dots)$ 
&  $[e_1,e_1,e_2]=(\alpha_2-\alpha_3)e_1$ 
& $[e_1,e_2,e_1]=(\alpha_3-\alpha_1)e_1$ \\&
& $[e_2,e_1,e_1]=(\alpha_1-\alpha_2)e_1$  
& $[e_1,e_2,e_2]=(\alpha_2-\alpha_1)e_2$  \\&
&  $[e_2,e_1,e_2]=(\alpha_1-\alpha_3)e_2$ 
& $[e_2,e_2,e_1]=(\alpha_3-\alpha_2)e_2$\\ \hline
\end{longtable}
}\noindent 
\vspace{+2mm}

Here  $(\alpha_1 -\alpha_2,\alpha_2- \alpha_3, \alpha_3-\alpha_1) \neq (0,0,0)$ and the triple $(\alpha_1,\alpha_2,\alpha_3)$ is determined (by the isomorphism class of $^1\mathbf{T}_{3}^{2}(\alpha_1,\alpha_2,\alpha_3)$) up to multiplication by a nonzero element of $\mathbb{C}$ and addition of an element of $\mathbb{C}$. Furthermore, $\epsilon\in\{0,1\}$, $\beta_1+\beta_2+\beta_3=\epsilon$, and if $\epsilon=0$, then the triple $(\beta_1,\beta_2,\beta_3)\neq (0,0,0)$ is determined (by the isomorphism class of $^1\mathbf{T}_{2}^{2}(\epsilon,\beta_1,\beta_2,\beta_3)$) up to multiplication by a nonzero element of $\mathbb{C}$.

\

In dimension $3$, we find the following  ternary algebras, up to isomorphism:

{\small

\begin{longtable}{|l|l|ll|}
\caption*{ }   \\

\hline \multicolumn{1}{|c|}{${\mathcal A}$} & \multicolumn{1}{c|}{ } & \multicolumn{2}{c|}{Multiplication table} \\ \hline 
\endfirsthead

 \multicolumn{4}{l}%
{{\bfseries  continued from previous page}} \\
\hline \multicolumn{1}{|c|}{${\mathcal A}$} & \multicolumn{1}{c|}{ } & \multicolumn{2}{c|}{Multiplication table} \\ \hline 
\endhead

\hline 

$^1\mathbf{T}_{1}^{3}$ & $p=(3,0,\dots)$  & $[e_1,e_1,e_1]=e_2$ & \\ \hline
$^1\mathbf{T}_{2}^{3}(\epsilon,\beta_1,\beta_2,\beta_3)$ & $p=(2,0,\dots)$ 
& $[e_1,e_1,e_1]=\epsilon e_1 $ 
& $[e_1,e_1,e_2]=\beta_3e_2 $ \\&
& $[e_1,e_2,e_1]=\beta_2e_2 $ 
& $[e_2,e_1,e_1]=\beta_1e_2 $ \\&
& $[e_1,e_1,e_3]=\beta_3e_3 $
& $[e_1,e_3,e_1]=\beta_2e_3$ \\ & 
& $[e_3,e_1,e_1]=\beta_1e_3$ & \\ \hline
$^1\mathbf{T}_{3}^{3}(\alpha_1,\alpha_2,\alpha_3)$ & $p=(1,1,0,\dots)$ 
&  $[e_1,e_1,e_2]=(\alpha_2-\alpha_3)e_1 $ 
& $[e_1,e_2,e_1]=(\alpha_3-\alpha_1)e_1 $ \\&
& $[e_2,e_1,e_1]=(\alpha_1-\alpha_2)e_1 $ 
& $[e_1,e_2,e_2]=(\alpha_2-\alpha_1)e_2$ \\&
&  $[e_2,e_1,e_2]=(\alpha_1-\alpha_3)e_2 $
& $[e_2,e_2,e_1]=(\alpha_3-\alpha_2)e_2$ \\ & 
& $[e_1,e_2,e_3]=\alpha_3e_3 $
& $[e_1,e_3,e_2]=-\alpha_2e_3$ \\&
& $[e_2,e_1,e_3]=-\alpha_3e_3 $ 
& $[e_2,e_3,e_1]=\alpha_2e_3 $ \\&
& $[e_3,e_1,e_2]=\alpha_1e_3 $ 
& $[e_3,e_2,e_1]=-\alpha_1e_3$\\ \hline
$^1\mathbf{T}_{4}^{3}(\alpha_1,\alpha_2,\alpha_3)$ & $p=(2,1,0,\dots)$ 
& $[e_1,e_1,e_2]=\alpha_3e_3 $
& $[e_1,e_2,e_1]=\alpha_2e_3 $ \\&
& $[e_2,e_1,e_1]=\alpha_1e_3$  & \\ \hline
\end{longtable}
}

\noindent Here $\alpha_1+\alpha_2+\alpha_3=0$ and the triple $(\alpha_1,\alpha_2,\alpha_3)\neq (0,0,0)$ is determined (by the isomorphism class of the corresponding algebra $^1\mathbf{T}_{3}^{3}(\alpha_1,\alpha_2,\alpha_3)$ or $^1\mathbf{T}_{4}^{3}(\alpha_1,\alpha_2,\alpha_3)$) up to multiplication by a nonzero element of $\mathbb{C}$. Furthermore, $\epsilon\in\{0,1\}$, $\beta_1+\beta_2+\beta_3=\epsilon$, and if $\epsilon=0$, then $(\beta_1,\beta_2,\beta_3)\neq (0,0,0)$ is determined (by the isomorphism class of $^1\mathbf{T}_{2}^{3}(\epsilon,\beta_1,\beta_2,\beta_3)$) up to multiplication by a nonzero element of $\mathbb{C}$.

In dimension $n\geq 4$, the classification is the following:

{\small
\begin{longtable}{|l|l|ll|}
\caption*{ }   \\

\hline \multicolumn{1}{|c|}{${\mathcal A}$} & \multicolumn{1}{c|}{ } & \multicolumn{2}{c|}{Multiplication table} \\ \hline 
\endfirsthead

 \multicolumn{4}{l}%
{{\bfseries  continued from previous page}} \\
\hline \multicolumn{1}{|c|}{${\mathcal A}$} & \multicolumn{1}{c|}{ } & \multicolumn{2}{c|}{Multiplication table} \\ \hline 
\endhead

\hline \multicolumn{4}{|r|}{{Continued on next page}} \\ \hline
\endfoot

\hline 
\endlastfoot
$^1\mathbf{T}_{1}^{n}$ & $p=(3,0,\dots)$ & $[e_1,e_1,e_1]=e_2$ & \\ \hline
$^1\mathbf{T}_{2}^{n}(\epsilon,\beta_1,\beta_2,\beta_3)$ & $p=(2,0,\dots)$ 
& $[e_1,e_1,e_1]=\epsilon e_1 $
& $[e_1,e_1,e_i]=\beta_3e_i $ \\&
& $[e_1,e_i,e_1]=\beta_2e_i $ 
& $[e_i,e_1,e_1]=\beta_1e_i$  \\ \hline
$^1\mathbf{T}_{3}^{n}(\alpha_1,\alpha_2,\alpha_3)$ & $p=(1,1,0,\dots)$ 
& $[e_1,e_1,e_2]=(\alpha_2-\alpha_3)e_1 $
& $[e_1,e_2,e_1]=(\alpha_3-\alpha_1)e_1 $ \\&
& $[e_2,e_1,e_1]=(\alpha_1-\alpha_2)e_1 $ 
&  $[e_1,e_2,e_2]=(\alpha_2-\alpha_1)e_2$ \\&
& $ [e_2,e_1,e_2]=(\alpha_1-\alpha_3)e_2$ 
& $[e_2,e_2,e_1]=(\alpha_3-\alpha_2)e_2$ \\ &
& $[e_1,e_2,e_j]=\alpha_3e_j$ 
& $[e_1,e_j,e_2]=-\alpha_2e_j$ \\&
& $[e_2,e_1,e_j]=-\alpha_3e_j $ 
& $[e_2,e_j,e_1]=\alpha_2e_j$ \\&
& $[e_j,e_1,e_2]=\alpha_1e_j $
& $[e_j,e_2,e_1]=-\alpha_1e_j$\\ \hline
$^1\mathbf{T}_{4}^{n}(\alpha_1,\alpha_2,\alpha_3)$ & $p=(2,1,0,\dots)$  
& $[e_1,e_1,e_2]=\alpha_3e_3 $
& $[e_1,e_2,e_1]=\alpha_2e_3$ \\&
& $[e_2,e_1,e_1]=\alpha_1e_3$ & \\ \hline
$^1\mathbf{T}_{5}^{n}$ & $p=(1,1,1,0,\dots)$ 
& $[e_1,e_2,e_3]=e_4 $
& $[e_1,e_3,e_2]=-e_4$ \\&
& $[e_2,e_1,e_3]=-e_4$ 
& $[e_2,e_3,e_1]=e_4$  \\&
& $[e_3,e_1,e_2]=e_4$ 
& $[e_3,e_2,e_1]=-e_4$\\ \hline
\end{longtable}
}
\vspace{+7mm}

\noindent Here $2\leq i\leq n$ and $3\leq j\leq n$. Regarding the coefficients, we have that $\alpha_1+\alpha_2+\alpha_3=0$, and $(\alpha_1,\alpha_2,\alpha_3)\neq (0,0,0)$ is determined up to multiplication by a nonzero element of $\mathbb{C}$ (by the isomorphism class of $^1\mathbf{T}_{3}^{n}(\alpha_1,\alpha_2,\alpha_3)$ or $^1\mathbf{T}_{4}^{n}(\alpha_1,\alpha_2,\alpha_3)$). Moreover, $\epsilon\in\{0,1\}$, $\beta_1+\beta_2+\beta_3=\epsilon$, and if $\epsilon=0$, then $(\beta_1,\beta_2,\beta_3)\neq (0,0,0)$ is determined up to multiplication by a nonzero element of $\mathbb{C}$ (by the isomorphism class of $^1\mathbf{T}_{2}^{n}(\epsilon,\beta_1,\beta_2,\beta_3)$).

\EditInfo{October 16, 2024}{October 16, 2024}{Adam Chapman, Mohamed Elhamdadi and Ivan Kaygorodov}
\end{document}